\documentclass[10pt]{book}%
\usepackage{amsmath}
\usepackage{amsfonts}
\usepackage{amssymb}
\usepackage{graphicx}
\usepackage{makeidx}%
\usepackage{titletoc}
\usepackage{titlesec}
\usepackage[square,numbers]{natbib}

\usepackage[switch,pagewise,mathlines]{lineno}   

\providecommand{\U}[1]{\protect
\rule{.1in}{.1in}}

\newtheorem{theorem}{Theorem}[section]

\newtheorem{corollary}[theorem]{Corollary}

\newtheorem{definition}[theorem]{Definition}
\newtheorem{example}[theorem]{Example}
\newtheorem{exercise}[theorem]{Exercise}
\newtheorem{lemma}[theorem]{Lemma}

\newtheorem{proposition}[theorem]{Proposition}
\newtheorem{remark}[theorem]{Remark}

\newenvironment{proof}[1][Proof]{\noindent \textbf{#1.} }{\hfill$\Box$}
\makeindex 

\newpagestyle{MainPgsStylex}{\sethead[\thepage][][\scriptsize{{Chap.\thechapter\quad\chaptertitle}}]
{\scriptsize{\S{\thesection\quad\sectiontitle}}}{}{\thepage}\setfoot[][][]{}{}{}\headrule}

\newpagestyle{appendix}{\sethead[\thepage][][\scriptsize{{Appendix }}]
{\scriptsize{\S{\thesection\quad\sectiontitle}}}{}{\thepage}\setfoot[][][]{}{}{}\headrule}

\newpagestyle{bibliography}{\sethead[\thepage][][\scriptsize{{Bibliography}}]
{\scriptsize{{Bibliography}}}{}{\thepage}\setfoot[][][]{}{}{}\headrule}

\newpagestyle{index}{\sethead[\thepage][][\scriptsize{{Index}}]
{\scriptsize{{Index}}}{}{\thepage}\setfoot[][][]{}{}{}\headrule}

\newpagestyle{notation}{\sethead[\thepage][][\scriptsize{{Index of\ Notation}}]
{\scriptsize{{Index of\
Notation}}}{}{\thepage}\setfoot[][][]{}{}{}\headrule}

\titleformat{\section}{\Large\bfseries}{\S \arabic{section}}{1em}{}
\renewcommand \thesection {\arabic{section}}
\renewcommand \thechapter {\Roman{chapter}}

\titlecontents{chapter}[0pt]{\vspace{1.0\baselineskip}\bfseries}
    {\chaptertitlename \ {\thecontentslabel}\quad }{}
    {\hspace{.5em}\titlerule*[8pt]{$\cdot$}\contentspage}

\titlecontents{section}[10pt]{}
    {{\S\thecontentslabel}\quad}{}
    {\hspace{.5em}\titlerule*[8pt]{$\cdot$}\contentspage}

\begin{document}

\title{\bf{Nonlinear Expectations and Stochastic Calculus under Uncertainty}
\newline\newline\normalsize{---with Robust Central Limit Theorem and G-Brownian
Motion}}
\medskip\medskip
\author{{Shige PENG}\\Institute of Mathematics\\Shandong University\\250100, Jinan,
China\\peng@sdu.edu.cn}
\date{Version: first edition}
\maketitle
\frontmatter


\chapter*{Preface}

This book is focused on the recent developments on problems of
probability model under uncertainty by using the notion of nonlinear
expectations and, in particular, sublinear expectations. Roughly
speaking, a nonlinear expectation $\mathbb{E}$ is a monotone and
constant preserving functional defined on a linear space of random
variables. We are particularly interested in
sublinear expectations, i.e., $\mathbb{E}[X+Y]\leq \mathbb{{E}%
}[X]+\mathbb{E}[Y]$ for all random variables $X$, $Y$ and
$\mathbb{{E}}[\lambda X]=\lambda \mathbb{E}[X]$ if $\lambda \geq0$.

A sublinear expectation $\mathbb{E}$ can be represented as the upper
expectation of a subset of linear expectations $\{E_{\theta}:\theta
\in \Theta \}$, i.e., $\mathbb{E}[X]=\sup_{\theta \in
\Theta}E_{\theta}[X]$. In most cases, this subset is often treated
as an uncertain model of probabilities $\{P_{\theta}:\theta \in
\Theta \}$ and the notion of sublinear expectation provides a robust
way to measure a risk loss $X$. In fact, the sublinear expectation
theory provides many rich, flexible and elegant tools.

A remarkable point of view is that we emphasize the
term ``expectation'' rather than  the well-accepted classical notion
``probability'' and its non-additive counterpart ``capacity''. A
technical reason is that in general the information contained in a
nonlinear expectation $\mathbb{E}$ will be lost if one consider only
its corresponding
``non-additive probability'' or ``capacity'' $\mathbb{{P}}(A)=\mathbb{E}%
[\mathbf{1}_{A}]$. Philosophically, the notion of expectation has
its direct meaning of ``mean'', ``average'' which is not necessary
to be derived from the corresponding ``relative frequency'' which is
the origin of the probability measure. For example, when a person
gets a sample $\{x_{1},\cdots,x_{N}\}$ from a random variable $X$,
he can directly use $\overline{X}=\frac{1}{N}\sum x_{i}$ to
calculate its mean. In general he uses
$\overline{\varphi(X)}=\frac{1}{N}\sum \varphi(x_{i})$ for the mean
of $\varphi(X)$. We will discuss in detail this issue after the
overview of our new law of large numbers (LLN) and central limit
theorem (CLT).

A theoretical foundation of the above expectation framework is our
new LLN and CLT under sublinear expectations. Classical LLN and CLT
have been widely used in probability theory, statistics, data
analysis as well as in many practical situations such as financial
pricing and risk management. They provide a strong and convincing
way to explain why in practice normal distributions are so widely
utilized. But often a serious problem is that, in general, the
``i.i.d.'' condition is difficult to be satisfied\textrm{. }In
practice, for the most real-time processes and data for which the
classical trials and samplings become impossible, the uncertainty of
probabilities and distributions can not be neglected. In fact the
abuse of normal distributions in finance and many other industrial
or commercial domains has been criticized.

Our new CLT does not need this strong ``i.i.d.'' assumption. Instead
of fixing a probability measure $P$, we introduce an uncertain
subset of probability measures $\{P_{\theta}:\theta \in \Theta \}$
and consider the corresponding sublinear expectation
$\mathbb{E}[X]=\sup_{\theta \in \Theta}E_{\theta}[X]$. Our main
assumptions are:\newline

(i) The distribution of $X_{i}$ is within a subset
of distributions $\{F_{\theta}(x):\theta \in \Theta \}$ with%
\[
\overline{\mu}=\mathbb{E}[X_{i}]\geq \underline{\mu}=-\mathbb{{E}%
}[-X_{i}];\
\]

(ii) Any realization of $X_{1},\cdots,X_{n}$ does not change the
distributional uncertainty of $X_{n+1}$.

Under $\mathbb{E}$, we call $X_{1},X_{2},\cdots$ to be identically
distributed if condition (i) is satisfied, and we call $X_{n+1}$ is
independent from $X_{1},\cdots,X_{n}$ if condition (ii) is
fulfilled. Mainly under the above weak ``i.i.d.'' assumptions, we
have proved that for each  continuous function $\varphi$ with linear
growth we have the following LLN:
\[
\lim_{n\rightarrow \infty}\mathbb{E}[\varphi(\frac{S_{n}}{n}%
)]=\sup_{\underline{\mu}\leq v\leq \overline{\mu}}\varphi(v).
\]
Namely, the uncertain subset of the distributions of $S_{n}/n$ is
approximately a subset of dirac measures $\{
\delta_{v}:\underline{\mu}\leq v\leq \overline{\mu}\}$.

In particular, if $\underline{\mu}=\overline{\mu}=0$, then $S_{n}/n$
converges in law to $0$. In this case, if we assume furthermore that
$\overline
{\sigma}^{2}=\mathbb{E}[X_{i}^{2}]$ and $\underline{\sigma}%
^{2}=-\mathbb{E}[-X_{i}^{2}]$, $i=1,2,\cdots$, then we have the
following generalization of the CLT:
\[
\lim_{n\rightarrow \infty}\mathbb{E}[\varphi(S_{n}/\sqrt{n}%
)]=\mathbb{E}[\varphi(X)].
\]
Here $X$ is called $G$-normal distributed and denoted by $N(\{0\}\times[\underline{\sigma}^2,\overline{\sigma}^2])$. The value $\mathbb{E}%
[\varphi(X)]$ can be calculated by defining $u(t,x):=\mathbb{E}%
[\varphi(x+\sqrt{t}X)]$ which solves the partial differential equation (PDE) $\partial_{t}%
u=G(u_{xx})$ with $G(a):=\frac{1}{2}(\overline{\sigma}^{2}a^{+}%
-\underline{\sigma}^{2}a^{-})$. Our results reveal a deep and
essential relation between the theory of probability and statistics
under uncertainty and second order fully nonlinear parabolic
equations (HJB equations). We have two interesting situations: when
$\varphi$ is a convex function, then
\[
\mathbb{E}[\varphi(X)]=\frac{1}{\sqrt{2\pi \overline{\sigma}^{2}}}%
\int_{-\infty}^{\infty}\varphi(x)\exp(-\frac{x^{2}}{2\overline{\sigma}^{2}%
})dx,
\]
but if $\varphi$ is a concave function, the above
$\overline{\sigma}^{2}$ must be replaced by
$\underline{\sigma}^{2}$. If $\underline{\sigma}=\overline
{\sigma}=\sigma$, then $N(\{0\} \times \lbrack \underline{\sigma}^{2}%
,\overline{\sigma}^{2}])=N(0,\sigma^{2})$ which is a classical
normal distribution.

This result provides a new way to explain a well-known puzzle: many
practitioners, e.g., traders and risk officials in financial markets
can widely use normal distributions without serious data analysis or
even with data inconsistence. In many typical situations\
$\mathbb{E}[\varphi(X)]$ can be calculated by using normal
distributions with careful choice of parameters, but it is also a
high risk calculation if the reasoning behind has not been
understood.

We call $N(\{0\} \times \lbrack
\underline{\sigma}^{2},\overline{\sigma}^{2}])$ the $G$-normal
distribution. This new type of sublinear distributions was first
introduced in Peng (2006)\cite{Peng2006a} (see also
\cite{Peng2007b}, \cite{Peng2007}, \cite{Peng2008},
\cite{Peng2008b}) for a new type of ``$G$-Brownian motion''
and the related calculus of It\^{o}'s type. The main motivations
were uncertainties in statistics, measures of risk and superhedging
in finance (see El Karoui, Peng and Quenez (1997)
\cite{EPQ},\ Artzner, Delbaen, Eber and
Heath (1999) \cite{ADEH2}, Chen and Epstein (2002)
\cite{CE}, F\"{o}llmer and Schied (2004)
\cite{F-Sch}). Fully nonlinear
super-hedging is also a possible domain of applications (see Avellaneda, Levy and Paras (1995) \cite{Avellaneda}%
, Lyons (1995) \cite{Lyons}, see also  Cheridito, Soner, Touzi and Victoir (2007) \cite{Touzi} where a new
BSDE approach was introduced).

Technically we introduce a new method to prove our CLT on a
sublinear expectation space. This proof is short since we have
borrowed a deep interior estimate of fully nonlinear partial
differential equation (PDE) in Krylov (1987) \cite{Krylov1}. In
fact the theory of fully nonlinear parabolic PDE plays an essential
role in deriving our new results of LLN and CLT. In the classical
situation the corresponding PDE becomes a heat equation which is
often hidden behind its heat kernel, i.e., the normal distribution.
In this book we use the powerful notion of viscosity solutions for
our nonlinear PDE initially introduced by Crandall and Lions
(1983) \cite{CrandallL}. This notion is specially useful when the
equation is degenerate. For reader's convenience, we provide an
introductory chapter in Appendix C. If readers are only interested
in the classical non-degenerate cases, the corresponding solutions
will become smooth (see the last section of Appendix C).

We define a sublinear expectation on the space of continuous paths
from $\mathbb{R}^{+}$ to $\mathbb{R}^{d}$ which is an analogue of
Wiener's law, by which a $G$-Brownian motion is formulated. Briefly
speaking, a $G$-Brownian motion $(B_{t})_{t\geq0}$ is a continuous
process with independent and stationary increments under a given
sublinear expectation $\mathbb{E}$.

$G$--Brownian motion has a very rich and interesting new structure
which non-trivially generalizes the classical one. We can establish
the related stochastic calculus, especially It\^{o}'s integrals and
the related quadratic variation process $\left \langle B\right
\rangle $. A very interesting new phenomenon of our $G$-Brownian
motion is that its quadratic variation process $\left \langle
B\right \rangle $ is also a continuous process with independent and
stationary increments, and thus can be still regarded as a Brownian
motion. The corresponding $G$-It\^{o}'s formula is obtained. We have
also established the existence and uniqueness of solutions to
stochastic differential equation under our stochastic calculus by
the same Picard iterations as in the classical situation.

New norms were introduced in the notion of $G$-expectation by which
the corresponding stochastic calculus becomes significantly more
flexible and powerful. Many interesting, attractive and challenging
problems are also automatically provided within this new framework.

In this book we adopt a novel method to present our $G$-Brownian
motion theory. In the first two chapters as well as the first two
sections of Chapter III, our sublinear expectations are only assumed
to be finitely sub-additive, instead of ``$\sigma$-sub-additive''.
This is just because all the related results obtained in this part
do not need the ``$\sigma$-sub-additive'' assumption, and readers
even need not to have the background of classical probability
theory. In fact, in the whole part of the first five chapters we
only use a very basic knowledge of functional analysis such as
Hahn-Banach Theorem (see Appendix A). A special situation is when
all the sublinear expectations in this book become linear. In this
case this book can be still considered as using a new and very
simple approach to teach the classical It\^o's stochastic calculus,
since this book does not need the knowledge of probability theory.
This is an important advantage to use expectation as our basic
notion.

The ``authentic probabilistic parts'', i.e., the pathwise analysis
of our $G$-Brownian motion and the corresponding random variables,
view as functions of $G$-Brownian path, is presented in Chapter VI.
Here just as the classical ``$P$-sure analysis'', we introduce
``{\^c}-sure analysis'' for $G$-capacity {\^c}. Readers who are not
interested in the deep parts of stochastic analysis of $G$-Brownian
motion theory do not need to read this chapter.

This book was based on the author's Lecture Notes \cite{Peng2007b}
for several series of lectures, for the 2nd Workshop “Stochastic
Equations and Related Topic” Jena, July 23--29, 2006; Graduate
Courses of Yantai Summer School in Finance, Yantai University, July
06--21, 2007; Graduate Courses of Wuhan Summer School, July 24--26,
2007; Mini-Course of Institute of Applied Mathematics, AMSS, April
16--18, 2007; Mini-course in Fudan University, May 2007 and August
2009; Graduate Courses of CSFI, Osaka University, May 15--June 13,
2007; Minerva Research Foundation Lectures of Columbia University in
Fall of 2008; Mini-Workshop of $G$-Brownian motion and
$G$-expectations in Weihai, July 2009, a series talks in Department
of Applied Mathematics, Hong Kong Polytechnic University,
November-December, 2009 and an intensive course in WCU Center for
Financial Engineering, Ajou University. The hospitalities and
encouragements of the above institutions and the enthusiasm of the
audiences are the main engine to realize this lecture notes. I thank
for many comments and suggestions given during those courses,
especially to Li Juan and Hu Mingshang. During the preparation of
this book, a special reading group was organized with members Hu
Mingshang, Li Xinpeng, Xu Xiaoming, Lin Yiqing, Su Chen, Wang Falei
and Yin Yue. They proposed very helpful suggestions for the revision
of the book. Hu Mingshang and Li Xinpeng have made a great effort
for the final edition. Their efforts are decisively important to the
realization of this book.

\tableofcontents

\mainmatter \pagestyle{MainPgsStylex}

%
%

\chapter{Sublinear Expectations and Risk Measures}

\label{ch1} The sublinear expectation is also called the upper expectation
or the upper prevision, and this notion is used in situations when the
probability models have uncertainty. In this chapter, we present the basic
notion of sublinear expectations and the corresponding sublinear expectation
spaces. We give the representation theorem of a sublinear expectation and
the notions of distributions and independence under the framework of
sublinear expectations. We also introduce a natural Banach norm of a
sublinear expectation in order to get the completion of a sublinear
expectation space which is a Banach space. As a fundamentally important
example, we introduce the notion of coherent risk measures in finance. A
large part of notions and results in this chapter will be throughout this
book.

\section{Sublinear Expectations and Sublinear Expectation Spaces}

Let $\Omega$ be a given set and let $\mathcal{H}$ \label{huah}be a linear
space of real valued functions defined on $\Omega$. In this book, we suppose
that $\mathcal{H}$ satisfies $c\in\mathcal{H}$ for each constant $c$ and $%
|X|\in\mathcal{H}$ if $X\in\mathcal{H}$. The space $\mathcal{H}$ can be
considered as the space of random variables.

\begin{definition}
{\label{Def-1 copy(1)} {A \textbf{Sublinear expectation }\label{sube}
\index{Sublinear expectation}$\mathbb{E}$ is a functional $\mathbb{{E} }:%
\mathcal{H}\rightarrow \mathbb{R}$ satisfying }}

\noindent{\textbf{\textup{(i)} Monotonicity: }%
\begin{equation*}
\mathbb{E}[X]\geq \mathbb{E}[Y] \ \
\text{if}\ X\geq Y.
\end{equation*}
}

\noindent \textbf{\textup{(ii)} Constant preserving:}
\begin{equation*}
\mathbb{E}[c]=c\ \ \ \text{for}\ c\in \mathbb{R}.
\end{equation*}

\noindent \textbf{\textup{(iii)} Sub-additivity: } For each $X,Y\in{\mathcal{%
H}}$,
\begin{equation*}
\mathbb{E}[X+Y]\leq \mathbb{E}[X]+\mathbb{E}[Y].
\end{equation*}

\noindent \textbf{\textup{(iv)} Positive homogeneity:}
\begin{equation*}
\mathbb{E}[\lambda X]=\lambda \mathbb{E}[X]\ \ \ \text{for}\ \lambda \geq0%
\text{.}
\end{equation*}

The triple $(\Omega ,\mathcal{H},\mathbb{E}\mathbb{)}$ is called a \textbf{%
sublinear expectation space}%
\index{Sublinear expectation space}.\label{subls} If \textup{(i)} and
\textup{(ii)} are satisfied, {{$\mathbb{E}$ is called a }}\textbf{nonlinear
expectation} and
\index{Nonlinear expectation} the triple $(\Omega ,\mathcal{H},\mathbb{E}%
\mathbb{)}$ is called a \textbf{nonlinear expectation space}
\index{Nonlinear expectation space}.
\end{definition}

\begin{definition}\label{DefI.1.2}
Let $\mathbb{E}_{1}$ and $\mathbb{E}_{2}$ be two nonlinear expectations
defined on $(\Omega ,\mathcal{H})$. $\mathbb{E}_{1}$ is said to be \textbf{%
dominated} by $\mathbb{E}_{2}$ if
\begin{equation}
\mathbb{E}_{1}[X]-\mathbb{E}_{1}[Y]\leq \mathbb{E}_{2}[X-Y]\ \ ~%
\text{for}\ X,Y\in \mathcal{H}.  \label{DominationI.1.2}
\end{equation}
\end{definition}

\begin{remark}
From \textup{(iii)}, a sublinear expectation is dominated by itself. In many
situations, \textup{(iii)} is also called the property of self-domination.
If the inequality in \textup{(iii)} becomes equality, then {{$\mathbb{E}$}}
is a linear expectation, i.e., {{$\mathbb{E}$}} is a linear functional
satisfying \textup{(i)} and \textup{(ii)}.
\end{remark}

\begin{remark}
\textup{(iii)+(iv)} is called \textbf{sublinearity}%
\index{Sublinearity}. This sublinearity implies

\noindent\textbf{\textup{(v)}} \textbf{Convexity}:
\begin{equation*}
\mathbb{E}[\alpha X+(1-\alpha)Y]\leq \alpha \mathbb{E}[X]+(1-\alpha)\mathbb{E%
}[Y]\ \
\text{for}\ \alpha \in \lbrack0,1].
\end{equation*}
If a nonlinear expectation $\mathbb{E}$ satisfies convexity, we call it a
\textbf{convex expectation}.%
\index{Convex expectation}

The properties \textup{(ii)}\textup{+}\textup{(iii)} implies

\noindent\textup{\textbf{(vi)}} \textbf{Cash translatability}:
\begin{equation*}
\mathbb{E}[X+c]=\mathbb{E}[X]+c\ \
\text{for} \ c\in\mathbb{R}.
\end{equation*}
In fact, we have
\begin{align*}
\mathbb{E}[X]+c=\mathbb{E}[X]-\mathbb{E}[-c]\leq \mathbb{E}[X+c]\leq \mathbb{%
E}[X]+\mathbb{E}[c]=\mathbb{E}[X]+c.
\end{align*}

For property \textup{(iv)}, an equivalence form is
\begin{equation*}
{\mathbb{E}}[\lambda X]=\lambda ^{+}{\mathbb{E}}[X]+\lambda ^{-}{\mathbb{E}}%
[-X]\ ~\text{for}\ \lambda \in \mathbb{R}.
\end{equation*}
\end{remark}

In this book, we will systematically study the sublinear expectation spaces.
In the following chapters, unless otherwise stated, we consider the
following sublinear expectation space $(\Omega,\mathcal{H},\mathbb{E}\mathbb{%
)}$: if $X_{1},\cdots,X_{n}\in \mathcal{H}$ then $\varphi(X_{1},%
\cdots,X_{n})\in \mathcal{H}$ for each $\varphi \in C_{l.Lip}(\mathbb{R}^{n})
$\label{cllip} where $C_{l.Lip}(\mathbb{R}^{n})$ denotes the linear space of
functions $\varphi$ satisfying
\begin{align*}
|\varphi(x)-\varphi(y)| & \leq C(1+|x|^{m}+|y|^{m})|x-y|\ \ \text{for}\
x,y\in \mathbb{R}^{n}\text{, \ } \\
\ & \text{some }C>0\text{, }m\in \mathbb{N}\text{ depending on }\varphi.
\end{align*}
In this case $X=(X_{1},\cdots,X_{n})$ is called an $n$-dimensional random
vector, denoted by $X\in \mathcal{H}^{n}$.

\begin{remark}
It is clear that if $X\in \mathcal{H}$ then $|X|$, $X^{m}\in \mathcal{H}$.
More generally, $\varphi(X)\psi(Y)\in \mathcal{H}$ if $X,Y\in \mathcal{H}$
and $\varphi,\psi \in C_{l.Lip}(\mathbb{R})$. In particular, if $X\in
\mathcal{H}$ then $\mathbb{E}[|X|^{n}]<\infty$ for each $n\in \mathbb{N}$.
\end{remark}

Here we use $C_{l.Lip}(\mathbb{R}^{n})$ in our framework only for some
convenience of techniques. In fact our essential requirement is that $%
\mathcal{H}$ contains all constants and, moreover, $X\in \mathcal{H}$
implies $\left \vert X\right \vert \in \mathcal{H}$. In general, $C_{l.Lip}(%
\mathbb{R}^{n})$ can be replaced by any one of the following spaces of
functions defined on $\mathbb{R}^{n}$.

\begin{itemize}
\item {\ $\mathbb{L}^{\infty}(\mathbb{R}^{n})$: the space of bounded
Borel-measurable functions; }

\item {\ $C_{b}(\mathbb{R}^{n})$: the space of bounded and continuous
functions; \label{cbr}}

\item {$C_{b}^{k}(\mathbb{R}^{n})$: the space of bounded and $k$-time
continuously differentiable functions with bounded derivatives of all orders
less than or equal to $k$;\label{cbk}}

\item {\ $C_{unif}(\mathbb{R}^{n})$: the space of bounded and uniformly
continuous functions; \label{cu}}

\item {\ $C_{b.Lip}(\mathbb{R}^{n})$: the space of bounded and Lipschitz
continuous functions; \label{cbl}}

\item {\ $L^{0}(\mathbb{R}^{n})$: the space of Borel measurable functions. %
\label{llll0}}
\end{itemize}

Next we give two examples of sublinear expectations.

\begin{example}
{\ In a game we select a ball from a box containing $W$ white, $B$ black and
$Y$ yellow balls. The owner of the box, who is the banker of the game, does
not tell us the exact numbers of $W,B$ and $Y$. He or she only informs us
that $W+B+Y=100$ and $W=B\in \lbrack20,25]$. Let $\xi$ be a random variable
defined by
\begin{equation*}
\xi=\left \{
\begin{array}{rcc}
1 &  & \text{if we get a white ball;} \\
0 &  & \text{if we get a yellow ball;} \\
-1 &  & \text{if we get a black ball.}%
\end{array}
\right.
\end{equation*}
Problem: how to measure a loss $X=\varphi(\xi)$ for a given function $\varphi
$ on $\mathbb{R}$. }

We know that the distribution of $\xi$ is%
\begin{equation*}
\left \{
\begin{array}{ccc}
-1 & 0 & 1 \\
\frac{p}{2} & 1-p & \frac{p}{2}%
\end{array}
\right \} \ \ \text{with uncertainty: $p\in$}[\underline{\mu} ,\overline{\mu}%
]=[0.4,0.5].
\end{equation*}
Thus the \textbf{robust expectation}%
\index{Robust expectation} of $X=\varphi(\xi)$ is
\begin{align*}
{\mathbb{E}}[\varphi(\xi)] & :=\sup_{P\in \mathcal{P}}E_{P}[\varphi (\xi)] \\
& =\sup_{p\in \lbrack
\underline{\mu},\overline{\mu}]}[\frac {p}{2}(\varphi(1)+\varphi(-1))+(1-p)%
\varphi(0)].
\end{align*}
Here, $\xi$ has distribution uncertainty.
\end{example}

\begin{example}
{\ A more general situation is that the banker of a game can choose among a
set of distributions }${\{F(\theta,A)}\}_{A\in \mathcal{B}(\mathbb{R}%
),\theta \in \Theta}${\ of a random variable }$\xi${\ . In this situation
the robust expectation of a risk position $\varphi(\xi)$ for some $\varphi
\in C_{l.Lip}(\mathbb{R})$ is}
\begin{equation*}
\mathbb{E}[\varphi(\xi)]:=\sup_{\theta \in \Theta}\int_{\mathbb{R}%
}\varphi(x)F(\theta,dx).
\end{equation*}
\end{example}

\begin{exercise}
Prove that a functional $\mathbb{E}$ satisfies sublinearity if and only if
it satisfies convexity and positive homogeneity.
\end{exercise}

\begin{exercise}
Suppose that all elements in $\mathcal{H}$ are bounded. Prove that the
strongest sublinear expectation on $\mathcal{H}$ is
\begin{equation*}
\mathbb{E}^{\infty}[X]:=X^{\ast}=\sup_{\omega \in \Omega}X(\omega).
\end{equation*}
Namely, all other sublinear expectations are dominated by $\mathbb{E }%
^{\infty}[\cdot]$.
\end{exercise}

\section{Representation of a Sublinear Expectation}

A sublinear expectation can be expressed as a supremum of linear
expectations.

\begin{theorem}
\label{t1} Let {$\mathbb{E}$} be a functional defined on a linear space $%
\mathcal{H}$ satisfying sub-additivity and positive homogeneity. Then there
exists a family of linear functionals $\{E_{\theta}:\theta \in \Theta \}$
defined on $\mathcal{H}$ such that
\begin{equation*}
{\mathbb{E}}[X]=\sup_{\theta \in \Theta}E_{\theta}[X]\ \ \text{for}\ X\in
\mathcal{H}
\end{equation*}
and, for each $X\in \mathcal{H}$, there exists $\theta_{X}\in \Theta$ such
that {$\mathbb{E}$}$[X]=E_{\theta_{X}}[X]$.

Furthermore, if $\mathbb{E}$ is a sublinear expectation, then the
corresponding $E_{\theta}$ is a linear expectation.
\end{theorem}

\begin{proof}
Let $\mathcal{Q}=\{E_{\theta}:\theta \in \Theta \}$ be the family of all
linear functionals dominated by ${\mathbb{E}}$, i.e., $E_{\theta}[X]\leq{%
\mathbb{E}}[X]$, for all $X\in \mathcal{H}$, $E_{\theta}\in \mathcal{Q}.$

We first prove that $\mathcal{Q}$ is non empty. For a given $X\in \mathcal{H}
$, we set $L=\{aX:a\in{\mathbb{R}}\}$ which is a subspace of $\mathcal{H}$.
We define $I:L\rightarrow${$\mathbb{R}$} by $I[aX]=a{\mathbb{E}}[X]$, $%
\forall a\in${$\mathbb{R}$}, then $I[\cdot]$ forms a linear functional on $%
\mathcal{H}$ and $I\leq{\mathbb{E}}$ on $L$. Since ${\mathbb{E}}[\cdot]$ is
sub-additive and positively homogeneous, by Hahn-Banach theorem (see
Appendix A), there exists a linear functional $E$ on $\mathcal{H}$ such that
$E=I$ on $L$ and $E\leq{\mathbb{E}}$ on $\mathcal{H}$. Thus $E$ is a linear
functional dominated by ${\mathbb{E}}$ such that ${\mathbb{E}}[X]=E[X]$.

We now define%
\begin{equation*}
{{\mathbb{E}}}_{\Theta}[X]:=\sup_{\theta \in \Theta}E_{\theta }[X]\ \ \text{%
for }X\in \mathcal{H}.
\end{equation*}
It is clear that ${\mathbb{E}}_{\Theta}={\mathbb{E}}$.

Furthermore, if $\mathbb{E}$ is a sublinear expectation, then we have that,
for each nonnegative element $X\in$ $\mathcal{H}$, $E[X]=-E[-X]\geq-{\mathbb{%
E}}[-X]\geq0.$ For each $c\in${$\mathbb{R}$}, $-E[c]=E[-c]\leq {\mathbb{E}}%
[-c]=-c$ and $E[c]\leq{\mathbb{E}}[c]=c$, so we get $E[c]=c$. Thus $E$ is a
linear expectation. The proof is complete.
\end{proof}

\begin{remark}
\label{r1} It is important to observe that the above linear expectation $%
E_{\theta}$ is only ``finitely additive''. A sufficient condition for the $%
\sigma$-additivity of $E_{\theta}$ is to assume that ${\mathbb{E}}%
[X_{i}]\rightarrow0$ for each sequence $\{X_{i}\}_{i=1}^{\infty}$ of $%
\mathcal{H}$ such that\ $X_{i}(\omega)\downarrow0$ for each $\omega$. In
this case, it is clear that $E_{\theta}[X_{i}]\rightarrow0$. Thus we can
apply the well-known Daniell-Stone Theorem (see Theorem \ref{TheoremDL} in
Appendix B) to find a $\sigma $-additive probability measure $P_{\theta}$ on
$(\Omega,\sigma(\mathcal{H}))$ such that%
\begin{equation*}
E_{\theta}[X]=\int_{\Omega}X(\omega)dP_{\theta},\ \ \ X\in \mathcal{H}\text{.%
}
\end{equation*}
The corresponding model uncertainty of probabilities is the subset $%
\{P_{\theta}:\theta \in \Theta \}$, and the corresponding uncertainty of
distributions for an $n$-dimensional random vector $X$ in $\mathcal{H}$ is $%
\{F_{X}(\theta,A):=P_{\theta}(X\in A):A\in \mathcal{B}({\mathbb{R}}^{n})\}$.
\end{remark}

In many situation, we may concern the probability uncertainty, and the
probability maybe only finitely additive. So next we will give another
version of the above representation theorem.

Let $\mathcal{P}_{f}$ be the collection of all finitely additive probability
measures on $(\Omega,\mathcal{F})$, we consider $\mathbb{L}%
_{0}^{\infty}(\Omega,\mathcal{F})$ the collection of risk positions with
finite values, which consists risk positions $X$ of the form
\begin{equation*}
X(\omega)=\sum_{i=1}^{N}x_{i}\mathbf{I}_{A_{i}}(\omega),\ x_{i}\in \mathbb{R}%
,\ A_{i}\in \mathcal{F},i=1,\cdots,N.
\end{equation*}
It is easy to check that, under the norm $\left \Vert \cdot \right
\Vert
_{\infty}$, $\mathbb{L}_{0}^{\infty}(\Omega,\mathcal{F})$ is dense in $%
\mathbb{L}^{\infty}(\Omega,\mathcal{F})$. For a fixed $Q\in \mathcal{P}_{f}$
and $X\in \mathbb{L}_{0}^{\infty}(\Omega,\mathcal{F})$ we define
\begin{equation*}
E_{Q}[X]=E_{Q}[\sum_{i=1}^{N}x_{i}\mathbf{I}_{A_{i}}(\omega)]:=%
\sum_{i=1}^{N}x_{i}Q(A_{i})=\int_{\Omega}X(\omega)Q(d\omega).
\end{equation*}
$E_{Q}:\mathbb{L}_{0}^{\infty}(\Omega,\mathcal{F})\rightarrow \mathbb{R}$ is
a linear functional. It is easy to check that $E_{Q}$ satisfies (i)
monotonicity and (ii) constant preserving. It is also continuous under $%
\left \Vert X\right \Vert _{\infty}$.
\begin{equation*}
|E_{Q}[X]|\leq \sup_{\omega \in \Omega}|X(\omega)|=\left \Vert X\right \Vert
_{\infty}\text{. }
\end{equation*}
Since $\mathbb{L}_{0}^{\infty}$ is dense in $\mathbb{L}^{\infty}$ we then
can extend $E_{Q}$ from $\mathbb{L}_{0}^{\infty}$ to a linear continuous
functional on $\mathbb{L}^{\infty}(\Omega,\mathcal{F})$.

\begin{proposition}
\label{prop1} The linear functional $E_{Q}[\cdot]:${$\mathbb{L}%
^{\infty}(\Omega,\mathcal{F})\rightarrow \mathbb{R}$ satisfies \textup{(i)}
and \textup{(ii)}. Inversely each linear functional }$\eta(\cdot):${$\mathbb{%
L}^{\infty}(\Omega,\mathcal{F}) \rightarrow \mathbb{R}$} satisfying \textup{%
(i)} and \textup{(ii)} induces a finitely additive probability measure via $%
Q_{\eta}(A)=\eta(\mathbf{I}_{A})$, $A\in \mathcal{F}$. The corresponding
expectation is $\eta$ itself
\begin{equation*}
\eta(X)=\int_{\Omega}X(\omega)Q_{\eta}(d\omega).
\end{equation*}
\end{proposition}

\begin{theorem}
A sublinear expectation ${\mathbb{E}}$ has the following representation:
there exists a subset $\mathcal{Q}\subset\mathcal{P}_f$, such that
\begin{equation*}
{\mathbb{E}}[X]=\sup_{Q\in\mathcal{Q}}E_Q[X]\ ~\text{for}\ X\in \mathcal{H}.
\end{equation*}
\end{theorem}

\begin{proof}
By Theorem \ref{t1}, we have
\begin{equation*}
{\mathbb{E}}[X]=\sup_{\theta \in \Theta}E_{\theta}[X]\ \ \text{for}\ X\in
\mathcal{H},
\end{equation*}
where $E_{\theta}$ is a linear expectation on $\mathcal{H}$ for fixed $%
\theta\in\Theta$.

We can define a new sublinear expectation on $\mathbb{L}^\infty(\Omega,%
\sigma(\mathcal{H}))$ by
\begin{equation*}
\tilde{\mathbb{E}}_{\theta}[X]:=\inf\{E_{\theta}[Y]; Y\geq X, Y\in \mathcal{H%
}\}.
\end{equation*}
It is not difficult to check that $\tilde{\mathbb{E}}_{\theta}$ is a
sublinear expectation on $\mathbb{L}^{\infty}(\Omega,\sigma(\mathcal{H}))$,
where $\sigma{(\mathcal{H})}$ is the smallest $\sigma$-algebra generated by $%
\mathcal{H}$. We also have $E_{\theta}\leq\tilde{\mathbb{E}}_{\theta}$ on $%
\mathcal{H}$, by Hahn-Banach theorem, $E_{\theta}$ can be extended from $%
\mathcal{H}$ to $\mathbb{L}^{\infty}(\Omega,\sigma(\mathcal{H}))$, by
Proposition \ref{prop1}, there exists $Q\in\mathcal{P}_f$, such that
\begin{equation*}
E_{\theta}[X]=E_Q[X]\ ~\text{for}\ X\in\mathcal{H}.
\end{equation*}
So there exists $\mathcal{Q}\subset\mathcal{P}_f$, such that
\begin{equation*}
{\mathbb{E}}[X]=\sup_{Q\in \mathcal{Q}}E_Q[X]\ \ \text{for}\ X\in \mathcal{H}%
.
\end{equation*}
\end{proof}

\begin{exercise}
Prove that $\tilde{\mathbb{E}}_{\theta}$ is a sublinear expectation.
\end{exercise}

\section{Distributions, Independence and Product Spaces}

We now give the notion of distributions of random variables under sublinear
expectations.

Let $X=(X_{1},\cdots ,X_{n})$ be a given $n$-dimensional random vector on a
nonlinear expectation space $(\Omega ,\mathcal{H},\mathbb{E})$. We define a
functional on $C_{l.Lip}(\mathbb{R}^{n})$ by
\begin{equation*}
\mathbb{F}_{X}[\varphi ]:=\mathbb{E}[\varphi (X)]:\varphi \in C_{l.Lip}(%
\mathbb{R}^{n})\rightarrow \mathbb{R}.
\end{equation*}%
The triple $(\mathbb{R}^{n},C_{l.Lip}(\mathbb{R}^{n}),\mathbb{F}_{X})$ forms
a nonlinear expectation space. $\mathbb{F}_{X}$ is called the \textbf{%
distribution}
\index{Distribution}of $X$ under $\mathbb{E}$. This notion is very useful
for a sublinear expectation space $\mathbb{E}$. In this case $\mathbb{F}_{X}$
is also a sublinear expectation. Furthermore we can prove that (see Remark %
\ref{r1}), there exists a family of probability measures $\{F_{X}^{\theta
}(\cdot )\}_{\theta \in \Theta }$ defined on $(\mathbb{R}^{n},\mathcal{B}(%
\mathbb{R}^{n}))$ such that%
\begin{equation*}
\mathbb{F}_{X}[\varphi ]=\sup_{\theta \in \Theta }\int_{\mathbb{R}%
^{n}}\varphi (x)F_{X}^{\theta }(dx),\ \
\text{for each }\varphi \in C_{b.Lip}(\mathbb{R}^{n}).
\end{equation*}%
Thus $\mathbb{F}_{X}[\cdot ]$ characterizes the uncertainty of the
distributions of $X$.

\begin{definition}
\label{d1} Let $X_{1}$ and $X_{2}$ be two $n$--dimensional random vectors
defined on {nonlinear expectation spaces }$(\Omega _{1},\mathcal{H}_{1},%
\mathbb{E}_{1})${\ and }$(\Omega _{2},\mathcal{H}_{2},\mathbb{E}_{2})$,
respectively. They are called \textbf{identically distributed},%
\index{Identically distributed}\label{iidd} denoted by $X_{1}\overset{d}{=}%
X_{2}$, if
\begin{equation*}
\mathbb{E}_{1}[\varphi (X_{1})]=\mathbb{E}_{2}[\varphi (X_{2})]\ \ \
\text{for}\ \varphi \in C_{l.Lip}(\mathbb{R}^{n}).
\end{equation*}%
It is clear that $X_{1}\overset{d}{=}X_{2}$ if and only if their
distributions coincide. We say that the distribution of $X_{1}$ is stronger
than that of $X_{2}$ if $\mathbb{E}_{1}[\varphi (X_{1})]\geq \mathbb{E}%
_{2}[\varphi (X_{2})]$, for each $\varphi \in C_{l.Lip}(\mathbb{R}^{n})$.
\end{definition}

\begin{remark}
In the case of sublinear expectations, $X_{1}\overset{d}{=}X_{2}$ implies
that the uncertainty subsets of distributions of $X_{1}$ and $X_{2}$ are the
same, e.g., in the framework of Remark \ref{r1},
\begin{equation*}
\{F_{X_{1}}(\theta_{1},\cdot):\theta_{1}\in \Theta_{1}\}=\{F_{X_{2}}(\theta
_{2},\cdot):\theta_{2}\in \Theta_{2}\}.
\end{equation*}
Similarly if the distribution of $X_{1}$ is stronger than that of $X_{2}$,
then
\begin{equation*}
\{F_{X_{1}}(\theta_{1},\cdot):\theta_{1}\in
\Theta_{1}\}\supset\{F_{X_{2}}(\theta _{2},\cdot):\theta_{2}\in
\Theta_{2}\}.
\end{equation*}
\end{remark}

The distribution of $X\in \mathcal{H}$ has the following four typical
parameters:
\begin{equation*}
\bar{\mu}:={\mathbb{E}}[X],\ \ \underline{\mu}:=-\mathbb{E}[-X],\ \ \ \ \ \
\ \ \bar{\sigma}^{2}:={\mathbb{E}}[X^{2}],\ \ \underline {\sigma}^{2}:=-{%
\mathbb{E}}[-X^{2}].\ \
\end{equation*}
The intervals $[\underline{\mu},\bar{\mu}]$ and $[\underline{\sigma}^{2},%
\bar{\sigma}^{2}]$ characterize the \textbf{mean-uncertainty}%
\index{Mean-uncertainty} and the \textbf{variance-uncertainty}
\index{Variance-uncertainty}of $X$ respectively.

A natural question is: can we find a family of distribution measures
to represent the above sublinear distribution of $X$? The answer is
affirmative:

\begin{lemma}\label{I-le3}
Let $(\Omega ,\mathcal{H},\mathbb{E})$ be a sublinear expectation
space. Let $X\in \mathcal{H}^{d}$ be given.  Then for each sequence
$\{\varphi _{n}\}_{n=1}^{\infty }\subset C_{l.Lip}(\mathbb{R}^{d})$
satisfying $\varphi _{n}\downarrow 0$, we have
$\mathbb{{E}}[\varphi _{n}(X)]\downarrow 0$.
\end{lemma}

\begin{proof}
For each fixed $N>0$,
\begin{equation*}
\varphi _{n}(x)\leq k_{n}^{N}+\varphi _{1}(x)\mathbf{I}_{[|x|>N]}\leq k_{n}^{N}+\frac{%
\varphi _{1}(x)|x|}{N}\text{\ for each }x\in \mathbb{R}^{d\times m},
\end{equation*}%
where $k_{n}^{N}=\max_{|x|\leq N}\varphi _{n}(x)$.  We then have
\begin{equation*}
\mathbb{{E}}[\varphi _{n}(X)]\leq k_{n}^{N}+\frac{1}{N^{\delta }}\mathbb{{E}}[\varphi _{1}(X)|X|^{\delta }].
\end{equation*}%
It follows from $\varphi _{n}\downarrow 0$ that $k_{n}^{N}\downarrow
0$. Thus we have $\lim_{n\rightarrow \infty
}\mathbb{{E}}[\varphi _{n}(X)]\leq
\frac{C}{N}\mathbb{{E}}[\varphi _{1}(X)|X|]$. Since $N$ can
be arbitrarily large, we get $\mathbb{{E}}[\varphi _{n}(X)]\downarrow 0$%
.
\end{proof}

\begin{lemma}\label{I-le4}
Let $(\Omega ,\mathcal{H},\mathbb{E})$ be a sublinear expectation
space and let $\mathbb{F}_{X}[\varphi ]:=\mathbb{E}[\varphi (X)]$ be
the sublinear distribution of $X\in \mathcal{H}^{d}$. Then there
exists a family of probability measures $\left\{ F_{\theta }\right\}
_{\theta \in \Theta }$ defined on
$(\mathbb{R}^{d},\mathcal{B}(\mathbb{R}^{d}))$ such that
\begin{equation}
\mathbb{F}_{X}[\varphi ]=\sup_{\theta \in \Theta }\int_{\mathbb{R}%
^{d}}\varphi (x)F_{\theta }(dx),\ \ \varphi \in
C_{l,Lip}(\mathbb{R}^{d}). \label{Distr}
\end{equation}
\end{lemma}

\begin{proof}
By the representation theorem, for the sublinear expectation $\mathbb{F}%
_{X}[\varphi ]$ defined on
$(\mathbb{R}^{d},C_{l.Lip}(\mathbb{R}^{n}))$, there exists a family
of linear expectations $\{f_{\theta }\}_{\theta \in \Theta }$ on
$(\mathbb{R}^{d},C_{l.Lip}(\mathbb{R}^{n}))$ such that
\begin{equation*}
\mathbb{F}_{X}[\varphi ]=\sup_{\theta \in \Theta }f_{\theta
}[\varphi ],\ \varphi \in C_{l.Lip}(\mathbb{R}^{n}).
\end{equation*}%
By the above lemma, for each sequence $\{\varphi _{n}\}_{n=1}^{\infty }$ in $%
C_{b.Lip}(\mathbb{R}^{n})$ such that $\varphi _{n}\downarrow 0$ on $\mathbb{R%
}^{d}$, $\mathbb{F}_{X}[\varphi _{n}]\downarrow 0$, thus $f_{\theta
}[\varphi _{n}]\downarrow 0$ for each $\theta \in \Theta $. It
follows from Daniell-Stone Theorem (see Theorem \ref{TheoremDL} in
Appendix B) that, for
each $\theta \in \Theta $, there exists a unique probability measure $%
F_{\theta }(\cdot )$ on $(\mathbb{R}^{d},\sigma (C_{b.Lip}(\mathbb{R}^{d}))=(%
\mathbb{R}^{d},\mathcal{B}(\mathbb{R}^{d}))$, such that $f_{\theta
}[\varphi
]=$  $\int_{\mathbb{R}^{d}}\varphi (x)F_{\theta }(dx)$. Thus we have (\ref{Distr}).
\end{proof}

\begin{remark}
The above lemma tells us that in fact the sublinear distribution $\mathbb{F}%
_{X}$ of $X$ characterizes the uncertainty of distribution of $X$
which is an subset of distributions $\left\{ F_{\theta }\right\}
_{\theta \in \Theta } $.
\end{remark}

The following property is very useful in our sublinear expectation theory.

\begin{proposition}
{\ \label{Prop-X+Y}Let }$(\Omega,\mathcal{H},\mathbb{E})$ be a sublinear
expectation space and {\ $X,Y{\ }$be two random variables such that $\mathbb{%
E}[Y]=-\mathbb{E}[-Y]$, i.e., }$Y$ has no mean-uncertainty.{\ Then we have}%
\begin{equation*}
\mathbb{E}[X+\alpha Y]=\mathbb{E}[X]+\alpha \mathbb{E}[Y]\ \ \
\text{for}\ \alpha \in \mathbb{R}.
\end{equation*}
{In particular, if $\mathbb{E}[Y]=\mathbb{E}[-Y]=0$, then $\mathbb{E}%
[X+\alpha Y]=\mathbb{E}[X]$. }
\end{proposition}

\begin{proof}
We have
\begin{equation*}
\mathbb{E}[\alpha Y]=\alpha ^{+}\mathbb{E}[Y]+\alpha ^{-}\mathbb{E}%
[-Y]=\alpha ^{+}\mathbb{E}[Y]-\alpha ^{-}\mathbb{E}[Y]=\alpha \mathbb{E}[Y]\
\text{for}\ \alpha \in \mathbb{R}.
\end{equation*}%
{\ Thus }%
\begin{equation*}
\mathbb{E}[X+\alpha Y]\leq \mathbb{E}[X]+\mathbb{E}[\alpha Y]=\mathbb{E}%
[X]+\alpha \mathbb{E}[Y]=\mathbb{E}[X]-\mathbb{E}[-\alpha Y]\leq \mathbb{E}%
[X+\alpha Y].
\end{equation*}%
\bigskip
\end{proof}

A more general form of the above proposition is:

\begin{proposition}
\label{PropI.3.4}We make the same assumptions as the previous
proposition. Let $\mathbb{\tilde{E}}$ be a nonlinear expectation on
$(\Omega ,\mathcal{H})
$ dominated by the sublinear expectation $\mathbb{E}$ in the sense of (\ref
{DominationI.1.2}). If{\ $\mathbb{E}[Y]=\mathbb{E}[-Y]$, then we have}%
\begin{equation}
\mathbb{\tilde{E}}[\alpha Y]=\alpha \mathbb{\tilde{E}}[Y]=\alpha \mathbb{E}%
[Y],\ \ \alpha \in \mathbb{R}  \label{eI.3.1}
\end{equation}%
as well as{\ }%
\begin{equation}
{\mathbb{\tilde{E}}[X+\alpha Y]=\mathbb{\tilde{E}}[X]+\alpha \mathbb{\tilde{E%
}}[Y],\ \ }X\in \mathcal{H},\ {\alpha \in \mathbb{R}.}
\label{eI.3.2}
\end{equation}%
{\ }In particular
\begin{equation}
{\mathbb{\tilde{E}}[X+c]=\mathbb{\tilde{E}}[X]+c},\text{ for }c\in \mathbb{R}%
.  \label{eI.3.3}
\end{equation}
\end{proposition}

\begin{proof}
We have%
\begin{equation*}
-\mathbb{\tilde{E}}[Y]=\mathbb{\tilde{E}}[0]-\mathbb{\tilde{E}}[Y]\leq
\mathbb{E}[-Y]=-\mathbb{E}[Y]\leq -\mathbb{\tilde{E}}[Y]
\end{equation*}%
and%
\begin{eqnarray*}
\mathbb{E}[Y] &=&-\mathbb{E}[-Y]\leq -\mathbb{\tilde{E}}[-Y] \\
&=&\mathbb{\tilde{E}}[0]-\mathbb{\tilde{E}}[-Y]\leq
\mathbb{{E}}[Y].
\end{eqnarray*}%
From the above relations we have $\mathbb{\tilde{E}}[Y]=\mathbb{{E}}[Y]=-%
\mathbb{\tilde{E}}[-Y]$ and thus (\ref{eI.3.1}). Still by the domination,%
\begin{align*}
\mathbb{\tilde{E}}[X+\alpha Y]-\mathbb{\tilde{E}}[X]&\leq
\mathbb{E}[\alpha
Y],\\
\mathbb{\tilde{E}}[X]-\mathbb{\tilde{E}}[X+\alpha Y]&\leq\mathbb{E}%
[-\alpha Y]=-\mathbb{E}[\alpha Y].
\end{align*}
Thus (\ref{eI.3.2}) holds.
\end{proof}

\begin{definition}
A sequence of $n$-dimensional random vectors $\left\{ \eta _{i}\right\}
_{i=1}^{\infty }$ defined on a sublinear expectation space $(\Omega ,%
\mathcal{H},\mathbb{E})$ is said to \textbf{converge in distribution}
\index{Converge in distribution}(or \textbf{converge in law})%
\index{Converge in law} under $\mathbb{E}$ if for each $\varphi \in
C_{b.Lip}(\mathbb{R}^{n})$, the sequence $\left\{ \mathbb{E}[\varphi (\eta
_{i})]\right\} _{i=1}^{\infty }$ converges.
\end{definition}

The following result is easy to check.

\begin{proposition}
Let $\left \{ \eta_{i}\right \} _{i=1}^{\infty}$ converge in law in the
above sense. Then the mapping $\mathbb{F}[\cdot]:C_{b.Lip}(\mathbb{R}%
^{n})\rightarrow \mathbb{R}$ defined by%
\begin{equation*}
\mathbb{F}[\varphi]:=\lim_{i\rightarrow \infty}\mathbb{E}[\varphi(\eta_{i})]%
\ \
\text{for}\ \varphi \in C_{b.Lip}(\mathbb{R}^{n})
\end{equation*}
is a sublinear expectation defined on $(\mathbb{R}^{n},C_{b.Lip}(\mathbb{R}%
^{n}))$.
\end{proposition}

The following notion of independence plays a key role in the
nonlinear expectation theory.

\begin{definition}
\label{d2} In a nonlinear expectation space $(\Omega,\mathcal{H},\mathbb{E})$%
, a random vector $Y\in \mathcal{H}^n$ is said to be \textbf{independent}
\index{Independent} from another random vector $X\in \mathcal{H}^m$ under $%
\mathbb{E}[\cdot]$ if for each test function $\varphi \in C_{l.Lip}(\mathbb{R%
}^{m+n})$ we have
\begin{equation*}
\mathbb{E}[\varphi(X,Y)]=\mathbb{E}[\mathbb{E}[\varphi(x,Y)]_{x=X}].
\end{equation*}
\end{definition}

\begin{remark}
In particular, for a sublinear expectation space
$(\Omega,\mathcal{H},\mathbb{E})$, $Y$ is
independent from $X$ means that the uncertainty of distributions $%
\{F_{Y}(\theta,\cdot):\theta \in \Theta \}$ of $Y$ does not change after the
realization of $X=x$. In other words, the ``conditional sublinear
expectation'' of $Y$ with respect to $X$ is $\mathbb{E}[\varphi(x,Y)]_{x=X}$%
. In the case of linear expectation, this notion of independence is just the
classical one.
\end{remark}

\begin{remark}
It is important to note that under sublinear expectations the condition
\textquotedblleft$Y$ is independent from $X$\textquotedblright \ does not
imply automatically that \textquotedblleft$X$ is independent from $Y$%
\textquotedblright.
\end{remark}

\begin{example}
We consider a case where $\mathbb{E}$ is a sublinear expectation and $X,Y\in
\mathcal{H}$ are identically distributed with $\mathbb{E}[X]=\mathbb{E}[-X]=0
$ and $%
\overline{\sigma}^{2}=\mathbb{E}[X^{2}]>\underline{\sigma}^{2}=-\mathbb{E}%
[-X^{2}]$. We also assume that $\mathbb{E}[|X|]=\mathbb{E}[X^{+}+X^{-}]>0$,
thus $\mathbb{\ {E}}[X^{+}]=\frac{1}{2}\mathbb{E}[|X|+X]=$$\frac{1}{2}\mathbb%
{E}[|X|]>0$. In the case where $Y$ is independent from $X$, we have%
\begin{equation*}
\mathbb{E}[XY^{2}]=\mathbb{E}[X^{+}\overline{\sigma}^{2}-X^{-}\underline{%
\sigma}^{2}]=(\overline{\sigma}^{2}-\underline{\sigma}^{2})\mathbb{E}%
[X^{+}]>0.
\end{equation*}
But if $X$ is independent from $Y$, we have%
\begin{equation*}
\mathbb{E}[XY^{2}]=0.
\end{equation*}
\end{example}

The independence property of two random vectors $X,Y$ involves only the
``joint distribution'' of $(X,Y)$. The following result tells us how to
construct random vectors with given ``marginal distributions'' and with a
specific direction of independence.

\begin{definition}
\label{dd1} Let $(\Omega_{i},\mathcal{H}_{i},\mathbb{E}_{i})$,
$i=1,2$ be two sublinear (resp. nonlinear) expectation spaces. We
denote
\begin{align*}
\mathcal{H}_{1}\otimes \mathcal{H}_{2} &
:=\{Z(\omega_{1},\omega_{2})=\varphi(X(\omega_{1}),Y(\omega_{2})):(%
\omega_{1},\omega_{2})\in \Omega _{1}\times \Omega_{2},\  \\
\ \ & \ (X,Y)\in\mathcal{H}_{1}^{m}\times\mathcal{H}_{2}^{n},\ \varphi \in
C_{l.Lip}(\mathbb{R}^{m+n}) \},\ \
\end{align*}
and, for each random variable of the above form $Z(\omega_{1},\omega
_{2})=\varphi(X(\omega_{1}),Y(\omega_{2}))$,
\begin{equation*}
(\mathbb{E}_{1}\otimes\mathbb{E}_{2})\mathbb{[}Z]:=\mathbb{E}_{1}\mathbb{[}%
\bar{\varphi}(X)],\ \ \text{where }\bar{\varphi}(x):=\mathbb{E}_{2}\mathbb{[}%
\varphi(x,Y)],\ x\in \mathbb{R}^{m}.
\end{equation*}
It is easy to check that the triple $(\Omega_{1}\times \Omega_{2},\mathcal{H}%
_{1}\otimes \mathcal{H}_{2},\mathbb{E}_{1}\otimes \mathbb{E}_{2})$ forms a
sublinear (resp. nonlinear) expectation space. We call it the \textbf{product space} \label%
{ppss}of sublinear (resp. nonlinear) expectation spaces
\index{Product space of sublinear expectation space} $(\Omega_{1},\mathcal{H}%
_{1},$ $\mathbb{E}_{1})$ and $(\Omega _{2},\mathcal{H}_{2},$ $\mathbb{E}_{2})
$. In this way, we can define the product space

\begin{equation*}
(\prod_{i=1}^n\Omega_i,\bigotimes_{i=1}^n\mathcal{H}_i,\bigotimes_{i=1}^n%
\mathbb{E}_i)
\end{equation*}

of given sublinear (resp. nonlinear) expectation spaces $(\Omega_{i},\mathcal{H}_{i},\mathbb{E}%
_{i})$, $i=1,2,\cdots,n$. In particular, when $(\Omega _{i},\mathcal{H}_{i},%
\mathbb{E}_{i})=(\Omega_{1},\mathcal{H}_{1},$ $\mathbb{E}_{1})$ we have the
product space of the form $(\Omega _{1}^{n},\mathcal{H}_{1}^{\otimes n},$ $%
\mathbb{E}_{1}^{\otimes n})$.
\end{definition}

Let $X,%
\bar{X}$ be two $n$-dimensional random vectors on a sublinear (resp.
nonlinear) expectation space $(\Omega,\mathcal{H},\mathbb{E})$.
$\bar{X}$ is called an independent
copy of $X$ if $\bar{X}\overset{d}{=}X$ and $\bar{X}$ is independent from $X$%
.

The following property is easy to check.

\begin{proposition}
\label{pp1} Let $X_{i}$ be an $n_{i}$-dimensional random vector on
sublinear (resp. nonlinear)
expectation space $(\Omega_{i},\mathcal{H}_{i},\mathbb{E}_{i})$ for $%
i=1,\cdots,n$, respectively. We denote
\begin{equation*}
Y_{i}(\omega_{1},\cdots,\omega_{n}):=X_{i}(\omega_{i}),\ \ i=1,\cdots,n.
\end{equation*}
Then $Y_{i}$, $i=1,\cdots,n$, are random vectors on $(\prod_{i=1}^n\Omega_i,%
\bigotimes_{i=1}^n\mathcal{H}_i,\bigotimes_{i=1}^n\mathbb{E}_i)$. Moreover
we have $Y_{i}\overset{d}{=}X_{i}$ and $Y_{i+1}$ is independent from $%
(Y_{1},\cdots,Y_{i})$, for each $i$.

Furthermore, if $(\Omega_{i},\mathcal{H}_{i},\mathbb{E}_{i})=(\Omega _{1},%
\mathcal{H}_{1},$ $\mathbb{E}_{1})$ and $X_{i}\overset{d}{=}X_{1}$, for all $%
i$, then we also have $Y_{i}\overset{d}{=}Y_{1}$. In this case $Y_{i}$ is
said to be an \textbf{independent copy}%
\index{Independent copy} of $Y_{1}$ for $i=2,\cdots,n$. \
\end{proposition}

\begin{remark}
In the above construction the integer $n$ can be also infinite. In this case
each random variable $X\in\bigotimes_{i=1}^\infty \mathcal{H}_{i}$ belongs
to $(\prod_{i=1}^k\Omega_i,\bigotimes_{i=1}^{k} \mathcal{H}%
_{i},\bigotimes_{i=1}^{k}\mathbb{E}_{i})$ for some positive integer $k<\infty
$ and
\begin{equation*}
\bigotimes_{i=1}^{\infty}\mathbb{E}_{i}[X]:=\bigotimes _{i=1}^{k} \mathbb{E}%
_{i}[X].
\end{equation*}
\end{remark}

\begin{remark}
{The situation \textquotedblleft$Y$ is independent from $X$%
\textquotedblright often appears when $Y$ occurs after $X$, thus a robust
expectation should take the information of $X$ into account. }
\end{remark}

\begin{exercise}
Suppose $X,Y\in \mathcal{H}^{d}$ and $Y$ is an independent copy of $X$.
Prove that for each $a\in \mathbb{R},b\in \mathbb{R}^{d}$,$a+\langle
b,Y\rangle$ is an independent copy of $a+\langle b,X\rangle$.
\end{exercise}

In a sublinear expectation space we have:

\begin{example}
{We consider a situation where two random variables $X$ and $Y$ in $\mathcal{%
H}$ are identically distributed and their common distribution is
\begin{equation*}
\mathbb{F}_{X}[\varphi]=\mathbb{F}_{Y}[\varphi]=\sup_{\theta \in
\Theta}\int_{\mathbb{R}}\varphi(y)F(\theta,dy)\ \ \text{for}\
\varphi \in C_{l.Lip}(\mathbb{R}),
\end{equation*}
where for each $\theta \in \Theta$, $\{F(\theta,A)\}_{A\in \mathcal{B}(%
\mathbb{R})}$ is a probability measure on $(\mathbb{R},\mathcal{B}(\mathbb{R}%
))$. In this case, "$Y$ is independent from $X$" means that the
joint distribution of $X$ and $Y$ is
\begin{align*}
\mathbb{F}_{X,Y}[\psi]=\sup_{\theta_{1}\in
\Theta}\int_{\mathbb{R}}\left[ \sup_{\theta_{2}\in
\Theta}\int_{\mathbb{R}}\psi(x,y)F(\theta _{2},dy)\right]
F(\theta_{1},dx)\ \ \text{for}\ \psi \in C_{l.Lip}(\mathbb{R}^{2}).
\end{align*}
}
\end{example}

\begin{exercise}
Let $(\Omega,\mathcal{H},\mathbb{E})$ be a sublinear expectation space.
Prove that if $\mathbb{E}[\varphi(X)]=\mathbb{E}[\varphi(Y)]$ for any $%
\varphi\in C_{b,Lip}$, then it still holds for any $\varphi\in C_{l,Lip}$.
That is, we can replace $\varphi\in C_{l,Lip}$ in Definition \ref{d1} by $%
\varphi\in C_{b,Lip}$.
\end{exercise}

\section{Completion of Sublinear Expectation Spaces}

\label{c1s5} Let $(\Omega,\mathcal{H},\mathbb{E})$ be a sublinear
expectation space. We have the following useful inequalities. {\ }

We first give the following well-known inequalities.

\begin{lemma}
{\ {F{or $r>0$ and $1<p,q<\infty$ with $\frac{1}{p}+\frac{1}{q}=1$, we have
\begin{align}
|a+b|^{r} & \leq \max \{1,2^{r-1}\}(|a|^{r}+|b|^{r})\ \ \text{for}\ a,b\in
\mathbb{R},  \label{ee04.3} \\
|ab| & \leq \frac{|a|^{p}}{p}+\frac{|b|^{q}}{q}.   \label{ee04.4}
\end{align}
} } }
\end{lemma}

\begin{proposition}
\label{ppp1} {\ {\ For each $X,Y\in$}}$\mathcal{H}${{, we have{\
\begin{align}
\mathbb{E}[|X+Y|^{r}] & \leq2^{r-1}(\mathbb{E}[|X|^{r}]+\mathbb{{E}[}%
|Y|^{r}]),  \label{ee04.5} \\
\mathbb{E}[|XY|] & \leq(\mathbb{E}[|X|^{p}])^{1/p}\cdot(\mathbb{E}%
[|Y|^{q}])^{1/q},  \label{ee04.6} \\
(\mathbb{E}[|X+Y|^{p}])^{1/p} & \leq(\mathbb{E}[|X|^{p}])^{1/p}+(\mathbb{E}%
[|Y|^{p}])^{1/p},   \label{ee04.7}
\end{align}
where {$r$}}}}${{{{\geq1}}}}${{{{{\ and $1<p,q<\infty$ with $\frac{1}{p}+%
\frac{1}{q}=1.$}} }}}

In particular, for $1\leq p<p^{\prime}$, we have $(\mathbb{E}%
[|X|^{p}])^{1/p}\leq(\mathbb{E}[|X|^{p^{\prime}}])^{1/p^{\prime}}.$
\end{proposition}

\begin{proof}
The inequality (\ref{ee04.5}) follows from (\ref{ee04.3}).

For the case $\mathbb{E}[|X|^{p}]\cdot \mathbb{E}[|Y|^{q}]>0$, we set
\begin{equation*}
\xi=\frac{X}{(\mathbb{E}[|X|^{p}])^{1/p}},\ \ \eta=\frac{Y}{(\mathbb{E}%
[|Y|^{q}])^{1/q}}.
\end{equation*}
By (\ref{ee04.4}) we have%
\begin{align*}
\mathbb{E}[|\xi \eta|] & \leq \mathbb{E}[\frac{|\xi|^{p}}{p}+\frac{|\eta|^{q}%
}{q}]\leq \mathbb{E}[\frac{|\xi|^{p}}{p}]+\mathbb{E}[\frac{|\eta|^{q}}{q}] \\
& =\frac{1}{p}+\frac{1}{q}=1.
\end{align*}
Thus (\ref{ee04.6}) follows.

For the case $\mathbb{E}[|X|^{p}]\cdot \mathbb{E}[|Y|^{q}]=0,$ we consider $%
\mathbb{E}[|X|^{p}]+\varepsilon$ and $\mathbb{E}[|Y|^{q}]+\varepsilon$ for $%
\varepsilon>0.$ Applying the above method and letting $\varepsilon
\rightarrow0,$ we get (\ref{ee04.6}).

We now prove (\ref{ee04.7}). We only consider the case ${\mathbb{E}}%
[|X+Y|^p]>0$.
\begin{align*}
\mathbb{E}[|X+Y|^{p}] & =\mathbb{E}[|X+Y|\cdot|X+Y|^{p-1}] \\
& \leq \mathbb{E}[|X|\cdot|X+Y|^{p-1}]+\mathbb{E}[|Y|\cdot |X+Y|^{p-1}] \\
& \leq(\mathbb{E}[|X|^{p}])^{1/p}\cdot(\mathbb{{E}[}|X+Y|^{(p-1)q}])^{1/q} \\
& \ \ \ +(\mathbb{E}[|Y|^{p}])^{1/p}\cdot(\mathbb{{E}[}%
|X+Y|^{(p-1)q}])^{1/q}.
\end{align*}
Since $(p-1)q=p$, we have (\ref{ee04.7}).

{By}(\ref{ee04.6}), it is easy to deduce that {{{$(\mathbb{E}%
[|X|^{p}])^{1/p}\leq (\mathbb{E}[|X|^{p^{\prime}}])^{1/p^{\prime}}$ for $%
1\leq p<p^{\prime}.$}}}
\end{proof}

\bigskip

{\ {For each fixed }$p\geq1$, we observe that $\mathcal{H}_{0}^{p}=\{X\in
\mathcal{H}$, $\mathbb{E}[|X|^{p}]=0\}$ is a linear subspace of $\mathcal{H}$%
. Taking $\mathcal{H}_{0}^{p}$ as our null space, we introduce the quotient
space $\mathcal{H}/\mathcal{H}_{0}^{p}$. Observing that, for every $\mathbf{%
\{}X\} \in \mathcal{H}/\mathcal{H}_{0}^{p}$ with a representation $X\in
\mathcal{H}$, we can define an expectation $\mathbb{E}\mathbf{[\{}X\}]:=%
\mathbb{E}[X]$ which is still a sublinear expectation. We set $\left \Vert
X\right \Vert _{p}:=(\mathbb{E}[|X|^{p}])^{\frac{1}{p}}$. By Proposition \ref%
{ppp1}, i{{t is easy to check that $\left
\Vert \cdot \right \Vert _{p}$
forms a Banach norm on {{$\mathcal{H}/\mathcal{H}_{0}^{p}$}}. We extend $%
\mathcal{H}/\mathcal{H}_{0}^{p}$ to its completion $\mathcal{\hat{H}}_{p}$
under\ this norm, then $(\mathcal{\hat{H}}_{p},\left \Vert \cdot \right
\Vert _{p})$ is a Banach space. \ In particular, when }}}$p=1,$ we denote it
by {{{$(\mathcal{\hat{H}},\left \Vert \cdot \right \Vert ).$}}}

{{{For each $X\in \mathcal{H}$, the mappings
\begin{equation*}
X^{+}(\omega):\mathcal{H\rightarrow H}\ \ \ \text{and \ \ }X^{-}(\omega):%
\mathcal{H\rightarrow H}
\end{equation*}
satisfy
\begin{equation*}
|X^{+}-Y^{+}|\leq|X-Y|\text{ \ \ and \ }\
|X^{-}-Y^{-}|=|(-X)^{+}-(-Y)^{+}|\leq|X-Y|.
\end{equation*}
Thus they are both contraction mappings under $\left \Vert \cdot \right
\Vert _{p}$ and can be continuously extended to the Banach space $(\mathcal{%
\hat{H}}_{p},\left \Vert \cdot \right \Vert _{p})$. } } }

{\ {\ {\ We can define the partial order \textquotedblleft$\geq$%
\textquotedblright \ in this Banach space. }}}

\begin{definition}
{\ {\ {\ An element $X$ in $(\mathcal{\hat{H}},\left \Vert \cdot
\right
\Vert )$ is said to be nonnegative, or $X\geq0$, $0\leq X$, if $%
X=X^{+}$. We also denote by $X\geq Y$, or $Y\leq X$, if $X-Y\geq0$. } } }
\end{definition}

{\ {\ {\ It is easy to check that $X\geq Y$ and $Y\geq X$ imply $X=Y$ on $(%
\mathcal{\hat{H}}_{p},\left \Vert \cdot \right \Vert _{p})$. }}}

For each {{{$X,Y\in \mathcal{H}$, note that }}}

\begin{equation*}
|{{{\mathbb{E}[X]-\mathbb{E}[Y]|\leq \mathbb{E}[|X-Y|]\leq||X-Y||}}}_{p}.
\end{equation*}
We then can define
\begin{definition}\label{DefI4.4}
 The {{{sublinear expectation $\mathbb{E}[\cdot]$ can be
continuously extended to $(\mathcal{\hat{H}}_{p},\left \Vert \cdot
\right \Vert _{p})$ on which it is still a sublinear expectation. We
still denote by $(\Omega,\mathcal{\hat{H}}_{p},\mathbb{E})$. }}}

Let $(\Omega,\mathcal{H},\mathbb{E}_{1})$ be a nonlinear expectation space. $%
\mathbb{E}_{1}$ is said to be \ dominated by $\mathbb{E}$ if

\begin{equation*}
\mathbb{E}_{1}[X]-\mathbb{E}_{1}[Y]\leq \mathbb{E}[X-Y]\ \ \ \ \text{for}\
X,Y\in \mathcal{H}.
\end{equation*}
From this we can easily deduce that $|{{{\mathbb{E}}}}_{1}{{{[X]-\mathbb{E}}}%
}_{1}{{{[Y]|\leq \mathbb{E}[|X-Y|],}}}$ thus the nonlinear {{{expectation $%
\mathbb{E}_{1}[\cdot]$ can be continuously extended to $(\mathcal{\hat{H}}%
_{p},\left \Vert \cdot \right \Vert _{p})$ on which it is still a
nonlinear expectation. We still denote by
$(\Omega,\mathcal{\hat{H}}_{p},\mathbb{E}_1)$.}}}
\end{definition}

\begin{remark}
It is important to note that $X_{1},\cdots,X_{n}\in \mathcal{\hat{H}}$ does
not imply $\varphi(X_{1},\cdots,X_{n})\in \mathcal{\hat{H}}$ for each $%
\varphi \in C_{l.Lip}(\mathbb{R}^{n}).$ Thus, when we talk about the notions
of distributions, independence and product spaces on $(\Omega,\mathcal{\hat{H%
}},\mathbb{E}),$ the space $C_{l.Lip}(\mathbb{R}^{n})$ is replaced by $%
C_{b.Lip}(\mathbb{R}^{n})$ unless otherwise stated.
\end{remark}

\begin{exercise}
Prove that the inequalities \textup{(\ref{ee04.5})},\textup{(\ref{ee04.6})},%
\textup{(\ref{ee04.7})} still hold for $(\Omega,\mathcal{\hat{H}},\mathbb{E}%
).$
\end{exercise}

\section{Coherent Measures of Risk}

Let the pair $(\Omega,\mathcal{H})$ be such that $\Omega$ is a set of
scenarios and $\mathcal{H}$ is the collection of all possible risk positions
in a financial market.

If $X\in \mathcal{H}$, then for each constant $c$, $X\vee c$, $X\wedge c$
are all in $\mathcal{H}$. One typical example in finance is that $X$ is the
tomorrow's price of a stock. In this case, any European call or put options
with strike price $K$ of forms $(S-K)^{+},\ \ (K-S)^{+}$ are in $\mathcal{H}$%
.

A risk supervisor is responsible for taking a rule to tell traders,
securities companies, banks or other institutions under his supervision,
which kind of risk positions is unacceptable and thus a minimum amount of
risk capitals should be deposited to make the positions acceptable. The
collection of acceptable positions is defined by
\begin{equation*}
\mathcal{{A}}=\{X\in \mathcal{H}: X\text{ is acceptable}\}.
\end{equation*}
This set has meaningful properties in economy.

\begin{definition}
\index{Coherent acceptable set}A set $\mathcal{A}$\label{mathcalA} is called
a \textbf{coherent acceptable set} if it satisfies \newline
\textbf{\textup{(i)} } {\textbf{Monotonicity:}
\begin{equation*}
X\in \mathcal{A},\;Y\geq X\; \;
\text{imply} \; \;Y\in \mathcal{A}.
\end{equation*}
\textbf{\textup{(ii)} }$0\in \mathcal{A}$ but $-1\not \in \mathcal{A}$. }%
\newline
\textbf{\textup{(iii)} }{\textbf{Positive homogeneity}
\begin{equation*}
X\in \mathcal{A}\; \; \text{implies} \; \; \lambda X\in
\mathcal{A}\; \; \text{for}\ \lambda \geq0.
\end{equation*}
\textbf{\textup{(iv)}} \textbf{Convexity:}
\begin{equation*}
X,Y\in \mathcal{A} \ \ \text{imply} \ \ \alpha X+(1-\alpha)Y\in \mathcal{A}%
\; \ \text{for}\ \alpha \in[0,1].
\end{equation*}
}
\end{definition}

\begin{remark}
{\textup{(iii)+}\textup{(iv)} imply \newline
\textbf{\textup{(v)}} \textbf{Sublinearity:}
\begin{equation*}
X,Y\in \mathcal{A}\Rightarrow \mu X+\nu Y\in \mathcal{A}\; \ \text{for}\ \mu
,\nu \geq0.
\end{equation*}
}
\end{remark}

\begin{remark}
If the set $\mathcal{A}$ only satisfies \textup{(i)},\textup{(ii)} and
\textup{(iv)}, then $\mathcal{A}$ is called a \textbf{convex acceptable set}%
\index{Convex acceptable set}.
\end{remark}

In this section we mainly study the coherent case. Once the rule of the
acceptable set is fixed, the minimum requirement of risk deposit is then
automatically determined.

\begin{definition}
Given a coherent acceptable set $\mathcal{A}$, the functional $\rho(\cdot)$
defined by
\begin{equation*}
\rho(X)=\rho_{\mathcal{A}}(X):=\inf \{m\in \mathbb{R}:\;m+X\in \mathcal{A}%
\},\ \ \ X\in \mathcal{H}
\end{equation*}
is called the \textbf{coherent risk measure}%
\index{Coherent risk measure} \label{rho}related to $\mathcal{A}$.
\end{definition}

It is easy to see that
\begin{equation*}
\rho(X+\rho(X))=0.
\end{equation*}
\begin{proposition}
{\ \textbf{\ }$\rho(\cdot)$ is a coherent risk measure satisfying four
properties:\newline
\newline
\textbf{\textup{(i)} Monotonicity:} If $X\geq Y$ then $\rho(X)\leq \rho(Y)$%
.\thinspace \ \thinspace \ \newline
\textbf{\textup{(ii)} Constant preserving:}$\mathbf{\ }\rho(1)=-\rho(-1)=-1.$
\newline
} {\textbf{\textup{(iii)} Sub-additivity:} For each $X,Y\in \mathcal{H}$,
\textbf{\ }$\rho(X+Y)\leq \rho(X)+\rho(Y).$ }\newline
\textbf{\textup{(iv)} Positive homogeneity: }$\rho(\lambda X)=\lambda
\rho(X)\ \
\text{for}\ \lambda \geq0.$
\end{proposition}

\begin{proof}
{\ (i), (ii) are obvious. }

We now prove (iii). Indeed,
\begin{align*}
\rho(X+Y)= & \inf \{m\in \mathbb{R}:\;m+(X+Y)\in \mathcal{A }\} \\
= & \inf \{m+n:m,n\in \mathbb{R},\;(m+X)+(n+Y)\in \mathcal{{A} }%
\} \\
\ \leq & \inf \{m\in \mathbb{R}:\;m+X\in \mathcal{{A} }\}+\inf
\{n\in \mathbb{R}:\;n+Y\in {\mathcal{A}} \} \\
= & \rho(X)+\rho(Y).
\end{align*}

To prove (iv), in fact the case ${\ \lambda=0}$ is trivial; when $\lambda>0$%
,
\begin{align*}
\rho(\lambda X) & =\inf \{m\in \mathbb{R}:\;m+\lambda X\in \mathcal{A }\} \
\ \ \ \ \ \  \\
& =\lambda \inf \{n\in \mathbb{R}:\;n+X\in \mathcal{A}\}=\lambda \rho(X),
\end{align*}
where $n=m/{\lambda}$.
\end{proof}

Obviously, if $\mathbb{E}$ is a sublinear expectation, we define $\rho(X):=%
\mathbb{E}[-X]$, then $\rho$ is a coherent risk measure. Inversely, if $\rho$
is a coherent risk measure, we define $\mathbb{E}[X]:=\rho(-X)$, then $%
\mathbb{E}$ is a sublinear expectation.

\begin{exercise}
Let ${\rho (\cdot )}$ be a coherent risk measure. Then we can inversely
define
\begin{equation*}
\mathcal{A}_{\rho }{:=\{X\in \mathcal{H}:\rho (X)\leq 0\}.}
\end{equation*}%
{Prove that }$\mathcal{A}_{\rho }$ is a coherent acceptable set.
\end{exercise}

\section*{Notes and Comments}

\addcontentsline{toc}{section}{Notes and Comments}

The sublinear expectation is also called the upper expectation (see Huber
(1981) \cite{Huber} in robust statistics), or the upper prevision in the
theory of imprecise probabilities (see Walley (1991) \cite{Walley} and a
rich literature provided in the Notes of this book). To our knowledge, the
Representation Theorem \ref{t1} was firstly obtained for the case where $%
\Omega$ is a finite set by \cite{Huber}, and this theorem was rediscovered
independently by Artzner, Delbaen, Eber and Heath (1999) \cite{ADEH2} and
then by Delbaen (2002) \cite{Delbaen} for the general $\Omega$. A typical
example of dynamic nonlinear expectation, called $g$--expectation (small $g$%
), was introduced in Peng (1997) \cite{Peng1997} in the framework of
backward stochastic differential equations. Readers are referred to Briand,
Coquet, Hu, M\'{e}min and Peng (2000) \cite{BCHMP1}, Chen (1998) \cite{Chen98}, Chen and
Epstein (2002) \cite{CE}, Chen, Kulperger and Jiang (2003) \cite{CKJ}, Chen and Peng (1998) \cite{CP} and (2000) \cite{CP1},
Coquet, Hu, M\'{e}min and Peng (2001) \cite{CHMP} (2002) \cite{CHMP3}
, Jiang (2004) \cite{Jiang}, Jiang and Chen (2004) \cite{JC,JC1}, Peng (1999) \cite{Peng1999}
and (2004) \cite{Peng2004c}, Peng and Xu (2003) \cite{PX2003} and Rosazza (2006) \cite{Roazza2003}
for the further development of this theory. It seems that the notions of
distributions and independence under nonlinear expectations were  new. We
think that these notions are perfectly adapted for the further development
of dynamic nonlinear expectations. For other types of the related notions of
distributions and independence under nonlinear expectations or non-additive
probabilities, we refer to the Notes of the book \cite{Walley} and the
references listed in Marinacci (1999) \cite{Marinacci} and Maccheroni and
Marinacci (2005) \cite{Marinacci1}. Coherent risk measures can be also
regarded as sublinear expectations defined on the space of risk positions in
financial market. This notion was firstly introduced in \cite{ADEH2}.
Readers can be referred also to the well-known book of F\"{o}llmer and
Schied (2004) \cite{F-Sch} for the systematical presentation of coherent risk
measures and convex risk measures. For the dynamic risk measure in
continuous time, see \cite{Peng2004c} or \cite{Roazza2003}, Barrieu and El
Karoui (2004) \cite{El-Bar} using $g$-expectations. Super-hedging and super
pricing (see El Karoui and Quenez (1995) \cite{EQ} and El Karoui, Peng and
Quenez (1997)  \cite{EPQ}) are also closely related to this formulation.

%
%
\chapter{Law of Large Numbers and Central Limit Theorem}
\label{ch2} In this chapter,  we first  introduce two types of
fundamentally important distributions, namely, maximal distribution
and $G$-normal distribution, in the theory of sublinear
expectations. The former corresponds to constants and the latter
corresponds to normal distribution in classical probability theory.
We then present the law of large numbers (LLN)  and central limit
theorem (CLT) under sublinear expectations. It is worth pointing out
that the limit in LLN is a maximal distribution and the limit in CLT
is a $G$-normal distribution.

\section{Maximal Distribution and $G$-normal Distribution}

We will firstly define a special type of very simple distributions
which are frequently used in practice, known as ``worst case risk
measure''.

\begin{definition}
(\textbf{maximal distribution}) \label{Prop-G1
copy(1)}\index{Maximal distribution}A $d$-dimensional random vector
$\eta=(\eta_{1},\cdots,\eta_{d})$ on a sublinear expectation space
$(\Omega,\mathcal{H},\mathbb{E})$ is called \textbf{maximal
distributed} if there exists a bounded, closed and convex subset
$\Gamma \subset \mathbb{R}^{d}$ such
that%
\[
\mathbb{E}[\varphi(\eta)]=\max_{y\in \Gamma}\varphi(y).
\]

\end{definition}

\begin{remark}
Here $\Gamma$ gives the degree of uncertainty of $\eta$. It is easy
to check that this maximal distributed random vector $\eta$
satisfies
\[
a\eta+b\bar{\eta}\overset{d}{=}(a+b)\eta\  \  \  \text{for
}a,b\geq0,\
\]
where $\bar{\eta}$ is an independent copy of $\eta$. We will see
later that in fact this relation characterizes a maximal
distribution. Maximal distribution is also called  ``worst case risk
measure'' in finance.
\end{remark}

\begin{remark}
When $d=1$ we have $\Gamma=[\underline{\mu},\overline{\mu}]$, where
$ \overline{\mu}=\mathbb{E}[\eta]\  \text{and} \
\underline{\mu}=-\mathbb{E }[-\eta]. $

The distribution of $\eta$ is
\[
\mathbb{F}_{\eta}[\varphi]=\mathbb{E}[\varphi(\eta)]=\sup
_{\underline{\mu}\leq y\leq \bar{\mu}}\varphi(y)\ \ \  \text{for} \
\varphi \in C_{l.Lip}(\mathbb{R}).
\]

\end{remark}

Recall a well-known characterization: $X\overset{d}{=}N(0,\Sigma)$ if and only if%

\begin{equation}
\label{ch2e1} aX+b\bar{X}\overset{d}{=}\sqrt{a^{2}+b^{2}}X\  \
\text{for }a,b\geq0,
\end{equation}
where $\bar{X}$ is an independent copy of $X$. The covariance matrix
$\Sigma$ is defined by $\Sigma=E[XX^{T}]$. We now consider the so
called $G$-normal distribution in probability model uncertainty
situation. The existence, uniqueness and characterization will be
given later.

\begin{definition}
(\textbf{$G$-normal distribution})%
\index{$G$-normal distribution}
\label{Def-Gnormal copy(2)} A $d$-dimensional random vector $X=(X_{1}%
,\cdots,X_{d})^{T}$ on a sublinear expectation space $(\Omega,\mathcal{H}%
,\mathbb{E})$ is called (centralized) \textbf{$G$-normal
distributed} if
\[
aX+b\bar{X}\overset{d}{=}\sqrt{a^{2}+b^{2}}X\  \  \  \text{for
}a,b\geq0,\
\]
where $\bar{X}$ is an independent copy of $X$.
\end{definition}

\begin{remark}
Noting that $\mathbb{E}[X+\bar{X}]=2\mathbb{E}[X]$ and
$\mathbb{E}[X+\bar{X}]=\mathbb{E}[\sqrt{2}X]=\sqrt{2}\mathbb{E}[X]$,
we then have $\mathbb{E}[X]=0.$ Similarly, we can prove that
$\mathbb{E}[-X]=0.$ Namely, $X$ has no mean-uncertainty.
\end{remark}

The following property is easy to prove by the definition.

\begin{proposition}
\label{GCD}Let $X$ be $G$-normal distributed. Then for each $A\in
\mathbb{R}^{m\times d}$, $AX$ is also $G$-normal distributed. In
particular,
for each $\mathbf{a}\in \mathbb{R}^{d},$ $\left \langle \mathbf{a}%
,X\right \rangle $ is a $1$-dimensional $G$-normal distributed
random variable, but its inverse is not true (see Exercise
\ref{ex1}).
\end{proposition}

We denote by $\mathbb{S}(d)$\label{sd} the collection of all
$d\times d$ symmetric matrices. Let $X$ be $G$-normal distributed
and $\eta$ be maximal distributed $d$-dimensional random  vectors on
$(\Omega,\mathcal{H},\mathbb{E})$. The following function is very
important to characterize their distributions:
\begin{equation}
G(p,A):=\mathbb{E}[\frac{1}{2}\left \langle AX,X\right \rangle
+\left \langle p,\eta \right \rangle ],\  \  \ (p,A)\in \mathbb{R}^{d}%
\times \mathbb{S}(d). \label{e313}%
\end{equation}
It is easy to check that $G$ is a sublinear function monotonic in
$A\in \mathbb{S}(d)$ in the following sense: for each $p,\bar{p}\in
\mathbb{R}^{d}$ and $A,\bar{A}\in \mathbb{S}(d)$%
\begin{equation}
\left \{
\begin{array}
[c]{rl}%
G(p+\bar{p},A+\bar{A})&\leq G(p,A)+G(\bar{p},\bar{A}),\\
G(\lambda p,\lambda A)&=\lambda G(p,A),\  \  \forall \lambda \geq0,\\
G(p,A) & \geq G(p,\bar{A}),\ \ \text{if}  \ A\geq \bar{A}.
\end{array}
\right.  \label{e314}%
\end{equation}
Clearly, $G$ is also a continuous function. By Theorem \ref{t1} in
Chapter \ref{ch1}, there exists a
bounded and closed subset $\Gamma \subset \mathbb{R}^{d}\times \mathbb{R}%
^{d\times d}$ such that
\begin{equation}
\label{ch2e2} G(p,A)=\sup_{(q,Q)\in
\Gamma}[\frac{1}{2}\mathrm{tr}[AQQ^{T}]+\left \langle p,q\right
\rangle ]\  \  \  \text{for}\ (p,A)\in \mathbb{R}^{d}\times
\mathbb{S}(d).
\end{equation}

We have the following result, which will be proved in the next
section.

\begin{proposition}
\label{Prop-Gnorm copy(1)} Let $G:\mathbb{R}^{d}\times \mathbb{S}%
(d)\rightarrow \mathbb{R}$ be a given sublinear and continuous
function, monotonic in $A\in \mathbb{S}(d)$ in the sense of
\textup{(\ref{e314})}. Then there exists a $G$-normal distributed
$d$-dimensional random vector $X$ and a maximal distributed
$d$-dimensional random vector $\eta$ on some sublinear expectation
space $(\Omega,\mathcal{H},\mathbb{E})$ satisfying
\textup{(\ref{e313})} and%
\begin{equation}
(aX+b\bar{X},a^{2}\eta+b^{2}\bar{\eta})\overset{d}{=}(\sqrt{a^{2}+b^{2}%
}X,(a^{2}+b^{2})\eta),\ \ \text{for}\ a,b\geq0, \label{e311}%
\end{equation}
where $(\bar{X},\bar{\eta})$ is an independent copy of $(X,\eta)$.
\end{proposition}

\begin{definition}
The pair $(X,\eta)$ satisfying \textup{(\ref{e311})} is called
\textbf{$G$-distributed}.\index{$G$-distributed}
\end{definition}

\begin{remark}
In fact, if the pair $(X,\eta)$ satisfies \textup{(\ref{e311})},
then
\[
aX+b\bar{X}\overset{d}{=}\sqrt{a^{2}+b^{2}}X,\ a\eta+b\bar{\eta}\overset{d}%
{=}(a+b)\eta\ \ \text{for}\ a,b\geq0.\  \
\]
Thus $X$ is $G$-normal and $\eta$ is maximal distributed.
\end{remark}

The above pair $(X,\eta)$ is characterized by the following
parabolic partial differential equation (PDE for short)
defined on $[0,\infty)\times \mathbb{R}^{d}\times \mathbb{R}^{d}:$%
\begin{equation}
\partial_{t}u-G(D_{y}u,D_{x}^{2}u)=0, \label{ee03}%
\end{equation}
with Cauchy condition$\  \ u|_{t=0}=\varphi$, where $G:\mathbb{R}^{d}%
\times \mathbb{S}(d)\rightarrow \mathbb{R}$ is defined by
(\ref{e313}) and
$D^{2}u=(\partial_{x_{i}x_{j}}^{2}u)_{i,j=1}^{d}$, $Du=(\partial_{x_{i}%
}u)_{i=1}^{d}$. The PDE (\ref{ee03}) is called a
\textbf{$G$-equation}\index{$G$-equation}.

In this book we will mainly use the notion of viscosity solution to
describe the solution of this PDE. For reader's convenience, we give
a systematical introduction of the notion of viscosity solution and
its related properties used in this book (see Appendix C, Section
1-3). It is worth to mention here that for the case where $G$ is
non-degenerate, the viscosity solution of the $G$-equation becomes a
classical $C^{1,2}$ solution (see Appendix C, Section 4). Readers
without knowledge of viscosity solutions can simply understand
solutions of the $G$-equation in the classical sense along the whole
book.

\begin{proposition}
\label{gG-P1 copy(1)}For the pair $(X,\eta)$ satisfying
\textup{(\ref{e311})} and a function $\varphi \in
C_{l.Lip}(\mathbb{R}^{d}\times \mathbb{R}^{d})$, we define
\[
u(t,x,y):=\mathbb{E}[\varphi(x+\sqrt{t}X,y+t\eta)],\ (t,x,y)\in
\lbrack0,\infty)\times \mathbb{R}^{d}\times \mathbb{R}^{d}.
\]
Then we have%
\begin{equation}
u(t+s,x,y)=\mathbb{E}[u(t,x+\sqrt{s}X,y+s\eta)],\  \ s\geq0. \label{e315}%
\end{equation}
We also have the estimates: for each $T>0,$ there exist constants
$C,k>0$ such that, for all $t,s\in \lbrack0,T]$ and
$x,\bar{x},y,\bar{y}\in \mathbb{R}^{d}$,
\begin{equation}
|u(t,x,y)-u(t,\bar{x},\bar{y})|\leq
C(1+|x|^{k}+|y|^{k}+|\bar{x}|^{k}+|\bar
{y}|^{k})(|x-\bar{x}|+|y-\bar{y}|) \label{e04}%
\end{equation}
and%
\begin{equation}
|u(t,x,y)-u(t+s,x,y)|\leq C(1+|x|^{k}+|y|^{k})(s+|s|^{1/2}). \label{e05}%
\end{equation}
Moreover, $u$ is the unique viscosity solution, continuous in the
sense of \textup{(\ref{e04})} and \textup{(\ref{e05})}, of the PDE
\textup{(\ref{ee03})}.
\end{proposition}

\begin{proof}
Since%
\begin{align*}
u(t,x,y)-u(t,\bar{x},\bar{y})  &  =\mathbb{E}[\varphi(x+\sqrt
{t}X,y+t\eta)]-\mathbb{E}[\varphi(\bar{x}+\sqrt{t}X,\bar{y}+t\eta)]\\
&  \leq
\mathbb{E}[\varphi(x+\sqrt{t}X,y+t\eta)-\varphi(\bar{x}+\sqrt
{t}X,\bar{y}+t\eta)]\\
&  \leq \mathbb{E}[C_{1}(1+|X|^{k}+|\eta|^{k}+|x|^{k}+|y|^{k}+|\bar
{x}|^{k}+|\bar{y}|^{k})]\\
&  \  \  \  \  \times(|x-\bar{x}|+|y-\bar{y}|)\\
&  \leq C(1+|x|^{k}+|y|^{k}+|\bar{x}|^{k}+|\bar{y}|^{k})(|x-\bar{x}%
|+|y-\bar{y}|),
\end{align*}
we have (\ref{e04}).

Let $(\bar{X},\bar{\eta})$ be an independent copy of $(X,\eta)$. {By
(\ref{e311}), }%
\begin{align*}
u(t+s,x,y)  &  =\mathbb{E}[\varphi(x+\sqrt{t+s}X,y+(t+s)\eta)]\\
& =\mathbb{E}[\varphi(x+\sqrt{s}X+\sqrt{t}\bar{X},y+s\eta+t\bar{\eta
})]\\
& =\mathbb{E}[\mathbb{E}[\varphi(x+\sqrt{s}\widetilde{x}+\sqrt
{t}\bar{X},y+s\widetilde{y}+t\bar{\eta})]_{(\widetilde{x},\widetilde
{y})=(X,\eta)}]\\
&  =\mathbb{E}[u(t,x+\sqrt{s}X,y+s\eta)],
\end{align*}
we thus obtain (\ref{e315}). From this and (\ref{e04}) it follows
that
\begin{align*}
&  u(t+s,x,y)-u(t,x,y)=\mathbb{E}[u(t,x+\sqrt{s}X,y+s\eta)-u(t,x,y)]\\
&  \leq \mathbb{E}[C_{1}(1+|x|^{k}+|y|^{k}+|X|^{k}+|\eta|^{k})(\sqrt
{s}|X|+s|\eta|)],
\end{align*}
thus we obtain (\ref{e05}).

Now, for a fixed $(t,x,y)\in(0,\infty)\times \mathbb{R}^{d}\times \mathbb{R}%
^{d}$, let $\psi \in C_{b}^{2,3}([0,\infty)\times
\mathbb{R}^{d}\times \mathbb{R}^{d})$ be such that $\psi \geq u$ and
$\psi(t,x,y)=u(t,x,y)$. By (\ref{e315}) and Taylor's expansion, it
follows that, for $\delta \in(0,t)$,
\begin{align*}
0  &  \leq \mathbb{E}[\psi(t-\delta,x+\sqrt{\delta}X,y+\delta \eta
)-\psi(t,x,y)]\\
&  \leq \bar{C}(\delta^{3/2}+\delta^{2})-\partial_{t}\psi(t,x,y)\delta \\
&  \ \ \ +\mathbb{E}[\left \langle D_{x}\psi(t,x,y),X\right \rangle
\sqrt
{\delta}+\left \langle D_{y}\psi(t,x,y),\eta \right \rangle \delta+\frac{1}%
{2}\left \langle D_{x}^{2}\psi(t,x,y)X,X\right \rangle \delta]\\
&  =-\partial_{t}\psi(t,x,y)\delta+\mathbb{E}[\left \langle D_{y}%
\psi(t,x,y),\eta \right \rangle +\frac{1}{2}\left \langle D_{x}^{2}%
\psi(t,x,y)X,X\right \rangle ]\delta+\bar{C}(\delta^{3/2}+\delta^{2})\\
&  =-\partial_{t}\psi(t,x,y)\delta+\delta G(D_{y}\psi,D_{x}^{2}\psi
)(t,x,y)+\bar{C}(\delta^{3/2}+\delta^{2}),
\end{align*}
from which it is easy to check that
\[
\lbrack \partial_{t}\psi-G(D_{y}\psi,D_{x}^{2}\psi)](t,x,y)\leq0.
\]
Thus $u$ is a viscosity subsolution of (\ref{ee03}). Similarly we
can prove that $u$ is a viscosity supersolution of (\ref{ee03}).
\end{proof}

\begin{corollary}
\label{gG-P1coro copy(1)} If both $(X,\eta)$ and
$(\bar{X},\bar{\eta})$
satisfy \textup{(\ref{e311})} {with the same }$G${, i.e.,}%
\[
G(p,A):=\mathbb{E}[\frac{1}{2}\left \langle AX,X\right \rangle
+\left \langle p,\eta \right \rangle ]=\mathbb{E}[\frac{1}{2}\left
\langle A\bar{X},\bar{X}\right \rangle +\left \langle
p,\bar{\eta}\right \rangle ]\ \  \ \text{for}\ (p,A)\in
\mathbb{R}^{d}\times \mathbb{S}(d),
\]
{ then }$(X,\eta)\overset{d}{=}(\bar{X},\bar{\eta})$. In particular,
$X\overset{d}{=}-X$.
\end{corollary}

\begin{proof}
For each $\varphi \in C_{l.Lip}(\mathbb{R}^{d}\times
\mathbb{R}^{d}),$ we set
\begin{align*}
u(t,x,y)  &  :=\mathbb{E}[\varphi(x+\sqrt{t}X,y+t\eta)],\  \\
\bar{u}(t,x,y)  &
:=\mathbb{E}[\varphi(x+\sqrt{t}\bar{X},y+t\bar{\eta })],\ (t,x,y)\in
\lbrack0,\infty)\times \mathbb{R}^{d}\times \mathbb{R}^{d}.
\end{align*}
By Proposition \ref{gG-P1 copy(1)}, both $u$ and $\bar{u}$ are
viscosity solutions of the $G$-equation (\ref{ee03}) with Cauchy
condition $u|_{t=0}=\bar{u}|_{t=0}=\varphi$. It follows from the
uniqueness of the
viscosity solution that $u\equiv \bar{u}$. In particular,%
\[
\mathbb{E}[\varphi(X,\eta)]=\mathbb{E}[\varphi(\bar{X},\bar{\eta
})]\text{.}%
\]
Thus $(X,\eta)\overset{d}{=}(\bar{X},\bar{\eta})$.
\end{proof}

\begin{corollary}
\label{Gv}Let $(X,\eta)$ satisfy \textup{(\ref{e311})}. For each
$\psi \in C_{l.Lip}(\mathbb{R}^{d})$ we define
\begin{equation}
v(t,x):=\mathbb{E}[\psi((x+\sqrt{t}X+t\eta)],\ (t,x)\in \lbrack
0,\infty)\times \mathbb{R}^{d}.\  \label{e320}%
\end{equation}
Then $v$ is the unique viscosity solution of the following parabolic PDE:%
\begin{equation}
\partial_{t}v-G(D_{x}v,D_{x}^{2}v)=0,\  \  \  \ v|_{t=0}=\psi. \label{e318}%
\end{equation}
Moreover, we have $v(t,x+y)\equiv u(t,x,y)$, where $u$ is the
solution of the PDE \textup{(\ref{ee03})} with initial condition
$u(t,x,y)|_{t=0}=\psi(x+y)$.
\end{corollary}

\begin{example}
Let $X$ be $G$-normal distributed. The distribution of $X$ is
characterized by
\[
u(t,x)=\mathbb{E}[\varphi(x+\sqrt{t}X)],\  \  \varphi \in C_{l.Lip}%
(\mathbb{R}^{d}).
\]
In particular, $\mathbb{E}[\varphi(X)]=u(1,0)$, where $u$ is the
solution of the following parabolic PDE defined on $[0,\infty)\times
\mathbb{R}^{d}:$%
\begin{equation}
\label{e03}
 \partial_{t}u-G(D^{2}u)=0,\ \ u|_{t=0}=\varphi,
\end{equation}
where $G=G_{X}(A):\mathbb{S}(d)\rightarrow \mathbb{R}$ is defined by
\[
G(A):=\frac{1}{2}\mathbb{E}[\left \langle AX,X\right \rangle ],\  \
\ A\in \mathbb{S}(d).
\]

The parabolic PDE \textup{(\ref{e03})} is called a \textbf{$G$-heat
equation}.\index{$G$-heat equation}

It is easy to check that $G$ is a sublinear function defined on $\mathbb{S}%
(d)$. By Theorem \ref{t1} in Chapter \ref{ch1}, there exists a
bounded, convex and closed subset $\Theta \subset \mathbb{S}(d)$
such that
\begin{equation}
\label{GaChII}
\frac{1}{2}\mathbb{E}[\left \langle AX,X\right \rangle ]=G(A)=\frac{1}%
{2}\sup_{Q\in \Theta}tr[AQ],\  \ A\in \mathbb{S}(d).
\end{equation}
Since $G(A)$ is monotonic: $G(A_{1})\geq G(A_{2})$, for $A_{1}\geq
A_{2}$, it
follows that%
\[
\Theta \subset \mathbb{S}_{+}(d)=\{ \theta \in \mathbb{S}(d):\theta
\geq 0\}=\{BB^{T}:B\in \mathbb{R}^{d\times d}\},
\]
where $\mathbb{R}^{d\times d}$\label{splusd} is the set of all
$d\times d$ matrices. If $\Theta$ is a singleton: $\Theta=\{Q\}$,
then $X$ is classical zero-mean normal distributed with covariance
$Q$. In general, $\Theta$ characterizes the covariance uncertainty
of $X$. We denote  $X\overset{d}{=}N(\{0\} \times \Theta)$ (Recall
equation \textup{(\ref{ch2e2})}, we can set $(q,Q)\in
\{0\}\times\Theta$).

When $d=1$, we have $X\overset{d}{=}N(\{0\} \times \lbrack
\underline{\sigma }^{2},\bar{\sigma}^{2}])$ (We also denoted by
$X\overset{d}{=}N(0,[\underline{\sigma}^2,\bar{\sigma}^2])$), where
$\bar{\sigma}^{2}=\mathbb{E}[X^{2}]$ and
$\underline{\sigma}^{2}=-\mathbb{E}[-X^{2}]$. The corresponding
$G$-heat equation is
\[
\partial_{t}u-\frac{1}{2}(\bar{\sigma}^{2}(\partial_{xx}^{2}u)^{+}%
-\underline{\sigma}^{2}(\partial_{xx}^{2}u)^{-})=0,~u|_{t=0}=\varphi.
\]
For the case $\underline{\sigma}^{2}>0$, this equation is also
called the Barenblatt equation.
\end{example}
In the following two typical situations, the calculation of
$\mathbb{E} [\varphi(X)]$ is very easy:

\begin{itemize}
\item For each \textbf{convex} function $\varphi$, we have
\[
\mathbb{E}[\varphi(X)]=\frac{1}{\sqrt{2\pi}}%
\int_{-\infty}^{\infty}\varphi(\overline{\sigma}^{2}y)\exp(-\frac{y^{2}}{2
})dy.
\]
Indeed, for each fixed $t\geq0$, it is easy to check that the
function {$u(t,x):=\mathbb{E}[\varphi(x+\sqrt{t}X)]$ is convex in
$x$:}
\begin{align*}
u(t,\alpha x+(1-\alpha)y)  &  =\mathbb{E}[{\varphi(\alpha
x+(1-\alpha
)y+\sqrt{t}X)]}\\
&  \leq{}{{\alpha}\mathbb{E}[\varphi(x+\sqrt{t}X)]+(1-{\alpha
)}\mathbb{E}[\varphi(x+\sqrt{t}X)]}\\
&  =\alpha u(t,x)+(1-\alpha)u(t,x).
\end{align*}
It follows that $(\partial_{xx}^{2}u)^{-}\equiv0$ and thus the above
$G$-heat
equation becomes {%
\[
\partial_{t}u=\frac{\overline{\sigma}^{2}}{2}\partial_{xx}^{2}u,\  \  \ u|_{t=0}%
=\varphi.\  \
\]
}

\item For each \textbf{concave} function $\varphi$, we have%
\[
\mathbb{E}[\varphi(X)]=\frac{1}{\sqrt{2\pi}}%
\int_{-\infty}^{\infty}\varphi(\underline{\sigma}^{2}y)\exp(-\frac{y^{2}}{2%
})dy.
\]

In particular,%
\[
\mathbb{E}[X]=\mathbb{E}[-X]=0,\  \  \mathbb{E}[X^{2}%
]=\overline{\sigma}^{2},\  \ -\mathbb{E}[-X^{2}]=\underline{\sigma}^{2}%
\]
and
\[
\mathbb{E}[X^{4}]=3\overline{\sigma}^{4},\ -\mathbb{E}%
[-X^{4}]=3\underline{\sigma}^{4}\ .\
\]

\end{itemize}

\begin{example}
\label{example116} Let $\eta$ be maximal distributed, the
distribution of $\eta$ is characterized
by the following parabolic PDE defined on $[0,\infty)\times \mathbb{R}^{d}:$%
\begin{align}
  \partial_{t}u-g(Du)=0,\label{e01}\ \ u|_{t=0}=\varphi,
\end{align}
where $g=g_{\eta}(p):\mathbb{R}^{d}\rightarrow \mathbb{R}$ is
defined by
\[
g_{\eta}(p):=\mathbb{E}[\left \langle p,\eta \right \rangle ],\  \ \
p\in \mathbb{R}^{d}.
\]

It is easy to check that $g_{\eta}$ is a sublinear function defined
on $\mathbb{R}^{d}$. By Theorem \ref{t1} in Chapter \ref{ch1}, there
exists a bounded, convex and closed subset
$\bar{\Theta}\subset R^{d}$ such that%
\begin{equation}
g(p)=\sup_{q\in \bar{\Theta}}\left \langle p,q\right \rangle ,\  \ \
\ p\in
\mathbb{R}^{d}. \label{e35}%
\end{equation}
By this characterization, we can prove that the distribution of
$\eta$ is given by
\begin{equation}
\mathbb{\hat{F}}_{\eta}[\varphi]=\mathbb{E}[\varphi(\eta)]=\sup
_{v\in \bar{\Theta}}\varphi(v)=\sup_{v\in \bar{\Theta}}\int_{\mathbb{R}^{d}%
}\varphi(x)\delta_{v}(dx),\  \  \varphi \in
C_{l.Lip}(\mathbb{R}^{d}),
\label{e003}%
\end{equation}
where $\delta_{v}$ is Dirac measure. Namely it is the maximal
distribution with the uncertainty subset of probabilities as Dirac
measures concentrated at $\bar{\Theta}$. We denote $\eta
\overset{d}{=}N(\bar{\Theta}\times \{0\})$ (Recall equation
\textup{(\ref{ch2e2})}, we can set $(q,Q)\in
\bar{\Theta}\times\{0\}$).

In particular, for $d=1,$
\[
g_{\eta}(p):=\mathbb{E}[p\eta]=\bar{\mu}p^{+}-\underline{\mu}%
p^{-},\  \  \ p\in \mathbb{R},
\]
where $\bar{\mu}=\mathbb{E}[\eta]$ and
$\underline{\mu}=-\mathbb{\hat {E}}[-\eta]$. The distribution of
$\eta$ is given by (\ref{e003}). We denote  $\eta
\overset{d}{=}N([\underline{\mu},\bar{\mu}]\times \{0\})$.
\end{example}

\begin{exercise}
\label{ex1} We consider $X=(X_{1},X_{2})$, where
$X_{1}\overset{d}{=}N(\{0\}\times[\underline
{\sigma}^{2},\overline{\sigma}^{2}])$ with
$\overline{\sigma}>\underline{\sigma}$, $X_{2}$ is an independent
copy of $X_{1}$. Show that

\textup{(1)} For each $a\in \mathbb{R}^{2}$, $\langle a,X\rangle$ is
a $1$-dimensional $G$-normal distributed random variable.

\textup{(2)} $X$ is not $G$-normal distributed.
\end{exercise}

\begin{exercise}
Let $X$ be $G$-normal distributed. For each $\varphi \in C_{l.Lip}%
(\mathbb{R}^{d}),$ we define a function%
\[
u(t,x):=\mathbb{E}[\varphi(x+\sqrt{t}X)],\ (t,x)\in \lbrack
0,\infty)\times \mathbb{R}^{d}.
\]
Show that $u$ is the unique viscosity solution of the PDE
(\ref{e03}) with Cauchy condition$\  \ u|_{t=0}=\varphi.$
\end{exercise}

\begin{exercise}
\label{ex2} Let $\eta$ be maximal distributed. For each $\varphi \in
C_{l.Lip}(\mathbb{R}^{d}),$ we define a function%
\[
u(t,y):=\mathbb{E}[\varphi(y+t\eta)],\ (t,y)\in \lbrack0,\infty
)\times \mathbb{R}^{d}.
\]
Show that $u$ is the unique viscosity solution of the PDE
(\ref{e01}) with Cauchy condition$\  \ u|_{t=0}=\varphi.$
\end{exercise}

\section{Existence of $G$-distributed Random Variables}

\label{c2s2} In this section, we give the proof of the existence of
G-distributed random variables, namely, the proof of Proposition
\ref{Prop-Gnorm copy(1)}.

Let $G:\mathbb{R}^{d}\times \mathbb{S}(d)\rightarrow \mathbb{R}$ be
a given sublinear function monotonic in $A\in \mathbb{S}(d)$ in the
sense of (\ref{e314}). We now construct a pair of $d$-dimensional
random vectors
$(X,\eta)$ on some sublinear expectation space $(\Omega,\mathcal{H}%
,\mathbb{E})$ satisfying (\ref{e313}) and (\ref{e311}).

For each $\varphi \in C_{l.Lip}(\mathbb{R}^{2d}),$ let
$u=u^{\varphi}$ be the unique viscosity solution of the $G$-equation
(\ref{ee03}) with $u^{\varphi}|_{t=0}=\varphi$. We take
$\widetilde{\Omega}=\mathbb{R}^{2d}$,
$\mathcal{\widetilde{H}}=C_{l.Lip}(\mathbb{R}^{2d})$ and
$\widetilde{\omega }=(x,y)\in \mathbb{R}^{2d}$. The corresponding
sublinear expectation
$\mathbb{\widetilde{E}}[\cdot]$ is defined by $\mathbb{\widetilde{E}}%
[\xi]=u^{\varphi}(1,0,0)$, for each $\xi \in
\mathcal{\widetilde{H}}$ of the form
$\xi(\widetilde{\omega})=(\varphi(x,y))_{(x,y)\in
\mathbb{R}^{2d}}\in C_{l.Lip}(\mathbb{R}^{2d})$. The monotonicity
and sub-additivity of $u^{\varphi}$ with respect to $\varphi$ are
known in the theory of viscosity solution. For reader's convenience
we provide a new and simple proof in Appendix C (see Corollary
\ref{Comparison} and Corollary \ref{Domination}). The
constant preserving and positive homogeneity of $\mathbb{\widetilde{E}}%
[\cdot]$ are easy to check. Thus the functional $\mathbb{\widetilde{E}%
}[\cdot]:\mathcal{\widetilde{H}}\rightarrow \mathbb{R}$ forms a
sublinear expectation.

We now consider a pair of $d$-dimensional random vectors
$(\widetilde {X},\widetilde{\eta})(\widetilde{\omega})=(x,y)$. We
have
\[
\mathbb{\widetilde{E}}[\varphi(\widetilde{X},\widetilde{\eta})]=u^{\varphi
}(1,0,0)\  \  \  \text{for}\ \varphi \in C_{l.Lip}(\mathbb{R}^{2d}).
\]
In particular, just setting $\varphi_{0}(x,y)=\frac{1}{2}\left
\langle Ax,x\right \rangle +\left \langle p,y\right \rangle $, we
can check that
\[
u^{\varphi_{0}}(t,x,y)=G(p,A)t+\frac{1}{2}\left \langle Ax,x\right
\rangle +\left \langle p,y\right \rangle .
\]
We thus have%
\[
\mathbb{\widetilde{E}}[\frac{1}{2}\left \langle
A\widetilde{X},\widetilde {X}\right \rangle +\left \langle
p,\widetilde{\eta}\right \rangle ]=u^{\varphi _{0}}(1,0,0)=G(p,A),\
\ (p,A)\in \mathbb{R}^{d}\times \mathbb{S}(d).
\]

We construct a product space%
\[
(\Omega,\mathcal{H},\mathbb{E})=(\widetilde{\Omega}\times \widetilde
{\Omega},\mathcal{\widetilde{H}\otimes
\widetilde{H}},\mathbb{\widetilde {E}\otimes \widetilde{E}}),
\]
and introduce two pairs of random vectors%
\[
(X,\eta)(\widetilde{\omega}_{1},\widetilde{\omega}_{2})=\widetilde{\omega}%
_{1},\ (\bar{X},\bar{\eta})(\widetilde{\omega}_{1},\widetilde{\omega}%
_{2})=\widetilde{\omega}_{2},\  \  \
(\widetilde{\omega}_{1},\widetilde{\omega }_{2})\in
\widetilde{\Omega}\times \widetilde{\Omega}.\
\]
By Proposition \ref{pp1} in Chapter \ref{ch1}, $(X,\eta)\overset{d}{=}%
(\widetilde{X},\widetilde{\eta})$ and $(\bar{X},\bar{\eta})$ is an
independent copy of $(X,\eta)$.

We now prove that the distribution of $(X,\eta)$ satisfies condition
(\ref{e311}). For each $\varphi \in C_{l.Lip}(\mathbb{R}^{2d})$ and
for each fixed $\lambda>0$, $(\bar{x},\bar{y})\in \mathbb{R}^{2d}$,
since the function
$v$ defined by $v(t,x,y):=u^{\varphi}(\lambda t,\bar{x}+\sqrt{\lambda}%
x,\bar{y}+\lambda y)$ solves exactly the same equation (\ref{ee03}),
but with Cauchy condition
\[
v|_{t=0}=\varphi(\bar{x}+\sqrt{\lambda}\times \cdot,\bar{y}+\lambda
\times \cdot).
\]
Thus
\[
\mathbb{E}[\varphi(\bar{x}+\sqrt{\lambda}X,\bar{y}+\lambda
\eta)]=v(1,0,0)=u^{\varphi}(\lambda,\bar{x},\bar{y}).\
\]

\bigskip

\bigskip By the definition of $\mathbb{E}$, for each $t>0$ and $s>0$,
\begin{align*}
\mathbb{E}[\varphi(\sqrt{t}X+\sqrt{s}\bar{X},t\eta+s\bar{\eta})] &
=\mathbb{E}[\mathbb{E}[\varphi(\sqrt{t}x+\sqrt{s}\bar{X}%
,ty+s\bar{\eta})]_{(x,y)=(X,\eta)}]\\
&  =\mathbb{E}[u^{\varphi}(s,\sqrt{t}X,t\eta)]=u^{u^{\varphi}%
(s,\cdot,\cdot)}(t,0,0)\\
&  =u^{\varphi}(t+s,0,0)\\
&  =\mathbb{E}[\varphi(\sqrt{t+s}X,(t+s)\eta)].
\end{align*}
Namely
$(\sqrt{t}X+\sqrt{s}\bar{X},t\eta+s\bar{\eta})\overset{d}{=}(\sqrt
{t+s}X,(t+s)\eta)$. Thus the distribution of $(X,\eta)$ satisfies
condition (\ref{e311}).
\begin{remark}
From now on, when we mention the sublinear expectation space
$(\Omega,\mathcal{H},\mathbb{E})$, we suppose that there exists a
pair of random vectors $(X,\eta)$ on
$(\Omega,\mathcal{H},\mathbb{E})$ such that $(X,\eta)$ is
G-distributed.
\end{remark}

\begin{exercise}
\label{exxee1} Prove that ${\mathbb{E}}[X^3]>0 $ for
$X\overset{d}{=}N(\{0\} \times \lbrack \underline{\sigma
}^{2},\bar{\sigma}^{2}])$ with $\underline{\sigma
}^{2}<\bar{\sigma}^{2}$.

It is worth to point that ${\mathbb{E}}[\varphi(X)]$ not always
equal to $\sup_{\underline{\sigma
}^{2}\leq\sigma\leq\bar{\sigma}^{2}} E_{\sigma}[\varphi(X)]$ for
$\varphi\in C_{l,Lip}(\mathbb{R})$, where $E_{\sigma}$ denotes the
linear expectation corresponding to the normal distributed density
function $N(0,\sigma^2)$.
\end{exercise}

\section{Law of Large Numbers and Central Limit Theorem}

\begin{theorem}
(\textbf{Law of large numbers})\index{Law of large numbers} {Let
$\left \{ Y_{i}\right \} _{i=1}^{\infty}$ be a sequence of
}$\mathbb{R}^{d}$-valued random variables on a sublinear expectation
space{ $($}$\Omega,${$\mathcal{H},\mathbb{E})$}. We assume that
$Y_{i+1}\overset{d}{=}Y_{i}${ and $Y_{i+1}$ is independent from
$\{Y_{1},\cdots,Y_{i}\}$ for each $i=1,2,\cdots$. Then the sequence
}$\{
\bar{S}_{n}\}_{n=1}^{\infty}$ defined by{ }%
\[
\bar{S}_{n}:=\frac{1}{n}\sum_{i=1}^{n}Y_{i}%
\]
{converges in law to a maximal distribution, i.e.,$\  \ $}%
\begin{equation}
\lim_{n\rightarrow \infty}\mathbb{E}[\varphi(\bar{S}_{n})]=\mathbb{
{E}}[\varphi(\eta)],\  \  \label{e325}%
\end{equation}
for all functions $\varphi \in C(\mathbb{R}^{d})$ satisfying linear
growth condition \textup{(}$|\varphi(x)|\leq C(1+|x|)$\textup{)},
{where }$\eta$ is {a maximal} distributed random vector and{ the
corresponding sublinear function
}$g:\mathbb{R}^{d}\rightarrow \mathbb{R}$ is defined by{ }%
\[
g(p):=\mathbb{E}[\left \langle p,Y_{1}\right \rangle ]{,}\  \ p\in
\mathbb{R}^{d}.
\]

\end{theorem}

\begin{remark}
When $d=1$, {the sequence }$\{ \bar{S}_{n}\}_{n=1}^{\infty}$
{converges in law
to }$N([\underline{\mu},\bar{\mu}]\times \{0\}),$ where $\bar{\mu}%
=\mathbb{E}[Y_{1}]$ and $\underline{\mu}=-\mathbb{E}[-Y_{1}]$. For
the general case, the sum $\frac{1}{n}\sum_{i=1}^{n}Y_{i}$
{converges in law to }$N(\bar{\Theta}\times \{0\})$, where
$\bar{\Theta}\subset \mathbb{R}^{d}$ is the bounded, convex and
closed subset defined in Example \ref{example116}. If we take in
particular $\varphi(y)=d_{\bar{\Theta}}(y)=\inf \{|x-y|:x\in
\bar{\Theta}\}$,
then by \textup{(\ref{e325})} we have the following generalized law of large numbers:%
\begin{equation}
\lim_{n\rightarrow \infty}\mathbb{E}[d_{\bar{\Theta}}(\frac{1}{n}%
\sum_{i=1}^{n}Y_{i})]=\sup_{\theta \in
\bar{\Theta}}d_{\bar{\Theta}}(\theta)=0.
\label{e319}%
\end{equation}
If $Y_{i}$ has no mean-uncertainty, or in other words,
$\bar{\Theta}$ is a singleton: $\bar{\Theta}=\{ \bar{\theta}\}$,
then \textup{(\ref{e319})} becomes
\[
\lim_{n\rightarrow \infty}\mathbb{E}[|\frac{1}{n}\sum_{i=1}^{n}Y_{i}%
-\bar{\theta}|]=0.
\]

\end{remark}

\begin{theorem}
{\textbf{(Central limit theorem with zero-mean)}\index{Central limit
theorem with zero-mean} Let $\left \{ X_{i}\right \}
_{i=1}^{\infty}$ be a sequence of }$\mathbb{R}^{d}$-valued random
variables on a sublinear expectation space{
$($}$\Omega,${$\mathcal{H},\mathbb{E})$}. We assume that
$X_{i+1}\overset{d}{=}X_{i}${ and $X_{i+1}$ is independent from
$\{X_{1},\cdots,X_{i}\}$ for each $i=1,2,\cdots$. We further assume
that
\[
\mathbb{E}[X_{1}]=\mathbb{E}[-X_{1}]=0.\  \
\]
Then the sequence }$\{ \bar{S}_{n}\}_{n=1}^{\infty}$ defined by{ }%
\[
\bar{S}_{n}:=\frac{1}{\sqrt{n}}\sum_{i=1}^{n}X_{i}%
\]
{converges in law to }$X${, i.e., }%
\[
\lim_{n\rightarrow \infty}\mathbb{E}[\varphi(\bar{S}_{n})]=\mathbb{
{E}}[\varphi(X)],\
\]
for all functions $\varphi \in C(\mathbb{R}^{d})$ satisfying linear
growth condition, {where }$X$ is {a }$G$-normal distributed random
vector and{ the corresponding sublinear function
}$G:\mathbb{S}(d)\rightarrow \mathbb{R}$
is defined by{ }%
\[
G(A):=\mathbb{E}[\frac{1}{2}\left \langle AX_{1},X_{1}\right \rangle
]{,\  \ }A\in \mathbb{S}(d).
\]

\end{theorem}

\begin{remark}
When $d=1$, {the sequence }$\{ \bar{S}_{n}\}_{n=1}^{\infty}$
{converges in law to }$N(\{0\}
\times[\underline{\sigma}^{2},\overline{\sigma}^{2}])$, where
$\overline{\sigma}^{2}=\mathbb{E}[X_{1}^{2}]$ and $\underline{\sigma
}^{2}=-\mathbb{E}[-X_{1}^{2}]$. In particular, if $\overline{\sigma}%
^{2}=\underline{\sigma}^{2}$, then it becomes a classical central
limit theorem.
\end{remark}

The following theorem is a nontrivial generalization of the above
two theorems.

\begin{theorem}
\label{CLT}{\textbf{(Central limit theorem with law of large
numbers)}\index{Central limit theorem with law of large numbers} Let
$\left \{
(X_{i},Y_{i})\right \}  _{i=1}^{\infty}$ be a sequence of }$\mathbb{R}%
^{d}\times \mathbb{R}^{d}$-valued random vectors on a sublinear
expectation space{ $($}$\Omega,${$\mathcal{H},\mathbb{E})$}. We
assume that $(X_{i+1},Y_{i+1})\overset{d}{=}(X_{i},Y_{i})${ and
$(X_{i+1},Y_{i+1})$ is independent from
$\{(X_{1},Y_{1}),\cdots,(X_{i},Y_{i})\}$ for each $i=1,2,\cdots$. We
further assume that
\[
\mathbb{E}[X_{1}]=\mathbb{E}[-X_{1}]=0.\
\]
Then the sequence }$\{ \bar{S}_{n}\}_{n=1}^{\infty}$ defined by{ }%
\[
\bar{S}_{n}:=\sum_{i=1}^{n}(\frac{X_{i}}{\sqrt{n}}+\frac{Y_{i}}{n})
\]
{converges in law to }$X+\eta$, i.e.,{
\begin{equation}
\lim_{n\rightarrow \infty}\mathbb{E}[\varphi(\bar{S}_{n})]=\mathbb{
{E}}[\varphi(X+\eta)],\  \  \  \  \label{e12}%
\end{equation}
}for all functions $\varphi \in C(\mathbb{R}^{d})$ satisfying a
linear growth condition, {where} the pair $(X,\eta)$ is
$G$-distributed.{ The corresponding sublinear function
}$G:\mathbb{R}^{d}\times \mathbb{S}(d)\rightarrow \mathbb{R}$ is defined by{ }%
\[
G(p,A):=\mathbb{E}[\left \langle p,Y_{1}\right \rangle +\frac{1}%
{2}\left \langle AX_{1},X_{1}\right \rangle ]{,\  \ }A\in \mathbb{S}%
(d),\  \ p\in \mathbb{R}^{d}.
\]
Thus $\mathbb{E}[\varphi(X+\eta)]$ can be calculated by Corollary
\ref{Gv}.
\end{theorem}

The following result is equivalent to the above central limit
theorem.

\begin{theorem}
\label{CLT1}We make the same assumptions as in Theorem \ref{CLT}.
Then for each function $\varphi \in C(\mathbb{R}^{d}\times
\mathbb{R}^{d})$ satisfying
linear growth condition, we have%
\[
\lim_{n\rightarrow \infty}\mathbb{E}[\varphi(\sum_{i=1}^{n}\frac{X_{i}%
}{\sqrt{n}},\sum_{i=1}^{n}\frac{Y_{i}}{n})]=\mathbb{E}[\varphi
(X,\eta)].
\]

\end{theorem}

\begin{proof}
It is easy to prove Theorem \ref{CLT} by Theorem \ref{CLT1}. To
prove Theorem \ref{CLT1} from Theorem \ref{CLT}, it suffices to
define a pair of $2d$-dimensional random vectors
\[
\bar{X}_{i}=(X_{i},0),\  \  \bar{Y}_{i}=(0,Y_{i})\  \  \text{for
}i=1,2,\cdots.
\]
We have
\begin{align*}
\lim_{n\rightarrow \infty}\mathbb{E}[\varphi(\sum_{i=1}^{n}\frac{X_{i}%
}{\sqrt{n}},\sum_{i=1}^{n}\frac{Y_{i}}{n})] &  =\lim_{n\rightarrow
\infty
}\mathbb{E}[\varphi(\sum_{i=1}^{n}(\frac{\bar{X}_{i}}{\sqrt{n}}%
+\frac{\bar{Y}_{i}}{n}))]=\mathbb{E}[\varphi(\bar{X}%
+\eta)]\\
&  =\mathbb{E}[\varphi(X,\eta)]
\end{align*}
with $\bar{X}=(X,0)$ and $\bar{\eta}=(0,\eta)$.
\end{proof}

To prove Theorem \ref{CLT}, we need the following norms to measure
the regularity of a given real functions $u$ defined on
$Q=[0,T]\times \mathbb{R}^{d}$:
\begin{align*}
\left \Vert u\right \Vert _{C^{0,0}(Q)} &  =\sup_{(t,x)\in Q}|u(t,x)|,\  \\
\left \Vert u\right \Vert _{C^{1,1}(Q)} &  =\left \Vert u\right
\Vert
_{C^{0,0}(Q)}+\left \Vert \partial_{t}u\right \Vert _{C^{0,0}(Q)}+\sum_{i=1}%
^{d}\left \Vert \partial_{x_{i}}u\right \Vert _{C^{0,0}(Q)},\\
\left \Vert u\right \Vert _{C^{1,2}(Q)} &  =\left \Vert u\right
\Vert _{C^{1,1}(Q)}+\sum_{i,j=1}^{d}\left \Vert
\partial_{x_{i}x_{j}}u\right \Vert _{C^{0,0}(Q)}.
\end{align*}
For given constants $\alpha,\beta \in(0,1)$, we denote%

\begin{align*}
\left \Vert u\right \Vert _{C^{\alpha,\beta}(Q)} &
=\sup_{\substack{x,y\in
\mathbb{R}^{d},\ x\not =y\\s,t\in \lbrack0,T],s\not =t}}\frac{|u(s,x)-u(t,y)|}%
{|r-s|^{\alpha}+|x-y|^{\beta}},\\
\left \Vert u\right \Vert _{C^{1+\alpha,1+\beta}(Q)} &  =\left \Vert
u\right \Vert _{C^{\alpha,\beta}(Q)}+\left \Vert \partial_{t}u\right
\Vert _{C^{\alpha,\beta }(Q)}+\sum_{i=1}^{d}\left \Vert
\partial_{x_{i}}u\right \Vert _{C^{\alpha,\beta
}(Q)},\\
\left \Vert u\right \Vert _{C^{1+\alpha,2+\beta}(Q)} &  =\left \Vert
u\right \Vert
_{C^{1+\alpha,1+\beta}(Q)}+\sum_{i,j=1}^{d}\left \Vert \partial_{x_{i}x_{j}%
}u\right \Vert _{C^{\alpha,\beta}(Q)}.
\end{align*}
If, for example, $\left \Vert u\right \Vert
_{C^{1+\alpha,2+\beta}(Q)}<\infty$, then $u$ is said to be a
$C^{1+\alpha,2+\beta}$-function on $Q$.

We need the following lemma.

\begin{lemma}
\label{Lem-CLT}We \ assume the same assumptions as in Theorem
\ref{CLT}. We further assume that there exists a constant $\beta>0$
such that, for each $A$,
$\bar{A}\in \mathbb{S}(d)$ with $A\geq \bar{A}$, we have {\ }%
\begin{equation}
\mathbb{E}[\left \langle AX_{1},X_{1}\right \rangle ]-\mathbb{{E}%
}[\left \langle \bar{A}X_{1},X_{1}\right \rangle ]\geq \beta
\text{tr}[A-\bar{A}].
\label{Ellip}%
\end{equation}
Then our main result \textup{(\ref{e12})} holds.
\end{lemma}

\begin{proof}
We first prove (\ref{e12}) for $\varphi \in
C_{b.Lip}(\mathbb{R}^{d})$. {For a
small but fixed $h>0$, let $V$ be the unique viscosity solution of}%
\begin{equation}
\partial_{t}V+G(DV,D^{2}V)=0,\ (t,x)\in \lbrack0,1+h)\times \mathbb{R}%
^{d}\text{,}\  \ V|_{t=1+h}=\varphi.\label{e14}%
\end{equation}
{Since }$(X,\eta)$ satisfies (\ref{e311}),{ we have }%
\begin{equation}
V(h,0)=\mathbb{E}[\varphi(X+\eta)],\  \ V(1+h,x)=\varphi
(x).\label{equ-h}%
\end{equation}
Since (\ref{e14}) is a uniformly parabolic PDE and $G$ is a convex
function, by the interior regularity of $V$ (see Appendix C), we
have
\[
\left \Vert V\right \Vert _{C^{1+\alpha/2,2+\alpha}([0,1]\times \mathbb{R}^{d}%
)}<\infty\  \text{for some }\alpha \in(0,1).
\]
We set $\delta=\frac{1}{n}$ and $S_{0}=0$. Then \
\begin{align*}
&  V(1,\bar{S}_{n})-V(0,0)=\sum_{i=0}^{n-1}\{V((i+1)\delta,\bar{S}%
_{i+1})-V(i\delta,\bar{S}_{i})\} \\
&  =\sum_{i=0}^{n-1}{\large
\{}[V((i+1)\delta,\bar{S}_{i+1})-V(i\delta,\bar
{S}_{i+1})]+[V(i\delta,\bar{S}_{i+1})-V(i\delta,\bar{S}_{i})]{\large \}}\\
&  =\sum_{i=0}^{n-1}\left \{  I_{\delta}^{i}+J_{\delta}^{i}\right \}
\end{align*}
with, by Taylor's expansion,%
\[
J_{\delta}^{i}=\partial_{t}V(i\delta,\bar{S}_{i})\delta+\frac{1}%
{2}\left \langle D^{2}V(i\delta,\bar{S}_{i})X_{i+1},X_{i+1}\right
\rangle
\delta+\left \langle DV(i\delta,\bar{S}_{i}),X_{i+1}\sqrt{\delta}+Y_{i+1}%
\delta \right \rangle
\]%
\begin{align*}
I_{\delta}^{i}&=\delta\int_{0}^{1}[\partial_{t}V((i+\beta)\delta,\bar{S}%
_{i+1})-\partial_{t}V(i\delta,\bar{S}_{i+1})]d\beta +[\partial
_{t}V(i\delta,\bar{S}_{i+1})-\partial_{t}V(i\delta,\bar{S}_{i})]\delta \\
& \ \ \ +\left \langle
D^{2}V(i\delta,\bar{S}_{i})X_{i+1},Y_{i+1}\right \rangle
\delta^{3/2}+\frac{1}{2}\left \langle D^{2}V(i\delta,\bar{S}_{i})Y_{i+1}%
,Y_{i+1}\right \rangle \delta^{2}\\
& \ \ \ +\int_{0}^{1}\int_{0}^{1}\left \langle \Theta_{\beta \gamma}^{i}(X_{i+1}%
\sqrt{\delta}+Y_{i+1}\delta),X_{i+1}\sqrt{\delta}+Y_{i+1}\delta
\right \rangle \gamma d\beta d\gamma
\end{align*}
with%
\[
\Theta_{\beta \gamma}^{i}=D^{2}V(i\delta,\bar{S}_{i}+\gamma \beta(X_{i+1}%
\sqrt{\delta}+Y_{i+1}\delta))-D^{2}V(i\delta,\bar{S}_{i}).
\]
{Thus
\begin{equation}
\mathbb{E}[\sum_{i=0}^{n-1}J_{\delta}^{i}]-\mathbb{E}[-\sum
_{i=0}^{n-1}I_{\delta}^{i}]\leq
\mathbb{E}[V(1,\bar{S}_{n})]-V(0,0)\leq
\mathbb{E}[\sum_{i=0}^{n-1}J_{\delta}^{i}]+\mathbb{E}[\sum
_{i=0}^{n-1}I_{\delta}^{i}].\label{c2ee15}%
\end{equation}
We now prove that $\mathbb{E}[\sum_{i=0}^{n-1}J_{\delta}^{i}]=0$.
For $J_{\delta}^{i}$, note that
\[
\mathbb{E}[\left \langle DV(i\delta,\bar{S}_{i}),X_{i+1}\sqrt{\delta
}\right \rangle ]=\mathbb{E}[-\left \langle DV(i\delta,\bar{S}%
_{i}),X_{i+1}\sqrt{\delta}\right \rangle ]=0,
\]
then, from the definition of the function }$G$,{ we have
\[
\mathbb{E}[J_{\delta}^{i}]=\mathbb{E}[\partial_{t}V(i\delta
,\bar{S}_{i})+G(DV(i\delta,\bar{S}_{i}),D^{2}V(i\delta,\bar{S}_{i}))]\delta.
\]
Combining the above two equalities with
$\partial_{t}V+G(DV,D^{2}V)=0$
as well as the independence of $(X_{i+1},Y_{i+1})$ from $\{(X_{1}%
,Y_{1}),\cdots,(X_{i},Y_{i})\}$, it follows that
\[
\mathbb{E}[\sum_{i=0}^{n-1}J_{\delta}^{i}]=\mathbb{E}[\sum
_{i=0}^{n-2}J_{\delta}^{i}]=\cdots=0.
\]
Thus (\ref{c2ee15}) can be rewritten as%
\[
-\mathbb{E}[-\sum_{i=0}^{n-1}I_{\delta}^{i}]\leq \mathbb{{E}%
}[V(1,\bar{S}_{n})]-V(0,0)\leq \mathbb{E}[\sum_{i=0}^{n-1}I_{\delta}%
^{i}].
\]
But since both $\partial_{t}V$ and $D^{2}V$ are uniformly
$\frac{\alpha}{2}$-h\"{o}lder continuous in $t$ and
$\alpha$-h\"{o}lder continuous in $x$
 on
$[0,1]\times$}$\mathbb{R}^{d}${, we then have }%
\[
{|I_{\delta}^{i}|\leq C\delta^{1+\alpha/2}(1+|X_{i+1}|^{2+\alpha}%
+|Y_{i+1}|^{2+\alpha}).}%
\]
{ It follows that
\[
\mathbb{E}[|I_{\delta}^{i}|]\leq C\delta^{1+\alpha/2}(1+\mathbb{{E}%
}[|X_{1}|^{2+\alpha}+|Y_{1}|^{2+\alpha}]).
\]
Thus%
\begin{align*}
-C(\frac{1}{n})^{\alpha/2}(1+\mathbb{E}[|X_{1}|^{2+\alpha}%
+|Y_{1}|^{2+\alpha}]) &  \leq \mathbb{E}[V(1,\bar{S}_{n})]-V(0,0)\\
&  \leq C(\frac{1}{n})^{\alpha/2}(1+\mathbb{E}[|X_{1}|^{2+\alpha}%
+|Y_{1}|^{2+\alpha}]).
\end{align*}
As $n\rightarrow \infty,$ we have
\begin{equation}
\lim_{n\rightarrow \infty}\mathbb{E}[V(1,\bar{S}_{n}%
)]=V(0,0).\label{equ-0}%
\end{equation}
On the other hand, for each $t,t^{\prime}\in \lbrack0,1+h]$ and
$x\in \mathbb{R}^{d}$, we have
\[
|V(t,x)-V(t^{\prime},x)|\leq
C(\sqrt{|t-t^{\prime}|}+|t-t^{\prime}|).
\]
Thus $|V(0,0)-V(h,0)|\leq C$}$(\sqrt{h}+h)$ and, by (\ref{equ-0}),
\[
|\mathbb{E}[V(1,\bar{S}_{n})]-\mathbb{E}[\varphi(\bar{S}%
_{n})]|=|\mathbb{E}[V(1,\bar{S}_{n})]-\mathbb{E}[V(1+h,\bar{S}%
_{n})]|\leq C(\sqrt{h}+h).
\]
It follows from (\ref{equ-h}) and (\ref{equ-0}) that
\[
\limsup_{n\rightarrow \infty}|\mathbb{E}[\varphi(\bar{S}_{n}%
)]-\mathbb{E}[\varphi(X+\eta)]|\leq2C(\sqrt{h}+h).
\]
Since $h$ can be arbitrarily small, we have
\[
\lim_{n\rightarrow \infty}\mathbb{E}[\varphi(\bar{S}_{n})]=\mathbb{
{E}}[\varphi(X+\eta)].
\]

\end{proof}

\begin{remark}
From the proof we can check that the main assumption of identical
distribution of $\{X_{i},Y_{i}\}_{i=1}^{\infty}$ can be weaken to
\[
\mathbb{E}[\left \langle p,Y_{i}\right \rangle +\frac{1}{2}\left
\langle
AX_{i},X_{i}\right \rangle ]=G(p,A),\  \ i=1,2,\cdots,\ p\in \mathbb{R}%
^{d},\ A\in \mathbb{S}(d).
\]
Another essential condition is $\mathbb{E}[|X_{i}|^{2+\delta
}]+\mathbb{E}[|Y_{i}|^{1+\delta}]\leq C$ for some $\delta>0$. We do
not
need the condition $\mathbb{E}[|X_{i}|^{n}]+\mathbb{E}[|Y_{i}%
|^{n}]$ $<\infty$ for each $n\in \mathbb{N}$.
\end{remark}

We now give the proof of Theorem \ref{CLT}.

{\noindent \textbf{Proof of Theorem \ref{CLT}}. For the case when
the uniform elliptic condition (\ref{Ellip}) does not hold, we first
introduce a perturbation to prove the above convergence for $\varphi
\in C_{b.Lip}(\mathbb{R}^{d})$. According to Definition \ref{dd1}
and Proposition \ref{pp1} in Chap I, we can construct a
sublinear expectation space $(\bar{\Omega},\mathcal{\bar{H}},\mathbb{\bar{E}%
)}$ and a sequence of three random vectors $\{(\bar{X}_{i},\bar{Y}_{i}%
,\bar{\kappa}_{i})\}_{i=1}^{\infty}$ such that, for each
$n=1,2,\cdots$,
$\{(\bar{X}_{i},\bar{Y}_{i})\}_{i=1}^{n}\overset{d}{=}\{(X_{i},Y_{i}%
)\}_{i=1}^{n}$ and
$(\bar{X}_{n+1},\bar{Y}_{n+1},\bar{\kappa}_{n+1})$ is independent
from $\{(\bar{X}_{i},\bar{Y}_{i},\bar{\kappa}_{i})\}_{i=1}^{n}$ and,
moreover,
\[
\mathbb{\bar{E}}[\psi(\bar{X}_{i},\bar{Y}_{i},\bar{\kappa}_{i})]=(2\pi)^{-d/2}\int_{\mathbb{R}^{d}}\mathbb{E}[\psi(X_{i}%
,Y_{i},x)]e^{-|x|^{2}/{2}}dx\  \  \text{for}\ \psi \in C_{l.Lip}(\mathbb{R}%
^{3\times d}).
\]
We then use the perturbation }$\bar{X}${$_{i}^{\varepsilon}=\bar{X}%
_{i}+\varepsilon \bar{\kappa}_{i}$ for a fixed $\varepsilon>0$. It
is easy to see that the sequence
$\{(\bar{X}_{i}^{\varepsilon},\bar{Y}_{i})\}_{i=1}^{\infty}$
satisfies all conditions in the above CLT, in particular,%
\[
G_{\varepsilon}(p,A):=\mathbb{\bar{E}}[\frac{1}{2}\left \langle A\bar{X}%
_{1}^{\varepsilon},\bar{X}_{1}^{\varepsilon}\right \rangle +\left
\langle p,\bar{Y}_{1}\right \rangle
]=G(p,A)+\frac{\varepsilon^{2}}{2}\text{tr}[A].
\]
Thus it is strictly elliptic. We then can apply Lemma \ref{Lem-CLT} to%
\[
\bar{S}_{n}^{\varepsilon}:=\sum_{i=1}^{n}(\frac{\bar{X}{_{i}^{\varepsilon}}%
}{\sqrt{n}}+\frac{\bar{Y}_{i}}{n})=\sum_{i=1}^{n}(\frac{\bar{X}{_{i}}}%
{\sqrt{n}}+\frac{\bar{Y}_{i}}{n})+\varepsilon J_{n},\  \ J_{n}=\sum_{i=1}%
^{n}\frac{\bar{\kappa}_{i}}{\sqrt{n}}%
\]
and obtain%
\[
\lim_{n\rightarrow
\infty}\mathbb{\bar{E}}[\varphi(\bar{S}_{n}^{\varepsilon
})]=\mathbb{\bar{E}}[\varphi(\bar{X}+\bar{\eta}+\varepsilon
\bar{\kappa})],
\]
where $((\bar{X},\bar{\kappa}),(\bar{\eta},0))$ is
$\bar{G}$-distributed under }$\mathbb{\bar{E}[\cdot]}${\ and
\[
{\bar{G}}(\bar{p},\bar{A}):=\mathbb{\bar{E}}[\frac{1}{2}\left
\langle \bar
{A}(\bar{X}_{1},\bar{\kappa}_{1})^{T},(\bar{X}_{1},\bar{\kappa}_{1}%
)^{T}\right \rangle +\left \langle \bar{p},(\bar{Y}_{1},0)^{T}\right
\rangle ],{\  \ }\bar{A}\in \mathbb{S}(2d),\  \  \bar{p}\in
\mathbb{R}^{2d}.
\]
By Proposition \ref{GCD}, it is easy to prove that $(\bar{X
}+\varepsilon \bar{\kappa},\bar{\eta})$ is
$G_{\varepsilon}$-distributed and
$(\bar{X},\bar{\eta})$ is $G$-distributed. But we have%
\begin{align*}
|\mathbb{E}[\varphi(\bar{S}_{n})]-\mathbb{\bar{E}}[\varphi(\bar{S}%
_{n}^{\varepsilon})]| &
=|\mathbb{\bar{E}}[\varphi(\bar{S}_{n}^{\varepsilon
}-\varepsilon J_{n})]-\mathbb{\bar{E}}[\varphi(\bar{S}_{n}^{\varepsilon})]|\\
&  \leq \varepsilon C\mathbb{\bar{E}}[|J_{n}|]\leq C'\varepsilon
\end{align*}
and similarly, }

{$|\mathbb{E}[\varphi(X+\eta)]-\mathbb{\bar{E}}[\varphi
(\bar{X}+\bar{\eta}+\varepsilon \bar{\kappa})]|={|\mathbb{\bar
{E}}[\varphi(}${$\bar{X}$}${+\bar{\eta})]-\mathbb{\bar{E}%
}[\varphi(}${$\bar{X}$}${+\bar{\eta}+\varepsilon \bar{\kappa}%
)]|}\leq C\varepsilon$. } {Since $\varepsilon$ can be arbitrarily
small, it follows that
\[
\lim_{n\rightarrow \infty}\mathbb{E}[\varphi(\bar{S}_{n})]=\mathbb{
{E}}[\varphi(X+\eta)]\  \  \text{for}\ \varphi \in C_{b.Lip}(\mathbb{R}%
^{d}).
\]
} On the other hand, it is easy to check that
$\sup_{n}\mathbb{{E}[}|\bar
{S}_{n}|^{2}]+\mathbb{E}[|X+\eta|^{2}]<\infty$. {We then can apply
the following lemma to prove that the above convergence holds for
$\varphi$}$\in${$C(\mathbb{R}^{d})$ with linear growth condition.
The proof is complete. \hfill$\Box$}

\begin{lemma}
Let $({\Omega},\mathcal{{H}},\mathbb{E})$ and $(\widetilde
{\Omega},\mathcal{\widetilde{H}},\mathbb{\widetilde{E}})$ be two
sublinear expectation spaces and let $Y_{n}\in \mathcal{{H}}$ and
$Y\in \mathcal{\widetilde{H}}$, $n=1,2,\cdots$, be given. We assume
that, for a
given $p\geq1,$ $\sup_{n}\mathbb{E}[|Y_{n}|^{p}]+\widetilde{\mathbb{E}%
}[|Y|^{p}]<\infty$. If the convergence $\lim_{n\rightarrow \infty}%
\mathbb{E}[\varphi(Y_{n})]=\widetilde{\mathbb{E}}[\varphi(Y)]$ holds
for each $\varphi \in C_{b.Lip}(\mathbb{R}^{d})$, then it also holds
for all functions $\varphi \in C(\mathbb{R}^{d})$ with the growth
condition $|\varphi(x)|\leq C(1+|x|^{p-1})$.
\end{lemma}

\begin{proof}
We first prove that the above convergence holds for $\varphi \in C_{b}%
(\mathbb{R}^{d})$ with a compact support. In this case, for each
$\varepsilon>0$, we can find a $\bar{\varphi}\in
C_{b.Lip}(\mathbb{R}^{d})$ such that $\sup_{x\in
\mathbb{R}^{d}}|\varphi(x)-\bar{\varphi}(x)|\leq
\frac{\varepsilon}{2}$. {We have } { } {
\begin{align*}
&  |\mathbb{E}[\varphi(Y_{n})]-\widetilde{\mathbb{E}}[\varphi
(Y)]|\leq|\mathbb{E}[\varphi(Y_{n})]-\mathbb{E}[\bar{\varphi
}(Y_{n})]|+|\widetilde{\mathbb{E}}[\varphi(Y)]-\widetilde{\mathbb{E}}%
[\bar{\varphi}(Y)]|\\
&  +|\mathbb{E}[\bar{\varphi}(Y_{n})]-\widetilde{\mathbb{E}}%
[\bar{\varphi}(Y)]|\leq \varepsilon+|\mathbb{E}[\bar{\varphi}%
(Y_{n})]-\widetilde{\mathbb{E}}[\bar{\varphi}(Y)]|.
\end{align*}
Thus }$\limsup_{n\rightarrow \infty}|\mathbb{E}[\varphi(Y_{n}%
)]-\widetilde{\mathbb{E}}[\varphi(Y)]|\leq \varepsilon$. The
convergence must hold since $\varepsilon$ can be arbitrarily small.

Now let $\varphi$ be an arbitrary $C(\mathbb{R}^{d})$-function with
growth condition $|\varphi(x)|\leq C(1+|x|^{p-1})$. For each $N>0$
we can find $\varphi_{1},\varphi_{2}\in C(\mathbb{R}^{d})$ such that
$\varphi=\varphi _{1}+\varphi_{2}$ where $\varphi_{1}$ has a compact
support and $\varphi _{2}(x)=0$ for $|x|\leq N$, and
$|\varphi_{2}(x)|\leq|\varphi(x)|$ for all $x$. It is clear that
\[
|\varphi_{2}(x)|\leq \frac{2C(1+|x|^{p})}{N}\  \  \text{for}\ x\in
\mathbb{R}^{d}.
\]
Thus
\begin{align*}
|\mathbb{E}[\varphi(Y_{n})]-\widetilde{\mathbb{E}}[\varphi(Y)]| &
=|\mathbb{E}[\varphi_{1}(Y_{n})+\varphi_{2}(Y_{n})]-\widetilde
{\mathbb{E}}[\varphi_{1}(Y)+\varphi_{2}(Y)]|\\
&  \leq|\mathbb{E}[\varphi_{1}(Y_{n})]-\widetilde{\mathbb{E}}%
[\varphi_{1}(Y)]|+\mathbb{E}[|\varphi_{2}(Y_{n})|]+\widetilde
{\mathbb{E}}[|\varphi_{2}(Y)|]\\
&  \leq|\mathbb{E}[\varphi_{1}(Y_{n})]-\widetilde{\mathbb{E}}%
[\varphi_{1}(Y)]|+\frac{2C}{N}(2+\mathbb{E}[|Y_{n}|^{p}]+\widetilde
{\mathbb{E}}[|Y|^{p}])\\
&  \leq|\mathbb{E}[\varphi_{1}(Y_{n})]-\widetilde{\mathbb{E}}%
[\varphi_{1}(Y)]|+\frac{\bar{C}}{N},
\end{align*}
where $\bar{C}=2C(2+\sup_{n}\mathbb{E}[|Y_{n}|^{p}]+\widetilde
{\mathbb{E}}[|Y|^{p}]).$ We thus have $\limsup_{n\rightarrow \infty
}|\mathbb{E}[\varphi(Y_{n})]-\widetilde{\mathbb{E}}[\varphi
(Y)]|\leq \frac{\bar{C}}{N}$. Since $N$ can be arbitrarily large,
$\mathbb{E}[\varphi(Y_{n})]$ must converge to $\widetilde{\mathbb{E}%
}[\varphi(Y)]$.
\end{proof}

\begin{exercise}
Let $X_{i}\in \mathcal{H}, i=1,2,\cdots$, be such that $X_{i+1}$ is
independent from $\{X_{1},\cdots,X_{i}\}$, for each $i=1,2,\cdots.$
We further assume that
\[
\mathbb{E}[X_{i}]=-\mathbb{E}[-X_{i}]=0,
\]
\[
\lim_{i\rightarrow
\infty}\mathbb{E}[X_{i}^{2}]=\overline{\sigma}^{2}< \infty,
\lim_{i\rightarrow \infty}-\mathbb{E}[-X_{i}^{2}]=\underline
{\sigma}^{2},
\]
\[
\mathbb{E}[|X_{i}|^{2+\delta}]\leq M\ ~\text{for~some~}\delta
>0~\text{and}~\text{a~constant}\ M.
\]

Prove that the sequence $\{ \bar{S}_{n}\}_{n=1}^{\infty}$ defined by
\[
\bar{S}_{n}=\frac{1}{\sqrt{n}}\sum_{i=1}^{n}X_{i}%
\]
converges in law to $X$, i.e.,
\[
\lim_{n\rightarrow \infty}\mathbb{E}[\varphi(\bar{S}_{n})]=
{\mathbb{E}}[\varphi(X)]\ \ \text{for}\ \varphi \in
C_{b,lip}(\mathbb{R}),
\]
where $X\sim N(\{0\} \times \lbrack \underline{\sigma}^{2}%
,\overline{\sigma}^{2}]).$

In particular, if $\overline{\sigma}^{2}=\underline{\sigma}^{2}$, it
becomes a classical central limit theorem.
\end{exercise}

\section*{Notes and Comments}
\addcontentsline{toc}{section}{Notes and Comments} The contents of
this chapter are mainly from Peng (2008) \cite{Peng2008b} (see also
Peng (2007) \cite{Peng2007a}).

The notion of $G$-normal distribution was firstly introduced by Peng
(2006) \cite{Peng2006a} for 1-dimensional case, and Peng (2008)
\cite{Peng2008} for multi-dimensional case. In the classical
situation, a distribution satisfying equation (\ref{ch2e1}) is said
to be stable (see L\'{e}vy (1925) \cite{Levy1} and (1965)
\cite{Levy2}). In this sense, our $G$-normal distribution can be
considered as the most typical stable distribution under the
framework of sublinear expectations.

Marinacci (1999) \cite{Marinacci} used different notions of
distributions and independence via capacity and the corresponding
Choquet expectation to obtain a law of large numbers and a central
limit theorem for non-additive probabilities (see also Maccheroni
and Marinacci (2005) \cite{Marinacci1} ). But since a sublinear
expectation can not be characterized by the corresponding capacity,
 our results can not be derived from theirs. In fact,  our results
show that the limit in CLT, under uncertainty, is a $G$-normal
distribution in which the distribution uncertainty is not just the
parameter of the classical normal distributions (see Exercise
\ref{exxee1}).

The notion of viscosity solutions plays a basic role in the
definition and properties of $G$-normal distribution and maximal
distribution. This notion was initially introduced by Crandall and
Lions (1983) \cite{CrandallL}. This is a fundamentally important
notion in the theory of nonlinear parabolic and elliptic PDEs.
Readers are referred to Crandall, Ishii and Lions (1992) \cite{CIL}
for rich references of the beautiful and powerful theory of
viscosity solutions. For books on the theory of viscosity solutions
and the related HJB equations, see Barles (1994) \cite{barles},
Fleming and Soner (1992) \cite{FS} as well as Yong and Zhou (1999)
\cite{Yong-Zhou}.

We note that, for the case when the uniform elliptic condition
holds, the viscosity solution (\ref{e320}) becomes a classical
$C^{1+\frac{\alpha}{2},2+\alpha}$-solution (see Krylov (1987)
\cite{Krylov1} and the recent works in Cabre and Caffarelli
 (1997) \cite{Caff1997} and Wang (1992) \cite{WangL}). In
1-dimensional situation, when $\underline{\sigma}^{2}>0$, the
$G$-equation becomes the following Barenblatt equation:
$$\partial_tu+\gamma|\partial_t u|=\triangle u,~ |\gamma|<1.$$
This equation was first introduced by Barenblatt (1979)
\cite{Barenblatt} (see also Avellaneda, Levy and Paras (1995)
\cite{Avellaneda}).

%
%
\chapter{$G$-Brownian Motion and It\^{o}'s Integral}
\label{ch3} The aim of this chapter is to introduce the concept of
$G$-Brownian motion, to study its properties and to construct
It\^{o}'s integral with respect to $G$-Brownian motion. We emphasize
here that our definition of $G$-Brownian motion is consistent with
the classical one in the sense that if there is no volatility
uncertainty. Our $G$-Brownian motion also has independent increments
with identical $G$-normal distributions. $G$-Brownian motion has a
very rich and interesting new structure which non-trivially
generalizes the classical one. We thus can establish the related
stochastic calculus, especially It\^o's integrals and the related
quadratic variation process. A very interesting new phenomenon of
our $G$-Brownian motion is that its quadratic process also has
independent increments which are identically distributed. The
corresponding $G$-It\^o's formula is obtained.

\section{$G$-Brownian Motion and its Characterization}

\begin{definition}
Let $(\Omega,\mathcal{H},\mathbb{{E})}$ be a sublinear expectation
space. $(X_{t})_{t\geq0}$ is called a $d$-dimensional
\textbf{stochastic process}\index{Stochastic process} if for each
$t\geq0$, $X_{t}$ is a $d$-dimensional random vector in
$\mathcal{H}$.
\end{definition}

Let $G(\cdot):\mathbb{S}(d)\rightarrow\mathbb{R}$ be a given
monotonic and sublinear function. By Theorem \ref{t1} in Chapter
\ref{ch1}, there exists a bounded, convex and closed subset $\Sigma
\subset \mathbb{S}_+(d)$ such that
\[
G(A)=\frac{1}{2}\sup_{B\in \Sigma}\left(  A,B\right)  ,\  \  \  \
A\in \mathbb{S}(d).
\]
By Section \ref{c2s2} in Chapter \ref{ch2}, we know that the
$G$-normal distribution $N(\{0\}\times\Sigma)$ exists.

We now give the definition of $G$-Brownian motion.

\begin{definition}
A $d$-dimensional process $(B_{t})_{t\geq0}$ on a sublinear
expectation space $(\Omega,\mathcal{H},\mathbb{E})$ is called a
$G$\textbf{--Brownian} \textbf{motion}
\label{bbb}\index{$G$-Brownian motion}if the following properties
are satisfied: \newline \textup{(i)} $B_{0}(\omega)=0$;\newline
\textup{(ii)} For each $t,s\geq0$, the increment $B_{t+s}-B_{t}$ is
$N(\{0\}\times s\Sigma)$-distributed and is independent from
$(B_{t_{1}},B_{t_{2}},\cdots,B_{t_{n}})$, for each $n\in \mathbb{N}$
and $0\leq t_{1}\leq \cdots \leq t_{n}\leq t$.
\end{definition}

\begin{remark}
We can prove that, for each $t_{0}>0$,
$(B_{t+t_{0}}-B_{t_{0}})_{t\geq0}$ is
a $G$-Brownian motion. For each $\lambda>0$, $(\lambda^{-\frac{1}{2}%
}B_{\lambda t})_{t\geq0}$ is also a $G$-Brownian motion. This is the
scaling property of $G$-Brownian motion, which is the same as that
of the classical Brownian motion.
\end{remark}

We will denote in the rest of this book
\[
B_{t}^{\mathbf{a}}=\langle \mathbf{a},B_{t}\rangle\  \  \ \text{for
each }\mathbf{a}=(a_{1},\cdots,a_{d})^{T}\in \mathbb{R}^{d}.
\]

By the above definition we have the following proposition which is
important in stochastic calculus.

\begin{proposition}
Let $(B_{t})_{t\geq0}$ be a $d$-dimensional $G$\textbf{-}Brownian
motion\textbf{ }on a sublinear expectation space $(\Omega,\mathcal{H}%
,\mathbb{E}).$ Then $(B_{t}^{\mathbf{a}})_{t\geq0}$ is a
$1$-dimensional $G_{\mathbf{a}}$\textbf{-}Brownian motion\textbf{
}for each $\mathbf{a\in }\mathbb{R}^{d}$, where
$G_{\mathbf{a}}(\alpha)=\frac{1}{2}(\sigma
_{\mathbf{aa}^{T}}^{2}\alpha^{+}-\sigma_{-\mathbf{aa}^{T}}^{2}\alpha^{-})$,
$\sigma_{\mathbf{aa}^{T}}^{2}=2G(\mathbf{aa}^{T})=\mathbb{E}[\langle \mathbf{a},B_1\rangle^2]$, $\sigma_{-\mathbf{aa}^{T}%
}^{2}=-2G(-\mathbf{aa}^{T})=-\mathbb{E}[-\langle
\mathbf{a},B_1\rangle^2].$

In particular, for each $t,s\geq0$, $B_{t+s}^{\mathbf{a}}-B_{t}^{\mathbf{a}%
}\overset{d}{=}N(\{0\}\times [s\sigma_{-\mathbf{aa}^{T}}^{2},s\sigma_{\mathbf{aa}^{T}%
}^{2}]).$
\end{proposition}

\begin{proposition}

For each convex function $\varphi$ , we have%
\[
\mathbb{E}[\varphi(B_{t+s}^{\mathbf{a}}-B_{t}^{\mathbf{a}})]=\frac
{1}{\sqrt{2\pi s\sigma_{\mathbf{aa}^{T}}^{2}}}\int_{-\infty}^{\infty}%
\varphi(x)\exp(-\frac{x^{2}}{2s\sigma_{\mathbf{aa}^{T}}^{2}})dx.
\]
For each concave function $\varphi$ and $\sigma^2_{-\mathbf{aa}^T}>0$, we have%
\[
\mathbb{E}[\varphi(B_{t+s}^{\mathbf{a}}-B_{t}^{\mathbf{a}})]=\frac
{1}{\sqrt{2\pi s\sigma_{-\mathbf{aa}^{T}}^{2}}}\int_{-\infty}^{\infty}%
\varphi(x)\exp(-\frac{x^{2}}{2s\sigma_{-\mathbf{aa}^{T}}^{2}})dx.
\]
In particular, we have%
\begin{align*}
\mathbb{E}[(B_{t}^{\mathbf{a}}-B_{s}^{\mathbf{a}})^{2}]  &
=\sigma_{\mathbf{aa}^{T}}^{2}(t-s),\  \  \  \mathbb{E}[(B_{t}^{\mathbf{a}%
}-B_{s}^{\mathbf{a}})^{4}]=3\sigma_{\mathbf{aa}^{T}}^{4}(t-s)^{2},\\
\mathbb{E}[-(B_{t}^{\mathbf{a}}-B_{s}^{\mathbf{a}})^{2}]  &
=-\sigma_{-\mathbf{aa}^{T}}^{2}(t-s),\  \  \  \mathbb{E}[-(B_{t}%
^{\mathbf{a}}-B_{s}^{\mathbf{a}})^{4}]=-3\sigma_{-\mathbf{aa}^{T}}%
^{4}(t-s)^{2}.
\end{align*}

\end{proposition}

The following theorem gives a characterization of $G$-Brownian
motion.

\begin{theorem}
Let $(B_{t})_{t\geq0}$ be a $d$-dimensional process defined on a
sublinear expectation space $(\Omega,\mathcal{H},\mathbb{E})${ such
that
\newline \textup{(i)} }$B_{0}(\omega)${$=0$;\newline \textup{(ii)} For each
$t,s\geq0$, $B${$_{t+s}-B_{t}$ and }$B${$_{s}$ are identically
distributed and
$B${$_{t+s}-B_{t}$ is }independent from $(B_{t_{1}},B_{t_{2}},\cdots,B_{t_{n}}%
)$, for each $n\in \mathbb{N}$ and $0\leq t_{1}\leq \cdots \leq
t_{n}\leq t$.
\newline \textup{(iii)} }}$\mathbb{E}${{$[B_{t}]=\mathbb{E}%
[-B_{t}]=0$ and
$\lim_{t\downarrow0}\mathbb{E}[|B_{t}|^{3}]t^{-1}=0$.
\newline Then }}$(B_{t})_{t\geq0}${{ is a $G$-Brownian motion with }%
}$G(A)=\frac{1}{2}\mathbb{E}[\langle AB_{1},B_{1}\rangle],\ A\in
\mathbb{S}(d)$.
\end{theorem}

\begin{proof}
{{\textbf{ }}We only need to prove that $B_{1}$ is }$G${-normal
distributed and }$B_{t}\overset{d}{=}\sqrt{t}B_{1}${. We first prove
that
\[
\mathbb{E}[\langle AB_{t},B_{t}\rangle]=2G(A)t,\  \ A\in
\mathbb{S}(d).
\]
For each given }$A\in \mathbb{S}(d)$, {we set $b(t)=$}$\mathbb{{E}
}[\langle AB_{t},B_{t}\rangle]${{. Then }$b(0)=0$ and $|b(t)|\leq|A|($%
}$\mathbb{E}${{$[|B_{t}|^{3}])^{2/3}\rightarrow0$ as }}$t\rightarrow
0${{.} {Since for each $t,s\geq0$,%
\begin{align*}
b(t+s)  &  =\mathbb{E}[\langle AB_{t+s},B_{t+s}\rangle]=\mathbb{\hat{E}%
}[\langle A(B_{t+s}-B_{s}+B_{s}),B_{t+s}-B_{s}+B_{s}\rangle]\\
&  =\mathbb{E}[\langle
A(B_{t+s}-B_{s}),(B_{t+s}-B_{s})\rangle+\langle
AB_{s},B_{s}\rangle+2\langle A(B_{t+s}-B_{s}),B_{s}\rangle]\\
&  =b(t)+b(s),
\end{align*}
we have }$b(t)=b(1)t=$}$2G(A)t${.{ } }

{ {We now prove that }$B_{1}$ is }$G${-normal distributed and }$B_{t}%
\overset{d}{=}\sqrt{t}B_{1}${. For this, we just need to prove that,
for each fixed $\varphi \in C_{b.Lip}(\mathbb{R}^{d}\mathbb{)}$, the
function
\[
u(t,x):=\mathbb{E}[\varphi(x+B_{t})],\  \ (t,x)\in \lbrack0,\infty
)\times \mathbb{R}^{d}%
\]
is the viscosity solution of the following $G$-heat equation:}%
\begin{equation}
\partial_{t}u-G(D^{2}u)=0,\  \ u|_{t=0}=\varphi. \label{G-heat-BM}%
\end{equation}

{ {We first prove that }$u$ is Lipschitz in $x$ and
$\frac{1}{2}$-H\"{o}lder continuous
in $t$. In fact, {for each fixed }$t$,{ }$u(t,\cdot)\in${$C_{b.Lip}%
(\mathbb{R}^{d}\mathbb{)}$ }since%
\begin{align*}
|u(t,x)-u(t,y)|&=
|\mathbb{E}[\varphi(x+B_{t})]-\mathbb{E}[\varphi(y+B_{t})]|\\
&
\leq \mathbb{E}[|\varphi(x+B_{t})-\varphi(y+B_{t})|]\\
&  \leq C|x-y|,
\end{align*}
where $C$ is Lipschitz constant of $\varphi$.

{For each $\delta \in \lbrack0,t]$, since $B_{t}-B_{\delta}$ is
independent from $B_{\delta}$, w}e also have
\begin{align*}
u(t,x)  & =\mathbb{E}[\mathbb{\varphi(}x+B_{\delta}+(B_{t}-B_{\delta
})]\\
&  =\mathbb{E}[\mathbb{E}[\varphi(y+(B_{t}-B_{\delta
}))]_{y=x+B_{\delta}}],
\end{align*}
hence%
\begin{equation}
u(t,x)=\mathbb{E}[u(t-\delta,x+B_{\delta})]. \label{Dyna}%
\end{equation}
{Thus}%
\begin{align*}
|u(t,x)-u(t-\delta,x)|  &  =|\mathbb{E}[u(t-\delta,x+B_{\delta
})-u(t-\delta,x)]|\\
&  \leq \mathbb{E}[|u(t-\delta,x+B_{\delta})-u(t-\delta,x)|]\\
&  \leq \mathbb{E}[C|B_{\delta}|]\leq C\sqrt{2G(I)}\sqrt{\delta}.
\end{align*}
To prove that $u$ is a viscosity solution of {(\ref{G-heat-BM})},
{we fix
$(t,x)\in(0,\infty)\times \mathbb{R}^{d}$} and{ let $v\in C_{b}^{2,3}%
([0,\infty)\times \mathbb{R}^{d})$ be such that $v\geq u$ and
$v(t,x)=u(t,x)$. From (\ref{Dyna}) we have
\[
v(t,x)=\mathbb{E}[u(t-\delta,x+B_{\delta})]\leq \mathbb{{E}%
}[v(t-\delta,x+B_{\delta})].
\]
Therefore by Taylor's expansion,
\begin{align*}
0  &  \leq \mathbb{E}[v(t-\delta,x+B_{\delta})-v(t,x)]\\
&  =\mathbb{E}[v(t-\delta,x+B_{\delta})-v(t,x+B_{\delta}%
)+(v(t,x+B_{\delta})-v(t,x))]\\
&  =\mathbb{E}[-\partial_{t}v(t,x)\delta+\langle Dv(t,x),B_{\delta
}\rangle+\frac{1}{2}\langle
D^{2}v(t,x)B_{\delta},B_{\delta}\rangle+I_{\delta
}]\\
&  \leq-\partial_{t}v(t,x)\delta+\frac{1}{2}\mathbb{E}[\langle
D^{2}v(t,x)B_{\delta},B_{\delta}\rangle]+\mathbb{E}[I_{\delta}]\\
&
=-\partial_{t}v(t,x)\delta+G(D^{2}v(t,x))\delta+\mathbb{E}[I_{\delta
}],
\end{align*}
where}%
\begin{align*}
I_{\delta}  &  =\int_{0}^{1}-[\partial_{t}v(t-\beta
\delta,x+B_{\delta
})-\partial_{t}v(t,x)]\delta d\beta \\
&  +\int_{0}^{1}\int_{0}^{1}\langle(D^{2}v(t,x+\alpha \beta B_{\delta}%
)-D^{2}v(t,x))B_{\delta},B_{\delta}\rangle \alpha d\beta d\alpha.
\end{align*}
With the assumption \textup{(iii)} we can check that $\lim_{\delta
\downarrow 0}\mathbb{E}[|I_{\delta}|]\delta^{-1}=0$, from which we
get $\partial_{t}v(t,x)-G(D^{2}v(t,x))\leq0$, hence $u$ is a
viscosity subsolution of {(\ref{G-heat-BM}). We can analogously
prove that }$u$ is a viscosity supersolution. Thus $u$ is a
viscosity solution and $(B_{t})_{t\geq0}$ is a $G$-Brownian motion.
The proof is complete. }
\end{proof}

\begin{exercise}
Let $B_t$ be a 1-dimensional Brownian motion, and $B_1\overset{d}{=}
N(\{0\}\times[\underline{\sigma}^2,\overline{\sigma}^2])$. Prove
that for each $m\in \mathbb{N}$,
\begin{equation*}
\hat{\mathbb{E}}[|B_t|^m]=\left\{\begin{array}{ll}
2(m-1)!!\overline{\sigma}^m t^{\frac{m}{2}}
/\sqrt{2\pi} &m~ \text{is~odd},\\
(m-1)!!\overline{\sigma}^m t^{\frac{m}{2}} &m
~\text{is~even}.\end{array}\right.
\end{equation*}

\end{exercise}

\section{{Existence of{ $G$-Brownian Motion}}}

In the rest of this book, we denote by
$\Omega=C_{0}^{d}(\mathbb{R}^{+})$ the
space of all $\mathbb{R}^{d}$--valued continuous paths $(\omega_{t}%
)_{t\in \mathbb{R}^{+}}$, with $\omega_{0}=0$, {{equipped with the
distance
\[
\rho(\omega^{1},\omega^{2}):=\sum_{i=1}^{\infty}2^{-i}[(\max_{t\in
\lbrack 0,i]}|\omega_{t}^{1}-\omega_{t}^{2}|)\wedge1].
\]
For each fixed }}$T\in \lbrack0,\infty),$ we set $\Omega_{T}:=\{
\omega _{\cdot \wedge T}:\omega \in \Omega \}.$ We will consider the
canonical process $B_{t}(\omega)=\omega_{t}$, $t\in
\lbrack0,\infty)$, for $\omega \in \Omega$.

For each fixed $T\in \lbrack0,\infty)$, we set
\[
L_{ip}(\Omega_{T}):=\{ \varphi(B_{t_{1}\wedge
T},\cdots,B_{t_{n}\wedge T}):n\in \mathbb{N},t_{1},\cdots,t_{n}\in
\lbrack0,\infty),\  \varphi \in
C_{l.Lip}(\mathbb{R}^{d\times n})\  \}.\text{ }%
\]
It is clear that
$L_{ip}(\Omega_{t})${$\subseteq$}$L_{ip}(\Omega_{T})${, for
$t\leq T$. We also set%
\[
L_{ip}(\Omega):=%
{\displaystyle \bigcup \limits_{n=1}^{\infty}}
L_{ip}(\Omega_{n}).
\]
}

\begin{remark}
It is clear that $C_{l.Lip}(\mathbb{R}^{d\times n})$,
$L_{ip}(\Omega_{T})$ and $L_{ip}(\Omega)${{ are vector lattices.
Moreover, note that $\varphi,\psi \in C_{l.Lip}(\mathbb{R}^{d\times
n})$ imply $\varphi
\cdot \psi \in C_{l.Lip}(\mathbb{R}^{d\times n})$, then $X$, $Y\in$}}%
$L_{ip}(\Omega_{T})${{ imply $X\cdot Y\in$}}$L_{ip}(\Omega_{T})${{.
}}In particular, for each $t\in \lbrack0,\infty)$, $B_{t}\in
L_{ip}(\Omega)$.
\end{remark}

Let $G(\cdot):\mathbb{S}(d)\rightarrow \mathbb{R}$ be a given
monotonic and sublinear function. In the following, we want to
construct a sublinear expectation on $(\Omega,L_{ip}(\Omega))$ such
that the canonical process {$(B_{t})_{t\geq0}$ is a $G$-Brownian
motion.} For this, we first construct a sequence of $d$-dimensional
random vectors $(${{$\xi_{i})_{i=1}^{\infty}$ }}on
a sublinear expectation space $(\widetilde{\Omega},\widetilde{\mathcal{H}%
},\widetilde{\mathbb{E}})$ such that {{$\xi_{i}$ is }}$G$-normal
distributed and {{$\xi_{i+1}$ is independent from
$(\xi_{1},\cdots,\xi_{i})$ for each $i=1,2,\cdots$. }}

We now introduce a sublinear expectation $\hat{\mathbb{E}}$ defined
on $L_{ip}(\Omega)$ via the following procedure: for each $X\in
L_{ip}(\Omega)$ with
\[
X=\varphi(B_{t_{1}}-B_{t_{0}},B_{t_{2}}-B_{t_{1}},\cdots,B_{t_{n}}-B_{t_{n-1}%
})
\]
{ for some {$\varphi \in C_{l.Lip}(\mathbb{R}^{d\times n})$ and $0=t_{0}%
<t_{1}<\cdots<t_{n}<\infty$, we set%
\[
\hat{\mathbb{E}}[\varphi(B_{t_{1}}-B_{t_{0}},B_{t_{2}}-B_{t_{1}}%
,\cdots,B_{t_{n}}-B_{t_{n-1}})]
\]%
\[
:=\widetilde{\mathbb{E}}[\varphi(\sqrt{t_{1}-t_{0}}\xi_{1},\cdots,\sqrt
{t_{n}-t_{n-1}}\xi_{n})].
\]
}}

{{The related conditional expectation of $X=\varphi(B_{t_{1}},B_{t_{2}%
}-B_{t_{1}},\cdots,B_{t_{n}}-B_{t_{n-1}})$ under $\Omega_{t_{j}}$ is
defined
by%
\begin{align}
\hat{\mathbb{E}}[X|{{\Omega_{t_{j}}}}]  &  =\hat{\mathbb{E}}%
[\varphi(B_{t_{1}},B_{t_{2}}-B_{t_{1}},\cdots,B_{t_{n}}-B_{t_{n-1}}%
)|{{\Omega_{t_{j}}}}]\label{Condition}\\
&  :=\psi(B_{t_{1}},\cdots,B_{t_{j}}-B_{t_{j-1}}),\nonumber
\end{align}
where}}%
\[
\psi(x_{1},\cdots,x_{j})=\widetilde{\mathbb{E}}[\varphi(x_{1},\cdots
,x_{j},\sqrt{t_{j+1}-t_{j}}\xi_{j+1},\cdots,\sqrt{t_{n}-t_{n-1}}\xi_{n})].
\]
{{It is easy to check that }}$\hat{\mathbb{E}}${{$[\cdot]$
consistently defines a sublinear expectation on }}$L_{ip}(\Omega)$
and {$(B_{t})_{t\geq0}$
is a $G$-Brownian motion}. Since $L_{ip}(\Omega_{T})${$\subseteq$}%
$L_{ip}(\Omega)$, $\hat{\mathbb{E}}${{$[\cdot]$ is also a sublinear
expectation on }}$L_{ip}(\Omega_{T})$.

\begin{definition}
{ { The {sublinear }expectation }}$\hat{\mathbb{E}}$$[\cdot]$%
{{$:L_{ip}(\Omega)\rightarrow \mathbb{R}$ defined through the above
procedure is called a
$G$\textbf{--expectation}.\index{$G$-expectation}\label{ge} The
corresponding canonical process
$(B_{t})_{t\geq0}$ on the sublinear expectation space $(\Omega,L_{ip}%
(\Omega),\hat{\mathbb{E}})$ is called a $G$--Brownian motion. } }
\end{definition}

In the rest of this book, when we talk about {{$G$--Brownian motion,
we mean that the canonical process $(B_{t})_{t\geq0}$ is under
$G$-expectation.}}

\begin{proposition}
\label{Prop-1-9-1 }{We list the properties of }$\hat{\mathbb{E}}%
${$[\cdot|\Omega_{t}]$ that hold for each $X,Y\in$}$L_{ip}(\Omega)${:}%
\newline \textup{(i) }{\textbf{ }{If $X\geq Y$, then }}$\hat{\mathbb{E}}%
${{$[X|\Omega_{t}]\geq \hat{\mathbb{E}}[Y|{\Omega}_{t}]$.\newline%
}\textup{(ii) }\textbf{
}}$\hat{\mathbb{E}}${$[\eta|\Omega_{t}]=\eta$, \ for
each }$t\in \lbrack0,\infty)$ and {{$\eta \in$}}$L_{ip}(\Omega_{t})$%
{{.\newline \textup{(iii) }\textbf{
}}}$\hat{\mathbb{E}}${{$[X|\Omega
_{t}]-\hat{\mathbb{E}}[Y|\Omega_{t}]\leq \hat{\mathbb{E}}%
[X-Y|\Omega_{t}].$\newline \textup{(iv) }{
}}}$\hat{\mathbb{E}}${{{$[\eta
X|\Omega_{t}]=\eta^{+}\hat{\mathbb{E}}[X|\Omega_{t}]+\eta^{-}%
\hat{\mathbb{E}}[-X|\Omega_{t}]$  for each $\eta \in L_{ip}(\Omega_{t}).$%
}}}\newline \textup{(v)} $
\hat{\mathbb{E}}{{[\hat{\mathbb{E}}[X|\Omega_{t}]|\Omega
_{s}]=\hat{\mathbb{E}}[X|\Omega_{t\wedge s}],\ }}\text{ in
particular},\
\hat{\mathbb{E}}[\hat{\mathbb{E}}[X|\Omega_{t}]]=\mathbb{\hat{E}
}[X]. $

For each $X\in L_{ip}(\Omega^{t})$, $\hat{\mathbb{E}}[X|{\Omega}%
_{t}]=\hat{\mathbb{E}}[X]$, where $L_{ip}(\Omega^{t})$ is the linear
space of random variables with the form
\begin{align*}
&  {\varphi(B_{t_{2}}-B_{t_{1}},B_{t_{3}}-B_{t_{2}},\cdots,B_{t_{n+1}%
}-B_{t_{n}}),}\\
\  &  \  \ n=1,2,\cdots,\  \varphi \in C_{l.Lip}(\mathbb{R}^{d\times n}%
),\ t_{1},\cdots,t_{n},t_{n+1}\in \lbrack t,\infty).
\end{align*}

\end{proposition}

\begin{remark}
{\textup{(ii)}} and {{\textup{(iii)}}} imply
\[
\hat{\mathbb{E}}{{[X+\eta|{\Omega_{t}}]=\hat{\mathbb{E}}[X|{\Omega
_{t}}]+\eta}}%
\ \ \text{for}\ \eta\in L_{ip}(\Omega_{t}).\]

\end{remark}

We now consider the completion of sublinear expectation space
$(\Omega ,L_{ip}(\Omega),\hat{\mathbb{E}}).$

We denote by $L_{G}^{p}(\Omega)$, $p\geq1$, the completion of
$L_{ip}(\Omega)$
under the norm $\left \Vert X\right \Vert _{p}:=(\hat{\mathbb{E}}%
[|X|^{p}])^{1/p}$. Similarly, we can define $L_{G}^{p}(\Omega_{T})$,
$L_{G}^{p}(\Omega_{T}^{t})$ and $L_{G}^{p}(\Omega^{t})$. It is clear
that for
each $0\leq t\leq T<\infty$, $L_{G}^{p}(\Omega_{t})\subseteq L_{G}^{p}%
(\Omega_{T})\subseteq L_{G}^{p}(\Omega)$.

According to Sec.\ref{c1s5} in Chap.\ref{ch1},
$\hat{\mathbb{E}}[\cdot]$  can be continuously extended to a
sublinear expectation on $(\Omega,L_{G}^{1}(\Omega))$ still denoted
by $\hat{\mathbb{E}}[\cdot]$. We now consider the extension of
conditional expectations. For each fixed $t\leq T$, the conditional
$G$-expectation ${{\hat{\mathbb{E}}[\cdot|{{\Omega}_{t}}]}}:L_{ip}%
(\Omega_{T})\rightarrow L_{ip}(\Omega_{t})$ is a continuous mapping
under $\left \Vert \cdot \right \Vert $. Indeed, we have
\[
{{\hat{\mathbb{E}}[X|{\Omega_{t}}]}}-{{\hat{\mathbb{E}}[Y|{\Omega_{t}%
}]}}\leq{{\hat{\mathbb{E}}[X-Y|{\Omega_{t}}]\leq \hat{\mathbb{E}}
[|X-Y||{{\Omega}_{t}}],}}%
\]
then
\[
|{{\hat{\mathbb{E}}[X|{\Omega_{t}}]}}-{{\hat{\mathbb{E}}[Y|{\Omega
_{t}}]|}}\leq{{\hat{\mathbb{E}}[|X-Y||{\Omega_{t}}].}}%
\]
We thus obtain
\[
\left \Vert {{\hat{\mathbb{E}}[X|{\Omega_{t}}]}}-{{\hat{\mathbb{E}}
[Y|{\Omega_{t}}]}}\right \Vert \leq \left \Vert X-Y\right \Vert .
\]
It follows that ${{\hat{\mathbb{E}}[\cdot|{\Omega_{t}}]}}$ can be
also extended as a continuous mapping
\[
{{\hat{\mathbb{E}}[\cdot|{\Omega_{t}}]}}:L_{G}^{1}(\Omega_{T})\rightarrow
L_{G}^{1}(\Omega_{t}).
\]
If the above $T$ is not fixed, then we can obtain
${{\hat{\mathbb{E}}
[\cdot|{\Omega_{t}}]}}:L_{G}^{1}(\Omega)\rightarrow
L_{G}^{1}(\Omega_{t})$.

\begin{remark}
The above proposition also holds for $X,Y\in L_{G}^{1}(\Omega)$. But
in {{\textup{(iv)}}}, $\eta \in L_{G}^{1}(\Omega_{t})$ should be
bounded, since $X,Y\in L_{G}^{1}(\Omega)$ does not imply $X\cdot
Y\in L_{G}^{1}(\Omega).$
\end{remark}

In particular, we have the following independence:
\[
{{\hat{\mathbb{E}}[X|{\Omega_{t}}]}}={{\hat{\mathbb{E}}[X}}%
],\medskip \  \  \forall X\in L_{G}^{1}(\Omega^{t}).
\]

We give the following definition similar to the classical one:

\begin{definition}
An $n$-dimensional random vector $Y\in(L_{G}^{1}({\Omega}))^{n}$ is
said to be independent from $\Omega_{t}$ for some given $t$ if for
each
$\varphi \in C_{b.Lip}(\mathbb{R}^{n})$ we have%
\[
{{\hat{\mathbb{E}}}}[\varphi(Y)|{{{\Omega_{t}}}}]={{\hat{\mathbb{E}}}%
}[\varphi(Y)].
\]

\end{definition}

\begin{remark}
Just as in the classical situation, the increments of $G$--Brownian
motion $(B_{t+s}-B_{t})_{s\geq0}$ are independent from $\Omega_{t}$.
\end{remark}

The following property is very useful.

\begin{proposition}
\label{d-E-x+y}Let $X,Y\in L_{G}^{1}(\Omega)$ be such that
${{\mathbb{\hat{E}
}[Y|{\Omega_{t}}}}]=-{{\hat{\mathbb{E}}[-Y|{\Omega_{t}}}}]$, for
some
$t\in \lbrack0,T]$. Then we have%
\[
\hat{\mathbb{E}}[X+Y|\Omega_{t}]=\hat{\mathbb{E}}[X|\Omega
_{t}]+\hat{\mathbb{E}}[Y|\Omega_{t}].
\]
In particular, if\  ${{\hat{\mathbb{E}}[Y|{\Omega_{t}}}}]={{\mathbb{\hat{E}%
}_{G}[-Y|{\Omega_{t}}}}]=0$, then $\hat{\mathbb{E}}[X+Y|\Omega
_{t}]=\hat{\mathbb{E}}[X|\Omega_{t}]$.
\end{proposition}

\begin{proof}
This follows from the following two inequalities:
\begin{align*}
&\mathbb{\hat{E}}%
[X+Y|\Omega_{t}]\leq \hat{\mathbb{E}}[X|\Omega_{t}]+\mathbb{\hat{E}%
}[Y|\Omega_{t}],\\
&\hat{\mathbb{E}}[X+Y|\Omega_{t}]\geq \hat{\mathbb{E}}[X|\Omega
_{t}]-{{\hat{\mathbb{E}}[-Y|{\Omega_{t}}}}]=\hat{\mathbb{E}}%
[X|\Omega_{t}]+\hat{\mathbb{E}}[Y|\Omega_{t}]\text{.}%
\end{align*}
\end{proof}

\begin{example}
For each fixed $\mathbf{a\in}\mathbb{R}^{d},s\leq t$, we have%
\begin{align*}
\hat{\mathbb{E}}[B_{t}^{\mathbf{a}}-B_{s}^{\mathbf{a}}|{\Omega}_{s}]
&
=0,\  \  \  \hat{\mathbb{E}}[-(B_{t}^{\mathbf{a}}-B_{s}^{\mathbf{a}}%
)|\Omega_{s}]=0,\\
\hat{\mathbb{E}}[(B_{t}^{\mathbf{a}}-B_{s}^{\mathbf{a}})^{2}|\Omega_{s}]
&  =\sigma_{\mathbf{aa}^{T}}^{2}(t-s),\  \  \  \hat{\mathbb{E}}%
[-(B_{t}^{\mathbf{a}}-B_{s}^{\mathbf{a}})^{2}|\Omega_{s}]=-\sigma
_{-\mathbf{aa}^{T}}^{2}(t-s),\\
\hat{\mathbb{E}}[(B_{t}^{\mathbf{a}}-B_{s}^{\mathbf{a}})^{4}|\Omega_{s}]
&  =3\sigma_{\mathbf{aa}^{T}}^{4}(t-s)^{2},\  \  \hat{\mathbb{E}}%
[-(B_{t}^{\mathbf{a}}-B_{s}^{\mathbf{a}})^{4}|\Omega_{s}]=-3\sigma
_{-\mathbf{aa}^{T}}^{4}(t-s)^{2},
\end{align*}
where $\sigma_{\mathbf{aa}^{T}}^{2}=2G(\mathbf{aa}^{T})$ and $\sigma
_{-\mathbf{aa}^{T}}^{2}=-2G(-\mathbf{aa}^{T}).$
\end{example}

\begin{example}
\label{eee1}For each $\mathbf{a\in}\mathbb{R}^{d},$ $n\in
\mathbb{N},$
$0\leq t\leq T,$ $X\in L_{G}^{1}(\Omega_{t})$ and $\varphi \in C_{l.Lip}%
(\mathbb{R)}$, we have%
\begin{align*}
\hat{\mathbb{E}}[X\varphi(B_{T}^{\mathbf{a}}-B_{t}^{\mathbf{a}}%
)|\Omega_{t}]  &  =X^{+}\hat{\mathbb{E}}[\varphi(B_{T}^{\mathbf{a}}%
-B_{t}^{\mathbf{a}})|\Omega_{t}]+X^{-}\hat{\mathbb{E}}[-\varphi
(B_{T}^{\mathbf{a}}-B_{t}^{\mathbf{a}})|\Omega_{t}]\\
&  =X^{+}\hat{\mathbb{E}}[\varphi(B_{T}^{\mathbf{a}}-B_{t}^{\mathbf{a}%
})]+X^{-}\hat{\mathbb{E}}[-\varphi(B_{T}^{\mathbf{a}}-B_{t}^{\mathbf{a}%
})].
\end{align*}
In particular, we have%
\[
\hat{\mathbb{E}}[X(B_{T}^{\mathbf{a}}-B_{t}^{\mathbf{a}})|{\Omega}%
_{t}]=X^{+}\hat{\mathbb{E}}[(B_{T}^{\mathbf{a}}-B_{t}^{\mathbf{a}}%
)]+X^{-}\hat{\mathbb{E}}[-(B_{T}^{\mathbf{a}}-B_{t}^{\mathbf{a}})]=0.
\]
This, together with Proposition \ref{d-E-x+y}, yields%
\[
\hat{\mathbb{E}}[Y+X(B_{T}^{\mathbf{a}}-B_{t}^{\mathbf{a}})|\Omega
_{t}]=\hat{\mathbb{E}}[Y|\Omega_{t}],\  \ Y\in L_{G}^{1}(\Omega).
\]
We also have
\begin{align*}
\hat{\mathbb{E}}[X(B_{T}^{\mathbf{a}}-B_{t}^{\mathbf{a}})^{2}|\Omega_{t}]
&  =X^{+}\hat{\mathbb{E}}[(B_{T}^{\mathbf{a}}-B_{t}^{\mathbf{a}}%
)^{2}]+X^{-}\hat{\mathbb{E}}[-(B_{T}^{\mathbf{a}}-B_{t}^{\mathbf{a}}%
)^{2}]\\
&  =[X^{+}\sigma_{\mathbf{aa}^{T}}^{2}-X^{-}\sigma_{-\mathbf{aa}^{T}}%
^{2}](T-t)
\end{align*}
and%
\begin{align*}
\hat{\mathbb{E}}[X(B_{T}^{\mathbf{a}}-B_{t}^{\mathbf{a}})^{2n-1}%
|\Omega_{t}]  &  =X^{+}\hat{\mathbb{E}}[(B_{T}^{\mathbf{a}}-B_{t}%
^{\mathbf{a}})^{2n-1}]+X^{-}\hat{\mathbb{E}}[-(B_{T}^{\mathbf{a}}%
-B_{t}^{\mathbf{a}})^{2n-1}]\\
&  =|X|\hat{\mathbb{E}}[(B_{T-t}^{\mathbf{a}})^{2n-1}].
\end{align*}

\end{example}

\begin{example}
\label{Exam-B2}Since
\[
\hat{\mathbb{E}}[2B_{s}^{\mathbf{a}}(B_{t}^{\mathbf{a}}-B_{s}^{\mathbf{a}%
})|\Omega_{s}]=\hat{\mathbb{E}}[-2B_{s}^{\mathbf{a}}(B_{t}^{\mathbf{a}%
}-B_{s}^{\mathbf{a}})|\Omega_{s}]=0,
\]
we have
\begin{align*}
\hat{\mathbb{E}}[(B_{t}^{\mathbf{a}})^{2}-(B_{s}^{\mathbf{a}})^{2}%
|\Omega_{s}]  &  =\hat{\mathbb{E}}[(B_{t}^{\mathbf{a}}-B_{s}^{\mathbf{a}%
}+B_{s}^{\mathbf{a}})^{2}-(B_{s}^{\mathbf{a}})^{2}|\Omega_{s}]\\
&  =\hat{\mathbb{E}}[(B_{t}^{\mathbf{a}}-B_{s}^{\mathbf{a}})^{2}%
+2(B_{t}^{\mathbf{a}}-B_{s}^{\mathbf{a}})B_{s}^{\mathbf{a}}|{\Omega}_{s}]\\
&  =\sigma_{\mathbf{aa}^{T}}^{2}(t-s).
\end{align*}
\end{example}

\begin{exercise}
Show that if $X\in Lip(\Omega_T)$ and
$\hat{\mathbb{E}}[X]=-\hat{\mathbb{E}}[-X]$, then
$\hat{\mathbb{E}}[X]={E}_P[X]$, where $P$ is a Wiener measure on
$\Omega$.
\end{exercise}

\begin{exercise}
\label{d-Exm-GBM-14a copy(1)} For each $s,t\geq0$, we set $B_{t}^{s}%
:=B_{t+s}-B_{s}$. Let $\eta=(\eta_{ij})_{i,j=1}^{d}\in
L_{G}^{1}(\Omega _{s};\mathbb{S}(d))$. Prove that
\begin{equation*}
\hat{\mathbb{E}}[\langle\eta
B_{t}^{s},B_{t}^{s}\rangle|\Omega_{s}]=2G(\eta)t.
\end{equation*}

\end{exercise}

\section{It{{{\^{o}'s Integral with $G$--Brownian Motion}}}}

\begin{definition}
\label{d-Def-4}For $T\in \mathbb{R}^{+},$ a partition $\pi_{T}$ of
$[0,T]$ is a finite ordered subset
$\pi_{T}=\{t_{0},t_{1},\cdots,t_{N}\}$ such that
$0=t_{0}<t_{1}<\cdots<t_{N}=T$.
\[
\mu(\pi_{T}):=\max \{|t_{i+1}-t_{i}|\ :\ i=0,1,\cdots,N-1\} \text{.}%
\]
We use $\pi_{T}^{N}=\{t_{0}^{N},t_{1}^{N},\cdots,t_{N}^{N}\}$ to
denote a sequence of partitions of $[0,T]$ such that
$\lim_{N\rightarrow \infty}\mu (\pi_{T}^{N})=0$.
\end{definition}

Let $p\geq1$ be fixed. We consider the following type of simple
processes: for
a given partition $\pi_{T}=\{t_{0},\cdots,t_{N}\}$ of $[0,T]$ we set%
\[
\eta_{t}(\omega)=\sum_{k=0}^{N-1}\xi_{k}(\omega)\mathbf{I}_{[t_{k},t_{k+1}%
)}(t),
\]
where $\xi_{k}\in L_{G}^{p}(\Omega_{t_{k}})$, $k=0,1,2,\cdots,N-1$
are given. The collection of these processes is denoted by
$M_{G}^{p,0}(0,T)$\label{mp0}.

\begin{definition}
\label{d-Def-5}For an $\eta \in M_{G}^{p,0}(0,T)$ with
$\eta_{t}(\omega)=\sum
_{k=0}^{N-1}\xi_{k}(\omega)\mathbf{I}_{[t_{k},t_{k+1})}(t)$, the
related \textbf{Bochner integral}\index{Bochner integral} is
\[
\int_{0}^{T}\eta_{t}(\omega)dt:=\sum_{k=0}^{N-1}\xi_{k}(\omega)(t_{k+1}%
-t_{k}).
\]

\end{definition}

For each $\eta \in M_{G}^{p,0}(0,T)$, we set
\[
\mathbb{\widetilde
{E}}_{T}[\eta]:=\frac{1}{T}\hat{\mathbb{E}}[\int_{0}^{T}\eta_{t}dt]
=\frac{1}{T}\hat{\mathbb{E}}[\sum_{k=0}^{N-1}\xi_{k}
(\omega)(t_{k+1}-t_{k})].
\]
It is easy to check that $\mathbb{\widetilde
{E}}_{T}:M_{G}^{p,0}(0,T)\rightarrow \mathbb{R}$ forms a sublinear
expectation. We then can introduce a natural norm $\Vert \eta \Vert
_{M_{G}^{p}(0,T)}$, under which, $M_{G}^{p,0}(0,T)$ can be extended
to $M_{G}^{p}(0,T)$ which is a Banach space.

\begin{definition}
For each $p\geq1$, we denote by $M_{G}^{p}(0,T)$ the completion of
$M_{G}\label{mgp}
^{p,0}(0,T)$ under the norm%
\[
{{\left \Vert \eta \right \Vert _{M_{G}^{p}(0,T)}:=\left \{ \hat{\mathbb{E}}[ {{\int_{0}%
^{T}|\eta_{t}|^{p}dt]}}\right \}  ^{1/p}}}.
\]

\end{definition}
It is clear that $M_{G}^{p}(0,T)\supset M_{G}^{q}(0,T)$ for $1\leq
p\leq q.$ We also use $M_G^p(0,T;\mathbb{R}^n)$ for all
$n$-dimensional stochastic processes
$\eta_t=(\eta_t^1,\cdots,\eta_t^n)$, $t\geq 0$ with $\eta_t^i\in
M_G^p(0,T)$, $i=1,2,\cdots,n$.

We now give the definition of It{{{\^{o}{{'s}} integral. For
simplicity, we first introduce It{{{\^{o}{{'s}} integral with
respect to }}}$1$-dimensional{{{ $G$--Brownian motion.}}}

Let $(B_{t})_{t\geq0}$ be a $1$-dimensional {{{$G$--Brownian motion with }}%
}$G(\alpha)=\frac{1}{2}(\bar{\sigma}^{2}\alpha^{+}-\underline{\sigma}%
^{2}\alpha^{-})$, where $0\leq \underline{\sigma}\leq
\bar{\sigma}<\infty.$

\begin{definition}
For each $\eta \in M_{G}^{2,0}(0,T)$ of the form
\[
\eta_{t}(\omega)=\sum_{j=0}^{N-1}\xi_{j}(\omega)\mathbf{I}_{[t_{j},t_{j+1}%
)}(t),
\]
we define
\[
I(\eta)=\int_{0}^{T}\eta_tdB_{t}:=\sum_{j=0}^{N-1}\xi_{j}(B_{t_{j+1}%
}-B_{t_{j}})\mathbf{.}%
\]

\end{definition}

\begin{lemma}
{ { { \label{bdd}The mapping $I:M_{G}^{2,0}(0,T)\rightarrow L_{G}^{2}%
(\Omega_{T})$ is a continuous linear mapping and thus can be
continuously extended to $I:M_{G}^{2}(0,T)\rightarrow
L_{G}^{2}(\Omega_{T})$. We have
\begin{align}
\hat{\mathbb{E}}[\int_{0}^{T}\eta_tdB_{t}]  &  =0,\  \  \label{e1}\\
\hat{\mathbb{E}}[(\int_{0}^{T}\eta_tdB_{t})^{2}]  &  \leq \bar{\sigma}%
^{2}\hat{\mathbb{E}}[\int_{0}^{T}\eta_t^{2}dt]. \label{e2}%
\end{align}
} } }
\end{lemma}

\begin{proof}
From Example \ref{eee1}, for each $j$,
\[
\hat{\mathbb{E}}\mathbf{[}\xi_{j}(B_{t_{j+1}}-B_{t_{j}})|{\Omega}_{t_{j}%
}]=\hat{\mathbb{E}}\mathbf{[-}\xi_{j}(B_{t_{j+1}}-B_{t_{j}})|\Omega
_{t_{j}}]=0.
\]
We have
\begin{align*}
\hat{\mathbb{E}}[\int_{0}^{T}\eta_tdB_{t}]  &  =\mathbb{\hat{E}}%
\mathbb{[}\int_{0}^{t_{N-1}}\eta_tdB_{t}+\xi_{N-1}(B_{t_{N}}-B_{t_{N-1}%
})]\\
&  =\hat{\mathbb{E}}\mathbb{[}\int_{0}^{t_{N-1}}\eta_tdB_{t}%
+\hat{\mathbb{E}}\mathbf{[}\xi_{N-1}(B_{t_{N}}-B_{t_{N-1}})|\Omega
_{t_{N-1}}]]\\
&  =\hat{\mathbb{E}}\mathbb{[}\int_{0}^{t_{N-1}}\eta_tdB_{t}].
\end{align*}
Then we can repeat this procedure to obtain (\ref{e1}).

We now give the proof of (\ref{e2}). Firstly, from Example
\ref{eee1}, we have
\begin{align*}
\hat{\mathbb{E}}[(\int_{0}^{T}\eta_tdB_{t})^{2}]  &  =\mathbb{\hat{E}%
}\mathbb{[}\left(  \int_{0}^{t_{N-1}}\eta_tdB_{t}+\xi_{N-1}(B_{t_{N}%
}-B_{t_{N-1}})\right)  ^{2}]\\
&=\hat{\mathbb{E}}\mathbb{[}\left(  \int_{0}^{t_{N-1}}\eta
_tdB_{t}\right)  ^{2}+\xi_{N-1}^{2}(B_{t_{N}}-B_{t_{N-1}})^{2}\\
&\ \ \ +2\left(  \int_{0}^{t_{N-1}}%
\eta_tdB_{t}\right)  \xi_{N-1}(B_{t_{N}}-B_{t_{N-1}})
]\\
&  =\hat{\mathbb{E}}\mathbb{[}\left(  \int_{0}^{t_{N-1}}\eta_tdB_{t}\right) ^{2}+\xi_{N-1}^{2}(B_{t_{N}}-B_{t_{N-1}})^{2}]\\
&=\cdots=\hat{\mathbb{E}}[\sum_{i=0}^{N-1}\xi_i^2(B_{t_{i+1}}-B_{t_{i}})^{2}].
\end{align*}
Then, for each $i=0,1,\cdots,N-1$, we have
\begin{align*}
&\hat{\mathbb{E}}[\xi_i^2(B_{t_{i+1}}-B_{t_i})^2-\overline{\sigma}^2\xi_i^2(t_{i+1}-t_i)]\\
=&\hat{\mathbb{E}}[\hat{\mathbb{E}}[\xi_i^2(B_{t_{i+1}}-B_{t_i})^2-\overline{\sigma}^2\xi_i^2(t_{i+1}-t_j)|\Omega_{t_i}]]\\
=&\hat{\mathbb{E}}[\overline{\sigma}^2\xi_i^2(t_{i+1}-t_i)-\overline{\sigma}^2\xi_i^2(t_{i+1}-t_i)]=0.
\end{align*}
Finally, we have
\begin{align*}
&\hat{\mathbb{E}}[(\int_{0}^{T}\eta_tdB_{t})^{2}]=\hat{\mathbb{E}}[\sum_{i=0}^{N-1}\xi_i^2(B_{t_{i+1}}-B_{t_{i}})^{2}]\\
\leq&
\hat{\mathbb{E}}[\sum_{i=0}^{N-1}\xi_i^2(B_{t_{i+1}}-B_{t_{i}})^{2}-\sum_{i=0}^{N-1}\overline{\sigma}^2\xi_i^2(t_{i+1}-t_i)]
+\hat{\mathbb{E}}[\sum_{i=0}^{N-1}\overline{\sigma}^2\xi_i^2(t_{i+1}-t_i)]\\
\leq&\sum_{i=0}^{N-1}\hat{\mathbb{E}}[\xi_i^2(B_{t_{i+1}}-B_{t_i})^2-\overline{\sigma}^2\xi_i^2(t_{i+1}-t_i)]+\hat{\mathbb{E}}[\sum_{i=0}^{N-1}\overline{\sigma}^2\xi_i^2(t_{i+1}-t_i)]\\
=&\hat{\mathbb{E}}[\sum_{i=0}^{N-1}\overline{\sigma}^2\xi_i^2(t_{i+1}-t_i)]=\bar{\sigma}%
^{2}\hat{\mathbb{E}}[\int_{0}^{T}\eta_t^{2}dt].
\end{align*}
\end{proof}

\begin{definition}
We define, for a fixed $\eta \in M_{G}^{2}(0,T)$, the stochastic
integral
\[
\int_{0}^{T}\eta_tdB_{t}:=I(\eta).
\]
\end{definition}
It is clear that (\ref{e1}) and (\ref{e2}) still hold for $\eta \in
M_{G}^{2}(0,T)$.

We list some main properties of It\^{o}'s integral of $G$--Brownian
motion. We denote, for some $0\leq s\leq t\leq T$,
\[
\int_{s}^{t}\eta_{u}dB_{u}:=\int_{0}^{T}\mathbf{I}_{[s,t]}(u)\eta_{u}dB_{u}.
\]
\begin{proposition}
{ { { \label{Prop-Integ}Let $\eta,\theta \in M_{G}^{2}(0,T)$ and let
$0\leq s\leq r\leq t\leq T$. Then we have\newline \textup{(i)}
$\int_{s}^{t}\eta
_{u}dB_{u}=\int_{s}^{r}\eta_{u}dB_{u}+\int_{r}^{t}\eta_{u}dB_{u}.$%
\newline \textup{(ii)} $\int_{s}^{t}(\alpha \eta_{u}+\theta_{u})dB_{u}%
=\alpha \int_{s}^{t}\eta_{u}dB_{u}+\int_{s}^{t}\theta_{u}dB_{u},\
$if$\  \alpha$ is bounded and in $L_{G}^{1}(\Omega_{s})$.}}}
{{{\textup{(iii)} $\hat{\mathbb{E}}\mathbb{[}X+\int_{r}^{T}\eta_{u}%
dB_{u}|{\Omega}_{s}]=\hat{\mathbb{E}}\mathbb{[}X|\Omega_{s}]$ \ \ \
for $X\in L_{G}^{1}(\Omega)$. } } }
\end{proposition}

We now consider the multi-dimensional case. Let $G(\cdot):\mathbb{S}%
(d)\rightarrow \mathbb{R}$ be a given monotonic and sublinear
function and let $(B_{t})_{t\geq0}$ be a $d$-dimensional
{{{$G$--Brownian motion. For each
fixed }}}$\mathbf{a\in}\mathbb{R}^{d}$, we still use $B_{t}^{\mathbf{a}%
}:=\langle \mathbf{a},B_{t}\rangle$. Then
$(B_{t}^{\mathbf{a}})_{t\geq0}$ is a
$1$-dimensional {{{$G_{\mathbf{a}}$--Brownian motion with }}}$G_{\mathbf{a}%
}(\alpha)=\frac{1}{2}(\sigma_{\mathbf{aa}^{T}}^{2}\alpha^{+}-\sigma
_{-\mathbf{aa}^{T}}^{2}\alpha^{-}),$ where $\sigma_{\mathbf{aa}^{T}}%
^{2}=2G(\mathbf{aa}^{T})$ and $\sigma_{-\mathbf{aa}^{T}}^{2}=-2G(-\mathbf{aa}%
^{T}).$ Similar to $1$-dimensional case, we can define
It{{{\^{o}{{'s}}
integral by}}}%
\[
I(\eta):=\int_{0}^{T}\eta_tdB_{t}^{\mathbf{a}},\  \  \text{for}\  \eta \in M_{G}%
^{2}(0,T).
\]
We still have, for each $\eta \in M_{G}^{2}(0,T)$,
\begin{align*}
\hat{\mathbb{E}}[\int_{0}^{T}\eta_tdB_{t}^{\mathbf{a}}]  &
=0,\  \ \\
\hat{\mathbb{E}}[(\int_{0}^{T}\eta_tdB_{t}^{\mathbf{a}})^{2}] & \leq
\sigma_{\mathbf{aa}^{T}}^{2}\hat{\mathbb{E}}[\int_{0}^{T}\eta_t^2dt]. %
\end{align*}
Furthermore, Proposition \ref{Prop-Integ} still holds for the
integral with respect to $B_t^{\mathbf{a}}$.

\begin{exercise}
Prove that, for a fixed $\eta\in M_G^2(0,T)$,
$$\underline{\sigma}^2\hat{\mathbb{E}}[\int_0^T\eta_t^2dt]\leq\hat{\mathbb{E}}[(\int_0^T\eta_tdB_t)^2]\leq\overline{\sigma}^2\hat{\mathbb{E}}[\int_0^T\eta_t^2dt],$$
where $\overline{\sigma}^2=\hat{\mathbb{E}}[B_1^2]$ and
$\underline{\sigma}^2=-\hat{\mathbb{E}}[-B_1^2]$.

\end{exercise}

\begin{exercise}
Prove that, for each $\eta\in M_G^p(0,T)$, we have
$$\hat{\mathbb{E}}[\int_0^T|\eta_t|^p dt]\leq \int_0^T\hat{\mathbb{E}}[|\eta_t|^p]dt.$$
\end{exercise}

\section{{{{Quadratic Variation Process of $G$--Brownian Motion}}}}

We first consider the quadratic variation process of $1$-dimensional
{{{$G$--Brownian motion }}}$(B_{t})_{t\geq0}$ with $B_{1}\overset{d}%
{=}N(\{0\}\times[\underline{\sigma}^{2},\bar{\sigma}^{2}])$. {{{Let
$\pi_{t}^{N}$, $N=1,2,\cdots$, be a sequence of partitions of
$[0,t]$. We consider} } }

{ { {
\begin{align*}
B_{t}^{2} &  =\sum_{j=0}^{N-1}(B_{t_{j+1}^{N}}^{2}-B_{t_{j}^{N}}^{2})\\
&
=\sum_{j=0}^{N-1}2B_{t_{j}^{N}}(B_{t_{j+1}^{N}}-B_{t_{j}^{N}})+\sum
_{j=0}^{N-1}(B_{t_{j+1}^{N}}-B_{t_{j}^{N}})^{2}.
\end{align*}
As $\mu(\pi_{t}^{N})\rightarrow0$, the first term of the right side
converges to $2\int_{0}^{t}B_{s}dB_{s}$ in
}}}$L_{G}^{2}(\Omega)${{{. The second term must be convergent. We
denote its limit by $\left \langle B\right \rangle _{t}$, i.e.,
\begin{equation}
\left \langle B\right \rangle
_{t}:=\lim_{\mu(\pi_{t}^{N})\rightarrow0}\sum
_{j=0}^{N-1}(B_{t_{j+1}^{N}}-B_{t_{j}^{N}})^{2}=B_{t}^{2}-2\int_{0}^{t}%
B_{s}dB_{s}.\label{quadra-def}%
\end{equation}
By the above construction, $( \left \langle B\right \rangle
_{t})_{t\geq0}$ is an increasing process with $\left \langle B\right
\rangle _{0}=0$. We call it the \textbf{quadratic variation
process}\index{Quadratic variation process}\label{qvb} of the
$G$--Brownian motion $B$. It characterizes the part of statistic
uncertainty of $G$--Brownian motion. It is important to keep in mind
that $\left \langle B\right \rangle _{t}$ is not a
deterministic process unless }}}$\underline{\sigma}${{{$=$}}}%
$\bar{\sigma}${{{, i.e., when $(B_t)_{t\geq 0}$ is a classical
Brownian motion. In fact we have the following lemma.} } }

\begin{lemma}
{ { { \label{Lem-Q1}For each $0\leq s\leq t<\infty,$ we have
\begin{align}
\hat{\mathbb{E}}[\left \langle B\right \rangle _{t}-\left \langle
B\right \rangle _{s}|\Omega_{s}] &  =\bar{\sigma}^{2}(t-s),\  \  \label{quadra}\\
\hat{\mathbb{E}}[-(\left \langle B\right \rangle _{t}-\left \langle
B\right \rangle _{s})|\Omega_{s}] &  =-\underline{\sigma}^{2}%
(t-s).\label{quadra1}%
\end{align}
} } }
\end{lemma}

\begin{proof}
{ { { By the definition of $\left \langle B\right \rangle $ and
Proposition \ref{Prop-Integ} (iii),
\begin{align*}
\hat{\mathbb{E}}[\left \langle B\right \rangle _{t}-\left \langle
B\right \rangle _{s}|\Omega_{s}]  &  =\hat{\mathbb{E}}[B_{t}^{2}-B_{s}%
^{2}-2\int_{s}^{t}B_{u}dB_{u}|\Omega_{s}]\\
&  =\hat{\mathbb{E}}[B_{t}^{2}-B_{s}^{2}|\Omega_{s}]=\bar{\sigma}%
^{2}(t-s).
\end{align*}
The last step follows from Example \ref{Exam-B2}. We then have
(\ref{quadra}). The equality (\ref{quadra1}) can be proved
analogously with the
consideration of }}}$\hat{\mathbb{E}}[-(B_{t}^{2}-B_{s}^{2})|\Omega_{s}%
]${{{$=-\underline{\sigma}^{2}(t-s)$. } } }
\end{proof}

{ {{A very interesting point of the quadratic variation process
$\left \langle B\right \rangle $ is, just like the $G$--Brownian
motion $B$ itself, the increment $\left \langle B\right \rangle
_{s+t}-\left \langle B\right \rangle _{s}$ is independent from
}}}$\Omega_{s}${{{ and identically distributed with $\left \langle
B\right \rangle _{t}$. In fact we have } } }

\begin{lemma}
{ { { {{\label{Lem-Qua2}}}For each fixed $s,$}}}$t\geq0$,{ {{$\left
\langle B\right \rangle _{s+t}-\left \langle B\right \rangle
_{s}$}}} is identically
distributed with {{{$\left \langle B\right \rangle _{t}$ and independent from }}%
}$\Omega_{s}${{{. } } }
\end{lemma}

\begin{proof}
{ { { The results follow directly from}}}%
\begin{align*}
\left \langle B\right \rangle _{s+t}-\left \langle B\right \rangle
_{s} &
=B_{s+t}^{2}-2\int_{0}^{s+t}B_{r}dB_{r}-[B_{s}^{2}-2\int_{0}^{s}B_{r}dB_{r}]\\
&  =(B_{s+t}-B_{s})^{2}-2\int_{s}^{s+t}(B_{r}-B_{s})d(B_{r}-B_{s})\\
&  =\left \langle B^{s}\right \rangle _{t},
\end{align*}
{{{where }}}$\left \langle B^{s}\right \rangle $ is the quadratic
variation process of {{{the {{$G$--}}Brownian motion
$B_{t}^{s}=B_{s+t}-B_{s}$, $t\geq 0$. }}}
\end{proof}

We now{ {{define the integral of a process $\eta \in M_{G}^{1}(0,T)$
with
respect to $\left \langle B\right \rangle $. We first define a mapping:%
\[
Q_{0,T}(\eta)=\int_{0}^{T}\eta_td\left \langle B\right \rangle
_{t}:=\sum _{j=0}^{N-1}\xi_{j}(\left \langle B\right \rangle
_{t_{j+1}}-\left \langle B\right \rangle
_{t_{j}}):M_{G}^{1,0}(0,T)\rightarrow L_{G}^{1}(\Omega_{T}).
\]
} } }

\begin{lemma}
{ { { \label{Lem-Q2}For each $\eta \in M_{G}^{1,0}(0,T)$,
\begin{equation}
\hat{\mathbb{E}}[|Q_{0,T}(\eta)|]\leq \bar{\sigma}^{2}\hat{\mathbb{E}}[\int_{0}%
^{T}|\eta_{t}|dt].\  \label{dA}%
\end{equation}
Thus
$Q_{0,T}:M_{G}^{1,0}(0,T)\rightarrow$}}}$L_{G}^{1}(\Omega_{T})${{{
is a continuous linear mapping. Consequently, $Q_{0,T}$ can be
uniquely extended to
}}}$M_{G}^{1}(0,T)${{{. We still denote this mapping \ by%
\[
\int_{0}^{T}\eta_td\left \langle B\right \rangle _{t}:=Q_{0,T}(\eta
)\  \  \text{for}\ \eta \in M_{G}^{1}(0,T)\text{.}%
\]
We still have
\begin{equation}
\hat{\mathbb{E}}[|\int_{0}^{T}\eta_td\left \langle B\right \rangle
_{t}|]\leq \bar{\sigma}^{2}\hat{\mathbb{E}}[\int_{0}^{T}|\eta
_{t}|dt]\  \  \text{for}\ \eta \in M_{G}^{1}(0,T)\text{.}\label{qua-ine}%
\end{equation}
} } }
\end{lemma}

\begin{proof}
Firstly, for each $j=1,\cdots,N-1$, we have
\begin{align*}
&\hat{\mathbb{E}}[|\xi_j|(\langle B\rangle_{t_{j+1}}-\langle
B\rangle_{t_j})-\overline{\sigma}^2|\xi_j|(t_{j+1}-t_j)]\\=&\hat{\mathbb{E}}[\hat{\mathbb{E}}[|\xi_j|(\langle
B\rangle_{t_{j+1}}-\langle B\rangle_{t_j})|\Omega_{t_j}]-\overline{\sigma}^2|\xi_j|(t_{j+1}-t_j)]\\
=&\hat{\mathbb{E}}[|\xi_j|\overline{\sigma}^2(t_{j+1}-t_j)-\overline{\sigma}^2|\xi_j|(t_{j+1}-t_j)]=0.
\end{align*}
Then (\ref{dA}) can be checked as follows:
\begin{align*}
&\hat{\mathbb{E}}[|\sum_{j=0}^{N-1}\xi_{j}(\left \langle B\right
\rangle _{t_{j+1}}-\left \langle B\right \rangle
_{t_{j}})|]\leq\hat{\mathbb{E}}[\sum_{j=0} ^{N-1}|\xi_{j}|\left
\langle B\right \rangle _{t_{j+1}}-\left \langle B\right \rangle
_{t_{j}}]\\
\leq &\hat{\mathbb{E}}[\sum_{j=0} ^{N-1}|\xi_{j}|[(\langle B\rangle
_{t_{j+1}}-\langle B\rangle
_{t_{j}})-\overline{\sigma}^2(t_{j+1}-t_j)]]+\hat{\mathbb{E}}[\overline{\sigma}^2\sum_{j=0}^{N-1}|\xi_j|(t_{j+1}-t_j)]\\
\leq& \sum_{j=0}^{N-1}\hat{\mathbb{E}}[|\xi_j|[(\langle
B\rangle_{t_{j+1}}-\langle
B\rangle_{t_j})-\overline{\sigma}^2(t_{j+1}-t_j)]]+\hat{\mathbb{E}}[\overline{\sigma}^2\sum_{j=0}^{N-1}|\xi_j|(t_{j+1}-t_j)]\\
=&\hat{\mathbb{E}}[\overline{\sigma}^2\sum_{j=0}^{N-1}|\xi_j|(t_{j+1}-t_j)]=\overline{\sigma}^2\hat{\mathbb{E}}[\int_{0}%
^{T}|\eta_{t}|dt].
\end{align*}

\end{proof}

\begin{proposition}
{ { { \label{Prop-temp}}}Let $0\leq s\leq t$, $\xi \in
L_{G}^{2}(\Omega_{s})$,
}$X\in${$L_{G}^{1}(\Omega)$. Then%
\begin{align*}
\hat{\mathbb{E}}[X+\xi(B_{t}^{2}-B_{s}^{2})]  &  =\mathbb{\hat{E}}%
[X+\xi(B_{t}-B_{s})^{2}]\\
&  =\hat{\mathbb{E}}[X+\xi(\left \langle B\right \rangle _{t}-\left
\langle B\right \rangle _{s})].
\end{align*}
}
\end{proposition}

\begin{proof}
{ { { By (\ref{quadra-def}) and Proposition \ref{Prop-Integ} (iii), we have%
\begin{align*}
\hat{\mathbb{E}}[X+\xi(B_{t}^{2}-B_{s}^{2})] &  =\mathbb{\hat{E}}%
[X+\xi(\left \langle B\right \rangle _{t}-\left \langle B\right
\rangle
_{s}+2\int_{s}^{t}B_{u}dB_{u})]\\
&  =\hat{\mathbb{E}}[X+\xi(\left \langle B\right \rangle _{t}-\left
\langle B\right \rangle _{s})].
\end{align*}
We also have
\begin{align*}
\hat{\mathbb{E}}[X+\xi(B_{t}^{2}-B_{s}^{2})] &  =\mathbb{\hat{E}}%
[X+\xi ((B_{t}-B_{s})^{2}+2(B_{t}-B_{s})B_{s})]\\
&  =\hat{\mathbb{E}}[X+\xi(B_{t}-B_{s})^{2}].
\end{align*}
} } }
\end{proof}

{ { { We have the following isometry. } } }

\begin{proposition}
{ { { Let $\eta \in M_{G}^{2}(0,T)$. Then
\begin{equation}
\hat{\mathbb{E}}[(\int_{0}^{T}\eta_tdB_{t})^{2}]=\mathbb{\hat{E}}%
[\int_{0}^{T}\eta^{2}_td\left \langle B\right \rangle _{t}].
\label{isometry}%
\end{equation}
} } }
\end{proposition}

\begin{proof}
{ { { We first consider $\eta \in M_{G}^{2,0}(0,T)$ of the form
\[
\eta_{t}(\omega)=\sum_{j=0}^{N-1}\xi_{j}(\omega)\mathbf{I}_{[t_{j},t_{j+1}%
)}(t)
\]
and then $\int_{0}^{T}\eta_tdB_{t}=\sum_{j=0}^{N-1}\xi_{j}(B_{t_{j+1}%
}-B_{t_{j}})$\textbf{.} From Proposition \ref{Prop-Integ}, we get
\[
\hat{\mathbb{E}}[X+2\xi_{j}(B_{t_{j+1}}-B_{t_{j}})\xi_{i}(B_{t_{i+1}%
}-B_{t_{i}})]=\hat{\mathbb{E}}[X]\ \text{ for }X\in{L_{G}^{1}(\Omega
)}\text{, }i\not =j.
\]
Thus%
\[
\hat{\mathbb{E}}[(\int_{0}^{T}\eta_tdB_{t})^{2}]=\mathbb{\hat{E}}%
\mathbb{[}( \sum_{j=0}^{N-1}\xi_{j}(B_{t_{j+1}}-B_{t_{j}}))
^{2}]=\hat{\mathbb{E}}\mathbb{[}\sum_{j=0}^{N-1}\xi_{j}^{2}(B_{t_{j+1}%
}-B_{t_{j}})^{2}].
\]
From this and Proposition \ref{Prop-temp}, it follows that
\[
\hat{\mathbb{E}}[(\int_{0}^{T}\eta_tdB_{t})^{2}]=\mathbb{\hat{E}}%
\mathbb{[}\sum_{j=0}^{N-1}\xi_{j}^{2}(\left \langle B\right \rangle
_{t_{j+1}}-\left \langle B\right \rangle _{t_{j}})]=\mathbb{\hat{E}}%
\mathbb{[}\int_{0}^{T}\eta^{2}_td\left \langle B\right \rangle
_{t}].
\]
Thus (\ref{isometry}) holds for $\eta \in M_{G}^{2,0}(0,T)$. We can
continuously extend the above equality to the case $\eta \in
M_{G}^{2}(0,T)$ and get (\ref{isometry}). }}}
\end{proof}

We now consider the multi-dimensional case. Let $(B_{t})_{t\geq0}$
be a $d$-dimensional {{{$G$--Brownian motion. For each fixed
}}}$\mathbf{a\in }\mathbb{R}^{d}$, $(B_{t}^{\mathbf{a}})_{t\geq0}$
is a $1$-dimensional {{{$G_{\mathbf{a}}$--Brownian motion}}}.
Similar to $1$-dimensional case, we
can define%
\begin{equation*}
\left \langle B^{\mathbf{a}}\right \rangle _{t}:=\lim_{\mu(\pi_{t}^{N}%
)\rightarrow0}\sum_{j=0}^{N-1}(B_{t_{j+1}^{N}}^{\mathbf{a}}-B_{t_{j}^{N}%
}^{\mathbf{a}})^{2}=(B_{t}^{\mathbf{a}})^{2}-2\int_{0}^{t}B_{s}^{\mathbf{a}%
}dB_{s}^{\mathbf{a}},
\end{equation*}
where $\left \langle B^{\mathbf{a}}\right \rangle $ is called the
\textbf{quadratic variation process} of $B^{\mathbf{a}}$. The above
results also hold for $\left \langle B^{\mathbf{a}}\right \rangle $.
In particular,
\begin{equation*}
\hat{\mathbb{E}}[|\int_{0}^{T}\eta_td\left \langle B^{\mathbf{a}%
}\right \rangle _{t}|]\leq \sigma_{\mathbf{aa}^{T}}^{2}\hat{\mathbb{E}}[\int_{0}^{T}%
|\eta_{t}|dt]\  \  \text{for}\ \eta \in M_{G}^{1}(0,T)
\end{equation*}
and%
\begin{equation*}
\hat{\mathbb{E}}[(\int_{0}^{T}\eta_tdB_{t}^{\mathbf{a}})^{2}%
]=\hat{\mathbb{E}}[\int_{0}^{T}\eta^{2}_td\left \langle B^{\mathbf{a}%
}\right \rangle _{t}]\  \  \text{for}\ \eta \in M_{G}^{2}(0,T).
\end{equation*}

Let $\mathbf{a}=(a_{1},\cdots,a_{d})^{T}$ and
$\mathbf{\bar{a}}=(\bar{a}_{1},\cdots,\bar{a}_{d})^{T}$ be two given
vectors in $\mathbb{R}^{d}$. We then have their quadratic variation
processes of $\left \langle B^{\mathbf{a}}\right \rangle $ and
$\left \langle B^{\mathbf{\bar
{a}}}\right \rangle $. We can define their \textbf{mutual variation process} \index{Mutual variation process}by%
\begin{align*}
\left \langle B^{\mathbf{a}},B^{\mathbf{\bar{a}}}\right \rangle _{t}
& :=\frac{1}{4}[\left \langle
B^{\mathbf{a}}+B^{\mathbf{\bar{a}}}\right \rangle
_{t}-\left \langle B^{\mathbf{a}}-B^{\mathbf{\bar{a}}}\right \rangle _{t}]\\
&  =\frac{1}{4}[\left \langle B^{\mathbf{a}+\mathbf{\bar{a}}}\right
\rangle _{t}-\left \langle B^{\mathbf{a}-\mathbf{\bar{a}}}\right
\rangle _{t}].
\end{align*}
Since $\left \langle B^{\mathbf{a}-\mathbf{\bar{a}}}\right \rangle
=\left \langle
B^{\mathbf{\bar{a}}-\mathbf{a}}\right \rangle =\left \langle -B^{\mathbf{a}%
-\mathbf{\bar{a}}}\right \rangle $, we see that $\left \langle B^{\mathbf{a}%
},B^{\mathbf{\bar{a}}}\right \rangle _{t}=\left \langle B^{\mathbf{\bar{a}}%
},B^{\mathbf{a}}\right \rangle _{t}$. In particular, we have $\left
\langle
B^{\mathbf{a}},B^{\mathbf{a}}\right \rangle =\left \langle B^{\mathbf{a}%
}\right \rangle $. Let $\pi_{t}^{N}$, $N=1,2,\cdots$, be a sequence
of
partitions of $[0,t]$. We observe that%
\[
\sum_{k=0}^{N-1}(B_{t_{k+1}^{N}}^{\mathbf{a}}-B_{t_{k}^{N}}^{\mathbf{a}%
})(B_{t_{k+1}^{N}}^{\mathbf{\bar{a}}}-B_{t_{k}^{N}}^{\mathbf{\bar{a}}}%
)=\frac{1}{4}\sum_{k=0}^{N-1}[(B_{t_{k+1}}^{\mathbf{a}+\mathbf{\bar{a}}%
}-B_{t_{k}}^{\mathbf{a}+\mathbf{\bar{a}}})^{2}-(B_{t_{k+1}}^{\mathbf{a}%
-\mathbf{\bar{a}}}-B_{t_{k}}^{\mathbf{a}-\mathbf{\bar{a}}})^{2}].
\]
Thus as $\mu(\pi_{t}^{N})\rightarrow0$ we have%
\[
\lim_{N\rightarrow \infty}\sum_{k=0}^{N-1}(B_{t_{k+1}^{N}}^{\mathbf{a}%
}-B_{t_{k}^{N}}^{\mathbf{a}})(B_{t_{k+1}^{N}}^{\mathbf{\bar{a}}}-B_{t_{k}^{N}%
}^{\mathbf{\bar{a}}})=\left \langle B^{\mathbf{a}},B^{\mathbf{\bar{a}}%
}\right \rangle _{t}.
\]
We also have
\begin{align*}
\left \langle B^{\mathbf{a}},B^{\mathbf{\bar{a}}}\right \rangle _{t}
& =\frac{1}{4}[\left \langle B^{\mathbf{a}+\mathbf{\bar{a}}}\right
\rangle
_{t}-\left \langle B^{\mathbf{a}-\mathbf{\bar{a}}}\right \rangle _{t}]\\
&  =\frac{1}{4}[(B_{t}^{\mathbf{a}+\mathbf{\bar{a}}})^{2}-2\int_{0}^{t}%
B_{s}^{\mathbf{a}+\mathbf{\bar{a}}}dB_{s}^{\mathbf{a}+\mathbf{\bar{a}}}%
-(B_{t}^{\mathbf{a}-\mathbf{\bar{a}}})^{2}+2\int_{0}^{t}B_{s}^{\mathbf{a}%
-\mathbf{\bar{a}}}dB_{s}^{\mathbf{a}-\mathbf{\bar{a}}}]\\
&  =B_{t}^{\mathbf{a}}B_{t}^{\mathbf{\bar{a}}}-\int_{0}^{t}B_{s}^{\mathbf{a}%
}dB_{s}^{\mathbf{\bar{a}}}-\int_{0}^{t}B_{s}^{\mathbf{\bar{a}}}dB_{s}%
^{\mathbf{a}}.
\end{align*}
Now for each $\eta \in M_{G}^{1}(0,T)$, we can consistently define
\[
\int_{0}^{T}\eta_{t}d\left \langle B^{\mathbf{a}},B^{\mathbf{\bar{a}}%
}\right \rangle _{t}=\frac{1}{4}\int_{0}^{T}\eta_{t}d\left \langle
B^{\mathbf{a}+\mathbf{\bar{a}}}\right \rangle _{t}-\frac{1}{4}\int_{0}^{T}%
\eta_{t}d\left \langle B^{\mathbf{a}-\mathbf{\bar{a}}}\right \rangle
_{t}.
\]

\begin{lemma}
\label{d-Lem-mutual}Let $\eta^{N}\in M_{G}^{2,0}(0,T)$,
$N=1,2,\cdots,$ be of the form
\[
\eta_{t}^{N}(\omega)=\sum_{k=0}^{N-1}\xi_{k}^{N}(\omega)\mathbf{I}_{[t_{k}%
^{N},t_{k+1}^{N})}(t)
\]
with $\mu(\pi_{T}^{N})\rightarrow0$ and $\eta^{N}\rightarrow \eta$
in $M_{G}^{2}(0,T)$, as $N\rightarrow \infty$. Then we have the
following
convergence in $L_{G}^{2}(\Omega_{T})$:%
\[
\sum_{k=0}^{N-1}\xi_{k}^{N}(B_{t_{k+1}^{N}}^{\mathbf{a}}-B_{t_{k}^{N}%
}^{\mathbf{a}})(B_{t_{k+1}^{N}}^{\mathbf{\bar{a}}}-B_{t_{k}^{N}}%
^{\mathbf{\bar{a}}})\rightarrow \int_{0}^{T}\eta_td\left \langle B^{\mathbf{a}%
},B^{\mathbf{\bar{a}}}\right \rangle _{t}.
\]

\end{lemma}

\begin{proof}
Since
\begin{align*}
\left \langle B^{\mathbf{a}},B^{\mathbf{\bar{a}}}\right \rangle _{t_{k+1}^{N}%
}-\left \langle B^{\mathbf{a}},B^{\mathbf{\bar{a}}}\right \rangle
_{t_{k}^{N}}
&  =(B_{t_{k+1}^{N}}^{\mathbf{a}}-B_{t_{k}^{N}}^{\mathbf{a}})(B_{t_{k+1}^{N}%
}^{\mathbf{\bar{a}}}-B_{t_{k}^{N}}^{\mathbf{\bar{a}}})\\
&\ \ \ \ -\int_{t_{k}^{N}}^{t_{k+1}^{N}}(B_{s}^{\mathbf{a}}-B_{t_{k}^{N}%
}^{\mathbf{a}})dB_{s}^{\mathbf{\bar{a}}}-\int_{t_{k}^{N}}^{t_{k+1}^{N}}%
(B_{s}^{\mathbf{\bar{a}}}-B_{t_{k}^{N}}^{\mathbf{\bar{a}}})dB_{s}^{\mathbf{a}%
},
\end{align*}
we only need to prove%
\[
\hat{\mathbb{E}}[\sum_{k=0}^{N-1}(\xi_{k}^{N})^{2}(\int_{t_{k}^{N}%
}^{t_{k+1}^{N}}(B_{s}^{\mathbf{a}}-B_{t_{k}^{N}}^{\mathbf{a}})dB_{s}%
^{\mathbf{\bar{a}}})^{2}]\rightarrow0.
\]
For each $k=1,\cdots,N-1$, we have
\begin{align*}
&\hat{\mathbb{E}}[(\xi_{k}^{N})^{2}(\int_{t_{k}^{N}
}^{t_{k+1}^{N}}(B_{s}^{\mathbf{a}}-B_{t_{k}^{N}}^{\mathbf{a}})dB_{s}
^{\mathbf{\bar{a}}})^{2}-C(\xi_{k}^{N})^{2}(t^N_{k+1}-t^N_k)^2]\\
=&\hat{\mathbb{E}}[\hat{\mathbb{E}}[(\xi_{k}^{N})^{2}(\int_{t_{k}^{N}
}^{t_{k+1}^{N}}(B_{s}^{\mathbf{a}}-B_{t_{k}^{N}}^{\mathbf{a}})dB_{s}
^{\mathbf{\bar{a}}})^{2}|\Omega_{t_k^N}]-C(\xi_{k}^{N})^{2}(t^N_{k+1}-t^N_k)^2]\\
\leq
&\hat{\mathbb{E}}[C(\xi_k^N)^2(t^N_{k+1}-t^N_k)^2-C(\xi_{k}^{N})^{2}(t^N_{k+1}-t^N_k)^2]=0,
\end{align*}
where
$C=\bar{\sigma}_{\mathbf{aa}^{T}}^{2}\bar{\sigma}_{\mathbf{\bar
{a}\bar{a}}^{T}}^{2}/{2}$.

Thus we have
\begin{align*}
&\hat{\mathbb{E}}[\sum_{k=0}^{N-1}(\xi_{k}^{N})^{2}(\int_{t_{k}^{N}%
}^{t_{k+1}^{N}}(B_{s}^{\mathbf{a}}-B_{t_{k}^{N}}^{\mathbf{a}})dB_{s}%
^{\mathbf{\bar{a}}})^{2}]\\
\leq&
\hat{\mathbb{E}}[\sum_{k=0}^{N-1}(\xi_{k}^{N})^{2}[(\int_{t_{k}^{N}
}^{t_{k+1}^{N}}(B_{s}^{\mathbf{a}}-B_{t_{k}^{N}}^{\mathbf{a}})dB_{s}
^{\mathbf{\bar{a}}})^{2}-C(t^N_{k+1}-t^N_k)^2]]\\&+\hat{\mathbb{E}}[\sum_{k=0}^{N-1}C(\xi_{k}^{N})^{2}(t^N_{k+1}-t^N_k)^2]\\
\leq&\sum_{k=0}^{N-1}
\hat{\mathbb{E}}[(\xi_{k}^{N})^{2}[(\int_{t_{k}^{N}
}^{t_{k+1}^{N}}(B_{s}^{\mathbf{a}}-B_{t_{k}^{N}}^{\mathbf{a}})dB_{s}
^{\mathbf{\bar{a}}})^{2}-C(t^N_{k+1}-t^N_k)^2]]\\&+\hat{\mathbb{E}}[\sum_{k=0}^{N-1}C(\xi_{k}^{N})^{2}(t^N_{k+1}-t^N_k)^2]\\
\leq
&\hat{\mathbb{E}}[\sum_{k=0}^{N-1}C(\xi_{k}^{N})^{2}(t^N_{k+1}-t^N_k)^2]
\leq  C\mu(\pi^N_T)\hat{\mathbb{E}}[\int_0^T|\eta^N_t|^2dt],
\end{align*}
As $\mu(\pi_{T}^{N})\rightarrow0$, the proof is complete.
\end{proof}

\begin{exercise}
Let $B_t$ be a 1-dimensional G-Brownian motion and $\varphi$ be a
bounded and Lipschitz function on $\mathbb{R}$. Show that
$$\lim_{N\rightarrow\infty}\hat{\mathbb{E}}[|\sum_{k=0}^{N-1}\varphi(B_{t_{k}^N})[(B_{t_{k+1}^N}-B_{t_{k}^N})^2-(\langle B\rangle_{t_{k+1}^N}-\langle B\rangle_{t_{k}^N})]|]=0,$$
where $t_k^N=kT/{N}, k=0,2,\cdots,N-1$.
\end{exercise}

\begin{exercise}
Prove that, for a fixed $\eta\in M_G^1(0,T)$,
$$\underline{\sigma}^2\hat{\mathbb{E}}[\int_0^T|\eta_t|dt]\leq\hat{\mathbb{E}}[\int_0^T|\eta_t|d\langle B\rangle_t]\leq\overline{\sigma}^2\hat{\mathbb{E}}[\int_0^T|\eta_t|dt],$$
where $\overline{\sigma}^2=\hat{\mathbb{E}}[B_1^2]$ and
$\underline{\sigma}^2=-\hat{\mathbb{E}}[-B_1^2]$.

\end{exercise}

\section{The Distribution of $\left \langle B\right \rangle$}

In this section, we first consider the $1$-dimensional
{{{$G$--Brownian motion
}}}$(B_{t})_{t\geq0}$ with $B_{1}\overset{d}{=}N(\{0\}\times[\underline{\sigma}%
^{2},\bar{\sigma}^{2}])$.

The quadratic variation process $\left \langle B\right \rangle $ of
$G$-Brownian motion $B$ is a very interesting process. We have seen
that the $G$-Brownian motion $B$ is a typical process with variance
uncertainty but without mean-uncertainty. In fact, $\langle
B\rangle$ is concentrated all uncertainty of the $G$-Brownian motion
$B$.  Moreover, $\left \langle B\right \rangle $ itself is a typical
process with mean-uncertainty. This fact will be applied to measure
the mean-uncertainty of risk positions.

\begin{lemma}
{{{We have
\begin{equation}
{{{\mathbb{\hat{E}}}}}{{{[\left \langle B\right \rangle _{t}^{2}]}}}%
\leq10\bar{\sigma}^{4}t^{2}. \label{Qua2}%
\end{equation}
}}}
\end{lemma}

\begin{proof}
Indeed,
\begin{align*}
{{{\mathbb{\hat{E}}}}}{{{[\left \langle B\right \rangle _{t}^{2}]}}}
&
=\hat{\mathbb{E}}[(B_{t}{}^{2}-2\int_{0}^{t}B_{u}dB_{u})^{2}]\\
& \leq2\hat{\mathbb{E}}[B_{t}^{4}]+8\hat{\mathbb{E}}[(\int_{0}
^{t}B_{u}dB_{u})^{2}]\\
& \leq6\bar{\sigma}^{4}t^{2}+8\bar{\sigma}^{2}\mathbb{\hat{E}
}[\int_{0}^{t}B_{u}{}^{2}du]\\
&  \leq6\bar{\sigma}^{4}t^{2}+8\bar{\sigma}^{2}\int_{0}^{t}\mathbb{\hat{E}%
}\mathbb{[}B_{u}{}^{2}]du\\
&  =10\bar{\sigma}^{4}t^{2}.
\end{align*}\end{proof}

\begin{proposition}
\label{c3p1} Let $(b_{t})_{t\geq0}$ be a process on a sublinear
expectation space $(\Omega,\mathcal{H},{{{\mathbb{\hat{E}}}}})$ such
that
\newline \textup{(i)} $b_{0}=0$;\newline \textup{(ii)} For each
$t,s\geq0$, $b${$_{t+s}-b_{t}$ is
identically distributed with }$b${$_{s}$ and independent from $(b_{t_{1}%
},b_{t_{2}},\cdots,b_{t_{n}})$ for each $n\in \mathbb{N}$ and $0\leq
t_{1},\cdots,t_{n}\leq t$;\newline \textup{(iii)}
}$\lim_{t\downarrow 0}{{{\mathbb{\hat{E}}}}}${$[b_{t}^{2}]t^{-1}=0$.
\newline Then }$b_{t}$ is
$N([\underline{\mu}t,\overline{\mu}t]\times\{0\})$-distributed{{
with }}$\overline{\mu
}={{{\mathbb{\hat{E}}}}}${{$[b_{1}]$ and $\underline{\mu}=-{{{\mathbb{\hat{E}%
}}}}[-b_{1}]$. } }
\end{proposition}

\begin{proof}
We first prove that
\[
{{{\mathbb{\hat{E}}}}}[b_{t}]=\overline{\mu}t\text{ \ and
}-{{{\mathbb{\hat {E}}}}}[-b_{t}]=\underline{\mu}t.
\]
We set $\varphi(t):={{{\mathbb{\hat{E}}}}}[b_{t}]$. Then
$\varphi(0)=0$ and
$\lim_{t\downarrow0}\varphi(t)=${$0$.} Since for each $t,s\geq0,$%
\begin{align*}
\varphi(t+s) &  ={{{\mathbb{\hat{E}}}}}[b_{t+s}]={{{\mathbb{\hat{E}}}}%
}[(b_{t+s}-b_{s})+b_{s}]\\
&  =\varphi(t)+\varphi(s).
\end{align*}
Thus $\varphi(t)$ is linear and uniformly continuous in $t,$ which
means that ${{{\mathbb{\hat{E}}}}}[b_{t}]=\overline{\mu}t$.
Similarly $-{{{\mathbb{\hat {E}}}}}[-b_{t}]=\underline{\mu}t$.

We now prove that $b_{t}$ is $N([\underline{\mu}t,\overline{\mu}%
t]\times\{0\})$-distributed. By Exercise \ref{ex2} in
Chap.\ref{ch2}, we just need to prove that for each fixed $\varphi
\in C_{b.Lip}(\mathbb{R)}$, the function
\[
u(t,x):={{{\mathbb{\hat{E}}}}}[\varphi(x+b_{t})],\  \ (t,x)\in
\lbrack
0,\infty)\times \mathbb{R}%
\]
is the viscosity solution of the following {parabolic PDE}:
\begin{equation}
\partial_{t}u-g(\partial_{x}u)=0,\  \  \  \ u|_{t=0}=\varphi \label{G-mean}%
\end{equation}
with $g(a)=\overline{\mu}a^{+}-\underline{\mu}a^{-}$.

We first prove that $u$ is {Lipschitz in $x$ and
$\frac{1}{2}$-H\"{o}lder continuous in $t$}. In fact, {for each
fixed }$t$,{ }$u(t,\cdot)\in${$C_{b.Lip}(\mathbb{R)}$ }since
\begin{align*}
|{{{\mathbb{\hat{E}}}}}[\varphi(x+b_{t})]-{{{\mathbb{\hat{E}}}}}%
[\varphi(y+b_{t})]|  &  \leq{{{\mathbb{\hat{E}}}}}[|\varphi(x+b_{t}%
)-\varphi(y+b_{t})|]\\
&  \leq C|x-y|.
\end{align*}
{For each $\delta \in \lbrack0,t]$, since $b_{t}-b_{\delta}$ is
independent from $b_{\delta}$, we  have
\begin{align*}
u(t,x)  &  ={{{\mathbb{\hat{E}}}}}[\mathbb{\varphi(}x+b_{\delta}%
+(b_{t}-b_{\delta})]\\
&
={{{\mathbb{\hat{E}}}}}[{{{\mathbb{\hat{E}}}}}[\varphi(y+(b_{t}-b_{\delta
}))]_{y=x+b_{\delta}}],
\end{align*}
hence}%
\begin{equation}
u(t,x)={{{\mathbb{\hat{E}}}}}[u(t-\delta,x+b_{\delta})]. \label{eq4.21}%
\end{equation}
{Thus}%
\begin{align*}
|u(t,x)-u(t-\delta,x)|  &
=|{{{\mathbb{\hat{E}}}}}[u(t-\delta,x+b_{\delta
})-u(t-\delta,x)]|\\
&  \leq{{{\mathbb{\hat{E}}}}}[|u(t-\delta,x+b_{\delta})-u(t-\delta,x)|]\\
&  \leq{{{\mathbb{\hat{E}}}}}[C|b_{\delta}|]\leq C_{1}\sqrt{\delta}.
\end{align*}

To prove that $u$ is a viscosity solution of the PDE (\ref{G-mean}),
we fix a
point $(t,x)\in(0,\infty)\times \mathbb{R}$ and let $v\in C_{b}^{2,2}%
([0,\infty)\times \mathbb{R})$ be such that $v\geq u$ and
$v(t,x)=u(t,x)$. From (\ref{eq4.21}), we have
\[
v(t,x)={{{\mathbb{\hat{E}}}}}[u(t-\delta,x+b_{\delta})]\leq{{{\mathbb{\hat{E}%
}}}}[v(t-\delta,x+b_{\delta})].
\]
Therefore, by Taylor's expansion,
\begin{align*}
0 &  \leq{{{\mathbb{\hat{E}}}}}[v(t-\delta,x+b_{\delta})-v(t,x)]\\
&  ={{{\mathbb{\hat{E}}}}}[v(t-\delta,x+b_{\delta})-v(t,x+b_{\delta
})+(v(t,x+b_{\delta})-v(t,x))]\\
&  ={{{\mathbb{\hat{E}}}}}[-\partial_{t}v(t,x)\delta+\partial_{x}%
v(t,x)b_{\delta}+I_{\delta}]\\
&  \leq-\partial_{t}v(t,x)\delta+{{{\mathbb{\hat{E}}}}}[\partial
_{x}v(t,x)b_{\delta}]+{{{\mathbb{\hat{E}}}}}[I_{\delta}]\\
&  =-\partial_{t}v(t,x)\delta+g(\partial_{x}v(t,x))\delta+{{{\mathbb{\hat{E}%
}}}}[I_{\delta}],
\end{align*}
where%
\begin{align*}
I_{\delta} &  =\delta \int_{0}^{1}[-\partial_{t}v(t-\beta
\delta,x+b_{\delta
})+\partial_{t}v(t,x)]d\beta \\
& \ \ \ +b_{\delta}\int_{0}^{1}[\partial_{x}v(t,x+\beta
b_{\delta})-\partial_{x}v(t,x)]d\beta .
\end{align*}
With the assumption that {$\lim_{t\downarrow0}$}${{{\mathbb{\hat{E}}}}}%
[b_{t}^{2}]t^{-1}=0,$ we can check that
\[
\lim_{\delta
\downarrow0}{{{\mathbb{\hat{E}}}}}[|I_{\delta}|]\delta^{-1}=0,
\]
from which we get $\partial_{t}v(t,x)-g(\partial_{x}v(t,x))\leq0$,
hence $u$ is a viscosity subsolution of (\ref{G-mean}). We can
analogously prove that $u$ is also a viscosity supersolution. It
follows that $b_{t}$ is $N([\underline
{\mu}t,\overline{\mu}t]\times\{0\})$-distributed. The proof is
complete.
\end{proof}

It is clear that $\left \langle B\right \rangle $ satisfies all the
conditions in the Proposition \ref{c3p1}, thus we immediately have

\begin{theorem}
$\left \langle B\right \rangle _{t}$ is
$N([\underline{\sigma}^{2}t,\bar{\sigma
}^{2}t]\times\{0\})$-distributed, i.e., for each $\varphi \in C_{l.Lip}(\mathbb{R)}$,%
\begin{equation}
\hat{\mathbb{E}}[\varphi(\left \langle B\right \rangle _{t})]=\sup
_{\underline{\sigma}^{2}\leq v\leq \bar{\sigma}^{2}}\varphi(vt).
\end{equation}

\end{theorem}

\begin{corollary}
For each $0\leq t\leq T<\infty$, we have
\[
\underline{\sigma}^{2}(T-t)\leq \left \langle B\right \rangle
_{T}-\left \langle B\right \rangle _{t}\leq \bar{\sigma}^{2}(T-t)\
~\text{in} \ L_G^1(\Omega).
\]

\end{corollary}

\begin{proof}
It is a direct consequence of
\[
\hat{\mathbb{E}}\mathbb{[(}\left \langle B\right \rangle _{T}-\left
\langle
B\right \rangle _{t}-\bar{\sigma}^{2}(T-t))^{+}]=\sup_{\underline{\sigma}%
^{2}\leq v\leq \bar{\sigma}^{2}}(v-\bar{\sigma}^{2})^{+}(T-t)=0
\]
and%
\[
\hat{\mathbb{E}}\mathbb{[(}\left \langle B\right \rangle _{T}-\left
\langle B\right \rangle
_{t}-\underline{\sigma}^{2}(T-t))^{-}]=\sup_{\underline{\sigma
}^{2}\leq v\leq
\bar{\sigma}^{2}}(v-\underline{\sigma}^{2})^{-}(T-t)=0.
\]

\end{proof}

\begin{corollary}
{ { { \label{Lem-Qua2 copy(1)} We have, for each }}}$t,s\geq0$,
$n\in \mathbb{N},${{{
\begin{equation}
\hat{\mathbb{E}}[\mathbb{(}\left \langle B\right \rangle
_{t+s}-\left \langle B\right \rangle
_{s})^{n}|\Omega_{s}]=\hat{\mathbb{E}}[\left \langle
B\right \rangle _{t}^{n}]=\bar{\sigma}^{2n}t^{n}%
\end{equation}
and}}}%
\begin{equation}
\hat{\mathbb{E}}[-\mathbb{(}\left \langle B\right \rangle _{t+s}%
-\left \langle B\right \rangle _{s})^{n}|\Omega_{s}]=\mathbb{\hat{E}}%
[-\left \langle B\right \rangle
_{t}^{n}]=-\underline{\sigma}^{2n}t^{n}.
\end{equation}

\end{corollary}

We now consider the multi-dimensional case. For notational
simplicity, we denote by $B^{i}:=B^{\mathbf{e}_{i}}$ the $i$-th
coordinate of the $G$--Brownian motion $B$, under a given
orthonormal basis $(\mathbf{e}_{1},\cdots,\mathbf{e}_{d})$ of
$\mathbb{R}^{d}$. We denote
\[
(\left \langle B\right \rangle _{t})_{ij}:=\left \langle
B^{i},B^{j}\right \rangle _{t}.
\]
Then $\left \langle B\right \rangle _{t}$, $t\geq0$, is an
$\mathbb{S}(d)$-valued
process. Since%

\[
\hat{\mathbb{E}}[\langle AB_{t},B_{t}\rangle]=2G(A)t\  \ \text{for}\
A\in \mathbb{S}(d),
\]
we have%
\begin{align*}
\hat{\mathbb{E}}[(\left \langle B\right \rangle _{t},A)] &
=\mathbb{\hat
{E}}[\sum_{i,j=1}^{d}a_{ij}\left \langle B^{i},B^{j}\right \rangle _{t}]\\
& =\hat{\mathbb{E}}[\sum_{i,j=1}^{d}a_{ij}(B_{t}^{i}B_{t}^{j}-\int
_{0}^{t}B_{s}^{i}dB_{s}^{j}-\int_{0}^{t}B_{s}^{j}dB_{s}^{i})]\\
&  =\hat{\mathbb{E}}[\sum_{i,j=1}^{d}a_{ij}B_{t}^{i}B_{t}^{j}%
]=2G(A)t\ \ \text{for} \ A\in \mathbb{S}(d),
\end{align*}
where $(a_{ij})_{i,j=1}^d=A$.

Now we set, for each $\varphi \in C_{l.Lip}(\mathbb{S}(d)),$
\[
v(t,X):=\hat{\mathbb{E}}[\varphi(X+\left \langle B\right \rangle
_{t})],\  \ (t,X)\in \lbrack0,\infty)\times \mathbb{S}(d).
\]
Let $\Sigma \subset \mathbb{S}_{+}(d)$ be the bounded, convex and
closed subset such that
\[
G(A)=\frac{1}{2}\sup_{B\in \Sigma}\left(  A,B\right)  ,\  \  \  \
A\in \mathbb{S}(d).
\]

\begin{proposition}
\label{p1} The function $v$ solves the following first order PDE:
\[
\partial_{t}v-2G(Dv)=0,\ v|_{t=0}=\varphi,\  \
\]
where $Dv=(\partial_{x_{ij}}v)_{i,j=1}^{d}$. We also have%
\[
v(t,X)=\sup_{\Lambda \in \Sigma}\varphi(X+t\Lambda).
\]

\end{proposition}

\noindent\textbf{Sketch of the Proof.} We have%
\begin{align*}
v(t+\delta,X) &  =\hat{\mathbb{E}}[\varphi(X+\left \langle B\right
\rangle _{\delta}+\left \langle B\right \rangle _{t+\delta}-\left
\langle B\right \rangle
_{\delta})]\\
&  =\hat{\mathbb{E}}[v(t,X+\left \langle B\right \rangle
_{\delta})].
\end{align*}
The rest part of the proof is similar to  the 1-dimensional case.
$\Box$

\begin{corollary}
We have%
\[
\left \langle B\right \rangle _{t}\in t\Sigma:=\{t\times
\gamma:\gamma \in \Sigma \},
\]
or equivalently, $d_{t\Sigma}(\left \langle B\right \rangle
_{t})=0$, where $d_{U}(X)=\inf \{ \sqrt{(X-Y,X-Y)}:Y\in U\}.$
\end{corollary}

\begin{proof}
Since
\[
\hat{\mathbb{E}}[d_{t\Sigma}(\left \langle B\right \rangle _{t}%
)]=\sup_{\Lambda \in \Sigma}d_{t\Sigma}(t\Lambda)=0\text{,}%
\]
it follows that $d_{t\Sigma}(\left \langle B\right \rangle _{t})=0$.
\end{proof}

\begin{exercise}
Complete the proof of Proposition \ref{p1}.
\end{exercise}

\section{$G$--It\^{o}'s Formula}

{{{ In this section, we give It\^{o}'s formula for a
\textquotedblleft$G$-It\^{o} process\textquotedblright \ $X$. For
simplicity, we first consider the case of the }}}function $\Phi$ is
sufficiently regular.{{{ }}}

\begin{lemma}
\label{d-Lem-26}Let $\Phi \in C^{2}(\mathbb{R}^{n})$ with $\partial_{x^{\nu}}%
\Phi,\  \partial_{x^{\mu}x^{\nu}}^{2}\Phi \in
C_{b.Lip}(\mathbb{R}^{n})$ for
$\mu,\nu=1,\cdots,n$. Let $s\in \lbrack0,T]$ be fixed and let $X=(X^{1}%
,\cdots,X^{n})^{T}$ be an $n$--dimensional process on $[s,T]$ of the
form
\[
X_{t}^{\nu}=X_{s}^{\nu}+\alpha^{\nu}(t-s)+\eta^{\nu ij}(\left
\langle B^{i},B^{j}\right \rangle _{t}-\left \langle
B^{i},B^{j}\right \rangle _{s})+\beta^{\nu j}(B_{t}^{j}-B_{s}^{j}),
\]
where, for $\nu=1,\cdots,n$, $i,j=1,\cdots,d$, $\alpha^{\nu}$,
$\eta^{\nu ij}$ and $\beta^{\nu j}$ are bounded elements in
$L_{G}^{2}(\Omega_{s})$ and $X_{s}=(X_{s}^{1},\cdots,X_{s}^{n})^{T}$
is a given random vector in $L_{G}^{2}(\Omega_{s})$. Then we have,
in $L_{G}^{2}(\Omega_{t})$,
\begin{align}
\Phi(X_{t})-\Phi(X_{s}) &
=\int_{s}^{t}\partial_{x^{\nu}}\Phi(X_{u})\beta^{\nu
j}dB_{u}^{j}+\int_{s}^{t}\partial_{x^{\nu}}\Phi(X_{u})\alpha^{\nu}%
du\label{d-B-Ito}\\
&\ \ \ +\int_{s}^{t}[\partial_{x^{\nu}}\Phi(X_{u})\eta^{\nu ij}+\frac{1}{2}%
\partial_{x^{\mu}x^{\nu}}^{2}\Phi(X_{u})\beta^{\mu i}\beta^{\nu j}]d\left \langle
B^{i},B^{j}\right \rangle _{u}.\nonumber
\end{align}
Here we use the \index{Einstein convention}, i.e., the above
repeated indices $\mu ,\nu$, $i$ and $j$  imply the summation.
\end{lemma}

\begin{proof}
For each positive integer $N$, we set $\delta=(t-s)/N$ and take the
partition
\[
\pi_{\lbrack
s,t]}^{N}=\{t_{0}^{N},t_{1}^{N},\cdots,t_{N}^{N}\}=\{s,s+\delta
,\cdots,s+N\delta=t\}.
\]
We have
\begin{align}
\Phi(X_{t})-\Phi(X_{s}) &  =\sum_{k=0}^{N-1}[\Phi(X_{t_{k+1}^{N}}%
)-\Phi(X_{t_{k}^{N}})]\label{d-Ito} \\
&  =\sum_{k=0}^{N-1}\{ \partial_{x^{\nu}}\Phi(X_{t_{k}^{N}})(X_{t_{k+1}^{N}}%
^{\nu}-X_{t_{k}^{N}}^{\nu})\nonumber \\
& \ \ \ +\frac{1}{2}[\partial_{x^{\mu}x^{\nu}}^{2}\Phi(X_{t_{k}^{N}})(X_{t_{k+1}^{N}%
}^{\mu}-X_{t_{k}^{N}}^{\mu})(X_{t_{k+1}^{N}}^{\nu}-X_{t_{k}^{N}}^{\nu}%
)+\eta_{k}^{N}]\},\nonumber
\end{align}
where
\[
\eta_{k}^{N}=[\partial_{x^{\mu}x^{\nu}}^{2}\Phi(X_{t_{k}^{N}}+\theta
_{k}(X_{t_{k+1}^{N}}-X_{t_{k}^{N}}))-\partial_{x^{\mu}x^{\nu}}^{2}\Phi
(X_{t_{k}^{N}})](X_{t_{k+1}^{N}}^{\mu}-X_{t_{k}^{N}}^{\mu})(X_{t_{k+1}^{N}%
}^{\nu}-X_{t_{k}^{N}}^{\nu})
\]
with $\theta_{k}\in \lbrack0,1]$. We have%
\begin{align*}
\hat{\mathbb{E}}[|\eta_{k}^{N}|^{2}] &  =\hat{\mathbb{E}}%
[|[\partial_{x^{\mu}x^{\nu}}^{2}\Phi(X_{t_{k}^{N}}+\theta_{k}(X_{t_{k+1}^{N}%
}-X_{t_{k}^{N}}))-\partial_{x^{\mu}x^{\nu}}^{2}\Phi(X_{t_{k}^{N}})]\\
&\ \ \
\times(X_{t_{k+1}^{N}}^{\mu}-X_{t_{k}^{N}}^{\mu})(X_{t_{k+1}^{N}}^{\nu
}-X_{t_{k}^{N}}^{\nu})|^{2}]\\
&\leq c\hat{\mathbb{E}}\mathbb{[}|X_{t_{k+1}^{N}}-X_{t_{k}^{N}}%
|^{6}]\leq C[\delta^{6}+\delta^{3}],
\end{align*}
where $c$ is the Lipschitz constant of $\{ \partial_{x^{\mu}x^{\nu}}^{2}%
\Phi \}_{\mu,\nu=1}^{n}${ }and $C$ is a constant independent of
$k${{{.}}} Thus
\[
\hat{\mathbb{E}}[|\sum_{k=0}^{N-1}\eta_{k}^{N}|^{2}]\leq N\sum_{k=0}%
^{N-1}\hat{\mathbb{E}}[|\eta_{k}^{N}|^{2}]\rightarrow0.
\]
The rest terms in the summation of the right side of (\ref{d-Ito})
are
$\xi_{t}^{N}+\zeta_{t}^{N}$ with%
\begin{align*}
\xi_{t}^{N} &  =\sum_{k=0}^{N-1}\{ \partial_{x^{\nu}}\Phi(X_{t_{k}^{N}}%
)[\alpha^{\nu}(t_{k+1}^{N}-t_{k}^{N})+\eta^{\nu ij}(\left \langle B^{i}%
,B^{j}\right \rangle _{t_{k+1}^{N}}-\left \langle B^{i},B^{j}\right
\rangle
_{t_{k}^{N}})\\
&\ \ \  +\beta^{\nu j}(B_{t_{k+1}^{N}}^{j}-B_{t_{k}^{N}}^{j})]+\frac{1}{2}%
\partial_{x^{\mu}x^{\nu}}^{2}\Phi(X_{t_{k}^{N}})\beta^{\mu i}\beta^{\nu
j}(B_{t_{k+1}^{N}}^{i}-B_{t_{k}^{N}}^{i})(B_{t_{k+1}^{N}}^{j}-B_{t_{k}^{N}%
}^{j})\}
\end{align*}
and
\begin{align*}
\zeta_{t}^{N} &  =\frac{1}{2}\sum_{k=0}^{N-1}\partial_{x^{\mu}x^{\nu}}^{2}%
\Phi(X_{t_{k}^{N}})\{[\alpha^{\mu}(t_{k+1}^{N}-t_{k}^{N})+\eta^{\mu
ij}(\left \langle B^{i},B^{j}\right \rangle _{t_{k+1}^{N}}-\left
\langle
B^{i},B^{j}\right \rangle _{t_{k}^{N}})]\\
&\ \ \  \times \lbrack \alpha^{\nu}(t_{k+1}^{N}-t_{k}^{N})+\eta^{\nu
lm}(\left \langle B^{l},B^{m}\right \rangle _{t_{k+1}^{N}}-\left
\langle B^{l},B^{m}\right \rangle
_{t_{k}^{N}})]\\
& \ \ \ +2[\alpha^{\mu}(t_{k+1}^{N}-t_{k}^{N})+\eta^{\mu ij}(\left
\langle B^{i},B^{j}\right \rangle _{t_{k+1}^{N}}-\left \langle
B^{i},B^{j}\right \rangle _{t_{k}^{N}})]\beta^{\nu
l}(B_{t_{k+1}^{N}}^{l}-B_{t_{k}^{N}}^{l})\}.
\end{align*}
We observe that, for each $u\in \lbrack t_{k}^{N},t_{k+1}^{N})$,
\begin{align*}
&  \hat{\mathbb{E}}[|\partial_{x^{\nu}}\Phi(X_{u})-\sum_{k=0}^{N-1}%
\partial_{x^{\nu}}\Phi(X_{t_{k}^{N}})\mathbf{I}_{[t_{k}^{N},t_{k+1}^{N}%
)}(u)|^{2}]\\
&  =\hat{\mathbb{E}}[|\partial_{x^{\nu}}\Phi(X_{u})-\partial_{x^{\nu}}%
\Phi(X_{t_{k}^{N}})|^{2}]\\
&  \leq c^{2}\hat{\mathbb{E}}[|X_{u}-X_{t_{k}^{N}}|^{2}]\leq
C[\delta+\delta^{2}],
\end{align*}
{{{where $c$ is the Lipschitz constant of }}}$\{
\partial_{x^{\nu}}\Phi \}_{\nu
=1}^{n}$ and $C$ is a constant independent of $k$. Thus $\sum_{k=0}%
^{N-1}\partial_{x^{\nu}}\Phi(X_{t_{k}^{N}})\mathbf{I}_{[t_{k}^{N},t_{k+1}^{N}%
)}(\cdot)$ converges to $\partial_{x^{\nu}}\Phi(X_{\cdot})$ in
$M_{G}^{2}(0,T)$. Similarly, $
\sum_{k=0}^{N-1}\partial_{x^{\mu}x^{\nu}}^{2}\Phi(X_{t_{k}^{N}})\mathbf{I}%
_{[t_{k}^{N},t_{k+1}^{N})}(\cdot)$ converges to $\partial_{x^{\mu}x^{\nu}}^{2}%
\Phi(X_{\cdot})\text{ in \ }M_{G}^{2}(0,T).$

From Lemma \ref{d-Lem-mutual} as well as the definitions of the
integrations of $dt$, $dB_{t}$ and $d\left \langle B\right \rangle
_{t}$, the limit of $\xi_{t}^{N}$ in $L_{G}^{2}(\Omega_{t})$ is just
the right hand side of (\ref{d-B-Ito}). By the next Remark we also
have $\zeta_{t}^{N}\rightarrow0$ in $L_{G}^{2}(\Omega_{t})$. We then
have proved (\ref{d-B-Ito}).
\end{proof}

\begin{remark}
To prove $\zeta_{t}^{N}\rightarrow0$ in $L_{G}^{2}(\Omega_{t})$, we
use
the following estimates: for $\psi^{N}\in M_{G}^{2,0}(0,T)$ with $\psi_{t}%
^{N}=\sum_{k=0}^{N-1}\xi_{t_{k}}^{N}\mathbf{I}_{[t_{k}^{N},t_{k+1}^{N})}(t)$,
and $\pi_{T}^{N}=\{t_{0}^{N},\cdots,t_{N}^{N}\}$ such that $\lim
_{N\rightarrow \infty}\mu(\pi_{T}^{N})=0$ and $\hat{\mathbb{E}}[\sum
_{k=0}^{N-1}|\xi_{t_{k}}^{N}|^{2}(t_{k+1}^{N}-t_{k}^{N})]\leq C$,
for all
$N=1,2,\cdots$, we have $\hat{\mathbb{E}}[|\sum_{k=0}^{N-1}\xi_{k}%
^{N}(t_{k+1}^{N}-t_{k}^{N})^{2}|^{2}]\rightarrow0$ and, for any
fixed
$\mathbf{a,\bar{a}\in}\mathbb{R}^{d}$,%
\begin{align*}
\hat{\mathbb{E}}[|\sum_{k=0}^{N-1}\xi_{k}^{N}(\left \langle B^{\mathbf{a}%
}\right \rangle _{t_{k+1}^{N}}-\left \langle B^{\mathbf{a}}\right
\rangle _{t_{k}^{N}})^{2}|^{2}] &  \leq
C\hat{\mathbb{E}}[\sum_{k=0}^{N-1}|\xi
_{k}^{N}|^{2}(\left \langle B^{\mathbf{a}}\right \rangle _{t_{k+1}^{N}%
}-\left \langle B^{\mathbf{a}}\right \rangle _{t_{k}^{N}})^{3}]\\
&  \leq C\hat{\mathbb{E}}[\sum_{k=0}^{N-1}|\xi_{k}^{N}|^{2}\sigma
_{\mathbf{aa}^{T}}^{6}(t_{k+1}^{N}-t_{k}^{N})^{3}]\rightarrow0,
\end{align*}%
\begin{align*}
&\hat{\mathbb{E}}[|\sum_{k=0}^{N-1}\xi_{k}^{N}(\left \langle
B^{\mathbf{a}}\right \rangle _{t_{k+1}^{N}}-\left \langle B^{\mathbf{a}%
}\right \rangle _{t_{k}^{N}})(t_{k+1}^{N}-t_{k}^{N})|^{2}]\\
\leq &C\hat{\mathbb{E}}[\sum_{k=0}^{N-1}|\xi_{k}^{N}|^{2}(t_{k+1}%
^{N}-t_{k}^{N})(\left \langle B^{\mathbf{a}}\right \rangle _{t_{k+1}^{N}%
}-\left \langle B^{\mathbf{a}}\right \rangle _{t_{k}^{N}})^{2}]\\
\leq& C\hat{\mathbb{E}}[\sum_{k=0}^{N-1}|\xi_{k}^{N}|^{2}\sigma
_{\mathbf{aa}^{T}}^{4}(t_{k+1}^{N}-t_{k}^{N})^{3}]\rightarrow0,
\end{align*}
as well as
\begin{align*}
&\hat{\mathbb{E}}[|\sum_{k=0}^{N-1}\xi_{k}^{N}(t_{k+1}^{N}-t_{k}%
^{N})(B_{t_{k+1}^{N}}^{\mathbf{a}}-B_{t_{k}^{N}}^{\mathbf{a}})|^{2}]\\
\leq&
C\hat{\mathbb{E}}[\sum_{k=0}^{N-1}|\xi_{k}^{N}|^{2}(t_{k+1}^{N}-t_{k}%
^{N})|B_{t_{k+1}^{N}}^{\mathbf{a}}-B_{t_{k}^{N}}^{\mathbf{a}}|^{2}]\\
\leq&C\hat{\mathbb{E}}[\sum_{k=0}^{N-1}|\xi_{k}^{N}|^{2}\sigma
_{\mathbf{aa}^{T}}^{2}(t_{k+1}^{N}-t_{k}^{N})^{2}]\rightarrow0\
\end{align*}
and%
\begin{align*}
&\hat{\mathbb{E}}[|\sum_{k=0}^{N-1}\xi_{k}^{N}(\left \langle
B^{\mathbf{a}}\right \rangle _{t_{k+1}^{N}}-\left \langle B^{\mathbf{a}%
}\right \rangle _{t_{k}^{N}})(B_{t_{k+1}^{N}}^{\mathbf{\bar{a}}}-B_{t_{k}^{N}%
}^{\mathbf{\bar{a}}})|^{2}]\\
\leq&C\hat{\mathbb{E}}[\sum_{k=0}^{N-1}|\xi_{k}^{N}|^{2}(\left
\langle
B^{\mathbf{a}}\right \rangle _{t_{k+1}^{N}}-\left \langle B^{\mathbf{a}%
}\right \rangle _{t_{k}^{N}})|B_{t_{k+1}^{N}}^{\mathbf{\bar{a}}}-B_{t_{k}^{N}%
}^{\mathbf{\bar{a}}}|^{2}]\\
\leq&C\hat{\mathbb{E}}[\sum_{k=0}^{N-1}|\xi_{k}^{N}|^{2}\sigma
_{\mathbf{aa}^{T}}^{2}\sigma_{\mathbf{\bar{a}\bar{a}}^{T}}^{2}(t_{k+1}%
^{N}-t_{k}^{N})^{2}]\rightarrow0.
\end{align*}
\endproof
\end{remark}

{ { { We now consider a general form of }}}$G$--It\^{o}'s
formula{{{.
Consider} } }%
\[
X_{t}^{\nu}=X_{0}^{\nu}+\int_{0}^{t}\alpha_{s}^{\nu}ds+\int_{0}^{t}\eta
_{s}^{\nu ij}d\left \langle B^{i},B^{j}\right \rangle
_{s}+\int_{0}^{t}\beta _{s}^{\nu j}dB_{s}^{j}.
\]

\begin{proposition}
\label{d-Prop-Ito}Let $\Phi \in C^{2}(\mathbb{R}^{n})$ with $\partial_{x^{\nu}%
}\Phi,\  \partial_{x^{\mu}x^{\nu}}^{2}\Phi \in
C_{b.Lip}(\mathbb{R}^{n})$ for $\mu,\nu=1,\cdots,n$. Let
$\alpha^{\nu}$, $\beta^{\nu j}$ and $\eta^{\nu ij}$,
$\nu=1,\cdots,n$, $i,j=1,\cdots,d$ be bounded processes in
$M_{G}^{2}(0,T)$.
Then for each $t\geq0$ we have, in $L_{G}^{2}(\Omega_{t})$%
\begin{align}
\Phi(X_{t})-\Phi(X_{s}) &
=\int_{s}^{t}\partial_{x^{\nu}}\Phi(X_{u})\beta _{u}^{\nu
j}dB_{u}^{j}+\int_{s}^{t}\partial_{x^{\nu}}\Phi(X_{u})\alpha_{u}^{\nu
}du\label{d-Ito-form1}\\
& \ \ \ +\int_{s}^{t}[\partial_{x^{\nu}}\Phi(X_{u})\eta_{u}^{\nu ij}+\frac{1}%
{2}\partial_{x^{\mu}x^{\nu}}^{2}\Phi(X_{u})\beta_{u}^{\mu
i}\beta_{u}^{\nu j}]d\left \langle B^{i},B^{j}\right \rangle _{u}.
\nonumber
\end{align}

\end{proposition}

\begin{proof}
We first consider the case where $\alpha$, $\eta$ and $\beta$ are
step
processes of the form%
\[
\eta_{t}(\omega)=\sum_{k=0}^{N-1}\xi_{k}(\omega)\mathbf{I}_{[t_{k},t_{k+1}%
)}(t).
\]
From the above lemma, it is clear that (\ref{d-Ito-form1}) holds
true. Now let
\[
X_{t}^{\nu,N}=X_{0}^{\nu}+\int_{0}^{t}\alpha_{s}^{\nu,N}ds+\int_{0}^{t}%
\eta_{s}^{\nu ij,N}d\left \langle B^{i},B^{j}\right \rangle _{s}+\int_{0}%
^{t}\beta_{s}^{\nu j,N}dB_{s}^{j},
\]
where $\alpha^{N}$, $\eta^{N}$ and $\beta^{N}$ are uniformly bounded
step processes that converge to $\alpha$, $\eta$ and $\beta$ in
$M_{G}^{2}(0,T)$ as $N\rightarrow \infty$, respectively. From Lemma
\ref{d-Lem-26},
\begin{align}
\Phi(X_{t}^{N})-\Phi(X_{s}^{N}) &
=\int_{s}^{t}\partial_{x^{\nu}}\Phi
(X_{u}^{N})\beta_{u}^{\nu j,N}dB_{u}^{j}+\int_{s}^{t}\partial_{x^{\nu}}%
\Phi(X_{u}^{N})\alpha_{u}^{\nu,N}du\label{d-N-Ito}\\
& \ \ \ +\int_{s}^{t}[\partial_{x^{\nu}}\Phi(X_{u}^{N})\eta_{u}^{\nu
ij,N}+\frac
{1}{2}\partial_{x^{\mu}x^{\nu}}^{2}\Phi(X_{u}^{N})\beta_{u}^{\mu
i,N}\beta _{u}^{\nu j,N}]d\left \langle B^{i},B^{j}\right \rangle
_{u}.\nonumber
\end{align}
Since%
\begin{align*}
&\hat{\mathbb{E}}\mathbb{[}|X_{t}^{\nu,N}-X_{t}^{\nu}|^{2}]\\
\leq&C\hat{\mathbb{E}}\mathbb{[}\int_{0}^{T}[(\alpha_{s}^{\nu,N}%
-\alpha_{s}^{\nu})^{2}+|\eta_{s}^{\nu,N}-\eta_{s}^{\nu}|^{2}+|\beta_{s}^{\nu,N}%
-\beta_{s}^{\nu}|^{2}]ds],
\end{align*}
where $C$ is a constant independent of $N$, we can prove that, in
$M_{G}^{2}(0,T)$,
\begin{align*}
\partial_{x^{\nu}}\Phi(X_{\cdot}^{N})\eta_{\cdot}^{\nu ij,N} &  \rightarrow
\partial_{x^{\nu}}\Phi(X_{\cdot})\eta_{\cdot}^{\nu ij},\\
\partial_{x^{\mu}x^{\nu}}^{2}\Phi(X_{\cdot}^{N})\beta_{\cdot}^{\mu i,N}%
\beta_{\cdot}^{\nu j,N} &  \rightarrow \partial_{x^{\mu}x^{\nu}}^{2}%
\Phi(X_{\cdot})\beta_{\cdot}^{\mu i}\beta_{\cdot}^{\nu j},\\
\partial_{x^{\nu}}\Phi(X_{\cdot}^{N})\alpha_{\cdot}^{\nu,N} &  \rightarrow
\partial_{x^{\nu}}\Phi(X_{\cdot})\alpha_{\cdot}^{\nu},\\
\partial_{x^{\nu}}\Phi(X_{\cdot}^{N})\beta_{\cdot}^{\nu j,N} &  \rightarrow
\partial_{x^{\nu}}\Phi(X_{\cdot})\beta_{\cdot}^{\nu j}.
\end{align*}
We then can pass to limit as $N\rightarrow\infty$ in both sides of
(\ref{d-N-Ito}) to get (\ref{d-Ito-form1}).
\end{proof}

In order to consider the general $\Phi$, we first prove a useful
inequality.

For the $G$-expectation $\hat{\mathbb{E}}$, we have the following
representation (see Chap.\ref{ch6}):
\begin{equation}
\hat{\mathbb{E}}[X]=\sup_{P\in \mathcal{P}}E_{P}[X]\  \ \text{for}\
X\in L_{G}^{1}(\Omega),
\end{equation}
where $\mathcal{P}$ is a weakly compact family of probability
measures on $(\Omega,\mathcal{B}(\Omega))$.

\begin{proposition}
\label{B-D-G} Let $\beta \in M_{G}^{p}(0,T)$ with $p\geq2$ and let
$\mathbf{a}\in
\mathbb{R}^{d}$ be fixed. Then we have $\int_{0}^{T}\beta_{t}dB_{t}%
^{\mathbf{a}}\in L_{G}^{p}(\Omega_{T})$ and
\begin{equation}
\label{ebdg}
\hat{\mathbb{E}}[|\int_{0}^{T}\beta_{t}dB_{t}^{\mathbf{a}}|^{p}]\leq
C_{p}\hat{\mathbb{E}}[|\int_{0}^{T}\beta_{t}^{2}d\langle B^{\mathbf{a}%
}\rangle_{t}|^{p/2}].
\end{equation}

\end{proposition}

\begin{proof}
It suffices to consider the case where $\beta$ is a step process of the form%
\[
\beta_{t}(\omega)=\sum_{k=0}^{N-1}\xi_{k}(\omega)\mathbf{I}_{[t_{k},t_{k+1}%
)}(t).
\]
For each $\xi \in L_{ip}(\Omega_{t})$ with $t\in \lbrack0,T]$, we have%
\[
\hat{\mathbb{E}}[\xi \int_{t}^{T}\beta_{s}dB_{s}^{\mathbf{a}}]=0.
\]
From this we can easily get $E_{P}[\xi \int_{t}^{T}\beta_{s}dB_{s}^{\mathbf{a}%
}]=0$ for each $P\in \mathcal{P}$, which implies that
$(\int_{0}^{t}\beta _{s}dB_{s}^{\mathbf{a}})_{t\in0,T]}$ is a
$P$-martingale. Similarly we can prove that
\[
M_{t}:=(\int_{0}^{t}\beta_{s}dB_{s}^{\mathbf{a}})^{2}-\int_{0}^{t}\beta
_{s}^{2}d\langle B^{\mathbf{a}}\rangle_{s},\  \ t\in \lbrack0,T]
\]
is a $P$-martingale for each $P\in \mathcal{P}$. By the
Burkholder-Davis-Gundy
inequalities, we have%
\[
E_{P}[|\int_{0}^{T}\beta_{t}dB_{t}^{\mathbf{a}}|^{p}]\leq
C_{p}E_{P}[|\int _{0}^{T}\beta_{t}^{2}d\langle
B^{\mathbf{a}}\rangle_{t}|^{p/2}]\leq
C_{p}\hat{\mathbb{E}}[|\int_{0}^{T}\beta_{t}^{2}d\langle B^{\mathbf{a}%
}\rangle_{t}|^{p/2}],
\]
where $C_{p}$ is a universal constant independent of $P$. Thus we
get (\ref{ebdg}).
\end{proof}

We now give the general $G$--It\^{o}'s formula.

\begin{theorem}\label{Thm6.5}
Let $\Phi$ be a $C^2$-function on $\mathbb{R}^n$ such that
$\partial_{x^{\mu}x^{\nu}}^{2}\Phi$ satisfy polynomial growth
condition for $\mu,\nu=1,\cdots,n$. Let $\alpha^{\nu}$, $\beta^{\nu
j}$ and $\eta^{\nu ij}$, $\nu=1,\cdots,n$, $i,j=1,\cdots,d$ be
bounded processes
in $M_{G}^{2}(0,T)$. Then for each $t\geq0$ we have in $L_{G}^{2}(\Omega_{t})$%
\begin{align}
\Phi(X_{t})-\Phi(X_{s}) &
=\int_{s}^{t}\partial_{x^{\nu}}\Phi(X_{u})\beta _{u}^{\nu
j}dB_{u}^{j}+\int_{s}^{t}\partial_{x^{\nu}}\Phi(X_{u})\alpha_{u}^{\nu
}du\label{e629}\\
& \ \ \ +\int_{s}^{t}[\partial_{x^{\nu}}\Phi(X_{u})\eta_{u}^{\nu ij}+\frac{1}%
{2}\partial_{x^{\mu}x^{\nu}}^{2}\Phi(X_{u})\beta_{u}^{\mu
i}\beta_{u}^{\nu j}]d\left \langle B^{i},B^{j}\right \rangle
_{u}.\nonumber
\end{align}

\end{theorem}

\begin{proof}
By the assumptions on $\Phi$, we can choose a sequence of functions
$\Phi
_{N}\in C_{0}^{2}(\mathbb{R}^{n})$ such that%
\[
|\Phi_{N}(x)-\Phi(x)|+|\partial_{x^{\nu}}\Phi_{N}(x)-\partial_{x^{\nu}}%
\Phi(x)|+|\partial_{x^{\mu}x^{\nu}}^{2}\Phi_{N}(x)-\partial_{x^{\mu}x^{\nu}%
}^{2}\Phi(x)|\leq \frac{C_{1}}{N}(1+|x|^{k}),
\]
where $C_{1}$ and $k$ are positive constants independent of $N$.
Obviously, $\Phi_{N}$ satisfies the conditions in Proposition
\ref{d-Prop-Ito}, therefore,
\begin{align}
\Phi_{N}(X_{t})-\Phi_{N}(X_{s}) &
=\int_{s}^{t}\partial_{x^{\nu}}\Phi _{N}(X_{u})\beta_{u}^{\nu
j}dB_{u}^{j}+\int_{s}^{t}\partial_{x^{v}}\Phi
_{N}(X_{u})\alpha_{u}^{\nu}du\label{e630}\\
& \ \ \ +\int_{s}^{t}[\partial_{x^{\nu}}\Phi_{N}(X_{u})\eta_{u}^{\nu
ij}+\frac
{1}{2}\partial_{x^{\mu}x^{\nu}}^{2}\Phi_{N}(X_{u})\beta_{u}^{\mu
i}\beta _{u}^{\nu j}]d\left \langle B^{i},B^{j}\right \rangle
_{u}.\nonumber
\end{align}
For each fixed $T>0$, by Proposition \ref{B-D-G}, there exists a
constant $C_{2}$ such that
\[
\hat{\mathbb{E}}[|X_{t}|^{2k}]\leq C_{2}\  \  \text{for}\ t\in
\lbrack0,T].
\]
Thus we can prove that $\Phi_{N}(X_{t})\rightarrow \Phi(X_{t})$ in $L_{G}%
^{2}(\Omega_{t})$ and in $M_{G}^{2}(0,T),$%
\begin{align*}
\partial_{x^{\nu}}\Phi_{N}(X_{\cdot})\eta_{\cdot}^{\nu ij} &  \rightarrow
\partial_{x^{\nu}}\Phi(X_{\cdot})\eta_{\cdot}^{\nu ij},\\
\partial_{x^{\mu}x^{\nu}}^{2}\Phi_{N}(X_{\cdot})\beta_{\cdot}^{\mu i}%
\beta_{\cdot}^{\nu j} &  \rightarrow
\partial_{x^{\mu}x^{\nu}}^{2}\Phi
(X_{\cdot})\beta_{\cdot}^{\mu i}\beta_{\cdot}^{\nu j},\\
\partial_{x^{\nu}}\Phi_{N}(X_{\cdot})\alpha_{\cdot}^{\nu} &  \rightarrow
\partial_{x^{\nu}}\Phi(X_{\cdot})\alpha_{\cdot}^{\nu},\\
\partial_{x^{\nu}}\Phi_{N}(X_{\cdot})\beta_{\cdot}^{\nu j} &  \rightarrow
\partial_{x^{\nu}}\Phi(X_{\cdot})\beta_{\cdot}^{\nu j}.
\end{align*}
We then can pass to limit as $N\rightarrow\infty$ in both sides of
(\ref{e630}) to get (\ref{e629}).
\end{proof}

\begin{corollary}
Let $\Phi$ be a polynomial and $\mathbf{a},\mathbf{a}^{\nu}\mathbf{\in}%
\mathbb{R}^{d}$ be fixed for $\nu=1,\cdots,n$. Then we have
\[
\Phi(X_{t})-\Phi(X_{s})=\int_{s}^{t}\partial_{x^{\nu}}\Phi(X_{u}%
)dB_{u}^{\mathbf{a}^{\nu}}+\frac{1}{2}\int_{s}^{t}\partial_{x^{\mu}x^{\nu}}%
^{2}\Phi(X_{u})d\left \langle B^{\mathbf{a}^{\mu}},B^{\mathbf{a}^{\nu}%
}\right \rangle _{u},
\]
where
$X_{t}=(B_{t}^{\mathbf{a}^{1}},\cdots,B_{t}^{\mathbf{a}^{n}})^{T}$.
In particular, we have, for $k=2,3,\cdots$,
\[
(B_{t}^{\mathbf{a}})^{k}=k\int_{0}^{t}(B_{s}^{\mathbf{a}})^{k-1}%
dB_{s}^{\mathbf{a}}+\frac{k(k-1)}{2}\int_{0}^{t}(B_{s}^{\mathbf{a}}%
)^{k-2}d\langle B^{\mathbf{a}}\rangle_{s}.
\]

\end{corollary}

If $\hat{\mathbb{E}}$ becomes a linear expectation, then the above
$G$--It\^{o}'s formula is the classical one.

\section{Generalized $G$-Brownian Motion}

Let $G:\mathbb{R}^{d}\times \mathbb{S}(d)\rightarrow \mathbb{R}$ be
a given continuous sublinear function monotonic in $A\in
\mathbb{S}(d)$. Then by Theorem \ref{t1} in Chap.\ref{ch1}, there
exists a bounded, convex and closed subset $\Sigma \subset \mathbb{R}^{d}%
\times \mathbb{S}_{+}(d)$ such that
\[
G(p,A)=\sup_{(q,B)\in \Sigma}[\frac{1}{2}\mathrm{tr}[AB]+\left
\langle p,q\right \rangle ]\  \  \  \text{for}\ (p,A)\in
\mathbb{R}^{d}\times \mathbb{S}(d).
\]
By Chapter \ref{ch2}, we know that there exists a pair of
$d$-dimensional random vectors $(X,Y)$ which is $G$-distributed.

We now give the definition of the generalized $G$-Brownian motion.

\begin{definition}
A $d$-dimensional process $(B_{t})_{t\geq0}$ on a sublinear
expectation space $(\Omega,\mathcal{H},\mathbb{\hat{E}})$ is called
a \textbf{generalized $G$-Brownian motion} \index{Generalized
$G$-Brownian motion} if the following properties are satisfied:
\newline \textup{(i)} $B_{0}(\omega)=0$;\newline \textup{(ii)} For each
$t,s\geq0$, the increment $B_{t+s}-B_{t}$ identically distributed
with $\sqrt{s}X+sY$ and is independent from
$(B_{t_{1}},B_{t_{2}},\cdots,B_{t_{n}})$, for each $n\in \mathbb{N}$
and $0\leq t_{1}\leq \cdots \leq t_{n}\leq t$, where $(X,Y)$ is
$G$-distributed.
\end{definition}

The following theorem gives a characterization of the generalized
$G$-Brownian motion.

\begin{theorem}
Let $(B_{t})_{t\geq0}$ be a $d$-dimensional process defined on a
sublinear expectation space $(\Omega,\mathcal{H},\mathbb{\hat{E}})${
such that
\newline \textup{(i)} }$B_{0}(\omega)${$=0$;\newline \textup{(ii)} For each
$t,s\geq0$, $B${$_{t+s}-B_{t}$ and }$B${$_{s}$ are identically
distributed and $B${$_{t+s}-B_{t}$ is }independent from $(B_{t_{1}},B_{t_{2}%
},\cdots,B_{t_{n}})$, for each $n\in \mathbb{N}$ and $0\leq
t_{1}\leq \cdots \leq
t_{n}\leq t$. \newline \textup{(iii)} $\lim_{t\downarrow0}\mathbb{\hat{E}%
}[|B_{t}|^{3}]t^{-1}=0$. \newline Then }}$(B_{t})_{t\geq0}${{ is a
generalized
$G$-Brownian motion with }}$G(p,A)=\lim_{\delta \downarrow0}\mathbb{\hat{E}%
}[\left \langle p,B_{\delta}\right \rangle +\frac{1}{2}\langle
AB_{\delta },B_{\delta}\rangle]\delta^{-1}\ \ \text{for}\ (p,A)\in
\mathbb{R}^{d}\times \mathbb{S}(d)$.
\end{theorem}

\begin{proof}
We first prove that $\lim_{\delta \downarrow0}\mathbb{\hat{E}}[\left
\langle
p,B_{\delta}\right \rangle +\frac{1}{2}\langle AB_{\delta},B_{\delta}%
\rangle]\delta^{-1}$ exists. For each fixed $(p,A)\in \mathbb{R}^{d}%
\times \mathbb{S}(d)$, we set%
\[
f(t):=\mathbb{\hat{E}}[\left \langle p,B_{t}\right \rangle
+\frac{1}{2}\langle AB_{t},B_{t}\rangle].
\]
Since
\[
|f(t+h)-f(t)|\leq \mathbb{\hat{E}}[(|p|+2|A||B_{t}|)|B_{t+h}-B_{t}%
|+|A||B_{t+h}-B_{t}|^{2}]\rightarrow0,
\]
we get that $f(t)$ is a continuous function. It is easy to prove that%
\[
\mathbb{\hat{E}}[\left \langle q,B_{t}\right \rangle ]=\mathbb{\hat{E}%
}[\left \langle q,B_{1}\right \rangle ]t\  \  \text{for}\ q\in
\mathbb{R}^{d}.
\]
Thus for each $t,s>0$,
\[
|f(t+s)-f(t)-f(s)|\leq C\mathbb{\hat{E}}[|B_{t}|]s,
\]
where $C=|A|\mathbb{\hat{E}}[|B_{1}|]$. By {{\textup{(iii)}}}, there
exists a constant $\delta_{0}>0$ such that
$\mathbb{\hat{E}}[|B_{t}|^{3}]\leq t$ for $t\leq \delta_{0}$. Thus
for each fixed $t>0$ and $N\in \mathbb{N}$ such that
$Nt\leq \delta_{0}$, we have%
\[
|f(Nt)-Nf(t)|\leq \frac{3}{4}C(Nt)^{4/3}.
\]
From this and the continuity of $f$, it is easy to show that $\lim
_{t\downarrow0}f(t)t^{-1}$ exists. Thus we can get $G(p,A)$ for each
$(p,A)\in \mathbb{R}^{d}\times \mathbb{S}(d)$. It is also easy to
check that $G$ is a continuous sublinear function monotonic in $A\in
\mathbb{S}(d)$.

{We only need to prove that, for each fixed $\varphi \in C_{b.Lip}%
(\mathbb{R}^{d}\mathbb{)}$, the function
\[
u(t,x):=\mathbb{\hat{E}}[\varphi(x+B_{t})],\  \ (t,x)\in
\lbrack0,\infty
)\times \mathbb{R}^{d}%
\]
is the viscosity solution of the following parabolic PDE}:
\begin{equation}
\label{e731}
\partial_{t}u-G(Du,D^{2}u)=0,\  \ u|_{t=0}=\varphi.
\end{equation}
{ {We first prove that }$u$ is Lipschitz in $x$ and
$\frac{1}{2}$-H\"{o}lder continuous
in $t$. In fact, {for each fixed }$t$,{ }$u(t,\cdot)\in${$C_{b.Lip}%
(\mathbb{R}^{d}\mathbb{)}$ }since%
\begin{align*}
|\mathbb{\hat{E}}[\varphi(x+B_{t})]-\mathbb{\hat{E}}[\varphi(y+B_{t})]|
&
\leq \mathbb{\hat{E}}[|\varphi(x+B_{t})-\varphi(y+B_{t})|]\\
&  \leq C|x-y|.
\end{align*}
{For each $\delta \in \lbrack0,t]$, since $B_{t}-B_{\delta}$ is
independent from $B_{\delta}$,}
\begin{align*}
u(t,x) &
=\mathbb{\hat{E}}[\mathbb{\varphi(}x+B_{\delta}+(B_{t}-B_{\delta
})]\\
&  =\mathbb{\hat{E}}[\mathbb{\hat{E}}[\varphi(y+(B_{t}-B_{\delta
}))]_{y=x+B_{\delta}}].
\end{align*}
Hence%
\begin{equation}
\label{e732} u(t,x)=\mathbb{\hat{E}}[u(t-\delta,x+B_{\delta})].
\end{equation}
{Thus}%
\begin{align*}
|u(t,x)-u(t-\delta,x)| &  =|\mathbb{\hat{E}}[u(t-\delta,x+B_{\delta
})-u(t-\delta,x)]|\\
&  \leq \mathbb{\hat{E}}[|u(t-\delta,x+B_{\delta})-u(t-\delta,x)|]\\
&  \leq \mathbb{\hat{E}}[C|B_{\delta}|]\leq
C\sqrt{G(0,I)+1}\sqrt{\delta}.
\end{align*}
To prove that $u$ is a viscosity solution of {(\ref{e731})}, {we fix
a
$(t,x)\in(0,\infty)\times \mathbb{R}^{d}$} and{ let $v\in C_{b}^{2,3}%
([0,\infty)\times \mathbb{R}^{d})$ be such that $v\geq u$ and
$v(t,x)=u(t,x)$. From (\ref{e732}), we have
\[
v(t,x)=\mathbb{\hat{E}}[u(t-\delta,x+B_{\delta})]\leq \mathbb{\hat{E}%
}[v(t-\delta,x+B_{\delta})].
\]
Therefore, by Taylor's expansion,
\begin{align*}
0 &  \leq \mathbb{\hat{E}}[v(t-\delta,x+B_{\delta})-v(t,x)]\\
&  =\mathbb{\hat{E}}[v(t-\delta,x+B_{\delta})-v(t,x+B_{\delta}%
)+(v(t,x+B_{\delta})-v(t,x))]\\
&  =\mathbb{\hat{E}}[-\partial_{t}v(t,x)\delta+\langle
Dv(t,x),B_{\delta }\rangle+\frac{1}{2}\langle
D^{2}v(t,x)B_{\delta},B_{\delta}\rangle+I_{\delta
}]\\
&  \leq-\partial_{t}v(t,x)\delta+\mathbb{\hat{E}}[\langle
Dv(t,x),B_{\delta }\rangle+\frac{1}{2}\langle
D^{2}v(t,x)B_{\delta},B_{\delta}\rangle
]+\mathbb{\hat{E}}[I_{\delta}],
\end{align*}
where}%
\begin{align*}
I_{\delta} &  =\int_{0}^{1}-[\partial_{t}v(t-\beta
\delta,x+B_{\delta
})-\partial_{t}v(t,x)]\delta d\beta \\
& \ \ \ +\int_{0}^{1}\int_{0}^{1}\langle(D^{2}v(t,x+\alpha \beta B_{\delta}%
)-D^{2}v(t,x))B_{\delta},B_{\delta}\rangle \alpha d\beta d\alpha.
\end{align*}
With the assumption \textup{(iii)} we can check that $\lim_{\delta
\downarrow 0}\mathbb{\hat{E}}[|I_{\delta}|]\delta^{-1}=0$, from
which we get $\partial_{t}v(t,x)-G(Dv(t,x),D^{2}v(t,x))\leq0$, hence
$u$ is a viscosity subsolution of {(\ref{e731}). We can analogously
prove that }$u$ is a
viscosity supersolution. Thus $u$ is a viscosity solution and $(B_{t}%
)_{t\geq0}$ is a generalized $G$-Brownian motion.}
\end{proof}

In many situations we are interested in a generalized
$2d$-dimensional
Brownian motion $(B_{t},b_{t})_{t\geq 0}$ such that $\mathbb{\hat{E}}%
[B_{t}]=-\mathbb{\hat{E}}[-B_{t}]=0$ and $\mathbb{\hat{E}}%
[|b_{t}|^{2}]/t\rightarrow 0$, as $t\downarrow 0$. In this case $B$ is in fact a $%
G$-Brownian motion defined on Definition 2.1 of Chapter \ref{ch2}.
Moreover
the process $b$ satisfies properties of Proposition 5.2. We define $u(t,x,y)=%
\mathbb{\hat{E}}[\varphi (x+B_{t},y+b_{t})]$.  By the above
proposition it follows that $u$ is the solution of the PDE%
\begin{equation*}
\partial _{t}u=G(D_{y}u,D_{xx}^{2}u),\ \ u|_{t=0}=\varphi \in C_{l.Lip}(%
\mathbb{R}^{2d}).
\end{equation*}%
where $G$ is a sublinear function of $(p,A)\in \mathbb{R}^{d}$,
defined by
\begin{equation*}
G(p,A):=\mathbb{\hat{E}}[\left\langle p,b_{t}\right\rangle
+\left\langle AB_{t},B_{t}\right\rangle ].
\end{equation*}%
Here $\left\langle \cdot ,\cdot \right\rangle =\left\langle \cdot
,\cdot \right\rangle _{\mathbb{R}^{d}}$.

\section{$\widetilde {G}$-Brownian Motion under a Nonlinear Expectation}

We can also define a $G$-Brownian motion on a nonlinear expectation space $%
(\Omega ,\mathcal{H},\mathbb{\widetilde {E}})$.

\begin{definition}
\label{defIII901}A $d$-dimensional process $(B_{t})_{t\geq 0}$ on a
nonlinear expectation space $(\Omega ,\mathcal{H},\mathbb{\widetilde
{E}})$ is called a (nonlinear)\textbf{\ $\widetilde {G}$-Brownian
motion} \index{Nonlinear G-Brownian motion} if the following
properties are satisfied: \newline \textup{(i)} $B_{0}(\omega
)=0$;\newline \textup{(ii)} For each $t,s\geq 0$, the increment
$B_{t+s}-B_{t}$
identically distributed with $B_{s}$ and is independent from $%
(B_{t_{1}},B_{t_{2}},\cdots ,B_{t_{n}})$, for each $n\in \mathbb{N}$ and $%
0\leq t_{1}\leq \cdots \leq t_{n}\leq t$;\newline
{{\textup{(iii)} $\lim_{t\downarrow 0}\mathbb{%
\hat{E}}[|B_{t}|^{3}]t^{-1}=0$.}}
\end{definition}

The following theorem gives a characterization of the nonlinear $\widetilde{G}$%
-Brownian motion, and give us the generator $\widetilde{G}$ of our $\widetilde{G}$%
-Brownian motion.

\begin{theorem}\label{ThmIII8.2}
Let $\mathbb{\widetilde{E}}$ be a nonlinear expectation and
$\mathbb{\hat{E}}$
be a sublinear expectation defined on $(\Omega ,\mathcal{H})$. let $\mathbb{%
\widetilde{E}}$ be dominated by $\mathbb{\hat{E}}$, namely%
\begin{equation*}
\mathbb{\widetilde{E}}[X]-\mathbb{\widetilde{E}}[Y]\leq
\mathbb{\hat{E}}[X-Y],\ \ X,Y\in \mathcal{H}.
\end{equation*}%
Let $(B_{t},b_{t})_{t\geq 0}$ be a given $\mathbb{R}^{2d}$--valued
$\widetilde{G}
$-Brownian motion on $(\Omega ,\mathcal{H},\mathbb{\widetilde{E}})$ such that $%
\mathbb{\hat{E}}[B_{t}]=\mathbb{\hat{E}}[-B_{t}]=0$ and
$\lim_{t\rightarrow 0}\mathbb{\hat{E}}[|b_{t}|^{2}]/t=0$. Then, {for
each fixed $\varphi \in
C_{b.Lip}(\mathbb{R}^{2d}\mathbb{)}$}, the function{%
\begin{equation*}
\tilde{u}(t,x,y):=\mathbb{\widetilde{E}}[\varphi
(x+B_{t},y+b_{t})],\ \ (t,x,y)\in \lbrack 0,\infty )\times
\mathbb{R}^{2d}
\end{equation*}%
is the viscosity solution of the following parabolic PDE}:
\begin{equation}
\partial _{t}\tilde{u}-\widetilde{G}(D_{y}\tilde{u},D_{x}^{2}\tilde{u})=0,\ \
u|_{t=0}=\varphi .  \label{eIII831}
\end{equation}%
where
\begin{equation*}
\widetilde{G}(p,A)=\mathbb{\widetilde{E}}[\left\langle p,b_{1}\right\rangle +\frac{1%
}{2}\langle AB_{1},B_{1}\rangle ],\ (p,A)\in \mathbb{R}^{d}\times \mathbb{S}%
(d).
\end{equation*}
\end{theorem}

\begin{remark}
\bigskip Let $G(p,A):=\mathbb{\hat{E}}[\left\langle p,b_{1}\right\rangle +%
\frac{1}{2}\langle AB_{1},B_{1}\rangle ]${. Then the function
}$\widetilde{G}$
is dominated by the sublinear function $G$ in the following sense:%
\begin{equation}
\widetilde{G}(p,A)-\widetilde{G}(p^{\prime },A^{\prime })\leq
G(p-p^{\prime
},A-A^{\prime }),\ \ (p,A),\ (p^{\prime },A^{\prime })\in \mathbb{R}%
^{d}\times \mathbb{S}(d).  \label{eIIIGdom}
\end{equation}
\end{remark}

\begin{proof}[Proof of Theorem \ref{ThmIII8.2}]
We set%
\begin{equation*}
f(t)=f_{A,t}(t):=\mathbb{\widetilde{E}}[\left\langle p,b_{t}\right\rangle +\frac{%
1}{2}\langle AB_{t},B_{t}\rangle ],\ t\geq 0.
\end{equation*}%
Since
\begin{equation*}
|f(t+h)-f(t)|\leq \mathbb{\hat{E}}%
[(|p|+2|A||B_{t}|)|B_{t+h}-B_{t}|+|A||B_{t+h}-B_{t}|^{2}]\rightarrow
0,
\end{equation*}%
we get that $f(t)$ is a continuous function. Since $\mathbb{\hat{E}}[B_{t}]=%
\mathbb{\hat{E}}[-B_{t}]=0$, it follows from Proposition
\ref{PropI.3.4}
that $\mathbb{\widetilde{E}}[X+\left\langle p,B_{t}\right\rangle ]=\mathbb{%
\widetilde{E}}[X]$ for each $X\in \mathcal{H}$ and $p\in
\mathbb{R}^{d}$. Thus
\begin{eqnarray*}
f(t+h) &=&\mathbb{\widetilde{E}}[\left\langle
p,b_{t+h}-b_{t}\right\rangle
+\left\langle p,b_{t}\right\rangle  \\
&&+\frac{1}{2}\langle AB_{t+h}-B_{t},B_{t+h}-B_{t}\rangle +\frac{1}{2}%
\left\langle AB_{t},B_{t}\right\rangle ] \\
&=&\mathbb{\widetilde{E}}[\left\langle p,b_{h}\right\rangle
+\frac{1}{2}\langle
AB_{h},B_{h}\rangle ]+\mathbb{\widetilde{E}}[\left\langle p,b_{t}\right\rangle +%
\frac{1}{2}\langle AB_{t},B_{t}\rangle ] \\
&=&f(t)+f(h).
\end{eqnarray*}%
It then follows that $f(t)=f(1)t=\widetilde{G}(A,p)t.${\ {We now
prove that the function }$u$ is Lipschitz in $x$ and uniformly
continuous in
$t$. In fact, {for each fixed }$t$,{\ }$u(t,\cdot )\in ${$C_{b.Lip}(\mathbb{R%
}^{d}\mathbb{)}$ }since%
\begin{align*}
& |\mathbb{\widetilde{E}}[\varphi
(x+B_{t},y+b_{t})]-\mathbb{\widetilde{E}}[\varphi
(x^{\prime }+B_{t},y^{\prime }+b_{t})]| \\
& \leq \mathbb{\hat{E}}[|\varphi (x+B_{t},y+b_{t})-\varphi
(x^{\prime }+B_{t},y^{\prime }+b_{t})|]\leq C(|x-x^{\prime
}|+|y-y^{\prime }|).
\end{align*}%
{For each $\delta \in \lbrack 0,t]$, since $(B_{t}-B_{\delta
},b_{t}-b_{\delta })$ is independent from {$(B_{\delta },b_{\delta
})$},}
\begin{align*}
\tilde{u}(t,x)& =\mathbb{\widetilde{E}}[\mathbb{\varphi
(}x+B_{\delta
}+(B_{t}-B_{\delta }),y+b_{\delta }+(b_{t}-b_{\delta })] \\
& =\mathbb{\widetilde{E}}[\mathbb{\widetilde{E}}[\varphi
(\bar{x}+(B_{t}-B_{\delta }),\bar{y}+(b_{t}-b_{\delta
}))]_{\bar{x}=x+B_{\delta },\bar{y}=y+b_{\delta }}].
\end{align*}%
Hence%
\begin{equation}
u(t,x)=\mathbb{\hat{E}}[u(t-\delta ,x+B_{\delta },y+b_{\delta })].
\label{eIII832}
\end{equation}%
{Thus}%
\begin{align*}
|\tilde{u}(t,x,y)-\tilde{u}(t-\delta ,x,y)|& =|\mathbb{\tilde{E}}[\tilde{u}%
(t-\delta ,x+B_{\delta },y+b_{\delta })-\tilde{u}(t-\delta ,x,y)]| \\
& \leq \mathbb{\hat{E}}[|\tilde{u}(t-\delta ,x+B_{\delta },y+b_{\delta })-%
\tilde{u}(t-\delta ,x,y)|] \\
& \leq \mathbb{\hat{E}}[C(|B_{\delta }|+|b_{\delta }|)].
\end{align*}%
It follows from (iii) of Definition \ref{defIII901} that  $u(t,x,y)$
is continuous in $t$ uniformly in $(t,x)\in [0,\infty)\times
\mathbb{R}^{2d}$.

To prove that $\tilde{u}$ is a viscosity solution of
{(\ref{eIII831})}, {we fix a $(t,x,y)\in (0,\infty )\times
\mathbb{R}^{2d}$} and{\ let $v\in
C_{b}^{2,3}([0,\infty )\times \mathbb{R}^{2d})$ be such that $v\geq u$ and $%
v(t,x,y)=u(t,x,y)$. From (\ref{eIII832}), we have
\begin{equation*}
v(t,x,y)=\mathbb{\widetilde{E}}[u(t-\delta ,x+B_{\delta
},y+b_{\delta })]\leq \mathbb{\widetilde{E}}[v(t-\delta ,x+B_{\delta
},y+b_{\delta })].
\end{equation*}%
Therefore, by Taylor's expansion,
\begin{align*}
0& \leq \mathbb{\widetilde{E}}[v(t-\delta ,x+B_{\delta },y+b_{\delta
})-v(t,x)]
\\
& =\mathbb{\widetilde{E}}[v(t-\delta ,x+B_{\delta },y+b_{\delta
})-v(t,x+B_{\delta },y+b_{\delta }) \\
& \ \ \ +(v(t,x+B_{\delta },y+b_{\delta })-v(t,x,y)] \\
& =\mathbb{\widetilde{E}}[-\partial _{t}v(t,x,y)\delta +\langle
D_{y}v(t,x,y),b_{\delta }\rangle +\langle \partial
_{x}v(t,x,y),B_{\delta }\rangle +\frac{1}{2}\langle
D_{xx}^{2}v(t,x,y)B_{\delta },B_{\delta
}\rangle +I_{\delta }] \\
& \leq -\partial _{t}v(t,x,y)\delta +\mathbb{\widetilde{E}}[\langle
D_{y}v(t,x,y),b_{\delta }\rangle +\frac{1}{2}\langle
D_{xx}^{2}v(t,x,y)B_{\delta },B_{\delta }\rangle ]+\mathbb{\hat{E}}%
[I_{\delta }],
\end{align*}%
where}%
\begin{align*}
I_{\delta }& =\int_{0}^{1}-[\partial _{t}v(t-\delta \gamma
,x+B_{\delta
},y+b_{\delta })-\partial _{t}v(t,x,y)]\delta d\gamma  \\
& +\int_{0}^{1}\left\langle \partial _{y}v(t,x+\gamma B_{\delta
},y+\gamma
b_{\delta })-\partial _{y}v(t,x,y),b_{\delta }\right\rangle d\gamma  \\
& +\int_{0}^{1}\left\langle \partial _{x}v(t,x,y+\gamma b_{\delta
})-\partial _{x}v(t,x,y),B_{\delta }\right\rangle d\gamma  \\
& +\int_{0}^{1}\int_{0}^{1}\langle (D_{xx}^{2}v(t,x+\alpha \gamma
B_{\delta },y+\gamma b_{\delta })-D_{xx}^{2}v(t,x,y))B_{\delta
},B_{\delta }\rangle \gamma d\gamma d\alpha .
\end{align*}%
With the assumption \textup{(iii)} we can check that $\lim_{\delta
\downarrow 0}\mathbb{\hat{E}}[|I_{\delta }|]\delta ^{-1}=0$, from
which we get $\partial _{t}v(t,x)-G(Dv(t,x),D^{2}v(t,x))\leq 0$,
hence $u$ is a viscosity subsolution of {(\ref{eIII831}). We can
analogously prove that }$u$ is a viscosity supersolution. Thus $u$
is a viscosity solution. }
\end{proof}

\section{Construction of ${\widetilde G}$-Brownian Motions under Nonlinear Expectation}

Let $G(\cdot ):\mathbb{R}^{d}\times \mathbb{S}(d)\rightarrow
\mathbb{R}$ be
a given sublinear function monotonic on $A\in \mathbb{S}(d)$ and $\widetilde{G}%
(\cdot ):\mathbb{R}^{d}\times \mathbb{S}(d)\rightarrow \mathbb{R}$
be a given function dominated by $G$ in the sense of
(\ref{eIIIGdom}). The
construction of a $\mathbb{R}^{2d}$-dimensional $\widetilde{G}$-Brownian motion $%
(B_{t},b_{t})_{t\geq 0}$ under a nonlinear expectation
$\mathbb{\widetilde{E}}$, dominated by a sublinear expectation
$\mathbb{\hat{E}}$ is based on a
similar approach introduced in Section 2. In fact we will see that by our construction $%
(B_{t},b_{t})_{t\geq 0}$ is also a $G$-Brownian motion of the
sublinear expectation $\mathbb{\hat{E}}$.

We denote by $\Omega =C_{0}^{2d}(\mathbb{R}^{+})$ the space of all $\mathbb{R%
}^{2d}$--valued continuous paths $(\omega _{t})_{t\in
\mathbb{R}^{+}}$. {For each fixed }$T\in \lbrack 0,\infty ),$ we set
$\Omega _{T}:=\{\omega _{\cdot
\wedge T}:\omega \in \Omega \}.$ We will consider the canonical process $%
(B_{t},b_{t})(\omega )=\omega _{t}$, $t\in \lbrack 0,\infty )$, for
$\omega \in \Omega $. We also follow section 2 to introduce the
spaces of random
variables $L_{ip}(\Omega _{T})$ and $L_{ip}(\Omega )$ so that to define $%
\mathbb{\hat{E}}$ and $\mathbb{\widetilde{E}}$ on $(\Omega
,L_{ip}(\Omega ))$.

To this purpose we first construct a sequence of $d$-dimensional
random vectors $(X_{i},\eta ${{$_{i})_{i=1}^{\infty }$ }}on a
sublinear expectation space $(\overline{\Omega
},\overline{\mathcal{H}},\overline{\mathbb{E}})$
such that $(X_{i},\eta ${{$_{i})$ is }}$G$-distributed and $(X_{i+1},\eta ${{%
$_{i+1})$ is independent from $((X_{1},\eta _{1}),\cdots
,(X_{i},\eta _{i}))$ for each $i=1,2,\cdots $. By the definition of
}}$G$-distribution the function
\begin{equation*}
u(t,x,y):=\mathbb{\hat{E}}[\varphi (x+\sqrt{t}X_{1},y+t\eta _{1})],\
\ t\geq 0,\ \ x,y\in \mathbb{R}^{d}
\end{equation*}%
is the viscosity solution of the following parabolic PDE, which is
the same as equation (\ref{ee03}) in Chap.II.
\begin{equation*}
\partial _{t}u-G(D_{y}u,D_{xx}^{2}u)=0,\ \ u|_{t=0}=\varphi \in C_{Lip}(%
\mathbb{R}^{2d}).
\end{equation*}%
We also consider the PDE (for the existence, uniqueness, comparison
and domination properties, see Theorem \ref{Com-G} in Appendix C).
\begin{equation*}
\partial _{t}\tilde{u}-\widetilde{G}(D_{y}\tilde{u},D_{xx}^{2}\tilde{u})=0,\ \ \tilde{u}%
|_{t=0}=\varphi \in C_{Lip}(\mathbb{R}^{2d}),
\end{equation*}%
and denote by $\widetilde{P}_{t}[\varphi ](x,y)=\tilde{u}(t,x,y)$.
Since $\widetilde{G}$ is dominated by $G$, it follows from the
domination theorem of viscosity
solutions, i.e., Theorem \ref{G-Tilde} in Appendix C, that, for each $%
\varphi ,\psi \in C_{b,Lip}(\mathbb{R}^{2d})$,
\begin{equation*}
\widetilde{P}_{t}[\varphi ](x,y)-\widetilde{P}_{t}[\psi ](x,y)\leq \overline{\mathbb{%
E}}[(\varphi -\psi )(x+\sqrt{t}X_{1},y+t\eta _{1})].
\end{equation*}

We now introduce a sublinear expectation $\hat{\mathbb{E}}$ and
a nonlinear $\mathbb{\widetilde{E}}$ defined on $L_{ip}(\Omega )$
via the following procedure: for each $X\in L_{ip}(\Omega )$ with
\begin{equation*}
X=\varphi (B_{t_{1}}-B_{t_{0}},b_{t_{1}}-b_{t_{0}},\cdots
,B_{t_{1}}-B_{t_{0}},b_{t_{n}}-b_{t_{n-1}})
\end{equation*}%
{\ for {$\varphi \in C_{l.Lip}(\mathbb{R}^{2d\times n})$ and $%
0=t_{0}<t_{1}<\cdots <t_{n}<\infty $, we set%
\begin{equation*}
\hat{\mathbb{E}}[\varphi
(B_{t_{1}}-B_{t_{0}},b_{t_{1}}-b_{t_{0}},\cdots
,B_{t_{n}}-B_{t_{n-1}},b_{t_{n}}-b_{t_{n-1}})]
\end{equation*}%
\begin{equation*}
:=\overline{\mathbb{E}}[\varphi
(\sqrt{t_{1}-t_{0}}X_{1},(t_{1}-t_{0})\eta _{1},\cdots
,\sqrt{t_{n}-t_{n-1}}X_{n},(t_{n}-t_{n-1})\eta _{n})].
\end{equation*}%
and}}%
\begin{equation*}
\mathbb{\widetilde{E}}[\varphi
(B_{t_{1}}-B_{t_{0}},b_{t_{1}}-b_{t_{0}},\cdots
,B_{t_{n}}-B_{t_{n-1}},b_{t_{n}}-b_{t_{n-1}})]=\varphi _{n}(0,0)
\end{equation*}%
where $\varphi _{n}\in C_{b.Lip}(\mathbb{R}^{2d})$ is defined
iteratively through
\begin{eqnarray*}
\varphi _{1}(x_{1},y_{1},\cdots ,x_{n-1},y_{n-1}) &=&\widetilde{P}%
_{t_{n}-t_{n-1}}[\varphi _{1}(x_{1},y_{1},\cdots
,x_{n-1},y_{n-1},\cdot
)](0,0), \\
&&\vdots  \\
\varphi _{n-1}(x_{1},y_{1}) &=&\widetilde{P}_{t_{2}-t_{1}}[\varphi
_{n-2}(x_{1},y_{1},\cdot )](0,0), \\
\varphi _{n}(x_{1},y_{1}) &=&\widetilde{P}_{t_{2}-t_{1}}[\varphi
_{n-1}(\cdot )](x_{1},y_{1}).
\end{eqnarray*}%
{The related conditional expectation of $X=$}$\varphi
(B_{t_{1}}-B_{t_{0}},b_{t_{1}}-b_{t_{0}},\cdots
,B_{t_{n}}-B_{t_{n-1}},b_{t_{n}}-b_{t_{n-1}})${{\ under $\Omega
_{t_{j}}$ is
defined by%
\begin{align}
\hat{\mathbb{E}}[X|{{\Omega _{t_{j}}}}]& =\hat{\mathbb{E}}[\varphi
(B_{t_{1}}-B_{t_{0}},b_{t_{1}}-b_{t_{0}},\cdots
,B_{t_{n}}-B_{t_{n-1}},b_{t_{n}}-b_{t_{n-1}}){\ }|{{\Omega
_{t_{j}}}}]
\label{Condition9} \\
& :=\psi (B_{t_{1}}-B_{t_{0}},b_{t_{1}}-b_{t_{0}},\cdots
,B_{t_{j}}-B_{t_{j-1}},b_{t_{j}}-b_{t_{j-1}}),  \notag
\end{align}%
where}}%
\begin{equation*}
\psi (x_{1},\cdots ,x_{j})=\overline{\mathbb{E}}[\varphi
(x_{1},\cdots
,x_{j},\sqrt{t_{j+1}-t_{j}}X_{j+1},(t_{1}-t_{0})\eta _{j+1},\cdots ,\sqrt{%
t_{n}-t_{n-1}}X_{n},(t_{1}-t_{0})\eta _{n})].
\end{equation*}%
Similarly
\begin{equation*}
\mathbb{\widetilde{E}}[X|{{\Omega _{t_{j}}}}]=\varphi
_{n-j}(B_{t_{1}}-B_{t_{0}},b_{t_{1}}-b_{t_{0}},\cdots
,B_{t_{j}}-B_{t_{j-1}},b_{t_{j}}-b_{t_{j-1}}).
\end{equation*}%
{It is easy to check that }$\hat{\mathbb{E}}${$[\cdot ]$ (resp. $\mathbb{%
\widetilde{E}}$}) {consistently defines a sublinear (resp.
nonlinear) expectation }and $\mathbb{\widetilde{E}}${$[\cdot ]$ on
}$(\Omega ,L_{ip}(\Omega
))$. Moreover {$(B_{t},b_{t})_{t\geq 0}$ is a $G$-Brownian motion under }$%
\mathbb{\hat{E}}$ and a $\widetilde{G}$-Brownian motion under
$\mathbb{\widetilde{E}} $.

\begin{proposition}
\label{Prop-1-9-2}{We also list the properties of }$\mathbb{\widetilde{E}}${$%
[\cdot |\Omega _{t}]$ that hold for each $X,Y\in $}$L_{ip}(\Omega )${:}%
\newline
\textup{(i) }{\textbf{\ }{If $X\geq Y$, then }}$\mathbb{\widetilde{E}}${{$%
[X|\Omega _{t}]\geq \mathbb{\widetilde{E}}[Y|{\Omega
}_{t}]$.\newline
}\textup{(ii) }\textbf{\ }}$\mathbb{\widetilde{E}}${$[X+\eta |\Omega _{t}]=$}$%
\mathbb{\widetilde{E}}$$[X|\Omega _{t}]+\eta ${, \ for each }$t\geq 0$ and {$%
\eta \in $}$L_{ip}(\Omega _{t})${{.\newline
\textup{(iii) }\textbf{\ }}}$\mathbb{\widetilde{E}}${{$[X|\Omega _{t}]-\mathbb{%
\widetilde{E}}[Y|\Omega _{t}]\leq \hat{\mathbb{E}}[X-Y|\Omega
_{t}].$\newline }}\textup{(iv)}
$\mathbb{\widetilde{E}}{{[\mathbb{\widetilde{E}}[X|\Omega
_{t}]|\Omega _{s}]=\mathbb{\widetilde{E}}[X|\Omega _{t\wedge s}],\
}}$ in
particular$,\ \mathbb{\widetilde{E}}[\mathbb{\widetilde{E}}[X|\Omega _{t}]]=\mathbb{%
\widetilde{E}}[X].$\newline
(v) For each $X\in L_{ip}(\Omega ^{t})$, $\mathbb{\widetilde{E}}[X|{\Omega }%
_{t}]=\mathbb{\widetilde{E}}[X]$, where $L_{ip}(\Omega ^{t})$ is the
linear space of random variables with the form
\begin{align*}
& {\varphi (W_{t_{2}}-W_{t_{1}},W_{t_{3}}-W_{t_{2}},\cdots
,W_{t_{n+1}}-W_{t_{n}}),} \\
\ & \ \ n=1,2,\cdots ,\ \varphi \in C_{l.Lip}(\mathbb{R}^{d\times
n}),\ t_{1},\cdots ,t_{n},t_{n+1}\in \lbrack t,\infty ).
\end{align*}
\end{proposition}

Since $\mathbb{\hat{E}}$ can be considered as a special nonlinear
expectation of $\mathbb{\widetilde{E}}$ dominated by its self, thus $\mathbb{%
\hat{E}}[\cdot |\Omega _{t}]$ also satisfies above properties
(i)--(v). Moreover

\begin{proposition}
\label{Prop-1-9-1}{The conditional sublinear expectation }$\hat{\mathbb{E}}${%
$[\cdot |\Omega _{t}]$ satisfies (i)-(v). Moreover }$\hat{\mathbb{E}}$$%
[\cdot |\Omega _{t}]$ itself is sublinear, i.e.,{\ }\newline
{{\textup{(vi) }\textbf{\ }}}$\hat{\mathbb{E}}${{$[X|\Omega _{t}]-\hat{%
\mathbb{E}}[Y|\Omega _{t}]\leq \hat{\mathbb{E}}[X-Y|\Omega _{t}],\ \ .$%
\newline
\textup{(vii) }{\ }}}$\hat{\mathbb{E}}${{{$[\eta X|\Omega _{t}]=\eta ^{+}%
\hat{\mathbb{E}}[X|\Omega _{t}]+\eta ^{-}\hat{\mathbb{E}}[-X|\Omega
_{t}]$ for each $\eta \in L_{ip}(\Omega _{t}).$}}}\newline
\end{proposition}

We now consider the completion of sublinear expectation space
$(\Omega ,L_{ip}(\Omega),\hat{\mathbb{E}}).$

We denote by $L_{G}^{p}(\Omega)$, $p\geq1$, the completion of $%
L_{ip}(\Omega) $ under the norm $\left \Vert X\right \Vert _{p}:=(\hat{%
\mathbb{E}}[|X|^{p}])^{1/p}$. Similarly, we can define $L_{G}^{p}(%
\Omega_{T}) $, $L_{G}^{p}(\Omega_{T}^{t})$ and
$L_{G}^{p}(\Omega^{t})$. It
is clear that for each $0\leq t\leq T<\infty$, $L_{G}^{p}(\Omega_{t})%
\subseteq L_{G}^{p}(\Omega_{T})\subseteq L_{G}^{p}(\Omega)$.

According to Sec.\ref{c1s5} in Chap.\ref{ch1},
$\hat{\mathbb{E}}[\cdot ]$ can be continuously extended to
$(\Omega,L_{G}^{1}(\Omega ))$. Moreover, since
$\mathbb{\widetilde{E}}$ is dominated by $\mathbb{\hat{E}}$, thus by
Definition \ref{DefI4.4} in Chap.I, $(\Omega ,L_{G}^{1}(\Omega ),\mathbb{%
\hat{E})}$ forms a sublinear expectation space and $(\Omega
,L_{G}^{1}(\Omega ),\mathbb{\widetilde{E})}$ forms a nonlinear
expectation space.

We now consider the extension of conditional $G$-expectation. For
each fixed
$t\leq T$, the conditional $G$-expectation ${{\hat{\mathbb{E}}[\cdot |{{%
\Omega }_{t}}]}}:L_{ip}(\Omega _{T})\rightarrow L_{ip}(\Omega _{t})$
is a continuous mapping under $\left\Vert \cdot \right\Vert $.
Indeed, we have
\begin{equation*}
{{\mathbb{\widetilde{E}}[X|{\Omega _{t}}]}}-{{\mathbb{\widetilde{E}}[Y|{\Omega _{t}}]%
}}\leq {{\hat{\mathbb{E}}[X-Y|{\Omega _{t}}]\leq \hat{\mathbb{E}}[|X-Y||{{%
\Omega }_{t}}],}}
\end{equation*}%
then
\begin{equation*}
|{{\mathbb{\widetilde{E}}[X|{\Omega _{t}}]}}-{{\mathbb{\widetilde{E}}[Y|{\Omega _{t}}%
]|}}\leq {{\hat{\mathbb{E}}[|X-Y||{\Omega _{t}}].}}
\end{equation*}%
We thus obtain
\begin{equation*}
\left\Vert {{\mathbb{\widetilde{E}}[X|{\Omega _{t}}]}}-{{\mathbb{\widetilde{E}}[Y|{%
\Omega _{t}}]}}\right\Vert \leq \left\Vert X-Y\right\Vert .
\end{equation*}%
It follows that ${{\mathbb{\widetilde{E}}[\cdot |{\Omega _{t}}]}}$
can be also extended as a continuous mapping
\begin{equation*}
{{\mathbb{\widetilde{E}}[\cdot |{\Omega _{t}}]}}:L_{G}^{1}(\Omega
_{T})\rightarrow L_{G}^{1}(\Omega _{t}).
\end{equation*}%
If the above $T$ is not fixed, then we can obtain ${{\mathbb{\widetilde{E}}%
[\cdot |{\Omega _{t}}]}}:L_{G}^{1}(\Omega )\rightarrow
L_{G}^{1}(\Omega _{t}) $.

\begin{remark}
The above proposition also holds for $X,Y\in L_{G}^{1}(\Omega)$. But in {{%
\textup{(iv)}}}, $\eta \in L_{G}^{1}(\Omega_{t})$ should be bounded, since $%
X,Y\in L_{G}^{1}(\Omega)$ does not imply $X\cdot Y\in
L_{G}^{1}(\Omega).$
\end{remark}

In particular, we have the following independence:
\begin{equation*}
{{\mathbb{\widetilde{E}}[X|{\Omega
_{t}}]}}={{\mathbb{\widetilde{E}}[X}}],\medskip \ \ \forall X\in
L_{G}^{1}(\Omega ^{t}).
\end{equation*}

We give the following definition similar to the classical one:

\begin{definition}
An $n$-dimensional random vector $Y\in (L_{G}^{1}({\Omega }))^{n}$
is said to be independent from $\Omega _{t}$ for some given $t$ if
for each $\varphi
\in C_{b.Lip}(\mathbb{R}^{n})$ we have%
\begin{equation*}
{{\mathbb{\widetilde{E}}}}[\varphi (Y)|{{{\Omega _{t}}}}]={{\mathbb{\widetilde{E}}}}%
[\varphi (Y)].
\end{equation*}
\end{definition}

\section*{Notes and Comments}
\addcontentsline{toc}{section}{Notes and Comments}

Bachelier (1900) \cite{Bachelier} proposed Brownian motion as a
model for fluctuations of the stock market, Einstein (1905)
\cite{Einstein} used Brownian motion to give experimental
confirmation of the atomic theory, and Wiener (1923) \cite{Wiener}
gave a mathematically rigorous construction of Brownian motion. Here
we follow Kolmogorov's idea (1956) \cite{K} to construct
$G$-Brownian motion by introducing infinite dimensional function
space and the corresponding family of infinite dimensional sublinear
distributions, instead of linear distributions in \cite{K}.

The notions of $G$-Brownian motion and the related stochastic
calculus of It\^o's type were firstly introduced by Peng (2006)
\cite{Peng2006a} for 1-dimensional case and then in (2008)
\cite{Peng2008} for multi-dimensional situation. It is very
interesting that Denis and Martini (2006) \cite{Denis-M} studied
super-pricing of contingent claims under model uncertainty of
volatility. They have introduced a norm on the space of continuous
paths $\Omega=C([0,T])$ which corresponds to our $L^2_G$-norm and
developed a stochastic integral. There is no notion of nonlinear
expectation and the related nonlinear distribution, such as
$G$-expectation, conditional $G$-expectation, the related $G$-normal
distribution and the notion of independence in their paper. But on
the other hand, powerful tools in capacity theory enable them to
obtain pathwise results for random variables and stochastic
processes through the language of ``quasi-surely'' (see e.g.
Dellacherie (1972) \cite{Del}, Dellacherie and Meyer (1978 and 1982)
\cite{DM}, Feyel and de La Pradelle (1989)
\cite{Feyel-DLP}) in place of ``almost surely'' in
classical probability theory.

The main motivations of $G$-Brownian motion were the pricing and risk
measures under volatility uncertainty in financial markets (see
Avellaneda, Levy and Paras (1995) \cite{Avellaneda} and Lyons (1995)
\cite{Lyons}). It was well-known that under volatility
uncertainty the corresponding uncertain probabilities are singular
from each other. This causes a serious problem for the related path
analysis to treat, e.g., path-dependent derivatives, under a
classical probability space. Our $G$-Brownian motion provides a
powerful tool to such type of problems.

Our new It\^o's calculus for $G$-Brownian motion is of course
inspired from It\^o's groundbreaking work since 1942 \cite{Ito} on
stochastic integration, stochastic differential equations and
stochastic calculus through interesting books cited in Chapter
\ref{ch5}. It\^o's formula given by Theorem \ref{Thm6.5} is from
Peng \cite{Peng2006a}, \cite{Peng2008}. Gao (2009)\cite{Gao} proved
a more general It\^o's formula for $G$-Brownian motion. An
interesting problem is: can we get an It\^o's formula in which the
conditions correspond the classical one? Recently Li and Peng have
solved this problem in \cite{Li-Peng2009}.

Using nonlinear Markovian semigroup known as Nisio's semigroup (see
Nisio (1976) \cite{Nisio2}), Peng (2005) \cite{Peng2005}
studied
 the processes with Markovian properties under a nonlinear
 expectation.

%
%
\chapter{$G$-martingales and Jensen's Inequality}
\label{ch5} In this chapter, we introduce the notion of
$G$-martingales and the related Jensen's inequality for a new type
of $G$-convex functions. Essentially different from the classical
situation, ``$M$ is a $G$-martingale'' does not imply that ``$-M$ is
a $G$-martingale''.

\section{The Notion of $G$-martingales}
We now give the notion of $G$--martingales.

\begin{definition}
A process $(M_{t})_{t\geq0}$ is called a
\textbf{$G$--martingale}\index{$G$-martingale} (respectively,
$G$\textbf{--supermartingale}\index{$G$-supermartingale},
$G$\textbf{--submartingale})\index{$G$-submartingale} if for each
$t\in \lbrack0,\infty)$, $M_{t}\in L_{G}^{1}(\Omega_{t})$ and for
each $s\in \lbrack0,t]$, we have
\[
\hat{\mathbb{E}}[M_{t}|\Omega_{s}]=M_{s}\  \  \
\text{(respectively,\ }\leq M_{s},\  \  \geq M_{s}).
\]

\end{definition}

\begin{example}
For each fixed $X\in L_{G}^{1}(\Omega)$, it is clear that $(\mathbb{\hat{E}%
}[X|\Omega_{t}])_{t\geq0}$ is a $G$--martingale.
\end{example}

\begin{example}
For each fixed $\mathbf{a}\in \mathbb{R}^{d}$, it is easy to check
that $(B_{t}^{\mathbf{a}})_{t\geq0}$ and
$(-B_{t}^{\mathbf{a}})_{t\geq0}$ are
$G$--martingales. The process $(\langle B^{\mathbf{a}}\rangle_{t}-\sigma_{\mathbf{aa}^{T}%
}^{2}t)_{t\geq0}$ is a $G$--martingale since
\begin{align*}
\hat{\mathbb{E}}[\langle B^{\mathbf{a}}\rangle_{t}-\sigma_{\mathbf{aa}%
^{T}}^{2}t|\Omega_{s}]  &  =\hat{\mathbb{E}}[\langle B^{\mathbf{a}}%
\rangle_{s}-\sigma_{\mathbf{aa}^{T}}^{2}t+(\langle
B^{\mathbf{a}}\rangle
_{t}-\langle B^{\mathbf{a}}\rangle_{s})|\Omega_{s}]\\
&  =\langle B^{\mathbf{a}}\rangle_{s}-\sigma_{\mathbf{aa}^{T}}^{2}%
t+\hat{\mathbb{E}}[\langle B^{\mathbf{a}}\rangle_{t}-\langle
B^{\mathbf{a}}\rangle_{s}]\\
&  =\langle B^{\mathbf{a}}\rangle_{s}-\sigma_{\mathbf{aa}^{T}}^{2}s.
\end{align*}
Similarly we can show that $(-(\langle B^{\mathbf{a}}\rangle_{t}%
-\sigma_{\mathbf{aa}^{T}}^{2}t))_{t\geq0}$ is a $G$--submartingale.
The process $((B_{t}^{\mathbf{a}})^{2})_{t\geq0}$ is a $G$--submartingale since%
\begin{align*}
\hat{\mathbb{E}}[(B_{t}^{\mathbf{a}})^{2}|\Omega_{s}]  &
=\mathbb{\hat
{E}}[(B_{s}^{\mathbf{a}})^{2}+(B_{t}^{\mathbf{a}}-B_{s}^{\mathbf{a}}%
)^{2}+2B_{s}^{\mathbf{a}}(B_{t}^{\mathbf{a}}-B_{s}^{\mathbf{a}})|\Omega_{s}]\\
&  =(B_{s}^{\mathbf{a}})^{2}+\hat{\mathbb{E}}[(B_{t}^{\mathbf{a}}%
-B_{s}^{\mathbf{a}})^{2}|\Omega_{s}]\\
&  =(B_{s}^{\mathbf{a}})^{2}+\sigma_{\mathbf{aa}^{T}}^{2}(t-s)\geq
(B_{s}^{\mathbf{a}})^{2}.
\end{align*}
Similarly we can prove that $((B_{t}^{\mathbf{a}})^{2}-\sigma_{\mathbf{aa}%
^{T}}^{2}t)_{t\geq0}$ and $((B_{t}^{\mathbf{a}})^{2}-\langle B^{\mathbf{a}%
}\rangle_{t})_{t\geq0}$ are $G$--martingales.
\end{example}

In general, we have the following important property.

\begin{proposition}
\label{ch5p1}
Let $M_{0}\in \mathbb{R}$, $\varphi=(\varphi^{j})_{j=1}^{d}\in M_{G}%
^{2}(0,T;\mathbb{R}^{d})$ and $\eta=(\eta^{ij})_{i,j=1}^{d}\in M_{G}%
^{1}(0,T;\mathbb{S}(d))$ be given and let%
\[
M_{t}=M_{0}+\int_{0}^{t}\varphi_{u}^{j}dB_{u}^{j}+\int_{0}^{t}\eta_{u}%
^{ij}d\left \langle B^{i},B^{j}\right \rangle
_{u}-\int_{0}^{t}2G(\eta _{u})du\ \ \text{for}\ t\in \lbrack0,T].
\]
Then $M$ is a $G$--martingale. Here we still use the \index{Einstein
convention}, i.e., the above repeated indices $i$ and $j$  imply the
summation.
\end{proposition}

\begin{proof}
Since $\hat{\mathbb{E}}[\int_{s}^{t}\varphi_{u}^{j}dB_{u}^{j}|\Omega
_{s}]=\hat{\mathbb{E}}[-\int_{s}^{t}\varphi_{u}^{j}dB_{u}^{j}|\Omega
_{s}]=0$, we only need to prove that%
\[
\bar{M}_{t}=\int_{0}^{t}\eta_{u}^{ij}d\left \langle
B^{i},B^{j}\right \rangle _{u}-\int_{0}^{t}2G(\eta_{u})du\ \
\text{for}\ t\in \lbrack0,T]
\]
is a $G$--martingale. It suffices to consider the case where $\eta
\in M_{G}^{1,0}(0,T;\mathbb{S}(d))$, i.e.,
\[
\eta_{t}=\sum_{k=0}^{N-1}\eta_{t_{k}}I_{[t_{k},t_{k+1})}(t).
\]
We have, for $s\in \lbrack t_{N-1},t_{N}]$,%
\begin{align*}
\hat{\mathbb{E}}[\bar{M}_{t}|\Omega_{s}]  &  =\bar{M}_{s}+\mathbb{\hat{E}%
}[(\eta_{t_{N-1}},\left \langle B\right \rangle _{t}-\left \langle
B\right \rangle _{s})-2G(\eta_{t_{N-1}})(t-s)|\Omega_{s}]\\
&  =\bar{M}_{s}+\hat{\mathbb{E}}[(A,\left \langle B\right \rangle
_{t}-\left \langle B\right \rangle _{s})]_{A=\eta_{t_{N-1}}}-2G(\eta_{t_{N-1}%
})(t-s)\\
&  =\bar{M}_{s}.
\end{align*}
Then we can repeat this procedure backwardly to prove the result for
$s\in \lbrack0,t_{N-1}]$.
\end{proof}

\begin{corollary}
Let $\eta \in M_{G}^{1}(0,T)$. Then for each fixed $\mathbf{a}\in \mathbb{R}%
^{d}$, we have%
\begin{equation}
\sigma_{-\mathbf{aa}^{T}}^{2}\hat{\mathbb{E}}[\int_{0}^{T}|\eta
_{t}|dt]\leq \hat{\mathbb{E}}[\int_{0}^{T}|\eta_{t}|d\langle B^{\mathbf{a}%
}\rangle_{t}]\leq \sigma_{\mathbf{aa}^{T}}^{2}\hat{\mathbb{E}}[\int_{0}%
^{T}|\eta_{t}|dt].
\end{equation}

\end{corollary}

\begin{proof}
For each $\xi \in M_{G}^{1}(0,T)$, by the above proposition, we have%
\[
\hat{\mathbb{E}}[\int_{0}^{T}\xi_{t}d\langle B^{\mathbf{a}}\rangle
_{t}-\int_{0}^{T}2G_{\mathbf{a}}(\xi_{t})dt]=0,
\]
where $G_{\mathbf{a}}(\alpha)=\frac{1}{2}(\sigma_{\mathbf{aa}^{T}}^{2}%
\alpha^{+}-\sigma_{-\mathbf{aa}^{T}}^{2}\alpha^{-})$. Letting
$\xi=|\eta|$ and
$\xi=-|\eta|$, we get%
\begin{align*}
\hat{\mathbb{E}}[\int_{0}^{T}|\eta_{t}|d\langle
B^{\mathbf{a}}\rangle
_{t}-\sigma_{\mathbf{aa}^{T}}^{2}\int_{0}^{T}|\eta_{t}|dt]  &  =0,\\
\hat{\mathbb{E}}[-\int_{0}^{T}|\eta_{t}|d\langle
B^{\mathbf{a}}\rangle
_{t}+\sigma_{-\mathbf{aa}^{T}}^{2}\int_{0}^{T}|\eta_{t}|dt]  &  =0.
\end{align*}
From the sub-additivity of $G$-expectation,  we can easily get the
result.
\end{proof}

\begin{remark}
It is worth to mention that for a $G$--martingale $M$, in general,
$-M$ is not a $G$--martingale. But in Proposition \ref{ch5p1}, when
$\eta \equiv0$, $-M$ is still a $G$--martingale.
\end{remark}

\begin{exercise}

\textup{(a)} Let $(M_{t})_{t\geq0}$ be a $G$--supermartingale. Show that $(-M_{t}%
)_{t\geq0}$ is a $G$--submartingale.

\textup{(b)} Find a $G$--submartingale $(M_{t})_{t\geq0}$ such that
$(-M_{t})_{t\geq0}$ is not a $G$--supermartingale.
\end{exercise}

\begin{exercise}

\textup{(a)} Let $(M_{t})_{t\geq0}$ and $(N_{t})_{t\geq0}$ be two $G$%
--supermartingales. Prove that $(M_{t}+N_{t})_{t\geq0}$ is a $G$%
--supermartingale.

\textup{(b)} Let $(M_{t})_{t\geq0}$ and $(-M_{t})_{t\geq0}$ be two
$G$--martingales. For each $G$--submartingale (respectively,
$G$--supermartingale) $(N_{t})_{t\geq0}$, prove that $(M_{t}
+N_{t})_{t\geq0}$ is a $G$--submartingale (respectively,
$G$--supermartingale).
\end{exercise}

\section{On $G$-martingale Representation Theorem}

How to give a $G$-martingale representation theorem is still a
largely open problem. Xu and Zhang (2009) \cite{Xu-Zhang} have obtained a
martingale representation for a special `symmetric' $G$-martingale
process. A more general situation have been proved by Soner, Touzi
and Zhang (preprint in private communications). Here we present the
formulation of this $G$-martingale representation theorem under a
very strong assumption.

In this section, we consider the generator
$G:\mathbb{S}(d)\rightarrow \mathbb{R}$ satisfying the uniformly
elliptic condition, i.e., there exists a $\beta>0$ such that, for
each $A,\bar{A}\in \mathbb{S}(d)$ with $A\geq \bar{A}$,
\[
G(A)-G(\bar{A})\geq \beta \mathrm{tr}[A-\bar{A}].
\]

For each $\xi=(\xi^{j})_{j=1}^{d}\in M_{G}^{2}(0,T;\mathbb{R}^{d})$
and $\eta=(\eta^{ij})_{i,j=1}^{d}\in M_{G}^{1}(0,T;\mathbb{S}(d))$,
we use the following notations
\[
\int_{0}^{T}\langle
\xi_{t},dB_{t}\rangle:=\sum_{j=1}^{d}\int_{0}^{T}\xi
_{t}^{j}dB_{t}^{j};\  \int_{0}^{T}(\eta_{t},d\langle B\rangle_{t}%
):=\sum_{i,j=1}^{d}\int_{0}^{T}\eta_{t}^{ij}d\left \langle B^{i},B^{j}%
\right \rangle _{t}.
\]

We first consider the representation of $\varphi(B_{T}-B_{t_{1}})$
for $0\leq t_{1}\leq T<\infty$.

\begin{lemma}
\label{ch5l1} Let $\xi=\varphi(B_{T}-B_{t_{1}})$, $\varphi \in
C_{b.Lip}(\mathbb{R}^{d})$. Then we have the following
representation:
\[
\xi=\hat{\mathbb{E}}[\xi]+\int_{t_{1}}^{T}\langle \beta_{t},dB_{t}%
\rangle+\int_{t_{1}}^{T}(\eta_{t},d\langle B\rangle_{t})-\int_{t_{1}}%
^{T}2G(\eta_{t})dt.
\]

\end{lemma}

\begin{proof}
We know that $u(t,x)=\hat{\mathbb{E}}[\varphi(x+B_{T}-B_{t})]$ is
the
solution of the following PDE:%
\[
\partial_{t}u+G(D^{2}u)=0\ \ (t,x)\in \lbrack0,T]\times
\mathbb{R}^{d},\ u(T,x)=\varphi(x).
\]
For each $\varepsilon>0$, by the interior regularity of $u$ (see
Appendix C), we have
\[
\left \Vert u\right \Vert _{C^{1+\alpha/2,2+\alpha}([0,T-\varepsilon
]\times \mathbb{R}^{d})}<\infty\  \text{for some }\alpha \in(0,1).
\]
Applying $G$-It\^{o}'s formula to $u(t,B_{t}-B_{t_{1}})$ on $[t_{1}%
,T-\varepsilon]$, since $Du(t,x)$ is uniformly bounded, letting
$\varepsilon \rightarrow0,$ we have
\begin{align*}
\xi &  =\hat{\mathbb{E}}[\xi]+\int_{t_{1}}^{T}\partial_{t}u(t,B_{t}%
-B_{t_{1}})dt+\int_{t_{1}}^{T}\langle Du(t,B_{t}-B_{t_{1}}),dB_{t}\rangle \\
& \ \ \
+\frac{1}{2}\int_{t_{1}}^{T}(D^{2}u(t,B_{t}-B_{t_{1}}),d\langle
B\rangle_{t})\\
&  =\hat{\mathbb{E}}[\xi]+\int_{t_{1}}^{T}\langle Du(t,B_{t}-B_{t_{1}%
}),dB_{t}\rangle+\frac{1}{2}\int_{t_{1}}^{T}(D^{2}u(t,B_{t}-B_{t_{1}%
}),d\langle B\rangle_{t})\\
& \ \ \ -\int_{t_{1}}^{T}G(D^{2}u(t,B_{t}-B_{t_{1}}))dt.
\end{align*}
\end{proof}

We now give the representation theorem of $\xi=\varphi(B_{t_{1}},B_{t_{2}%
}-B_{t_{1}},\cdots,B_{t_{N}}-B_{t_{N-1}}).$

\begin{theorem}
\label{ch5t1} Let
$\xi=\varphi(B_{t_{1}},B_{t_{2}}-B_{t_{1}},\cdots,B_{t_{N}}-B_{t_{N-1}})$,
$\varphi \in C_{b.Lip}(\mathbb{R}^{d\times N})$, $0\leq t_{1}<t_{2}%
<\cdots<t_{N}=T<\infty$. Then we have the following representation:%
\[
\xi=\hat{\mathbb{E}}[\xi]+\int_{0}^{T}\langle
\beta_{t},dB_{t}\rangle +\int_{0}^{T}(\eta_{t},d\langle
B\rangle_{t})-\int_{0}^{T}2G(\eta_{t})dt.
\]

\end{theorem}

\begin{proof}
We only need to prove the case
$\xi=\varphi(B_{t_{1}},B_{T}-B_{t_{1}})$. We set, for each $(x,y)\in
\mathbb{R}^{2d},$
\[
u(t,x,y)=\hat{\mathbb{E}}[\varphi(x,y+B_{T}-B_{t})];\  \varphi
_{1}(x)=\hat{\mathbb{E}}[\varphi(x,B_{T}-B_{t_{1}})].
\]
For each $x\in \mathbb{R}^{d}$, we denote $\bar{\xi}=\varphi(x,B_{T}%
-B_{t_{1}})$. By Lemma \ref{ch5l1}, we have%
\begin{align*}
\bar{\xi}  &  =\varphi_{1}(x)+\int_{t_{1}}^{T}\langle D_{y}u(t,x,B_{t}%
-B_{t_{1}}),dB_{t}\rangle+\frac{1}{2}\int_{t_{1}}^{T}(D_{y}^{2}u(t,x,B_{t}%
-B_{t_{1}}),d\langle B\rangle_{t})\\
& \ \ \ -\int_{t_{1}}^{T}G(D_{y}^{2}u(t,x,B_{t}-B_{t_{1}}))dt.
\end{align*}
By the definitions of the integrations of $dt$, $dB_{t}$ and
$d\langle B\rangle_{t}$, we can replace $x$ by $B_{t_{1}}$ and get
\begin{align*}
\xi &  =\varphi_{1}(B_{t_{1}})+\int_{t_{1}}^{T}\langle D_{y}u(t,B_{t_{1}%
},B_{t}-B_{t_{1}}),dB_{t}\rangle \\
& \ \ \ +\frac{1}{2}\int_{t_{1}}^{T}(D_{y}^{2}u(t,B_{t_{1}},B_{t}-B_{t_{1}%
}),d\langle B\rangle_{t})-\int_{t_{1}}^{T}G(D_{y}^{2}u(t,B_{t_{1}}%
,B_{t}-B_{t_{1}}))dt.
\end{align*}
Applying Lemma \ref{ch5l1} to $\varphi_{1}(B_{t_{1}})$, we complete
the proof.
\end{proof}

We then immediately have the following $G$-martingale representation
theorem.

\begin{theorem}
Let $(M_{t})_{t\in \lbrack0,T]}$ be a $G$-martingale with
$M_{T}=\varphi
(B_{t_{1}},B_{t_{2}}-B_{t_{1}},\cdots,$$B_{t_{N}}-B_{t_{N-1}})$,
$\varphi \in C_{b.Lip}(\mathbb{R}^{d\times N})$, $0\leq
t_{1}<t_{2}<\cdots<t_{N}=T<\infty$. Then \
\[
M_{t}=\hat{\mathbb{E}}[M_{T}]+\int_{0}^{t}\langle \beta_{s},dB_{s}%
\rangle+\int_{0}^{t}(\eta_{s},d\langle
B\rangle_{s})-\int_{0}^{t}2G(\eta _{s})ds,\ t\leq T.
\]

\end{theorem}

\begin{proof}
For $M_{T}$, by Theorem \ref{ch5t1}, we have
\[
M_{T}=\hat{\mathbb{E}}[M_{T}]+\int_{0}^{T}\langle \beta_{s},dB_{s}%
\rangle+\int_{0}^{T}(\eta_{s},d\langle
B\rangle_{s})-\int_{0}^{T}2G(\eta _{s})ds.
\]
Taking the conditional $G$-expectation on both sides of the above
equality and by Proposition \ref{ch5p1}, we obtain the result.
\end{proof}

\section{$G$--convexity and Jensen's Inequality for $G$--expectations}

A very interesting question is whether the well--known Jensen's
inequality still holds for $G$--expectations. 

First, we give a new notion of convexity.
\begin{definition}
A continuous function $h:\mathbb{R\rightarrow R}$ is called $G$%
\textbf{--convex}\index{$G$-convex} if for each bounded $\xi \in
L_{G}^{1}(\Omega)$, the
following Jensen's inequality holds:%
\[
\hat{\mathbb{E}}[h(\xi)]\geq h(\hat{\mathbb{E}}[\xi]).
\]

\end{definition}

In this section, we mainly consider $C^{2}$-functions.

\begin{proposition}
\label{GG-p1} Let $h\in C^{2}(\mathbb{R})$. Then the following
statements are
equivalent: \newline \textbf{\textup{(i)}} The function $h$ is $G$%
--convex.\newline \textbf{\textup{(ii)}} For each bounded $\xi \in L_{G}%
^{1}(\Omega)$, the following Jensen's inequality holds:
\[
\hat{\mathbb{E}}[h(\xi)|\Omega_{t}]\geq h(\hat{\mathbb{E}}[\xi
|\Omega_{t}])\ \ \text{for}\  t\geq0.
\]
\textbf{\textup{(iii) }}For each $\varphi \in
C_{b}^{2}(\mathbb{R}^{d})$, the
following Jensen's inequality holds:%
\[
\hat{\mathbb{E}}[h(\varphi(B_{t}))]\geq h(\hat{\mathbb{E}}%
[\varphi(B_{t})])\ \ \text{for}\ t\geq0.
\]
\newline \textbf{\textup{(iv) }}The following condition holds for each
$(y,z,A)\in \mathbb{R}\times \mathbb{R}^{d}\times \mathbb{S}(d)$:
\begin{equation}
G(h^{\prime}(y)A+h^{\prime \prime}(y)zz^{T})-h^{\prime}(y)G(A)\geq
0.\  \label{d-G-conv}%
\end{equation}

\end{proposition}

To prove the above proposition, we need the following lemmas.

\begin{lemma}
\label{ch5l2} Let $\Phi:\mathbb{R}^{d}\rightarrow \mathbb{S}(d)$ be
continuous with
polynomial growth. Then%
\begin{equation}
\lim_{\delta \downarrow0}\hat{\mathbb{E}}[\int_{t}^{t+\delta}(\Phi
(B_{s}),d\langle B\rangle_{s})]\delta^{-1}=2\hat{\mathbb{E}}[G(\Phi
(B_{t}))].
\end{equation}

\end{lemma}

\begin{proof}
If $\Phi$ is a Lipschitz function, it is easy to prove that
\[
\hat{\mathbb{E}}[|\int_{t}^{t+\delta}(\Phi(B_{s})-\Phi(B_{t}),d\langle
B\rangle_{s})|]\leq C_{1}\delta^{3/2},
\]
where $C_{1}$ is a constant independent of $\delta$. Thus
\begin{align*}
\lim_{\delta \downarrow0}\hat{\mathbb{E}}[\int_{t}^{t+\delta}(\Phi
(B_{s}),d\langle B\rangle_{s})]\delta^{-1}  &  =\lim_{\delta
\downarrow 0}\hat{\mathbb{E}}[(\Phi(B_{t}),\langle
B\rangle_{t+\delta}-\langle
B\rangle_{s})]\delta^{-1}\\
&  =2\hat{\mathbb{E}}[G(\Phi(B_{t}))].
\end{align*}
Otherwise, we can choose a sequence of Lipschitz functions $\Phi
_{N}:\mathbb{R}^{d}\rightarrow \mathbb{S}(d)$ such that%
\[
|\Phi_{N}(x)-\Phi(x)|\leq \frac{C_{2}}{N}(1+|x|^{k}),
\]
where $C_{2}$ and $k$ are positive constants independent of $N$. It
is easy to show that
\[
\hat{\mathbb{E}}[|\int_{t}^{t+\delta}(\Phi(B_{s})-\Phi_{N}(B_{s}),d\langle
B\rangle_{s})|]\leq \frac{C}{N}\delta
\]
and%
\[
\hat{\mathbb{E}}[|G(\Phi(B_{t}))-G(\Phi_{N}(B_{t}))|]\leq
\frac{C}{N},
\]
where $C$ is a universal constant. Thus%
\begin{align*}
&  |\hat{\mathbb{E}}[\int_{t}^{t+\delta}(\Phi(B_{s}),d\langle
B\rangle
_{s})]\delta^{-1}-2\hat{\mathbb{E}}[G(\Phi(B_{t}))]|\\
\leq&|\hat{\mathbb{E}}[\int_{t}^{t+\delta}(\Phi_{N}(B_{s}),d\langle
B\rangle_{s})]\delta^{-1}-2\hat{\mathbb{E}}[G(\Phi_{N}(B_{t}))]|+\frac
{3C}{N}.
\end{align*}
Then we have
\[
\limsup_{\delta \downarrow0}|\hat{\mathbb{E}}[\int_{t}^{t+\delta}%
(\Phi(B_{s}),d\langle B\rangle_{s})]\delta^{-1}-2\hat{\mathbb{E}}%
[G(\Phi(B_{t}))]|\leq \frac{3C}{N}.
\]
Since $N$ can be arbitrarily large, we complete the proof.
\end{proof}

\begin{lemma}
\label{ch5l3} Let $\Psi$ be a $C^{2}$-function on $\mathbb{R}^d$
such that $D^{2}\Psi$ satisfy polynomial growth condition.
Then we have%
\begin{equation}
\lim_{\delta \downarrow0}(\hat{\mathbb{E}}[\Psi(B_{\delta})]-\Psi
(0))\delta^{-1}=G(D^{2}\Psi(0)).
\end{equation}

\end{lemma}

\begin{proof}
Applying $G$-It\^{o}'s formula to $\Psi(B_{\delta})$, we get%
\[
\Psi(B_{\delta})=\Psi(0)+\int_{0}^{\delta}\langle D\Psi(B_{s}),dB_{s}%
\rangle+\frac{1}{2}\int_{0}^{\delta}(D^{2}\Psi(B_{s}),d\langle
B\rangle_{s}).
\]
Thus we have%
\[
\hat{\mathbb{E}}[\Psi(B_{\delta})]-\Psi(0)=\frac{1}{2}\mathbb{\hat{E}}%
_{G}[\int_{0}^{\delta}(D^{2}\Psi(B_{s}),d\langle B\rangle_{s})].
\]
By Lemma \ref{ch5l2}, we obtain the result.
\end{proof}

\begin{lemma}
\label{ch5l4} Let $h\in C^{2}(\mathbb{R})$ and satisfy
\textup{(\ref{d-G-conv})}. For each $\varphi \in
C_{b.Lip}(\mathbb{R}^{d})$, let $u(t,x)$ be the solution of the
$G$-heat equation:%
\begin{equation}
\label{ch5e1}
\partial_{t}u-G(D^{2}u)=0\ \ \ (t,x)\in \lbrack
0,\infty)\times \mathbb{R}^{d}, \ u(0,x)=\varphi(x).
\end{equation}
Then $\tilde{u}(t,x):=h(u(t,x))$ is a viscosity subsolution of
$G$-heat equation \textup{(\ref{ch5e1})} with initial condition
$\tilde{u}(0,x)=h(\varphi(x))$.
\end{lemma}

\begin{proof}
For each $\varepsilon>0$, we denote by $u^{\varepsilon}$ the
solution of the
following PDE:%
\[
\partial_{t}u^{\varepsilon}-G_{\varepsilon}(D^{2}u^{\varepsilon}%
)=0\ \ \ (t,x)\in \lbrack0,\infty )\times \mathbb{R}^{d},  \
u^{\varepsilon}(0,x)=\varphi(x),
\]
where $G_{\varepsilon}(A):=G(A)+\varepsilon \mathrm{tr}[A]$. Since
$G_{\varepsilon}$ satisfies the uniformly elliptic condition, by
Appendix C, we have $u^{\varepsilon}\in C^{1,2}((0,\infty)\times
\mathbb{R}^{d})$. By
simple calculation, we have%
\[
\partial_{t}h(u^{\varepsilon})=h^{\prime}(u^{\varepsilon})\partial
_{t}u^{\varepsilon}=h^{\prime}(u^{\varepsilon})G_{\varepsilon}(D^{2}%
u^{\varepsilon})
\]
and
\[
\partial_{t}h(u^{\varepsilon})-G_{\varepsilon}(D^{2}h(u^{\varepsilon
}))=f_{\varepsilon}(t,x),\ h(u^{\varepsilon}(0,x))=h(\varphi(x)),
\]
where%
\[
f_{\varepsilon}(t,x)=h^{\prime}(u^{\varepsilon})G(D^{2}u^{\varepsilon
})-G(D^{2}h(u^{\varepsilon}))-\varepsilon h^{\prime
\prime}(u^{\varepsilon })|Du^{\varepsilon}|^{2}.
\]
Since $h$ is $G$--convex, it follows that
$f_{\varepsilon}\leq-\varepsilon h^{\prime
\prime}(u^{\varepsilon})|Du^{\varepsilon}|^{2}$. We can also deduce
that $|Du^{\varepsilon}|$ is uniformly bounded by the Lipschitz
constant of $\varphi$. It is easy to show that $u^{\varepsilon}$
uniformly converges to $u$ as $\varepsilon \rightarrow0$. Thus
$h(u^{\varepsilon})$ uniformly converges to $h(u)$ and $h^{\prime
\prime}(u^{\varepsilon})$ is uniformly
bounded. Then we get%
\[
\partial_{t}h(u^{\varepsilon})-G_{\varepsilon}(D^{2}h(u^{\varepsilon}))\leq
C\varepsilon,\ h(u^{\varepsilon}(0,x))=h(\varphi(x)),
\]
where $C$ is a constant independent of $\varepsilon$. By Appendix C,
we conclude that $h(u)$ is a viscosity subsolution.
\end{proof}

\textbf{Proof of Proposition \ref{GG-p1}. }Obviously \textbf{\textup{(ii)}%
}$\Longrightarrow$\textbf{\textup{(i)}}$\Longrightarrow$\textbf{\textup{(iii)}%
}. We now prove \textbf{\textup{(iii)}}$\Longrightarrow$\textbf{\textup{(ii)}%
}. For  $\xi \in L_{G}^{1}(\Omega)$ of the form%
\[
\xi=\varphi(B_{t_{1}},B_{t_{2}}-B_{t_{1}},\cdots,B_{t_{n}}-B_{t_{n-1}}),
\]
where $\varphi \in C_{b}^{2}(\mathbb{R}^{d\times n})$, $0\leq
t_{1}\leq \cdots \leq t_{n}<\infty$, by the definitions of
$\hat{\mathbb{E}}[\cdot]$
and $\hat{\mathbb{E}}[\cdot|\Omega_{t}]$, we have%
\[
\hat{\mathbb{E}}[h(\xi)|\Omega_{t}]\geq h(\hat{\mathbb{E}}[\xi
|\Omega_{t}]),\ t\geq0.
\]
We then can extend this Jensen's inequality, under the norm $||\cdot
||=\hat{\mathbb{E}}[|\cdot|]$, to each bounded $\xi \in
L_{G}^{1}(\Omega )$.\newline
\textbf{\textup{(iii)}}$\Longrightarrow$\textbf{\textup{(iv)}}: for
each $\varphi \in C_{b}^{2}(\mathbb{R}^{d})$, we have $\mathbb{\hat{E}}%
[h(\varphi(B_{t}))]\geq h(\hat{\mathbb{E}}[\varphi(B_{t})])$ for
each
$t\geq0$. By Lemma \ref{ch5l3}, we know that%
\[
\lim_{\delta \downarrow0}(\hat{\mathbb{E}}[\varphi(B_{\delta}%
)]-\varphi(0))\delta^{-1}=G(D^{2}\varphi(0))
\]
and%
\[
\lim_{\delta \downarrow0}(\hat{\mathbb{E}}[h(\varphi(B_{\delta
}))]-h(\varphi(0)))\delta^{-1}=G(D^{2}h(\varphi)(0)).
\]
Thus we get%
\[
G(D^{2}h(\varphi)(0))\geq h^{\prime}(\varphi(0))G(D^{2}\varphi(0)).
\]
For each $(y,z,A)\in \mathbb{R\times R}^{d}\times \mathbb{S}(d)$, we
can choose a $\varphi \in C_{b}^{2}(\mathbb{R}^{d})$ such that
$(\varphi(0),D\varphi (0),D^{2}\varphi(0))=$ $(y,z,A)$. Thus we
obtain \textbf{\textup{(iv)}}.\newline
\textbf{\textup{(iv)}}$\Longrightarrow$\textbf{\textup{(iii)}}: for
each $\varphi \in C_{b}^{2}(\mathbb{R}^{d})$,
$u(t,x)=\hat{\mathbb{E}}[\varphi(x+B_{t})]$ (respectively, $\bar{u}%
(t,x)=\hat{\mathbb{E}}[h(\varphi(x+B_{t}))]$) solves the $G$-heat
equation (\ref{ch5e1}). By Lemma \ref{ch5l4}, $h(u)$ is a viscosity
subsolution of $G$-heat equation (\ref{ch5e1}). It follows from the
maximum principle that $h(u(t,x))\leq \bar{u}(t,x)$. In particular,
\textbf{\textup{(iii)}} holds. {$\Box$}

\begin{remark}
In fact, \textbf{\textup{(i)}}$\Longleftrightarrow$\textbf{\textup{(ii)}%
}$\Longleftrightarrow$\textbf{\textup{(iii)}} still hold without the
assumption $h\in C^{2}(\mathbb{R}).$
\end{remark}

\begin{proposition}
Let $h$ be a $G$--convex function and $X\in L_{G}^{1}(\Omega)$ be
bounded. Then $Y_{t}=h(\hat{\mathbb{E}}[X|\Omega_{t}])$, $t\geq0$,
is a $G$--submartingale.
\end{proposition}

\begin{proof}
For each $s\leq t$,
\[
\hat{\mathbb{E}}[Y_{t}|\Omega_{s}]=\hat{\mathbb{E}}[h(\mathbb{\hat{E}%
}[X|\Omega_{t}])|\Omega_{s}]\geq h(\hat{\mathbb{E}}[X|\Omega
_{s}])=Y_{s}\text{.}%
\]

\end{proof}

\begin{exercise}
Suppose that $G$ satisfies the uniformly elliptic condition and $h\in C^{2}(\mathbb{R}%
)$. Show that $h$ is $G$-convex if and only if $h$ is convex.
\end{exercise}

\section*{Notes and Comments}
\addcontentsline{toc}{section}{Notes and Comments} This chapter is
mainly from Peng (2007) \cite{Peng2007b}.

Peng (1997) \cite{Peng1997} introduced a filtration consistent (or
time consistent, or dynamic) nonlinear expectation, called
$g$-expectation, via BSDE, and then in (1999) \cite{Peng1999} for
some basic properties of the $g$-martingale such as nonlinear
Doob-Meyer decomposition theorem, see also Briand, Coquet,
Hu, M\'{e}min  and Peng (2000) \cite{BCHMP1}, Chen,
Kulperger and Jiang  (2003) \cite{CKJ},  Chen and Peng (1998) \cite{CP} and (2000) \cite{CP1}, Coquet, Hu,
M\'{e}min and Peng (2001) \cite{CHMP}, and (2002)
\cite{CHMP3}, Peng (1999) \cite{Peng1999}, (2004)
\cite{Peng2004c},  Peng and Xu (2003) \cite{PX2003},
Rosazza (2006) \cite{Roazza2003}. Our conjecture is that all
properties obtained for $g$-martingales must has its correspondence
for $G$-martingale. But this conjecture is still far from being
complete. Here we present some properties of $G$-martingales.

The problem $G$-martingale representation theorem has been raised as
a problem in Peng (2007) \cite{Peng2007b}. In Section 2, we only
give a result with very regular random variables. Some very
interesting developments to this important problem can be found in
Soner, Touzi and Zhang (2009) \cite{Soner} and Song (2009)
\cite{Song2}.

Under the framework of $g$-expectation, Chen, Kulperger and Jiang (2003)
 \cite{CKJ}, Hu (2005) \cite{huying}, Jiang and Chen (2004)
\cite{JC} investigate the Jensen's inequality for
$g$-expectation. Recently, Jia and Peng (2007) \cite{JiaPeng}
introduced the notion of $g$-convex function and obtained many
interesting properties. Certainly a $G$-convex function concerns
fully nonlinear situations.

%
%

\chapter{Stochastic Differential Equations}
\label{ch4}

In this chapter, we consider the stochastic differential equations
and backward stochastic differential equations driven by
$G$-Brownian motion. The conditions and proofs of existence and
uniqueness of a stochastic differential equation is similar to the
classical situation. However the corresponding problems for backward
stochastic differential equations are not that easy, many are still
open. We only give partial results to this direction.

\section{Stochastic Differential Equations}
In this chapter, we denote by $\bar{M}_{G}^{p}(0,T;\mathbb{R}^{n})$,
$p\geq1$\label{mbargp}, the completion
of $M_{G}^{p,0}(0,T;\mathbb{R}^{n})$ under the norm $(\int_{0}^{T}%
\hat{\mathbb{E}}[|\eta_{t}|^{p}]dt)^{1/p}$. It is not hard to prove
that $\bar{M}_{G}^{p}(0,T;\mathbb{R}^{n})\subseteq
M_{G}^p(0,T;\mathbb{R}^{n})$. We consider all the problems in the
space $\bar{M}_{G}^{p}(0,T;\mathbb{R}^{n})$, and the sublinear
expectation space $(\Omega,\mathcal{H},\hat{\mathbb{E}})$ is fixed.

We consider the following SDE driven by a $d$-dimensional
$G$-Brownian motion:
\begin{equation}
X_{t}=X_{0}+\int_{0}^{t}b(s,X_{s})ds+\int_{0}^{t}h_{ij}(s,X_{s})d\left
\langle
B^{i},B^{j}\right \rangle _{s}+\int_{0}^{t}\sigma_{j}(s,X_{s})dB_{s}^{j}%
,\ t\in \lbrack0,T], \label{d-SDE}%
\end{equation}
where the initial condition $X_{0}\in \mathbb{R}^{n}$ is a given constant, and $b,h_{ij}%
,\sigma_{j}$ are given functions satisfying $b(\cdot,x)$,
$h_{ij}(\cdot,x)$, $\sigma_{j}(\cdot,x)\in
\bar{M}_{G}^{2}(0,T;\mathbb{R}^{n})$ for each $x\in \mathbb{R}^{n}$
and the Lipschitz condition, i.e., $|\phi(t,x)-\phi
(t,x^{\prime})|\leq K|x-x^{\prime}|$, for each $t\in \lbrack0,T]$,
$x$, $x^{\prime}\in \mathbb{R}^{n}$, $\phi=b$, $h_{ij}$ and
$\sigma_{j}$, respectively. Here the horizon $[0,T]$ can be
arbitrarily large. The solution is a process $X\in
\bar{M}_{G}^{2}(0,T;\mathbb{R}^{n})$ satisfying the SDE
(\ref{d-SDE}).

We first introduce the following mapping on a fixed interval $[0,T]$:%
\[
\Lambda_{\cdot}:\bar{M}_{G}^{2}(0,T;\mathbb{R}^{n})\rightarrow \bar
{M}_{G}^{2}(0,T;\mathbb{R}^{n})\  \
\]
by setting $\Lambda_{t}$, $t\in \lbrack0,T]$, with
\[
\Lambda_{t}(Y)=X_{0}+\int_{0}^{t}b(s,Y_{s})ds+\int_{0}^{t}h_{ij}%
(s,Y_{s})d\left \langle B^{i},B^{j}\right \rangle
_{s}+\int_{0}^{t}\sigma _{j}(s,Y_{s})dB_{s}^{j}.
\]

We immediately have the following lemma.

\begin{lemma}
\label{ch4l1} For each $Y,Y^{\prime}\in
\bar{M}_{G}^{2}(0,T;\mathbb{R}^{n}),$ we have the
following estimate:%
\begin{equation}
\label{ch4-1}
\hat{\mathbb{E}}[|\Lambda_{t}(Y)-\Lambda_{t}(Y^{\prime})|^{2}]\leq
C\int_{0}^{t}\hat{\mathbb{E}}[|Y_{s}-Y_{s}^{\prime}|^{2}]ds,\ t\in
\lbrack0,T],
\end{equation}
where the constant $C$ depends only on the Lipschitz constant $K$.
\end{lemma}
%

We now prove that SDE (\ref{d-SDE}) has a unique solution. By
multiplying $e^{-2Ct}$ on both sides of (\ref{ch4-1}) and
integrating them on
$[0,T]$, it follows that%
\begin{align*}
\int_{0}^{T}\hat{\mathbb{E}}[|\Lambda_{t}(Y)-\Lambda_{t}(Y^{\prime}%
)|^{2}]e^{-2Ct}dt  &  \leq C\int_{0}^{T}e^{-2Ct}\int_{0}^{t}\mathbb{\hat{E}%
}_{G}[|Y_{s}-Y_{s}^{\prime}|^{2}]dsdt\\
&  =C\int_{0}^{T}\int_{s}^{T}e^{-2Ct}dt\hat{\mathbb{E}}[|Y_{s}%
-Y_{s}^{\prime}|^{2}]ds\\
&  =\frac{1}{2}\int_{0}^{T}(e^{-2Cs}-e^{-2CT})\hat{\mathbb{E}}%
[|Y_{s}-Y_{s}^{\prime}|^{2}]ds.
\end{align*}
We then have
\begin{equation}
\label{BSDE-1}
\int_{0}^{T}\hat{\mathbb{E}}[|\Lambda_{t}(Y)-\Lambda_{t}(Y^{\prime}%
)|^{2}]e^{-2Ct}dt\leq \frac{1}{2}\int_{0}^{T}\hat{\mathbb{E}}[|Y_{t}%
-Y_{t}^{\prime}|^{2}]e^{-2Ct}dt.
\end{equation}
We observe that the following two norms are equivalent on
$\bar{M}_{G} ^{2}(0,T;\mathbb{R}^{n})$, i.e.,
\[
(\int_{0}^{T}\hat{\mathbb{E}}[|Y_{t}|^{2}]dt)^{1/2}\thicksim
(\int_{0} ^{T}\hat{\mathbb{E}}[|Y_{t}|^{2}]e^{-2Ct}dt)^{1/2}.
\]
From (\ref{BSDE-1}) we can obtain that $\Lambda(Y)$ is a contraction
mapping. Consequently, we have the following theorem.

\begin{theorem}
\label{BSDE-th1} There exists a unique solution $X\in
\bar{M}_{G}^{2}(0,T;\mathbb{R}^{n})$ of the stochastic differential
equation \textup{(\ref{d-SDE})}.
\end{theorem}

We now consider the following linear SDE. For simplicity, we assume that $d=1$ and $n=1.$%
\begin{equation}
\label{BSDE-2}
X_{t}=X_{0}+\int_{0}^{t}(b_{s}X_{s}+\tilde{b}_{s})ds+\int_{0}^{t}(h_{s}%
X_{s}+\tilde{h}_{s})d\langle B\rangle_{s}+\int_{0}^{t}(\sigma_{s}X_{s}%
+\tilde{\sigma}_{s})dB_{s},\ t\in \lbrack0,T],
\end{equation}
where $X_{0}\in \mathbb{R}$ is given, $b_{.},h_{.},\sigma_{.}$ are
given
bounded processes in $\bar{M}_{G}^{2}(0,T;\mathbb{R})$ and $\tilde{b}%
_{.},\tilde{h}_{.},\tilde{\sigma}_{.}$ are given processes in $\bar{M}_{G}%
^{2}(0,T;\mathbb{R})$. By Theorem \ref{BSDE-th1}, we know that the
linear SDE (\ref{BSDE-2}) has a unique solution.

\begin{remark}
\label{ch4p1} The solution of the linear SDE \textup{(\ref{BSDE-2})}
is
\begin{equation*}
X_{t}=\Gamma_{t}^{-1}(X_{0}+\int_{0}^{t}\tilde{b}_{s}\Gamma_{s}ds+\int_{0}%
^{t}(\tilde{h}_{s}-\sigma_{s}\tilde{\sigma}_{s})\Gamma_{s}d\langle
B\rangle_{s}+\int_{0}^{t}\tilde{\sigma}_{s}\Gamma_{s}dB_{s}),\ t\in
\lbrack0,T],
\end{equation*}
where $\Gamma_{t}=\exp(-\int_{0}^{t}b_{s}ds-\int_{0}^{t}(h_{s}-\frac{1}%
{2}\sigma_{s}^{2})d\langle
B\rangle_{s}-\int_{0}^{t}\sigma_{s}dB_{s})$.

In particular, if $b_{.},h_{.},\sigma_{.}$ are constants and $\tilde{b}%
_{.},\tilde{h}_{.},\tilde{\sigma}_{.}$ are zero, then $X$ is a
geometric $G$-Brownian motion.
\end{remark}

\begin{definition}
We call $X$ is a \textbf{geometric $G$-Brownian
motion}\index{Geometric $G$-Brownian motion} if
\begin{equation}
X_{t}=\exp(\alpha t+\beta \langle B\rangle_{t}+\gamma B_{t}),
\end{equation}
where $\alpha,\beta,\gamma$ are constants.
\end{definition}

\begin{exercise}
Prove that $\bar{M}_{G}^{p}(0,T;\mathbb{R}^{n})\subseteq
M_{G}^p(0,T;\mathbb{R}^{n})$.
\end{exercise}
\begin{exercise}
Complete the proof of Lemma \ref{ch4l1}.
\end{exercise}

\section{Backward Stochastic Differential Equations}

We consider the following type of BSDE:
\begin{equation}
Y_{t}=\hat{\mathbb{E}}[\xi+\int_{t}^{T}f(s,Y_{s})ds+\int_{t}^{T}%
h_{ij}(s,Y_{s})d\left \langle B^{i},B^{j}\right \rangle _{s}|\Omega
_{t}],\  \ t\in \lbrack0,T], \label{BSDE}%
\end{equation}
where $\xi \in L_{G}^{1}(\Omega_{T};\mathbb{R}^{n})$ is given, and
$f,h_{ij}$ are given functions satisfying $f(\cdot,y)$,
$h_{ij}(\cdot,y)\in \bar{M}_{G}^{1}(0,T;\mathbb{R}^{n})$ for each
$y\in \mathbb{R}^{n}$ and the Lipschitz condition, i.e.,
$|\phi(t,y)-\phi(t,y^{\prime})|\leq K|y-y^{\prime}|$, for each $t\in
\lbrack0,T]$, $y$, $y^{\prime}\in \mathbb{R}^{n}$, $\phi=f$ and
$h_{ij}$, respectively. The solution is a process $Y\in \bar{M}_{G}%
^{1}(0,T;\mathbb{R}^{n})$ satisfying the above BSDE.

We first introduce the following mapping on a fixed interval $[0,T]$:%
\[
\Lambda_{\cdot}:\bar{M}_{G}^{1}(0,T;\mathbb{R}^{n})\rightarrow \bar
{M}_{G}^{1}(0,T;\mathbb{R}^{n})\  \
\]
by setting $\Lambda_{t}$, $t\in \lbrack0,T]$, with
\[
\Lambda_{t}(Y)=\hat{\mathbb{E}}[\xi+\int_{t}^{T}f(s,Y_{s})ds+\int_{t}%
^{T}h_{ij}(s,Y_{s})d\left \langle B^{i},B^{j}\right \rangle
_{s}|\Omega_{t}].
\]

We immediately have

\begin{lemma}
\label{ch4l2} \label{BSDE-l1} For each $Y,Y^{\prime}\in
\bar{M}_{G}^{1}(0,T;\mathbb{R}^{n}),$ we have the
following estimate:%
\begin{equation}
\label{ch4-2}
\hat{\mathbb{E}}[|\Lambda_{t}(Y)-\Lambda_{t}(Y^{\prime})|]\leq C\int
_{t}^{T}\hat{\mathbb{E}}[|Y_{s}-Y_{s}^{\prime}|]ds,\ t\in
\lbrack0,T],
\end{equation}
where the constant $C$ depends only on the Lipschitz constant $K$.
\end{lemma}
%

We now prove that BSDE (\ref{BSDE}) has a unique solution. By
multiplying $e^{2Ct}$ on both sides of (\ref{ch4-2}) and integrating
them on
$[0,T]$, it follows that%
\begin{align}
\int_{0}^{T}\hat{\mathbb{E}}[|\Lambda_{t}(Y)-\Lambda_{t}(Y^{\prime
})|]e^{2Ct}dt  &  \leq C\int_{0}^{T}\int_{t}^{T}\hat{\mathbb{E}}%
[|Y_{s}-Y_{s}^{\prime}|]e^{2Ct}dsdt\nonumber\\
&  =C\int_{0}^{T}\hat{\mathbb{E}}[|Y_{s}-Y_{s}^{\prime}|]\int_{0}%
^{s}e^{2Ct}dtds\nonumber\\
&  =\frac{1}{2}\int_{0}^{T}\hat{\mathbb{E}}[|Y_{s}-Y_{s}^{\prime
}|](e^{2Cs}-1)ds\nonumber\\
&  \leq \frac{1}{2}\int_{0}^{T}\hat{\mathbb{E}}[|Y_{s}-Y_{s}^{\prime
}|]e^{2Cs}ds.\label{BSDE-3}
\end{align}
We observe that the following two norms are equivalent on
$\bar{M}_{G} ^{1}(0,T;\mathbb{R}^{n})$, i.e.,
\[
\int_{0}^{T}\hat{\mathbb{E}}[|Y_{t}|]dt\thicksim \int_{0}^{T}%
\hat{\mathbb{E}}[|Y_{t}|]e^{2Ct}dt.
\]
From (\ref{BSDE-3}), we can obtain that $\Lambda(Y)$ is a
contraction mapping. Consequently, we have the following theorem.

\begin{theorem}
There exists a unique solution $(Y_{t})_{t\in \lbrack0,T]}\in \bar{M}_{G}%
^{1}(0,T;\mathbb{R}^{n})$ of the backward stochastic differential
equation \textup{(\ref{BSDE})}.
\end{theorem}
Let $Y^{v}$, $v=1,2$, be the solutions of the following BSDE:
\[
Y_{t}^{v}=\hat{\mathbb{E}}[\xi^{v}+\int_{t}^{T}(f(s,Y_{s}^{v})+\varphi
_{s}^{v})ds+\int_{t}^{T}(h_{ij}(s,Y_{s}^{v})+\psi_{s}^{ij,v})d\left
\langle B^{i},B^{j}\right \rangle _{s}|\Omega_{t}].
\]
Then the following estimate holds.
\begin{proposition}
We have%
\begin{equation}
\hat{\mathbb{E}}[|Y_{t}^{1}-Y_{t}^{2}|]\leq Ce^{C(T-t)}(\mathbb{\hat{E}%
}[|\xi^{1}-\xi^{2}|]+\int_{t}^{T}\hat{\mathbb{E}}[|\varphi_{s}%
^{1}-\varphi_{s}^{2}|+|\psi_{s}^{ij,1}-\psi_{s}^{ij,2}|]ds),
\end{equation}
where the constant $C$ depends only on the Lipschitz constant $K$.
\end{proposition}

\begin{proof}
Similar to Lemma \ref{BSDE-l1}, we have%
\begin{align*}
\hat{\mathbb{E}}[|Y_{t}^{1}-Y_{t}^{2}|]  &  \leq C(\int_{t}^{T}%
\hat{\mathbb{E}}[|Y_{s}^{1}-Y_{s}^{2}|]ds+\hat{\mathbb{E}}[|\xi
^{1}-\xi^{2}|]\\
& \ \ \ +\int_{t}^{T}\hat{\mathbb{E}}[|\varphi_{s}^{1}-\varphi_{s}^{2}%
|+|\psi_{s}^{ij,1}-\psi_{s}^{ij,2}|]ds).
\end{align*}
By the Gronwall inequality (see Exercise \ref{ch4e1}), we conclude
the result.
\end{proof}

\begin{remark}
In particular, if $\xi^{2}=0$, $\varphi_{s}^{2}=-f(s,0)$, $\psi_{s}%
^{ij,2}=-h_{ij}(s,0)$, $\varphi_{s}^{1}=0$, $\psi_{s}^{ij,1}=0$, we
obtain the estimate of the solution of the BSDE. Let $Y$ be the
solution of the BSDE \textup{(\ref{BSDE})}. Then
\begin{equation}
\label{BSDE-7} \hat{\mathbb{E}}[|Y_{t}|]\leq
Ce^{C(T-t)}(\hat{\mathbb{E}}[|\xi
|]+\int_{t}^{T}\hat{\mathbb{E}}[|f(s,0)|+|h_{ij}(s,0)|]ds),
\end{equation}
where the constant $C$ depends only on the Lipschitz constant $K$.
\end{remark}

\begin{exercise}
\label{ch4e1} (The Gronwall inequality) Let $u(t)$ be a nonnegative
function such that
$$u(t)\leq C+A\int_0^tu(s)ds\ ~\text{for}~0\leq t\leq T,$$
where $C$ and $A$ are constants. Prove that $u(t)\leq Ce^{At}$ for
$0\leq t\leq T$.
\end{exercise}

\begin{exercise}
For each $\xi \in L_{G}^{1}(\Omega_{T};\mathbb{R}^{n})$. Show that
the process $(\hat{\mathbb{E}}[\xi|\Omega_{t}])_{t\in \lbrack0,T]}$
belongs to $\bar {M}_{G}^{1}(0,T;\mathbb{R}^{n})$.
\end{exercise}

\begin{exercise}
Complete the proof of Lemma \ref{ch4l2}.
\end{exercise}

\section{Nonlinear Feynman-Kac Formula}

Consider the following SDE:%
\begin{equation}
\label{BSDE-8} \left\{\begin{aligned}dX_{s}^{t,\xi}&
=b(X_{s}^{t,\xi})ds+h_{ij}(X_{s}^{t,\xi})d\left \langle
B^{i},B^{j}\right \rangle _{s}+\sigma_{j}(X_{s}^{t,\xi})dB_{s}^{j}%
,\ s\in \lbrack t,T],\\
X_{t}^{t,\xi}&=\xi,\end{aligned}\right.
\end{equation}
where $\xi \in L_{G}^{2}(\Omega_{t};\mathbb{R}^{n})$ is given and $b$, $h_{ij}%
$, $\sigma_{j}:\mathbb{R}^{n}\rightarrow \mathbb{R}^{n}$ are given
Lipschitz functions, i.e., $|\phi(x)-\phi(x^{\prime})|\leq
K|x-x^{\prime}|$, for each $x$, $x^{\prime}\in \mathbb{R}^{n}$,
$\phi=b$, $h_{ij}$ and $\sigma_{j}$.

\noindent We then consider associated BSDE:%
\begin{equation}
Y_{s}^{t,\xi}=\hat{\mathbb{E}}[\Phi(X_{T}^{t,\xi})+\int_{s}^{T}%
f(X_{r}^{t,\xi},Y_{r}^{t,\xi})dr+\int_{s}^{T}g_{ij}(X_{r}^{t,\xi},Y_{r}%
^{t,\xi})d\left \langle B^{i},B^{j}\right \rangle _{r}|\Omega_{s}],
\end{equation}
where $\Phi:\mathbb{R}^{n}\rightarrow \mathbb{R}$ is a given
Lipschitz function and $f$, $g_{ij}:\mathbb{R}^{n}\times
\mathbb{R}\rightarrow \mathbb{R}$ are given Lipschitz functions,
i.e., $|\phi(x,y)-\phi(x^{\prime},y^{\prime})|\leq
K(|x-x^{\prime}|+|y-y^{\prime}|)$, for each $x$, $x^{\prime}\in \mathbb{R}^{n}%
$, $y$, $y^{\prime}\in \mathbb{R}$, $\phi=f$ and $g_{ij}$.

We have the following estimates:

\begin{proposition}
\label{ch4p2} For each $\xi$, $\xi^{\prime}\in
L_{G}^{2}(\Omega_{t};\mathbb{R}^{n})$, we
have, for each $s\in \lbrack t,T],$\bigskip%
\begin{equation}
\hat{\mathbb{E}}[|X_{s}^{t,\xi}-X_{s}^{t,\xi^{\prime}}|^{2}|\Omega
_{t}]\leq C|\xi-\xi^{\prime}|^{2}%
\end{equation}
and%
\begin{equation}
\label{BSDE-5} \hat{\mathbb{E}}[|X_{s}^{t,\xi}|^{2}|\Omega_{t}]\leq
C(1+|\xi|^{2}),
\end{equation}
where the constant $C$ depends only on the Lipschitz constant $K$.
\end{proposition}

\begin{proof}
It is easy to obtain%
\[
\hat{\mathbb{E}}[|X_{s}^{t,\xi}-X_{s}^{t,\xi^{\prime}}|^{2}|\Omega
_{t}]\leq C_{1}(|\xi-\xi^{\prime}|^{2}+\int_{t}^{s}\hat{\mathbb{E}}%
[|X_{r}^{t,\xi}-X_{r}^{t,\xi^{\prime}}|^{2}|\Omega_{t}]dr).
\]
By the
Gronwall inequality, we obtain%
\[
\hat{\mathbb{E}}[|X_{s}^{t,\xi}-X_{s}^{t,\xi^{\prime}}|^{2}|\Omega
_{t}]\leq C_{1}e^{C_{1}T}|\xi-\xi^{\prime}|^{2}.
\]
Similarly, we can get (\ref{BSDE-5}).
\end{proof}

\begin{corollary}
For each $\xi \in L_{G}^{2}(\Omega_{t};\mathbb{R}^{n})$, we have
\begin{equation}
\label{BSDE-12}
\hat{\mathbb{E}}[|X_{t+\delta}^{t,\xi}-\xi|^{2}|\Omega_{t}]\leq
C(1+|\xi|^{2})\delta \ \ \text{for} \ \delta\in [0,T-t],
\end{equation}
where the constant $C$ depends only on the Lipschitz constant $K$.
\end{corollary}

\begin{proof}
It is easy to obtain
\[
\hat{\mathbb{E}}[|X_{t+\delta}^{t,\xi}-\xi|^{2}|\Omega_{t}]\leq C_{1}%
\int_{t}^{t+\delta}(1+\hat{\mathbb{E}}[|X_{s}^{t,\xi}|^{2}|\Omega
_{t}])ds.
\]
By Proposition \ref{ch4p2}, we obtain the result.
\end{proof}

\begin{proposition}
For each $\xi$, $\xi^{\prime}\in
L_{G}^{2}(\Omega_{t};\mathbb{R}^{n})$, we
have%
\begin{equation}
\label{BSDE-13} |Y_{t}^{t,\xi}-Y_{t}^{t,\xi^{\prime}}|\leq
C|\xi-\xi^{\prime}|
\end{equation}
and%
\begin{equation}
\label{BSDE-14} |Y_{t}^{t,\xi}|\leq C(1+|\xi|),
\end{equation}
where the constant $C$ depends only on the Lipschitz constant $K$.
\end{proposition}

\begin{proof}
For each $s\in \lbrack0,T]$, it is easy to check that%
\[
|Y_{s}^{t,\xi}-Y_{s}^{t,\xi^{\prime}}|\leq C_{1}\hat{\mathbb{E}}%
[|X_{T}^{t,\xi}-X_{T}^{t,\xi^{\prime}}|+\int_{s}^{T}(|X_{r}^{t,\xi}%
-X_{r}^{t,\xi^{\prime}}|+|Y_{r}^{t,\xi}-Y_{r}^{t,\xi^{\prime}}|)dr|\Omega
_{s}].
\]
Since%
\[
\hat{\mathbb{E}}[|X_{s}^{t,\xi}-X_{s}^{t,\xi^{\prime}}||\Omega_{t}%
]\leq(\hat{\mathbb{E}}[|X_{s}^{t,\xi}-X_{s}^{t,\xi^{\prime}}|^{2}%
|\Omega_{t}])^{1/2},
\]
we have
\[
\hat{\mathbb{E}}[|Y_{s}^{t,\xi}-Y_{s}^{t,\xi^{\prime}}||\Omega_{t}]\leq
C_{2}(|\xi-\xi^{\prime}|+\int_{s}^{T}\hat{\mathbb{E}}[|Y_{r}^{t,\xi}%
-Y_{r}^{t,\xi^{\prime}}||\Omega_{t}]dr).
\]
By the Gronwall inequality, we obtain (\ref{BSDE-13}). Similarly we
can get (\ref{BSDE-14}).
\end{proof}

We are more interested in the case when $\xi=x\in \mathbb{R}^{n}$. Define%
\begin{equation}
u(t,x):=Y_{t}^{t,x},\  \ (t,x)\in \lbrack0,T]\times \mathbb{R}^{n}.
\end{equation}
By the above proposition, we immediately have the following estimates:%
\begin{equation}
\label{BSDE-16} |u(t,x)-u(t,x^{\prime})|\leq C|x-x^{\prime}|,
\end{equation}%
\begin{equation}
\label{BSDE-17} |u(t,x)|\leq C(1+|x|),
\end{equation}
where the constant $C$ depends only on the Lipschitz constant $K$.

\begin{remark}
It is important to note that $u(t,x)$ is a deterministic function of
$(t,x)$, because $X_{s}^{t,x}$ and $Y_{s}^{t,x}$ are independent
from $\Omega_{t}$.
\end{remark}

\begin{theorem}
\label{ch4t2}
For each $\xi \in L_{G}^{2}(\Omega_{t};\mathbb{R}^{n})$, we have%
\begin{equation}
u(t,\xi)=Y_{t}^{t,\xi}.
\end{equation}

\end{theorem}


\begin{proposition}
We have, for $\delta \in \lbrack0,T-t],$
\begin{equation}
\label{BSDE-19}
u(t,x)=\hat{\mathbb{E}}[u(t+\delta,X_{t+\delta}^{t,x})+\int_{t}^{t+\delta
}f(X_{r}^{t,x},Y_{r}^{t,x})dr+\int_{t}^{t+\delta}g_{ij}(X_{r}^{t,x}%
,Y_{r}^{t,x})d\left \langle B^{i},B^{j}\right \rangle _{r}].
\end{equation}

\end{proposition}

\begin{proof}
Since $X_{s}^{t,x}=X_{s}^{t+\delta,X_{t+\delta}^{t,x}}$ for $s\in
\lbrack t+\delta,T]$, we get
$Y_{t+\delta}^{t,x}=Y_{t+\delta}^{t+\delta,X_{t+\delta }^{t,x}}$. By
Theorem \ref{ch4t2}, we have $Y_{t+\delta}^{t,x}=u(t+\delta
,X_{t+\delta}^{t,x})$, which implies the result.
\end{proof}

For each $A\in \mathbb{S}(n)$, $p\in \mathbb{R}^{n}$, $r\in
\mathbb{R}$, we set
\begin{equation*}
F(A,p,r,x):=G(B(A,p,r,x))+\langle p,b(x)\rangle+f(x,r),
\end{equation*}
where $B(A,p,r,x)$ is a $d\times d$ symmetric matrix with%
\begin{equation*}
B_{ij}(A,p,r,x):=\langle A\sigma_{i}(x),\sigma_{j}(x)\rangle+\langle
p,h_{ij}(x)+h_{ji}(x)\rangle+g_{ij}(x,r)+g_{ji}(x,r).
\end{equation*}

\begin{theorem}
$u(t,x)$ is a viscosity solution of the following PDE:%
\begin{equation}
\label{BSDE-22} \left \{
\begin{array}
[c]{l}%
\partial_{t}u+F(D^{2}u,Du,u,x)=0,\\
u(T,x)=\Phi(x).
\end{array}
\right.
\end{equation}

\end{theorem}

\begin{proof}
We first show that $u$ is a continuous function. By (\ref{BSDE-16})
we know that $u$ is a Lipschitz function in $x$. It follows from
(\ref{BSDE-7}) and (\ref{BSDE-5}) that for $s\in \lbrack t,T],$
$\hat{\mathbb{E}}[|Y_{s}^{t,x}|]\leq C(1+|x|)$. Noting
(\ref{BSDE-12}) and (\ref{BSDE-19}), we get
$|u(t,x)-u(t+\delta,x)|\leq C(1+|x|)(\delta
^{1/2}+\delta)$ for $\delta \in \lbrack0,T-t]$. Thus $u$ is $\frac{1}{2}%
$-H\"{o}lder continuous in $t$, which implies that $u$ is a
continuous function. We can also show, that for each $p\geq2$,
\begin{equation}
\label{BSDE-23} \hat{\mathbb{E}}[|X_{t+\delta}^{t,x}-x|^{p}]\leq
C(1+|x|^{p})\delta^{p/2},
\end{equation}
Now for
fixed $(t,x)\in(0,T)\times \mathbb{R}^{n}$, let $\psi \in C_{b}^{2,3}%
([0,T]\times \mathbb{R}^{n})$ be such that $\psi \geq u$ and
$\psi(t,x)=u(t,x)$. By (\ref{BSDE-19}), (\ref{BSDE-23}) and Taylor's
expansion, it follows that, for $\delta
\in(0,T-t),$%
\begin{align*}
0  & \leq \hat{\mathbb{E}}[\psi(t+\delta,X_{t+\delta}^{t,x})-\psi
(t,x)+\int_{t}^{t+\delta}f(X_{r}^{t,x},Y_{r}^{t,x})dr\\
&\ \ \ \ +\int_{t}^{t+\delta }g_{ij}(X_{r}^{t,x},Y_{r}^{t,x})d\left
\langle B^{i},B^{j}\right \rangle
_{r}]\\
& \leq \frac{1}{2}\hat{\mathbb{E}}[(B(D^{2}\psi(t,x),D\psi(t,x),\psi
(t,x),x),\langle B\rangle_{t+\delta}-\langle B\rangle_{t})]\\
& \ \ \ +(\partial_{t}\psi(t,x)+\langle
D\psi(t,x),b(x)\rangle+f(x,\psi
(t,x)))\delta+C(1+|x|+|x|^{2}+|x|^{3})\delta^{3/2}\\
&
\leq(\partial_{t}\psi(t,x)+F(D^{2}\psi(t,x),D\psi(t,x),\psi(t,x),x))\delta
+C(1+|x|+|x|^{2}+|x|^{3})\delta^{3/2},
\end{align*}
then it is easy to check that%
\[
\partial_{t}\psi(t,x)+F(D^{2}\psi(t,x),D\psi(t,x),\psi(t,x),x)\geq0.
\]
Thus $u$ is a viscosity subsolution of (\ref{BSDE-22}). Similarly we
can prove that $u$ is a viscosity supersolution of (\ref{BSDE-22}).
\end{proof}

\begin{example}
Let $B=(B^{1},B^{2})$ be a $2$-dimensional $G$-Brownian motion with
\[
G(A)=G_{1}(a_{11})+G_{2}(a_{22}),
\]
where%
\[
G_{i}(a)=\frac{1}{2}(\overline{\sigma}_{i}^{2}a^{+}-\underline{\sigma}_{i}%
^{2}a^{-}),\  \ i=1,2.
\]
In this case, we consider the following $1$-dimensional SDE:%
\[
dX_{s}^{t,x}=\mu X_{s}^{t,x}ds+\nu X_{s}^{t,x}d\left \langle
B^{1}\right \rangle _{s}+\sigma X_{s}^{t,x}dB_{s}^{2},\  \
X_{t}^{t,x}=x,
\]
where $\mu$, $\nu$ and $\sigma$ are constants.

The corresponding function $u$ is defined by\newline%
\[
u(t,x):=\hat{\mathbb{E}}[\varphi(X_{T}^{t,x})].
\]
Then%
\[
u(t,x)=\hat{\mathbb{E}}[u(t+\delta,X_{t+\delta}^{t,x})]
\]
and $u$ is the viscosity solution of the following PDE:
\[
\partial_{t}u+\mu x\partial_{x}u+2G_{1}(\nu x\partial_{x}u)+\sigma^{2}%
x^{2}G_{2}(\partial_{xx}^{2}u)=0,\ u(T,x)=\varphi(x).
\]

\end{example}

\begin{exercise}
For each $\xi \in L_{G}^{p}(\Omega_{t};\mathbb{R}^{n})$ with
$p\geq2$, show that SDE \textup{(\ref{BSDE-8})} has a unique
solution in $\bar{M}_{G}^{p}(t,T;\mathbb{R}^{n})$. Furthermore, show
that the following estimates hold.
\begin{equation*}
\hat{\mathbb{E}}[|X_{s}^{t,x}-X_{s}^{t,x^{\prime}}|^{p}]\leq
C|x-x^{\prime}|^{p},
\end{equation*}%
\begin{equation*}
\hat{\mathbb{E}}[|X_{s}^{t,x}|^{p}]\leq C(1+|x|^{p}),
\end{equation*}%
\begin{equation*}
\hat{\mathbb{E}}[|X_{t+\delta}^{t,x}-x|^{p}]\leq
C(1+|x|^{p})\delta^{p/2}.
\end{equation*}

\end{exercise}

\section*{Notes and Comments}
\addcontentsline{toc}{section}{Notes and Comments} This chapter is
mainly from Peng (2007) \cite{Peng2007b}.

There are many excellent books on It\^{o}'s stochastic calculus and
stochastic differential equations founded by It\^o's original paper
\cite{Ito}, as well as on martingale theory. Readers are referred to
Chung and Williams (1990) \cite{Chu-Will}, Dellacherie and Meyer
(1978 and 1982) \cite{DM}, He, Wang and Yan (1992) \cite{HWY},
It\^{o} and McKean (1965) \cite{Ito-McKean}, Ikeda and Watanabe
(1981) \cite{IW}, Kallenberg (2002) \cite{Kallenberg}, Karatzas and
Shreve (1988) \cite{KSh}, \O ksendal (1998) \cite{Oksendal}, Protter
(1990) \cite{Protter}, Revuz and Yor (1999)\cite{RM} and Yong and
Zhou (1999) \cite{Yong-Zhou}.

Linear backward stochastic differential equation (BSDE) was first
introduced by Bismut in (1973) \cite{Bis1973} and (1978)
\cite{Bis1978}. Bensoussan developed this approach in (1981)
\cite{Ben1981} and (1982) \cite{Ben1982}. The
 existence and uniqueness theorem of a general nonlinear BSDE, was
obtained in 1990 in Pardoux and Peng \cite{Pardoux-Peng}. The
present version of the proof was based on El Karoui, Peng and Quenez
(1997) \cite{EPQ}, which is also a very good survey on BSDE theory
and its applications, specially in finance. Comparison theorem of
BSDEs was obtained in Peng (1992) \cite{Peng1992} for the case when
$g$ is a $C^1$-function and then in \cite{EPQ} when $g$ is
Lipschitz. Nonlinear Feynman-Kac formula for BSDE was introduced by
Peng (1991) \cite{Peng1992c} and (1992) \cite{Peng1992b}. Here we obtain
the corresponding Feynman-Kac formula under the framework of
$G$-expectation. We also refer to Yong and Zhou (1999)
\cite{Yong-Zhou}, as well as in Peng (1997) \cite{Peng1997b} (in
1997, in Chinese) and (2004) \cite{Peng2004a} for systematic
presentations of BSDE theory. For contributions in the developments
of this theory, readers can be referred to the literatures listing
in the Notes and Comments in Chap. \ref{ch1}.

%
%
\newcommand{\1}{{\bf 1}}
\newcommand{\lun}{\mathbb{L}^{1}}
\newcommand{\lunb}{\mathbb{L}^{1}_b}
\newcommand{\lunc}{\mathbb{L}^{1}_c}
\newcommand{\lp}{\mathbb{L}^{p}}
\newcommand{\lpb}{\mathbb{L}^{p}_b}
\newcommand{\lpc}{\mathbb{L}^{p}_c}
\newcommand{\linf}{\mathbb{L}^{\infty}}
\newcommand{\linfb}{\mathbb{L}^{\infty}_b}
\newcommand{\linfc}{\mathbb{L}^{\infty}_c}
\newcommand{\N}{\mathbb{N}}
\newcommand{\Ne}{\N^{\ast}}
\newcommand{\R}{\mathbb{R}}
\newcommand{\Z}{\mathbb{Z}}
\newcommand{\C}{\mathbb{C}}
\newcommand{\E}{\mathbb{E}}
\newcommand{\AAA}{{\cal A}}
\newcommand{\BB}{{\cal B}}
\newcommand{\CC}{{\cal C}}
\newcommand{\DD}{{\cal D}}
\newcommand{\EE}{{\cal E}}
\newcommand{\FF}{{\cal F}}
\newcommand{\GG}{{\cal G}}
\newcommand{\HH}{{\cal H}}
\newcommand{\II}{{\cal I}}
\newcommand{\JJ}{{\cal J}}
\newcommand{\KK}{{\cal K}}
\newcommand{\LL}{{\cal L}}
\newcommand{\NN}{{\cal N}}
\newcommand{\MM}{{\cal M}}
\newcommand{\OO}{{\cal O}}
\newcommand{\PP}{{\cal P}}
\newcommand{\QQ}{{\cal Q}}
\newcommand{\RR}{{\cal R}}
\newcommand{\TT}{{\cal T}}
\newcommand{\UU}{{\cal U}}
\newcommand{\VV}{{\cal V}}
\newcommand{\WW}{{\cal W}}
\newcommand{\XX}{{\cal X}}
\newcommand{\YY}{{\cal Y}}
\newcommand{\ZZ}{{\cal Z}}
\newcommand{\tom}{\tilde{\Omega}}
\newcommand{\tp}{\tilde{P}}
\newcommand{\tB}{\tilde{B}}
\newcommand{\dis}{\displaystyle}
\newcommand{\lc}{{\cal L}^1 (\Om )}
\newcommand{\omu}{\overline{\mu}}
\newcommand{\umu}{\underline{\mu}}
\newcommand{\ep}{\epsilon}
\newcommand{\rgt}{\rightarrow}
\newcommand{\ced}{ \begin{flushright}$\Box$\end{flushright}}
\newcommand{\mb}[1]{\makebox{#1}}
\newcommand{\p}{^{\prime}}
\newcommand{\vphi}{\varphi}

\newcommand{\ET}[1]{\tilde{E}_P (#1 \mid \FF_t )}
\newcommand{\Om}{\Omega}
\newcommand{\tO}{\tilde{\Omega}}
\newcommand{\tF}{\tilde{\FF}}
\newcommand{\SSS}{{\cal S}}
\newcommand{\tP}{\widetilde{\PP}}
\newcommand{\cb}{C_b (\Omega )}
\newcommand{\ct}{\tilde{c}}
\newcommand{\mut}{\tilde{\mu}}
\setlength{\parindent}{0cm}
\newcommand{\wt}{\tilde{\omega}}
\newcommand{\so}{\sigma_0}
\newcommand{\Etx}{\E_{t,x}}
\newcommand{\vep}{\varepsilon}
\newcommand{\tX}{\tilde{X}}

\chapter{Capacity and Quasi-Surely Analysis for $G$-Brownian Paths} \label{ch6}
In this chapter, we first present a general framework for an upper
expectation defined on a metric space $(\Omega,
\mathcal{B}({\Omega}))$ and the  corresponding capacity to introduce
the quasi-surely analysis. The results are important for us to
obtain the pathwise analysis for $G$-Brownian motion.

\section{Integration Theory associated to an Upper Probability}

Let $\Omega$ be a complete separable metric space equipped with the
distance $d$, $\mathcal{B}(\Omega)$ the Borel $\sigma$-algebra of
$\Omega$ and $\mathcal{M}$ the collection of all probability
measures on $(\Omega ,\mathcal{B}(\Omega))$.

\begin{itemize}
\item $L^{0}(\Omega)$\label{l0omega}: the space of all $\mathcal{B}(\Omega)$-measurable real functions;

\item $B_{b}(\Omega)$\label{bbomega}: all bounded functions in $L^{0}(\Omega)$;

\item $C_{b}(\Omega)$\label{cbomega}: all continuous functions in $B_{b}(\Omega)$.
\end{itemize}

All along this section, we consider a given subset
$\mathcal{P}\subseteq \mathcal{M}$.

\subsection{Capacity associated to $\mathcal{P}$\  \ }

We denote
\[
c(A):=\sup_{P\in \mathcal{P}} P(A),\  \  \  \ A\in
\mathcal{B}(\Omega).
\]

One can easily verify the following theorem.

\begin{theorem}
The set function $c(\cdot)$ is a Choquet capacity, i.e. (see
\cite{Choquet,Del}),
\begin{enumerate}
\item $0\leq c(A)\leq1,\  \  \forall A\subset \Omega$.

\item If $A\subset B$, then $c(A)\leq c(B)$.

\item If $(A_{n})_{n=1}^{\infty}$ is a
sequence in $\mathcal{B}(\Omega)$, then $c(\cup A_{n})\leq \sum c(A_{n}%
)$.

\item If $(A_{n})_{n=1}^{\infty}$ is an increasing
sequence in $\mathcal{B}(\Omega)$: $A_{n}\uparrow A=\cup A_{n}$,
then $c(\cup A_{n})=\lim _{n\rightarrow \infty}c(A_{n})$.
\end{enumerate}
\end{theorem}

Furthermore, we have

\begin{theorem}
For each $A\in \mathcal{B}(\Omega)$, we have
\[
c(A) =\sup \{ c(K):\ K\makebox{ compact }\ K\subset A\}.
\]
\end{theorem}

\begin{proof} It is simply because
\[
c(A) =\sup_{P\in \mathcal{P}}\sup_{\substack{K\,
\text{compact}\\K\subset A}}P(K)=\sup_{\substack{K\,
\text{compact}\\K\subset A}}\sup_{P\in
\mathcal{P}}P(K)=\sup_{\substack{K\, \text{compact}\\K\subset
A}}c(K).
\]
\end{proof}

\begin{definition}
We use the standard capacity-related vocabulary: a set $A$ is
\textbf{polar}\index{Polar} if $c(A)=0$ and a property holds
``\textbf{quasi-surely}''\index{Quasi-surely} (q.s.)\label{qs} if it
holds outside a polar set.
\end{definition}

\begin{remark}
In other words, $A\in \mathcal{B} (\Omega)$ is polar if and only if
$P(A)=0$ for any $P\in \mathcal{P}$.
\end{remark}

We also have in a trivial way a Borel-Cantelli Lemma.

\begin{lemma}
{\label{BorelC}} Let $(A_{n} )_{n\in \mathbb{N}}$ be a sequence of
Borel sets such that
\[
\sum_{n=1}^{\infty} c(A_{n} )<\infty.
\]
Then $\limsup_{n\rightarrow\infty}A_{n}$ is polar .
\end{lemma}

\begin{proof} Applying the Borel-Cantelli Lemma under each probability
$P\in \mathcal{P}$.
\end{proof}

The following theorem is Prokhorov's theorem.

\begin{theorem} \label{newth6}
$\mathcal{P}$ is relatively compact if and only if for each
$\varepsilon >0$, there exists a compact set $K$ such that
$c(K^{c})<\varepsilon$.
\end{theorem}

The following two lemmas can be found in \cite{HuSt}.

\begin{lemma}
\label{Lemma1}$\mathcal{P}$ is relatively compact if and only if for
each sequence of closed sets $F_{n}\downarrow \emptyset$, we have
$c(F_{n})\downarrow0$.
\end{lemma}

\begin{proof} We outline the
proof for the convenience of readers.\\
\textquotedblleft$\Longrightarrow$\textquotedblright \ part: It
follows from Theorem \ref{newth6} that for each fixed
$\varepsilon>0$, there exists a compact set $K$ such that
$c(K^{c})<\varepsilon$. Note that $F_{n}\cap K \downarrow
\emptyset$, then there exists an $N>0$ such that $F_{n}\cap K
=\emptyset$ for $n \geq N$, which implies
$\lim_{n}c(F_n)<\varepsilon$. Since $\varepsilon$ can be arbitrarily
small, we obtain $c(F_{n})\downarrow0$.
\newline \textquotedblleft$\Longleftarrow$\textquotedblright \ part:
For each $\varepsilon>0$, let $(A_{i}^{k})_{i=1}^{\infty}$ be a
sequence of open balls of radius $1/k$ covering $\Omega$. Observe
that $(\cup_{i=1}^{n}A_{i}^{k})^{c}\downarrow \emptyset$, then there
exists an $n_k$ such that
$c((\cup_{i=1}^{n_k}A_{i}^{k})^{c})<\varepsilon 2^{-k}$. Set
$K=\overline{\cap_{k=1}^{\infty}\cup_{i=1}^{n_k}A_{i}^{k}}$. It is
easy to check that $K$ is compact and $c(K^{c})<\varepsilon$. Thus
by Theorem \ref{newth6}, $\mathcal{P}$ is relatively compact.
\end{proof}

\begin{lemma}\label{newle8}
Let $\mathcal{P}$ be weakly compact. Then for each sequence of
closed sets $F_{n}\downarrow F$, we have $c(F_{n})\downarrow c(F)$.
\end{lemma}

\begin{proof} We outline the
proof for the convenience of readers. For each fixed
$\varepsilon>0$, by the definition of $c(F_n)$, there exists a $P_n
\in \mathcal{P}$ such that $P_n(F_n) \geq c(F_n)-\varepsilon$. Since
$\mathcal{P}$ is weakly compact, there exist $P_{n_k}$ and $P \in
\mathcal{P}$ such that $P_{n_k}$ converge weakly to $P$. Thus
\[
P(F_m) \geq \limsup_{k\rightarrow\infty}P_{n_k}(F_m) \geq
\limsup_{k\rightarrow\infty}P_{n_k}(F_{n_k}) \geq
\lim_{n\rightarrow\infty}c(F_n)-\varepsilon.
\]
Letting $m\rightarrow\infty$, we get $P(F)\geq
\lim_{n\rightarrow\infty}c(F_n)-\varepsilon$, which yields
$c(F_{n})\downarrow c(F)$.
\end{proof}

Following \cite{HuSt} (see also \cite{Delbaen,FoSch}) the upper
expectation of $\mathcal{P}$ is defined as follows: for each $X\in
L^{0}(\Omega)$ such that $E_{P}[X]$ exists for each $P\in
\mathcal{P}$,
\[
\mathbb{E} [X]=\mathbb{E} ^{\mathcal{P}}[X]:=\sup_{P\in \mathcal{P}}%
E_{P}[X].
\]

It is easy to verify

\begin{theorem}\label{THM2}
The upper expectation $\mathbb{E} [\cdot]$ of the family
$\mathcal{P}$  is a sublinear expectation on $B_{b}(\Omega)$ as well
as on $C_{b}(\Omega)$, i.e.,
\begin{enumerate}
\item for all $X,Y$ in $B_{b}(\Omega)$, $X\geq
Y\Longrightarrow \mathbb{E} [X]\geq \mathbb{E} [Y]$.

\item for all $X,Y$ in $B_{b}(\Omega)$, $\mathbb{E} [X+Y]\leq
\mathbb{E} [X]+\mathbb{E} [Y]$.

\item for all $\lambda \geq0$, $X\in B_{b}(\Omega)$,
 $\mathbb{E} [\lambda X]=\lambda \mathbb{E} [X]$.

\item for all $c\in \mathbb{R}$, $X\in
B_{b}(\Omega)$ , $\mathbb{E} [X+c]=\mathbb{E} [X]+c$.

\end{enumerate}
\end{theorem}

Moreover, it is also easy to check

\begin{theorem}
We have
\begin{enumerate}
\item Let $\mathbb{E}[X_n]$ and $\mathbb{E}[\sum _{n=1}^{\infty}X_n]$ be finite. Then
$\mathbb{E}[\sum _{n=1}^{\infty}X_n] \leq \sum _{n=1}^{\infty}
\mathbb{E}[X_n].$

\item Let $X_n \uparrow X$ and $\mathbb{E}[X_n]$, $\mathbb{E}[X]$ be finite. Then
$\mathbb{E}[X_n] \uparrow \mathbb{E}[X]$.
\end{enumerate}
\end{theorem}

\begin{definition}
The functional $\mathbb{E} [\cdot]$ is said to be \textbf{regular}\index{Regular} if for each $\{X_{n}%
\}_{n=1}^{\infty}$ in $C_{b}(\Omega)$ such that $X_{n}\downarrow0$ \
on $\Omega$, we have $\mathbb{E} [X_{n}]\downarrow0$.
\end{definition}

Similar to Lemma \ref{Lemma1} we have:

\begin{theorem}
\label{Thm2}$\mathbb{\mathbb{E}}[\cdot]$ is regular if and only if
$\mathcal{P}$ is relatively compact.
\end{theorem}

\begin{proof}
\textquotedblleft$\Longrightarrow$\textquotedblright \ part: For
each sequence of closed subsets $F_{n}\downarrow \emptyset$ such
that $F_{n}$, $n=1,2,\cdots$, are non-empty (otherwise the proof is
trivial), there exists $\{g_{n}\}_{n=1}^{\infty}\subset
C_{b}(\Omega)$ satisfying
\[
0\leq g_{n}\leq1,\  \ g_{n}=1\text{ on }F_{n}\text{ and
}g_{n}=0\text{ on }\{ \omega \in \Omega:d(\omega,F_{n})\geq
\frac{1}{n}\}.
\]
We set $f_{n}=\wedge_{i=1}^{n}g_{i}$, it is clear that $f_{n}\in
C_{b}(\Omega)$ and $\mathbf{1}_{F_{n}}\leq f_{n}\downarrow0$. \
$\mathbb{E} [\cdot]$ is regular implies $\mathbb{E}
[f_{n}]\downarrow0$ and thus $c(F_{n})\downarrow 0$. It follows from
Lemma \ref{Lemma1} that $\mathcal{P}$ is relatively compact.
\newline \textquotedblleft$\Longleftarrow$\textquotedblright \ part:
For each $\left \{  X_{n}\right \}  _{n=1}^{\infty}\subset
C_{b}(\Omega)$ such
that $X_{n}\downarrow0$, we have%
\[
\mathbb{E} [X_{n}]=\sup_{P\in \mathcal{P}}E_{P}[X_{n}]=\sup_{P\in \mathcal{P}}%
\int_{0}^{\infty}P(\{X_{n}\geq t\})dt\leq
\int_{0}^{\infty}c(\{X_{n}\geq t\})dt.
\]
For each fixed $t>0$, $\{X_{n}\geq t\}$ is a closed subset and
$\{X_{n}\geq t\} \downarrow \emptyset$ as $n\uparrow \infty$. By
Lemma \ref{Lemma1}, $c(\{X_{n}\geq t\})\downarrow0$ and thus
$\int_{0}^{\infty}c(\{X_{n}\geq t\})dt\downarrow0$. Consequently
$\mathbb{E} [X_{n}]\downarrow0$.\end{proof}


\subsection{Functional spaces}

We set, for $p>0$,

\begin{itemize}
\item $\mathcal{L}^{p}:=\{X\in L^{0}(\Omega):\mathbb{E} [|X|^{p}]=\sup
_{P\in \mathcal{P}}E_{P}[|X|^{p}]<\infty \}$;\  \

\item $\mathcal{N}^{p}:=\{X\in L^{0}(\Omega):\mathbb{E}
[|X|^{p}]=0\}$;

\item $\mathcal{N}:=\{X\in L^{0}(\Omega):X=0$, $c$-q.s.$\}$.
\end{itemize}

It is seen that $\mathcal{L}^{p}$ and $\mathcal{N}^{p}$ are linear
spaces and $\mathcal{N}^{p}=\mathcal{N}$, for each $p>0$.

We denote $\mathbb{L}^{p}:=\mathcal{L}^{p}/\mathcal{N}$. As usual,
we do not take care about the distinction between classes and their
representatives.\newline

\begin{lemma}
{\label{markov}} Let $X\in \mathbb{L}^{p}$. Then for each $\alpha>0$
\[
c(\{| X | >\alpha \})\leq \displaystyle
\frac{\mathbb{E}[|X|^{p}]}{\alpha^{p}}.
\]

\end{lemma}

\begin{proof} Just apply Markov inequality under each $P\in
\mathcal{P}$.
\end{proof}

Similar to the classical results, we get the following proposition
and the proof is omitted which is similar to the classical
arguments.
\begin{proposition}
\label{Prop3}We have
\begin{enumerate}
\item For each $p\geq1$, $\mathbb{L}^{p}$ is a Banach
space under the norm $\left \Vert X\right \Vert _{p}:=\left(  \mathbb{{E}%
}[|X|^{p}]\right)  ^{\frac{1}{p}}$.

\item For each $p<1$,
$\mathbb{L}^{p}$
is a complete metric space under the distance \newline $d(X,Y):= \mathbb{{E}%
}[|X-Y|^{p}]  $.
\end{enumerate}
\end{proposition}

\  \

We set%
\begin{align*}
\mathcal{L}^{\infty}  &  :=\{X\in L^{0}(\Omega):\exists \text{ a
constant
}M\text{, s.t. }|X|\leq M,\  \  \text{q.s.}\}; \\
\mathbb{L}^{\infty}  &  :=\mathcal{L}^{\infty}/\mathcal{N}.
\end{align*}

\begin{proposition}
\label{Prop4}Under the norm
\[
\left \Vert X\right \Vert _{\infty}:=\inf \left \{  M\geq0:|X|\leq
M,\  \  \text{q.s.}\right \}  ,
\]
$\mathbb{L}^{\infty}$ is a Banach space.
\end{proposition}

\begin{proof}From $\left \{  \left \vert X\right \vert >\left \Vert X\right
\Vert _{\infty}\right \}  =\cup_{n=1}^{\infty}\left \{ \left \vert
X\right \vert \geq \left \Vert X\right \Vert
_{\infty}+\frac{1}{n}\right \}  $ we know that $\left \vert X\right
\vert \leq \left \Vert X\right \Vert _{\infty}$, q.s., then it is
easy to check that $\left \Vert \cdot \right \Vert_{\infty} $ is a
norm. The proof of the completeness of $\mathbb{L}^{\infty}$ is
similar to the classical result. \end{proof}

With respect to the distance defined on $\mathbb{L}^{p}$, $p>0$, we
denote by

\begin{itemize}
\item $\mathbb{L}^{p}_{b}$\label{Lpb} the completion of $B_{b}(\Omega)$.

\item $\mathbb{L}^{p}_{c}$\label{Lpc} the completion of $C_{b}(\Omega)$.
\end{itemize}

By Proposition \ref{Prop3}, we have%
\[
\mathbb{L}^{p}_{c}\subset \mathbb{L}^{p}_{b}\subset \mathbb{L}^{p},\
\  \ p>0.
\]

The following Proposition is obvious and the proof is left to the
reader.

\begin{proposition}
We have
\begin{enumerate}
\item Let $p,q>1$, $\frac{1}{p}+\frac{1}{q}=1$. Then $X\in \mathbb{L}^{p}$ and
$Y\in \mathbb{L}^{q}$ implies
\[
XY\in \mathbb{L}^{1}\text{ and }\mathbb{E} [|XY|]\leq \left(
\mathbb{E} [|X|^{p}]\right)  ^{\frac{1}{p}}\left(
\mathbb{E}[|Y|^{q}]\right)
^{\frac{1}{q}};%
\]
\newline Moreover $X\in \mathbb{L}^{p}_{c}$ and $Y\in \mathbb{L}^{q}_{c}$  implies
$XY\in \mathbb{L}^{1}_{c}$;

\item $\mathbb{L}^{p_{1}}\subset \mathbb{L}^{p_{2}}$,\ $\mathbb{L}^{p_{1}}%
_{b}\subset \mathbb{L}^{p_{2}}_{b}$, $\mathbb{L}^{p_{1}}_{c}\subset
\mathbb{L}^{p_{2}}_{c}$, $0<p_{2}\leq p_{1}\leq \infty$;

\item $\left \Vert X\right \Vert _{p}\uparrow \left \Vert X\right \Vert _{\infty}$,
for each $X\in \mathbb{L}^{\infty}$.


\end{enumerate}
\end{proposition}

\begin{proposition}
Let $p\in(0,\infty]$ and $(X_{n} )$ be a sequence in
$\mathbb{L}^{p}$ which converges to $X$ in $\mathbb{L}^{p}$. Then
there exists a subsequence $(X_{n_{k}})$ which converges to $X$
quasi-surely in the sense that it converges to $X$ outside a polar
set.
\end{proposition}

\begin{proof}Let us assume $p\in (0,\infty)$, the case
$p=\infty$ is obvious since the convergence in $\linf$ implies the
convergence in $\lp$ for all $p$.\\
One can extract a subsequence $(X_{n_k})$ such that
\[ \E [|X-X_{n_k}|^p]\leq 1/k^{p+2},\  \  \ k\in \N .\]
We set for all $k$
\[ A_k =\{ |X-X_{n_k}|>1/k \},\]
then as a consequence of the Markov property (Lemma \ref{markov})
and the Borel-Cantelli Lemma \ref{BorelC},
$c(\overline{\lim}_{k\rightarrow \infty}A_k )=0$. As it is clear
that on $(\overline{\lim}_{k\rightarrow \infty}A_k )^c$, $(X_{n_k})$
converges to $X$, the proposition is proved. \end{proof}

We now give a description of $\mathbb{L}^{p}_{b}$.

\begin{proposition}
\label{Prop5}For each $p>0$,%
\[
\mathbb{L}^{p}_{b}=\{X\in \mathbb{L}^{p}:\lim_{n\rightarrow
\infty}\mathbb{E} [|X|^{p}\mathbf{1}_{\{|X|>n\}}]=0\}.
\]

\end{proposition}

\begin{proof}We denote $J_{p}=\{X\in
\mathbb{L}^{p}:\lim_{n\rightarrow \infty}\E
[|X|^{p}\mathbf{1}_{\{|X|>n\}}]=0\}$. For each $X\in
J_{p}$ let $X_{n}=(X\wedge n)\vee(-n)\in B_{b}(\Omega)$. We have%
\[
\E [|X-X_{n}|^{p}]\leq \E [|X|^{p}\mathbf{1}%
_{\{|X|>n\}}]\rightarrow0\text{, as }n\rightarrow \infty \text{. }%
\]
Thus $X\in \mathbb{L}^{p}_{b}$.

On the other hand, for each $X\in \mathbb{L}^{p}_{b}$, we can find a
sequence $\left \{ Y_{n}\right \}
_{n=1}^{\infty}$ in $B_{b}(\Omega)$ such that $\E [|X-Y_{n}%
|^p]\rightarrow0$. Let $y_{n}=\sup_{\omega \in
\Omega}|Y_{n}(\omega)|$ and $X_{n}=(X\wedge y_{n})\vee (-y_{n})$.
Since $|X-X_{n}|\leq|X-Y_{n}|$, we have $\E [|X-X_{n}|^p
]\rightarrow0$. This clearly implies that for any sequence
$(\alpha_n )$
tending to $\infty$, $\lim_{n\rightarrow \infty} \E [|X-(X\wedge \alpha_{n})\vee (-\alpha_{n})|^p ]=0$.\\
Now we have, for all $n\in \N$,
\begin{align*}
\E [|X|^{p}\mathbf{1}_{\{|X|>n\}}]  &  =\mathbb{{E}%
}[(|X|-n+n)^{p}\mathbf{1}_{\{|X|>n\}}]\\
&  \leq(1\vee2^{p-1})\left( \E [(|X|-n)^{p}\mathbf{1}_{\{|X|>n\}}%
]+n^{p}c(|X|>n)\right).
\end{align*}
The first term of the right hand side tends to $0$ since%
\[
\E [(|X|-n)^{p}\mathbf{1}_{\{|X|>n\}}]= \E %
[|X-(X\wedge {n})\vee {(-n)}|^{p}]\rightarrow0.
\]
For the second term, since%
\[
\frac{n^{p}}{2^{p}}\mathbf{1}_{\{|X|>n\}}\leq(|X|-\frac{n}{2})^{p}%
\mathbf{1}_{\{|X|>n\}}\leq(|X|-\frac{n}{2})^{p}\mathbf{1}_{\{|X|>\frac{n}%
{2}\}},
\]
we have
\[
\frac{n^{p}}{2^{p}}c(|X|>n)=\frac{n^{p}}{2^{p}}\E [\mathbf{1}%
_{\{|X|>n\}}]\leq \mathbb{E}[(|X|-\frac{n}{2}%
)^{p}\mathbf{1}_{\{|X|>\frac{n}{2}\}}]\rightarrow0.
\]
Consequently $X\in J_{p}$.
\end{proof}

\begin{proposition}
\label{Prop12}Let $X\in \mathbb{L}^{1}_{b}$. Then for each
$\varepsilon>0$, there exists a $\delta>0$, such that for all $A\in
\mathcal{B}(\Omega)$ with $c(A)\leq \delta$, we have
$\mathbb{E}[|X|\mathbf{1}_{A}]\leq \varepsilon$.
\end{proposition}

\begin{proof}For each $\varepsilon>0$, by Proposition \ref{Prop5}, there exists an $N>0$ such that
$\E [|X|\mathbf{1}_{\{|X|>N\}}]\leq \frac{\varepsilon}{2}$. Take $\delta=\frac{\varepsilon}{2N}$. Then for a subset $A\in \mathcal{B}%
(\Omega)$ with $c(A)\leq \delta$, we have%
\begin{align*}
\E [|X|\mathbf{1}_{A}]  &  \leq \E [|X|\mathbf{1}%
_{A}\mathbf{1}_{\{|X|>N\}}]+\E [|X|\mathbf{1}_{A}\mathbf{1}%
_{\{|X|\leq N\}}]\\
&  \leq \E [|X|\mathbf{1}_{\{|X|>N\}}]+Nc(A)\leq \varepsilon
\text{.}%
\end{align*}
\end{proof}

It is important to note that not every element in $\mathbb{L}^{p}$
satisfies the condition $\lim_{n\rightarrow \infty}\mathbb{E}
[|X|^{p}\mathbf{1}_{\{|X|>n\}}]=0$. We give the following two
counterexamples to show that $\mathbb{L}^{1}$ and
$\mathbb{L}^{1}_{b}$ are different spaces even under the case that
$\mathcal{P}$ is weakly compact.

\begin{example}
\label{Exm2}Let $\Omega=\mathbb{N}$, $\mathcal{P}=\{P_{n}:n\in
\mathbb{N}\}$ where $P_{1}(\{1\})=1$ and
$P_{n}(\{1\})=1-\frac{1}{n}$, $P_{n}(\{n\})=\frac{1}{n}$, for
$n=2,3,\cdots$. $\mathcal{P}$ is weakly compact. We consider a
function $X$ on $\mathbb{N}$ defined by $X(n)=n$, $n\in \mathbb{N}$.
We have $\mathbb{E} [|X|]=2$ but $\mathbb{E}
[|X|\mathbf{1}_{\{|X|>n\}}]=1\not \rightarrow 0$. In this case,
$X\in \mathbb{L}^{1}$ but $X\not \in \mathbb{L}^{1}_{b}$.
\end{example}

\begin{example}
\label{Exm3}Let $\Omega=\mathbb{N}$, $\mathcal{P}=\{P_{n}:n\in
\mathbb{N}\}$ where $P_{1}(\{1\})=1$ and
$P_{n}(\{1\})=1-\frac{1}{n^{2}}$,  $P_{n}(\{kn\})=\frac{1}{n^{3}}$,
$k=1,2,\ldots,n$,for $n=2,3,\cdots$. $\mathcal{P}$ is weakly
compact. We consider a function $X$ on $\mathbb{N}$
defined by $X(n)=n$, $n\in \mathbb{N}$. We have $\mathbb{E}%
[|X|]=\frac{25}{16}$ and $n\mathbb{E}[\mathbf{1}_{\{|X|\geq
n\}}]=\frac {1}{n}\rightarrow0$, but
$\mathbb{E}[|X|\mathbf{1}_{\{|X|\geq n\}}]=\frac
{1}{2}+\frac{1}{2n}\not \rightarrow 0$. In this case, $X$ is in $\mathbb{L}%
^{1}$, continuous and $n\mathbb{E}[\mathbf{1}_{\{|X|\geq
n\}}]\rightarrow0$, but it is not in $\mathbb{L}^{1}_{b}$.
\end{example}

\subsection{Properties of elements in $\mathbb{L}^{p}_{c}$}

\begin{definition}
A mapping $X$ on $\Omega$ with values in a topological space is said
to be \textrm{quasi-continuous} (q.c.) if
\[
\forall \varepsilon>0 ,\; \makebox{there exists an open set } O\
\makebox{ with}\ c(O)<\varepsilon \makebox{ such that } X|_{O^{c}}
\makebox{ is continuous}.
\]

\end{definition}

\begin{definition}
We say that $X:\Omega \rightarrow \mathbb{R}$ has a quasi-continuous
version if there exists a quasi-continuous function $Y:\Omega
\rightarrow \mathbb{R}$ with $X=Y$ q.s..
\end{definition}

%
\begin{proposition}
{\label{qc}} Let $p>0$. Then each element in $\mathbb{L}^{p}_{c}$
has a quasi-continuous version.
\end{proposition}

\begin{proof} Let $(X_n )$ be a Cauchy sequence in $C_b (\Omega)$ for the distance on $\lp$.
Let us choose a subsequence $(X_{n_k} )_{k\geq 1}$ such that
\[ \E [|X_{n_{k+1}}-X_{n_k}|^p]\leq 2^{-2k},\  \  \  \forall k\geq 1,\]
and set for all $k$,
\[ A_k =\bigcup_{i=k}^{\infty} \{ |X_{n_{i+1}}-X_{n_i}|> 2^{-i/p}\}.\]
Thanks to the subadditivity property and the Markov inequality, we
have
\[ c(A_k )\leq \sum_{i=k}^{\infty} c(|X_{n_{i+1}}-X_{n_i}|>
2^{-i/p})\leq \sum_{i=k}^{\infty} 2^{-i}=2^{-k+1}.\] As a
consequence, $\lim_{k\rightarrow \infty} c(A_k )=0$, so the Borel
set $A=\bigcap_{k=1}^{\infty}A_{k}$ is polar.\\
As each $X_{n_k}$ is continuous, for all $k\geq 1$, $A_k$ is an open
set. Moreover, for all $k$, $(X_{n_i})$ converges uniformly on
$A_k^c$ so that the limit is continuous on each $A_k ^c$. This
yields the result.
\end{proof}

The following theorem gives a concrete characterization of the space
$\mathbb{L}^{p}_{c}$.

\begin{theorem}
\label{Thm8}For each $p>0$,%
\[
\mathbb{L}^{p}_{c}=\{X\in \mathbb{L}^{p} :X\text{ has a quasi-continuous version, }%
\lim_{n\rightarrow \infty}\mathbb{\mathbb{E} }[|X|^{p}\mathbf{1}_{\{|X|>n\}}%
]=0\}.
\]

\end{theorem}

\begin{proof}We denote
\[
J_{p}=\{X\in \mathbb{L}^{p}:X\text{ has a quasi-continuous version,
}\lim_{n\rightarrow \infty}\E [|X|^{p}\mathbf{1}_{\{|X|>n\}}]=0\}.
\]
Let $X\in \lpc$, we know by Proposition \ref{qc} that $X$ has a
quasi-continuous version. Since $X\in \mathbb{L}^{p}_b$,
we have by Proposition \ref{Prop5} that $\lim_{n\rightarrow \infty}%
\E [|X|^{p}\mathbf{1}_{\{|X|>n\}}]=0$. Thus $X\in J_{p}$.
\newline On the other hand, let $X\in J_{p}$ be quasi-continuous. Define $Y_{n}=(X\wedge n)\vee(-n)$ for all
$n\in \N$. As $\E [|X|^{p}\mathbf{1}_{\{|X|>n\}}]\rightarrow0$, we
have $\E [|X-Y_{n}|^{p}]\rightarrow0$. \\
Moreover, for all $n\in \N$, as  $Y_{n}$ is quasi-continuous , there
exists a closed set $F_{n}$ such that
$c(F_{n}^{c})<\frac{1}{n^{p+1}}$ and $Y_{n}$ is continuous on
$F_{n}$. It follows from Tietze's extension theorem that there
exists $Z_{n}\in C_{b}(\Omega)$ such that
\[
|Z_{n}|\leq n\text{ and }Z_{n}=Y_{n}\text{ \ on }F_{n}.
\]
We then have%
\[
\E [|Y_{n}-Z_{n}|^{p}]\leq(2n)^{p}c(F_{n}^{c})\leq \frac{(2n)^{p}%
}{n^{p+1}}.
\]
So $\mathbb{E}[|X-Z_{n}|^{p}]\leq(1\vee2^{p-1})(\mathbb{E}%
[|X-Y_{n}|^{p}]+\mathbb{E}[|Y_{n}-Z_{n}|^{p}])\ \rightarrow0,$ and
$X\in \lpc$.
\end{proof}

We give the following example to show that $\mathbb{L}^{p}_{c}$ is
different from $\mathbb{L}^{p}_{b}$ even under the case that
$\mathcal{P}$ is weakly compact.

\begin{example}
\label{Exm1}Let $\Omega=[0,1]$, $\mathcal{P}=\{ \delta_{x}:x\in
\lbrack0,1]\}$ is weakly compact. It is seen that
$\mathbb{L}^{p}_{c}=C_{b}(\Omega)$ which is different from
$\mathbb{L}^{p}_{b}$.
\end{example}

We denote $\mathbb{L}^{\infty}_{c}:=\{X\in \mathbb{L}^{\infty}:X$
has a quasi-continuous version$\}$, we have

\begin{proposition}
\label{Prop9}$\mathbb{L}^{\infty}_{c}$ is a closed linear subspace
of $\mathbb{L}^{\infty}$.
\end{proposition}

\begin{proof} For each Cauchy sequence $\left \{  X_{n}\right \}
_{n=1}^{\infty}$ of $\mathbb{L}^{\infty}_{c}$ under $\left \Vert
\cdot \right \Vert _{\infty}$, we can find a subsequence $\left \{  X_{n_{i}%
}\right \}  _{i=1}^{\infty}$ such that $\left \Vert X_{n_{i+1}}-X_{n_{i}%
}\right \Vert _{\infty}\leq 2^{-i}$. We may further assume that each
$X_n$ is quasi-continuous. Then it is easy to prove that for each
$\varepsilon>0$, there exists an open set $G$ such that
$c(G)<\varepsilon$ and $\left \vert X_{n_{i+1}}-X_{n_{i}}\right
\vert \leq 2^{-i}$ for all $i\geq 1$ on $G^{c}$, which implies that
the limit belongs to $\mathbb{L}^{\infty}_{c}$.
\end{proof}

As an application of Theorem \ref{Thm8}, we can easily get the
following results.

\begin{proposition}
\label{Prop10}Assume that $X:\Omega \rightarrow \mathbb{R}$ has a
quasi-continuous version and that there exists a function $f:\mathbb{R}^{+}%
\rightarrow \mathbb{R}^{+}$ satisfying $\lim_{t\rightarrow
\infty}\frac {f(t)}{t^{p}}=\infty$ and $\mathbb{E} [f(|X|)]<\infty$.
Then $X\in \mathbb{L}^{p}_{c}$.
\end{proposition}

\begin{proof} For each $\varepsilon>0$, there exists an $N>0$ such that
$\frac{f(t)}{t^{p}}\geq \frac{1}{\varepsilon}$, for all $t\geq N$. Thus%
\[
\E [|X|^{p}\mathbf{1}_{\{|X|>N\}}]\leq \varepsilon \mathbb{E}%
[f(|X|)\mathbf{1}_{\{|X|>N\}}]\leq \varepsilon \E [f(|X|)]\text{.}%
\]
Hence $\lim_{N\rightarrow \infty}\E [|X|^{p}\mathbf{1}_{\{|X|>N\}}%
]=0$. From Theorem \ref{Thm8} we infer $X\in$ $\lpc$. \end{proof}

\begin{lemma}
\label{Lemma14}Let $\left \{  P_{n}\right \}  _{n=1}^{\infty}\subset
\mathcal{P}$ converge weakly to $P\in \mathcal{P}$. Then for each
$X\in \mathbb{L}^{1}_{c}$, we have $E_{P_{n}}[X]\rightarrow
E_{P}[X]$.
\end{lemma}

\begin{proof} We may assume that $X$ is quasi-continuous, otherwise we can
consider its quasi-continuous version which does not change the
value $E_{Q}$ for each $Q \in \mathcal{P}$. For each
$\varepsilon>0$, there exists an $N>0$ such that $\E
[|X|\mathbf{1}_{\left \{ |X|>N\right \} }]<\frac{\varepsilon }{2}$.
Set $X_{N}=(X\wedge N)\vee (-N)$. We can find an open subset $G$
such that $c(G)<\frac{\varepsilon}{4N}$ and $X_{N}$ is continuous on
$G^{c}$. By Tietze's extension theorem, there exists $Y\in
C_{b}(\Omega)$ such that $\left \vert Y\right \vert \leq N$ and
$Y=X_{N}$ on $G^{c}$. Obviously, for each $Q \in \mathcal{P}$,
\begin{align*}
|E_{Q}[X]-E_{Q}[Y]|  & \leq E_{Q}[|X-X_{N}|]+E_{Q}[|X_{N}-Y|]\\
& \leq \frac{\varepsilon}{2}+2N\frac{\varepsilon}{4N}=\varepsilon \text{. }%
\end{align*}
It then follows that
\[
\limsup_{n\rightarrow \infty}E_{P_{n}}[X]\leq \lim_{n\rightarrow \infty}E_{P_{n}%
}[Y]+\varepsilon=E_{P}[Y]+\varepsilon \leq E_{P}[X]+2\varepsilon,
\]
and similarly $\liminf_{n\rightarrow \infty}E_{P_{n}}[X]\geq E_{P}%
[X]-2\varepsilon$. Since $\varepsilon$ can be arbitrarily small, we
then have $E_{P_{n}}[X]\rightarrow E_{P}[X]$.
\end{proof}

\begin{remark}
For continuous $X$, the above lemma is Lemma 3.8.7 in \cite{bog}.
\end{remark}

Now we give an extension of Theorem \ref{Thm2}.
\begin{theorem}
\label{Thm15}Let $\mathcal{P}$ be weakly compact and let $\left \{
X_{n}\right \}  _{n=1}^{\infty}\subset \mathbb{L}
^{1}_{c}$ be such that $X_{n}\downarrow X$, q.s.. Then $\mathbb{E}%
[X_{n}]\downarrow \mathbb{E}[X]$.
\end{theorem}

\begin{remark}
It is important to note that $X$ does not necessarily belong to $\mathbb{L}%
^{1}_{c}$.
\end{remark}

\begin{proof}For the case $\E [X]>-\infty$, if there exists a
$\delta>0$ such that $\E [X_{n}]>\E [X]+\delta$,
$n=1,2,\cdots$, we then can find a $P_{n}\in \mathcal{P}$ such that $E_{P_{n}%
}[X_{n}]>\E [X]+\delta-\frac{1}{n}$, $n=1,2,\cdots$. Since
$\mathcal{P}$ is weakly compact, we then can find a subsequence
$\left \{ P_{n_{i}}\right \}  _{i=1}^{\infty}$ that converges weakly
to some $P\in \mathcal{P}$. From which it follows that
\begin{align*}
E_{P}[X_{n_{i}}]  & =\lim_{j\rightarrow \infty}E_{P_{n_{j}}}[X_{n_{i}}%
]\geq \limsup_{j\rightarrow \infty}E_{P_{n_{j}}}[X_{n_{j}}]
\\&\geq \limsup _{j\rightarrow \infty}\{ \E
[X]+\delta-\frac{1}{n_{j}}\} =\E [X]+\delta \text{,\  \
}i=1,2,\cdots.
\end{align*}
Thus   $E_{P}[X]\geq \E [X]+\delta$. This contradicts  the
definition of $\E [\cdot]$. The proof for the case $\E [X]=-\infty$
is analogous.\end{proof}

We immediately have the following corollary.
\begin{corollary}\label{newco33}
Let $\mathcal{P}$ be weakly compact and let $\left \{  X_{n}\right
\} _{n=1}^{\infty}$ be a sequence in $\mathbb{L}^{1}_{c}$
decreasingly converging to $0$ q.s.. Then $\mathbb{E}
[X_{n}]\downarrow0$.
\end{corollary}

\subsection{Kolmogorov's criterion}

\begin{definition}
Let $I$ be a set of indices, $(X_{t})_{t\in I}$ and $(Y_{t} )_{t\in
I}$ be two processes indexed by $I$ . We say that $Y$ is \textrm{a
quasi-modification} of $X$ if for all $t\in I$, $X_{t} =Y_{t} $
q.s..
\end{definition}

\begin{remark}
In the above definition, quasi-modification is also called
modification in some papers.
\end{remark}

We now give a Kolmogorov criterion for a process indexed by
$\mathbb{R}^{d}$ with $d\in \mathbb{N}$.

\begin{theorem}\label{ch6t128}
Let $p>0$ and $(X_{t} )_{t\in[0,1]^{d}}$ be a process such that for
all $t\in[0,1]^{d}$, $X_{t}$ belongs to $\mathbb{L}^{p}$ . Assume
that there exist positive constants $c$ and $\varepsilon$ such that
\[
\mathbb{E} [|X_{t} -X_{s} |^{p}]\leq c|t-s|^{d+\varepsilon}.
\]
Then $X$ admits a modification $\tilde{X}$ such that
\[
\mathbb{E}\left[ \left(  \sup_{s\neq t} \displaystyle
\frac{|\tilde{X}_{t} -\tilde{X}_{s}|}{|t-s|^{\alpha}}\right)
^{p}\right] <\infty,
\]
for every $\alpha \in[0,\varepsilon/p)$. As a consequence, paths of
$\tilde{X}$ are quasi-surely H\"oder continuous of order $\alpha$
for every $\alpha<\varepsilon/p$ in the sense that there exists a
Borel set $N$ of capacity $0$ such that for all $w\in N^{c}$, the
map $t\rightarrow \tilde{X} (w)$ is H\"oder continuous of order
$\alpha$ for every $\alpha <\varepsilon/p$. Moreover, if
$X_t\in\mathbb{L}_c^p$ for each $t$, then we also have
$\tilde{X}_t\in\mathbb{L}_c^p$.
\end{theorem}

\begin{proof} Let $D$ be the set of dyadic points in $[0,1]^d$:
\[ D=\left \{ (\frac{i_1}{2^n} ,\cdots ,\frac{i_d}{2^n}); \ n\in \N , i_1
,\cdots ,i_d \in \{0,1,\cdots ,2^n \} \right \}.\] Let $\alpha \in
[0,\vep/p)$. We set
\[ M=\sup_{s,t\in D, s\neq t} \displaystyle \frac{|X_t
-X_s|}{|t-s|^{\alpha}}.\] Thanks to the classical Kolmogorov's
criterion (see Revuz-Yor \cite{RM}), we know that for any $P\in
\mathcal{P}$, $E_P [M^{p}]$ is finite and uniformly bounded with
respect to $P$ so that
\[ \E [M^{p}]=\sup_{P\in \PP}E_P
[M^{p}]<\infty.\] As a consequence, the map $t\rightarrow X_t $ is
uniformly continuous on $D$ quasi-surely and so we can define
\[ \forall t\in [0,1]^d ,\  \tX_t =\lim_{s\rightarrow t , s\in D} X_s
.\] It is now clear that $\tX$ satisfies the enounced properties.
\end{proof}

\section{$G$-expectation as an Upper Expectation}

In the following sections of this Chapter, let
$\Omega=C_{0}^{d}(\mathbb{R}^{+})$ denote the space of all
$\mathbb{R}^{d}-$valued continuous functions $(\omega_{t})_{t\in
\mathbb{R}^{+}}$, with $\omega_{0}=0$, equipped with the distance
\[
\rho(\omega^{1},\omega^{2}):=\sum_{i=1}^{\infty}2^{-i}[(\max_{t\in
\lbrack 0,i]}|\omega_{t}^{1}-\omega_{t}^{2}|)\wedge1],
\]
and let$\  \bar{\Omega}=(\mathbb{R}^{d})^{[0,\infty)}$ denote the
space of all $\mathbb{R}^{d}-$valued functions
$(\bar{\omega}_{t})_{t\in \mathbb{R}^{+}}$. Let
$\mathcal{B}(\Omega)$ denote the $\sigma$-algebra generated by all
open sets and let $\mathcal{B}(\bar{\Omega})$ denote the
$\sigma$-algebra generated by all finite dimensional cylinder sets.
The corresponding canonical process is $B_{t}(\omega)=\omega_{t}$
(respectively, $\bar{B}_{t}(\bar{\omega})=\bar{\omega}
_{t}$), $t\in \lbrack0,\infty)$ for $\omega \in \Omega$\ (respectively,\ $\bar{\omega}%
\in \bar{\Omega}$). The spaces of Lipschitzian cylinder functions on
$\Omega$
and $\bar{\Omega}$ are denoted respectively by%
\[
L_{ip}(\Omega):=\{
\varphi(B_{t_{1}},B_{t_{2}},\cdots,B_{t_{n}}):\forall
n\geq1,t_{1},\cdots,t_{n}\in \lbrack0,\infty),\forall \varphi \in C_{l.Lip}%
(\mathbb{R}^{d\times n})\},
\]
\[
L_{ip}(\bar{\Omega}):=\{
\varphi(\bar{B}_{t_{1}},\bar{B}_{t_{2}},\cdots
,\bar{B}_{t_{n}}):\forall n\geq1,t_{1},\cdots,t_{n}\in
\lbrack0,\infty ),\forall \varphi \in C_{l.Lip}(\mathbb{R}^{d\times
n})\}.
\]

Let $G(\cdot):\mathbb{S}(d)\rightarrow \mathbb{R}$ be a given
continuous monotonic and sublinear function. Following Sec.2 in
Chap.\ref{ch3}, we can construct the corresponding $G$-expectation
$\hat{\mathbb{E}}$ on $(\Omega,L_{ip}(\Omega))$. Due to the natural
correspondence of $L_{ip}(\bar{\Omega})$ and $L_{ip}(\Omega)$, we
also construct a sublinear expectation $\mathbb{\bar{E}}$ on
$(\bar{\Omega },L_{ip}(\bar{\Omega}))$ such that
$(\bar{B}_{t}(\bar{\omega}))_{t\geq0}$ is a $G$-Brownian motion.

The main objective of this section is to find a weakly compact
family of ($\sigma$-additive) probability measures on
$(\Omega,\mathcal{B}(\Omega))$ to represent $G$-expectation
$\hat{\mathbb{E}}$. The following lemmas are a variety of Lemma
\ref{I-le3} and \ref{I-le4}.

\begin{lemma}
\label{le3} Let $0\leq t_{1}<t_{2}<\cdots<t_{m}<\infty$ and $\{
\varphi _{n}\}_{n=1}^{\infty}\subset C_{l.Lip}(\mathbb{R}^{d\times
m})$ satisfy $\varphi_{n}\downarrow0$. Then
$\mathbb{\bar{E}}[\varphi_{n}(\bar
{B}_{t_{1}},\bar{B}_{t_{2}},\cdots,\bar{B}_{t_{m}})]\downarrow0$.
\end{lemma}

We denote  $\mathcal{T}:=\{
\underline{t}=(t_{1},\ldots,t_{m}):\forall m\in \mathbb{N},0\leq
t_{1}<t_{2}<\cdots<t_{m}<\infty \}.$

\begin{lemma}
\label{le4} Let $E$ be a finitely additive linear expectation
dominated by $\mathbb{\bar{E}}$ on $L_{ip}(\bar{\Omega})$. Then
there exists a unique probability measure $Q$ on
$(\bar{\Omega},\mathcal{B}(\bar{\Omega}))$ such that $E[X]=E_{Q}[X]$
for each $X\in L_{ip}(\bar{\Omega})$.
\end{lemma}

\begin{proof}
For each fixed $\underline{t}=(t_{1},\ldots,t_{m})\in \mathcal{T}$,
by Lemma \ref{le3}, for each sequence $\{
\varphi_{n}\}_{n=1}^{\infty}\subset C_{l.Lip}(\mathbb{R}^{d\times
m})$ satisfying $\varphi_{n}\downarrow0$, we
have $E[\varphi_{n}(\bar{B}_{t_{1}},\bar{B}_{t_{2}},\cdots,\bar{B}_{t_{m}%
})]\downarrow0$. By Daniell-Stone's theorem (see Appendix B), there
exists a unique probability
measure $Q_{\underline{t}}$ on $(\mathbb{R}^{d\times m},\mathcal{B}%
(\mathbb{R}^{d\times m}))$ such that $E_{Q_{\underline{t}}}[\varphi
]=E[\varphi(\bar{B}_{t_{1}},\bar{B}_{t_{2}},\cdots,\bar{B}_{t_{m}})]$
for each $\varphi \in C_{l.Lip}(\mathbb{R}^{d\times m})$. Thus we
get a family of
finite dimensional distributions $\{Q_{\underline{t}}:\underline{t}%
\in \mathcal{T}\}$. It is easy to check that
$\{Q_{\underline{t}}:\underline{t}\in \mathcal{T}\}$ is consistent.
Then by Kolmogorov's consistent theorem, there exists a probability
measure $Q$ on
$(\bar{\Omega},\mathcal{B}(\bar{\Omega}))$ such that $\{Q_{\underline{t}%
}:\underline{t}\in \mathcal{T}\}$ is the finite dimensional
distributions of $Q$. Assume that there exists another probability
measure $\bar{Q}$ satisfying the condition, by Daniell-Stone's
theorem, $Q$ and $\bar{Q}$ have the same finite-dimensional
distributions. Then by monotone class theorem, $Q=\bar{Q}$. The
proof is complete.
\end{proof}

\begin{lemma}
\label{le5} There exists a family of probability measures
$\mathcal{P}_{e}$ on $(\bar{\Omega},\mathcal{B}(\bar{\Omega}))$ such
that
\[
\mathbb{\bar{E}}[X]=\max_{Q\in \mathcal{P}_{e}}E_{Q}[X],\quad
\text{for}\ X\in L_{ip}(\bar{\Omega}).
\]

\end{lemma}

\begin{proof}
By the representation theorem of sublinear expectation and Lemma
\ref{le4}, it is easy to get the result.
\end{proof}

\  \  \

For this $\mathcal{P}_{e}$, we define the associated capacity:
\[
\tilde{c}(A):=\sup_{Q\in \mathcal{P}_{\! \!e}}Q(A),\quad A\in \mathcal{B}%
(\bar{\Omega}),
\]
and the upper expectation for each
$\mathcal{B}(\bar{\Omega})$-measurable real function $X$ which makes
the following definition meaningful:
\[
\mathbb{\tilde{E}}[X]:=\sup_{Q\in \mathcal{P}_{\! \!e}}E_{Q}[X].
\]

\begin{theorem}
\label{le6} For $(\bar{B})_{t\geq 0}$ , there exists a continuous
modification $(\tilde{B})_{t\geq 0}$ of $\bar{B}$ (i.e., $\tilde{c}
(\{ \tilde{B}_{t}\not =\bar{B}_{t}\})=0$, for each $t\geq0$) such
that $\tilde{B}_{0}=0$.
\end{theorem}

\begin{proof}
By Lemma \ref{le5}, we know that
$\mathbb{\bar{E}}=\mathbb{\tilde{E}}$ on $L_{ip}(\bar{\Omega})$. On
the other hand, we have
\[
\mathbb{\tilde{E}}[|\bar{B}_{t}-\bar{B}_{s}|^{4}]=\mathbb{\bar{E}}
[|\bar{B}_{t}-\bar{B}_{s}|^{4}]=d|t-s|^{2}\ \  \text{for}\  s,t\in
\lbrack0,\infty),
\]
where $d$ is a constant depending only on $G$. By Theorem \ref{ch6t128}, there exists a continuous modification $\tilde{B}$ of $\bar{B}%
$. Since $\tilde{c}(\{ \bar{B}_{0}\not =0\})=0$, we can set
$\tilde{B}_{0}=0$. The proof is complete.
\end{proof}

\  \  \

For each $Q\in \mathcal{P}_{\! \!e}$, let $Q\circ \tilde{B}^{-1}$
denote the probability measure on $(\Omega,\mathcal{B}(\Omega))$
induced by $\tilde{B}$
with respect to $Q$. We denote  $\mathcal{P}_{1}=\{Q\circ \tilde{B}^{-1}%
:Q\in \mathcal{P}_{\! \!e}\}$. By Lemma \ref{le6}, we get
\[
\mathbb{\tilde{E}}[|\tilde{B}_{t}-\tilde{B}_{s}|^{4}]=\mathbb{\tilde{E}}%
[|\bar{B}_{t}-\bar{B}_{s}|^{4}]=d|t-s|^{2},\forall s,t\in
\lbrack0,\infty).
\]
Applying the well-known result of moment criterion for tightness of
Kolmogorov-Chentsov's type (see Appendix B), we conclude that
$\mathcal{P}_{1}$ is tight. We denote by
$\mathcal{P}=\overline{\mathcal{P}}_{1}$ the closure of
$\mathcal{P}_{1}$ under the topology of weak convergence, then
$\mathcal{P}$ is weakly compact.

Now, we give the representation of $G$-expectation.

\begin{theorem}
\label{Gt34} For each continuous monotonic and sublinear function
$G:\mathbb{S}(d)\rightarrow \mathbb{R}$, let $\hat{\mathbb{E}}$ be
the corresponding $G$-expectation on $(\Omega,L_{ip}(\Omega))$. Then
there exists a weakly compact family of probability measures
$\mathcal{P}$ on
$(\Omega,\mathcal{B}(\Omega))$ such that%

\[
\hat{\mathbb{E}}[X]=\max_{P\in \mathcal{P}}E_{P}[X]\quad \text{for}\
X\in L_{ip}(\Omega).
\]

\end{theorem}

\begin{proof}
By Lemma \ref{le5} and Lemma \ref{le6}, we have
\[
\hat{\mathbb{E}}[X]=\max_{P\in \mathcal{P}_{1}}E_{P}[X]\quad
\text{for}\  X\in L_{ip}(\Omega).
\]
For each $X\in L_{ip}(\Omega)$, by Lemma \ref{le3}, we get $\mathbb{\hat{E}%
}[|X-(X\wedge N)\vee(-N)|]\downarrow0$ as $N\rightarrow \infty$.
Noting that $\mathcal{P}=\overline{\mathcal{P}}_{1}$, by the
definition of weak convergence, we get the result.
\end{proof}

\begin{remark}
In fact, we can construct the family $\mathcal{P}$ in a more
explicit way: Let $(W_{t})_{t\geq0}=(W_{t}^{i})_{i=1,t\geq0}^{d}$ be
a $d$-dimensional Brownian motion in this space. The filtration
generated by $W$ is denoted by $\mathcal{F}_{t}^{W}$. Now let
$\Gamma$ be the bounded, closed and convex subset in
$\mathbb{R}^{d\times d}$ such that
\[
G(A)=\sup_{\gamma \in \Gamma}\text{tr}[A\gamma \gamma^{T}],\  \ A\in
\mathbb{S}(d),
\]
(see (\ref{GaChII}) in Ch. II) and $\mathcal{A}_{\Gamma}$ the
collection of all $\Theta$-valued
$(\mathcal{F}_{t}^{W})_{t\geq0}$-adapted process $[0,\infty)$. We
denote
\[
B_{t}^{\gamma}:=\int_{0}^{T}\gamma_{s}dW_{s},\ t\geq0,\  \  \gamma
\in \mathcal{A}_{\Gamma}.
\]
and $\mathcal{P}_{0}$ the collection of probability measures on the
canonical space $(\Omega,\mathcal{B}(\Omega))$ induced by
$\{B^{\gamma}:\gamma \in \mathcal{A}_{\Gamma}\}$. Then
$\mathcal{P}=\overline{\mathcal{P}}_{0}$ (see \cite{DHP} for
details).
\end{remark}

\section{$G$-capacity and Paths of $G$-Brownian Motion}

According to Theorem \ref{Gt34}, we obtain a weakly compact family
of probability measures $\mathcal{P}$ on
$(\Omega,\mathcal{B}(\Omega))$ to represent $G$-expectation
$\hat{\mathbb{E}}[\cdot]$. For this $\mathcal{P}$, we define the
associated $G$-capacity:
\[
\hat{c}(A):=\sup_{P\in \mathcal{P}}P(A),\quad A\in
\mathcal{B}(\Omega)
\]
and upper expectation for each $X\in L^{0}(\Omega)$ which makes the
following definition meaningful:
\[
\bar{\mathbb{E}}[X]:=\sup_{P\in \mathcal{P}}E_{P}[X].
\]
By Theorem \ref{Gt34}, we know that
$\bar{\mathbb{E}}=\hat{\mathbb{E}}$ on $L_{ip}(\Omega)$, thus the
$\hat{\mathbb{E}}[|\cdot|]$-completion and the
$\bar{\mathbb{E}}[|\cdot|]$-completion of $L_{ip}(\Omega)$ are the
same.

For each $T>0$, we also denote by $\Omega_{T}=C_{0}^{d}([0,T])$
equipped with the distance
\[
\rho(\omega^{1},\omega^{2})=\left \Vert \omega^{1}-\omega^{2} \right
\Vert _{C_{0}^{d}([0,T])}:=\max_{0\leq t\leq
T}|\omega^{1}_{t}-\omega^{2}_{t}|.
\]

We now prove that $L_{G}^{1}(\Omega)=\mathbb{L}_{c}^{1}$, where $\mathbb{L}%
_{c}^{1}$ is defined in Sec.1. First, we need the following
classical approximation lemma.

\begin{lemma}
\label{le10} For each $X\in C_{b}(\Omega)$ and $n=1,2,\cdots$, we
denote
\[
X^{(n)}(\omega):=\inf_{\omega^{\prime}\in \Omega}\{X(\omega^{\prime
})+n\left \Vert \omega-\omega^{\prime}\right \Vert
_{C_{0}^{d}([0,n])}\} \quad \text{for}\ \omega \in \Omega.
\]
Then the sequence $\{X^{(n)}\}_{n=1}^{\infty}$ satisfies:

\begin{enumerate}
\item[\textup{(i)}] $-M\leq X^{(n)}\leq X^{(n+1)}\leq \cdots \leq X$, $M=\sup_{\omega \in
\Omega}|X(\omega)|;$\

\item[\textup{(ii)}] $|X^{(n)}(\omega_{1})-X^{(n)}(\omega_{2})|\leq n\left \Vert \omega
_{1}-\omega_{2}\right \Vert _{C_{0}^{d}([0,n])}\  \  \textup{for}\ \omega_{1}%
,\omega_{2}\in \Omega;$

\item[\textup{(iii)}] $X^{(n)}(\omega)\uparrow X(\omega)\  \  \textup{for}\ \omega \in \Omega.$
\end{enumerate}
\end{lemma}

\begin{proof}
\textup{(i)} is obvious.

For \textup{(ii)}, we have
\[%
\begin{array}
[c]{l}%
\qquad \quad X^{(n)}(\omega_{1})-X^{(n)}(\omega_{2})\\
\qquad \leq \sup_{\omega^{\prime}\in \Omega}\{[X(\omega^{\prime
})+n\left \Vert \omega_{1}-\omega^{\prime}\right \Vert _{C_{0}^{d}%
([0,n])}]-[X(\omega^{\prime})+n\left \Vert \omega_{2}-\omega^{\prime
}\right \Vert _{C_{0}^{d}([0,n])}]\} \\
 \qquad \leq n\left \Vert \omega_{1}-\omega_{2}\right
\Vert
_{C_{0}^{d}([0,n])}%
\end{array}
\]
and, symmetrically, $X^{(n)}(\omega_{2})-X^{(n)}(\omega_{1})\leq
n\left \Vert \omega_{1}-\omega_{2}\right \Vert _{C_{0}^{d}([0,n])}$.
Thus \textup{(ii)} follows.

We now prove \textup{(iii)}. For each fixed $\omega \in \Omega$, let $\omega_{n}%
\in \Omega$ be such that
\[
X(\omega_{n})+n\left \Vert \omega-\omega_{n}\right \Vert
_{C_{0}^{d}([0,n])}\leq X^{(n)}(\omega)+\frac{1}{n}.
\]
It is clear that $n\left \Vert \omega-\omega_{n}\right \Vert _{C_{0}^{d}%
([0,n])}\leq2M+1$ or $\left \Vert \omega-\omega_{n}\right \Vert _{C_{0}%
^{d}([0,n])}\leq \frac{2M+1}{n}$. Since $X\in C_{b}(\Omega)$, we get
$X(\omega_{n})\rightarrow X(\omega)$ as $n\rightarrow \infty$. We
have
\[
X(\omega)\geq X^{(n)}(\omega)\geq X(\omega_{n})+n\left \Vert
\omega-\omega _{n}\right \Vert _{C_{0}^{d}([0,n])}-\frac{1}{n},
\]
thus%
\[
n\left \Vert \omega-\omega_{n}\right \Vert _{C_{0}^{d}([0,n])}\leq
|X(\omega)-X(\omega_{n})|+\frac{1}{n}.
\]
We also have%
\begin{align*}
X(\omega_{n})-X(\omega)+n\left \Vert \omega-\omega_{n}\right \Vert _{C_{0}%
^{d}([0,n])} &  \geq X^{(n)}(\omega)-X(\omega)\\
&  \geq X(\omega_{n})-X(\omega)+n\left \Vert \omega-\omega_{n}\right
\Vert _{C_{0}^{d}([0,n])}-\frac{1}{n}.
\end{align*}
From the above two relations, we obtain%
\begin{align*}
|X^{(n)}(\omega)-X(\omega)| &  \leq|X(\omega_{n})-X(\omega)|+n\left
\Vert
\omega-\omega_{n}\right \Vert _{C_{0}^{d}([0,n])}+\frac{1}{n}\\
&  \leq2(|X(\omega_{n})-X(\omega)|+\frac{1}{n})\rightarrow0\  \text{as}%
\ n\rightarrow \infty.
\end{align*}
Thus \textup{(iii)} is obtained.
\end{proof}

\begin{proposition}
\label{pr11} For each $X\in C_{b}(\Omega)$ and $\varepsilon>0$,
there exists $Y\in L_{ip}(\Omega)$ such that
$\bar{\mathbb{E}}[|Y-X|]\leq \varepsilon$.
\end{proposition}

\begin{proof}
We denote $M=\sup_{\omega \in \Omega}|X(\omega)|$. By Theorem
\ref{Thm2} and Lemma \ref{le10}, we can find $\mu>0$, $T>0$ and
$\bar{X}\in C_{b}(\Omega_{T})$ such that
$\bar{\mathbb{E}}[|X-\bar{X}|]<\varepsilon/3$, $\sup_{\omega \in
\Omega}|\bar {X}(\omega)|\leq M$ and
\[
|\bar{X}(\omega)-\bar{X}(\omega^{\prime})|\leq \mu \left \Vert
\omega -\omega^{\prime}\right \Vert _{C_{0}^{d}([0,T])}\  \
\textup{for}\ \omega ,\omega^{\prime}\in \Omega.
\]
Now for each positive integer $n$, we introduce a mapping $\omega^{(n)}%
(\omega):\Omega \rightarrow\Omega$:
\[
\omega^{(n)}(\omega)(t)=\sum_{k=0}^{n-1}\frac{\mathbf{1}_{[t_{k}^{n}%
,t_{k+1}^{n})}(t)}{t_{k+1}^{n}-t_{k}^{n}}[(t_{k+1}^{n}-t)\omega(t_{k}%
^{n})+(t-t_{k}^{n})\omega(t_{k+1}^{n})]+\mathbf{1}_{[T,\infty)}(t)\omega(t),\
\]
where $t_{k}^{n}=\frac{kT}{n},\ k=0,1,\cdots,n$. We set $\bar{X}^{(n)}%
(\omega):=\bar{X}(\omega^{(n)}(\omega))$, then
\begin{align*}
|\bar{X}^{(n)}(\omega)-\bar{X}^{(n)}(\omega^{\prime})| &  \leq \mu
\sup
_{t\in \lbrack0,T]}|\omega^{(n)}(\omega)(t)-\omega^{(n)}(\omega^{\prime})(t)|\\
&  =\mu \sup_{k\in \lbrack0,\cdots,n]}|\omega(t_{k}^{n})-\omega^{\prime}%
(t_{k}^{n})|.
\end{align*}
We now choose a compact subset $K\subset \Omega$ such that $\mathbb{\bar{E}%
}[\mathbf{1}_{K^{C}}]\leq \varepsilon/6M$. Since $\sup_{\omega \in
K}\sup _{t\in
\lbrack0,T]}|\omega(t)-\omega^{(n)}(\omega)(t)|\rightarrow0$, as
$n\rightarrow \infty$, we can choose a sufficiently large $n_{0}$
such that
\begin{align*}
\sup_{\omega \in K}|\bar{X}(\omega)-\bar{X}^{(n_{0})}(\omega)| &
=\sup
_{\omega \in K}|\bar{X}(\omega)-\bar{X}(\omega^{(n_{0})}(\omega))|\\
&  \leq \mu \sup_{\omega \in K}\sup_{t\in \lbrack0,T]}|\omega(t)-\omega^{(n_{0}%
)}(\omega)(t)|\\
&  <\varepsilon/3.
\end{align*}
Set $Y:=\bar{X}^{(n_{0})}$, it follows that
\begin{align*}
\bar{\mathbb{E}}[|X-Y|] &  \leq \bar{\mathbb{E}}[|X-\bar{X}|]+\mathbb{\bar{E}%
}[|\bar{X}-\bar{X}^{(n_{0})}|]\\
&  \leq
\bar{\mathbb{E}}[|X-\bar{X}|]+\bar{\mathbb{E}}[\mathbf{1}_{K}|\bar
{X}-\bar{X}^{(n_{0})}|]+2M\bar{\mathbb{E}}[\mathbf{1}_{K^{C}}]\\
&  <\varepsilon.
\end{align*}
The proof is complete.
\end{proof}

\  \  \

By Proposition \ref{pr11}, we can easily get $L_{G}^{1}(\Omega)=\mathbb{L}%
_{c}^{1}$. Furthermore, we can get
$L_{G}^{p}(\Omega)=\mathbb{L}_{c}^{p}$, $\forall p>0$.

Thus, we obtain a pathwise description of $L_{G}^{p}(\Omega)$ for
each $p>0$:
\[
L_{G}^{p}(\Omega)=\{X\in L^{0}(\Omega):X\  \text{has a  quasi-continuous  version and }%
\lim_{n\rightarrow
\infty}\bar{\mathbb{E}}[|X|^{p}I_{\{|X|>n\}}]=0\}.
\]
Furthermore, $\bar{\mathbb{E}}[X]=\hat{\mathbb{E}}[X]$, for each
$X\in L_{G}^{1}(\Omega)$.

\begin{exercise}
Show that, for each $p>0$,
\[
L_{G}^{p}(\Omega_{T})=\{X\in L^{0}(\Omega_{T}):X\  \text{has a
quasi-continuous version and }\lim_{n\rightarrow
\infty}\bar{\mathbb{E}}[|X|^{p}I_{\{|X|>n\}}]=0\}.
\]
\newline
\end{exercise}

\section*{Notes and Comments}
\addcontentsline{toc}{section}{Notes and Comments}

The results of this chapter for $G$-Brownian motions were mainly
obtained by Denis, Hu and Peng (2008) \cite{DHP} (see also Denis and
Martini (2006) \cite{Denis-M} and the related comments after Chapter
III). Hu and Peng  (2009) \cite{HP} then have introduced an
intrinsic and simple approach. This approach can be regarded as a
combination and extension of the original Brownian motion
construction approach of Kolmogorov (for more general stochastic
processes) and a sort of cylinder Lipschitz functions technique
already introduced in Chap. \ref{ch3}. Section 1 is from \cite{DHP}
and Theorem \ref{Gt34} is firstly obtained in \cite{DHP}, whereas
contents of Sections 2 and 3 are mainly from \cite{HP}.

Choquet capacity was first introduced by Choquet (1953)
\cite{Choquet}, see also Dellacherie (1972) \cite{Del} and the
references therein for more properties. The capacitability of
Choquet capacity was first studied by Choquet \cite{Choquet} under
2-alternating case, see Dellacherie and Meyer (1978 and 1982)
\cite{DM}, Huber and Strassen (1972) \cite{HuSt} and the references
therein for more general case. It seems that the notion of upper
expectations was first discussed by Huber (1981) \cite{Huber} in
robust statistics. Recently, it was rediscovered in mathematical
finance, especially in risk measure, see Delbaen (1992, 2002)
\cite{Delbaen1,Delbaen}, F\"{o}llmer and Schied (2002, 2004)
\cite{FoSch} and etc..

\appendix
 \setcounter{tocdepth}{0}
\titlecontents{chapter}[0pt]{\vspace{1.0\baselineskip}\bfseries}
    {Appendix {\thecontentslabel}\quad}{}
    {\hspace{.5em}\titlerule*[8pt]{$\cdot$}\contentspage}


%
%

\chapter*{Appendix\  A\newline Preliminaries in Functional Analysis}\label{A}
\addcontentsline{toc}{chapter}{Appendix A\ Preliminaries in
Functional Analysis} \pagestyle{appendix}
\section{Completion of Normed Linear Spaces}
In this section, we suppose $\mathcal{H}$ is a linear space under
the norm $\|\cdot\|$.

\begin{definition}
$\{x_n\}\in\mathcal{H}$ is a {\textbf{Cauchy sequence}}\index{Cauchy
sequence}, if $\{x_n\}$ satisfies {\textbf{Cauchy's convergence
condition}}\index{Cauchy's convergence condition}:
$$\lim_{n, m\rightarrow \infty}\|x_n-x_m\|=0.$$
\end{definition}

\begin{definition}
A normed linear space $\mathcal{H}$ is called a {\textbf{Banach
space}} \index{Banach space}if it is
{\textbf{complete}}\index{Complete}, i.e., if every Cauchy sequence
$\{x_n\}$ of $\mathcal{H}$ converges strongly to a point $x_\infty$
of $\mathcal{H}$:
$$\lim_{n\rightarrow \infty}\|x_n-x_\infty\|=0.$$
Such a limit point $x_\infty$, if exists, is uniquely determined
because of the triangle inequality $\|x-x^\prime\|\leq
\|x-x_n\|+\|x_n-x^\prime\|$.
\end{definition}

The completeness of a Banach space plays an important role in
functional analysis. We introduce the following theorem of
completion.

\begin{theorem}
Let $\mathcal{H}$ be a normed linear space which is not complete.
Then $\mathcal{H}$ is isomorphic and isometric to a dense linear
subspace of a Banach-space $\tilde{\mathcal{H}}$, i.e., there exists
a one-to-one correspondence $x\leftrightarrow \tilde{x}$ of
$\mathcal{H}$ onto a dense linear subspace of $\tilde{\mathcal{H}}$
such that
$$\widetilde{x+y}=\tilde{x}+\tilde{y},\ \widetilde{\alpha x}=\alpha \tilde{x},\ \|\tilde{x}\|=\|x\|.$$
The space $\tilde{\mathcal{H}}$ is uniquely determined up to
isometric isomorphism.
\end{theorem}

For a proof see Yosida \cite{Yosida} (1980, p.56).

\section{The Hahn-Banach Extension Theorem}

\begin{definition}
Let $T_1$ and $T_2$ be two linear operators with domains $D(T_1)$
and $D(T_2)$ both contained in a linear space $\mathcal{H}$, and the
ranges $R(T_1)$ and $R(T_2)$ both contained in a linear space
$\mathcal{M}$. Then $T_1=T_2$ if and only if $D(T_1)=D(T_2)$ and
$T_1 x=T_2 x$ for all $x\in D(T_1)$. If $D(T_1)\subseteq D(T_2)$ and
$T_1 x= T_2 x$ for all $x\in D(T_1)$, then $T_2$ is called an
{\textbf{extension}}\index{Extension} of $T_1$, or $T_1$ is called a
\textbf{restriction}\index{Restriction} of $T_2$.
\end{definition}

\begin{theorem} {(\textbf{Hahn-Banach extension theorem in real linear spaces})}\index{Hahn-Banach extension theorem} Let
$\mathcal{H}$ be a real linear space and let $p(x)$ be a real-valued
function defined on $\mathcal{H}$ satisfying the following
conditions:
\begin{align*}
&p(x+y)\leq p(x)+p(y)\ \text{(subadditivity)};\\
&p(\alpha x)=\alpha p(x)\ \text{for}\ \alpha \geq 0~
\text{(positive~homogeneity)}.
\end{align*}
Let $\mathcal{L}$ be a real linear subspace of $\mathcal{H}$ and
$f_0$ be a real-valued linear functional defined on $\mathcal{L}:$
$$f_0(\alpha x+\beta y)=\alpha f_0(x)+\beta f_0(y)\ for\ x, y\in \mathcal{L}\ and\ \alpha, \beta\in\mathbb{R}.$$
Let $f_0$ satisfy $f_0(x)\leq p(x)$ on $\mathcal{L}$. Then there
exists a real-valued linear functional $F$ defined on $\mathcal{H}$
such that

\textup{(i)} $F$ is an extension of $f_0$, i.e., $F(x)=f_0(x)$ for
all $x\in \mathcal{L}$.

\textup{(ii)} $F(x)\leq p(x)$ \ \ \text{for}\ $x\in\mathcal{H}$.
\end{theorem}

For a proof see Yosida \cite{Yosida} (1980, p.102).

\begin{theorem} {(\textbf{Hahn-Banach extension theorem in normed linear spaces})} Let $\mathcal{H}$
be a normed linear space under the norm $\|\cdot\|$, $\mathcal{L}$
be a linear subspace of $\mathcal{H}$ and let $f_1$ be a continuous
linear functional defined on $\mathcal{L}$. Then there exists a
continuous linear functional $f$, defined on $\mathcal{H}$, such
that

\textup{(i)} $f$ is an extension of $f_1$.

\textup{(ii)} $\|f_1\|=\|f\|$.
\end{theorem}

For a proof see for example Yosida \cite{Yosida} (1980, p.106).

\section{Dini's Theorem and Tietze's Extension Theorem}

\begin{theorem}(\textbf{Dini's theorem}\index{Dini's theorem})
Let $\mathcal{H}$ be a compact topological space. If a monotone
sequence of bounded continuous functions converges pointwise to a continuous
function, then it also converges uniformly.
\end{theorem}

\begin{theorem}(\textbf{Tietze's extension theorem}) Let $\mathcal{L}$ be a closed subset of a normal space $\mathcal{H}$ and let
$f:\mathcal{L}\rightarrow \mathbb{R}$ be a continuous function. Then
there exists a continuous extension of $f$ to all of $\mathcal{H}$
with values in $\mathbb{R}$.
\end{theorem}

\chapter*{Appendix\  B\newline Preliminaries in Probability Theory}\label{B}
\addcontentsline{toc}{chapter}{Appendix B\ Preliminaries in
Probability Theory} \setcounter{section}{0}

\section{Kolmogorov's Extension Theorem}

Let $X$ be a random variable with values in $\mathbb{R}^n$ defined
on a probability space $(\Omega, \mathcal{F}, P)$. Denote by
$\mathcal{B}$ the Borel $\sigma$-algebra on $\mathbb{R}^n.$ We
define $X$'s law of distribution $P_X$ and its expectation $E_P$
with respect to $P$ as follows respectively:
$$P_X(B):=P(\omega: X(\omega)\in B);\ E_P[X]:=\int_{-\infty}^{+\infty} x P(dx),$$
where $B\in \mathcal{B}$.

In fact, we have $P_X(B)=E_P[{\mathbf{I}}_B(X)]$.

Now let $\{X_t\}_{t\in T}$ be a stochastic process with values in
$\mathbb{R}^n$ defined on a probability space $(\Omega, \mathcal{F},
P)$, where the parameter space $T$ is usually the halfline $[0,
+\infty)$.

\begin{definition}
The \textbf{finite dimensional distributions} \index{Finite
dimensional distributions}of the process $\{X_t\}_{t\in T}$ are the
measures $\mu_{t_1,\cdots, t_k}$ defined on $\mathbb{R}^{nk}, k=1,
2, \cdots$, by
$$\mu_{t_1,\cdots, t_k}(B_1\times \cdots \times B_k):=P[X_{t_1}\in B_1, \cdots, X_{t_k}\in B_k],\ t_i\in T,\ i=1, 2, \cdots, k,$$
where $B_i\in \mathcal{B}, i=1, 2, \cdots, k$.
\end{definition}

The family of all finite-dimensional distributions determine many
(but not all) important properties of the process $\{X_t\}_{t\in
T}$.

Conversely, given a family $\{\nu_{t_1,\cdots, t_k}: t_i\in T, i=1,
2, \cdots, k, k\in \mathbb{N}\}$ of probability measures on
$\mathbb{R}^{nk}$, it is important to be able to construct a
stochastic process $(Y_t)_{t\in T}$ with $\nu_{t_1,\cdots, t_k}$ as
its finite-dimensional distributions. The following famous theorem
states that this can be done provided that $\{\nu_{t_1,\cdots,
t_k}\}$ satisfy two natural consistency conditions.

\begin{theorem} {(\textbf{Kolmogorov's extension theorem})}\index{Kolmogorov's extension theorem}
For all $t_1$, $t_2$, $\cdots$, $t_k$, $k\in \mathbb{N}$, let
$\nu_{t_1,\cdots, t_k}$ be probability measures on $\mathbb{R}^{nk}$
such that
$$\nu_{t_{\pi(1)},\cdots, t_{\pi(k)}}(B_1\times \cdots \times B_k)=\nu_{t_1,\cdots, t_k}(B_{\pi^{-1}(1)}\times \cdots \times B_{\pi^{-1}(k)})$$
for all permutations $\pi$ on $\{1, 2, \cdots, k\}$ and
$$\nu_{t_1,\cdots, t_k}(B_1\times \cdots \times B_k)=\nu_{t_1,\cdots, t_k, t_{k+1}, \cdots, t_{k+m}}(B_1\times \cdots \times B_k\times \mathbb{R}^n \times \cdots \times \mathbb{R}^n)$$
for all $m\in \mathbb{N}$, where the set on the right hand side has
a total of $k+m$ factors.

Then there exists a probability space $(\Omega, \mathcal{F}, P)$ and
a stochastic process $(X_t)$ on $\Omega$, $X_t: \Omega \rightarrow
\mathbb{R}^n$, such that
$$\nu_{t_1,\cdots, t_k}(B_1\times \cdots \times B_k)=P[X_{t_1}\in B_1, \cdots, X_{t_k}\in B_k]$$
for all $t_i\in T$ and all Borel sets $B_i$, $i=1, 2, \cdots, k$,
$k\in \mathbb{N}.$
\end{theorem}

For a proof see Kolmogorov \cite{K} (1956, p.29).

\section{Kolmogorov's Criterion}

\begin{definition}
Suppose that $(X_t)$ and $(Y_t)$ are two stochastic processes
defined on $(\Omega, \mathcal{F}, P)$. Then we say that $(X_t)$ is a
{\textbf{version}}\index{Version} of (or a
{\textbf{modification}}\index{Modification} of) $(Y_t)$ if
$$P(\{\omega: X_t(\omega)=Y_t(\omega)\})=1\ \text{for~all}~ t.$$
\end{definition}

\begin{theorem} {(\textbf{Kolmogorov's continuity
criterion})}\index{Kolmogorov's continuity criterion} Suppose that
the process $X=\{X_t\}_{t\geq 0}$ satisfies the following condition:
for all $T> 0$ there exist positive constants $\alpha$, $\beta$, $D$
such that
$$E[|X_t-X_s|^\alpha]\leq D|t-s|^{1+\beta},\ \ \ \ 0\leq s, t\leq T.$$
Then there exists a continuous version of $X$.
\end{theorem}

For a proof see  Stroock and Varadhan \cite{SV} (1979, p.51).

Let $E$ be a metric space and $\mathcal{B}$ be the Borel
$\sigma$-algebra on $E$. We recall a few facts about the weak
convergence of probability measures on $(E, \mathcal{B})$. If $P$ is
such a measure, we say that a subset $A$ of $E$ is a $P$-continuity
set if $P(\partial A)=0$, where $\partial A$ is the boundary of $A$.

\begin{proposition}
For probability measures $P_n (n\in \mathbb{N})$ and $P$, the
following conditions are equivalent: \textup{\item{(i)}} For every
bounded continuous function $f$ on $E$,
$$\lim_{n\rightarrow\infty} \int f dP_n= \int f dP;$$
\textup{\item{(ii)}} For every bounded uniformly continuous function
$f$ on $E$,
$$\lim_{n\rightarrow\infty} \int f d P_n=\int f d P;$$
\textup{\item{(iii)}} For every closed subset $F$ of $E$,
$\limsup_{n\rightarrow\infty} P_n(F)\leq P(F)$; \textup{\item{(iv)}}
For every open subset $G$ of $E$, $\liminf_{n\rightarrow\infty}
P_n(G)\geq P(G);$ \textup{\item{(v)}} For every $P$-continuity set
$A$, $\lim_{n\rightarrow\infty} P_n(A)=P(A).$
\end{proposition}

\begin{definition}
If $P_n$ and $P$ satisfy the equivalent conditions of the preceding
proposition, we say that $(P_n)$ {\textbf{converges
weakly}}\index{Converge weakly} to $P$.
\end{definition}

Now let $\pi$ be a family of probability measures on $(E,
\mathcal{B}).$

\begin{definition}
A family $\pi$ is {\textbf{weakly relatively compact}}\index{Weakly
relatively compact} if every sequence of $\pi$ contains a weakly
convergent subsequence.
\end{definition}

\begin{definition}
A family $\pi$ is {\textbf{tight}}\index{Tight} if for every
$\varepsilon\in (0, 1)$, there exists a compact set $K_\varepsilon$
such that
$$P(K_\varepsilon)\geq 1-\varepsilon\ \ \ \text{for\ every}\ P\in \pi. $$
\end{definition}

With this definition, we have the following theorem.

\begin{theorem} ({\textbf{Prokhorov's criterion}})\index{Prokhorov's criterion}
If a family $\pi$ is tight, then it is weakly relatively compact. If
$E$ is a Polish space (i.e., a separable completely metrizable
topological space), then a weakly relatively compact family is
tight.
\end{theorem}

\begin{definition}
If $(X_n)_{n\in \mathbb{N}}$ and $X$ are random variables taking
their values in a metric space $E$, we say that $(X_n)$
\textbf{converges in distribution}\index{Converge in distribution}
or {\textbf{converges in law}}\index{Converge in law} to $X$ if
their laws $P_{X_n}$ converge weakly to the law $P_X$ of $X$.
\end{definition}

We stress the fact that the $(X_n)$ and $X$ need not be defined on
the same probability space.

\begin{theorem} ({\textbf{Kolmogorov's criterion for weak
compactness}})\index{Kolmogorov's criterion for weak compactness}
Let $\{X^n\}$ be a sequence of $\mathbb{R}^d$-valued continuous
processes defined on probability spaces $(\Omega^n, \mathcal{F}^n,
P^n)$ such that \textup{\item{(i)}} the family $\{P^n_{X_0^n}\}$ of
initial laws is tight in $\mathbb{R}^d$. \textup{\item{(ii)}} there
exist three strictly positive constants $\alpha$, $\beta$, $\gamma$
such that for each $s, t\in \mathbb{R}_+$ and each $n$,
$$E_{P^n}[|X_s^n-X_t^n|^\alpha]\leq \beta |s-t|^{\gamma+1},$$
then the set $(P^n_{X^n})$ of the laws of the $(X_n)$ is weakly
relatively compact.
\end{theorem}

For the proof see Revuz and Yor \cite{RM} (1999, p.517)

\section{Daniell-Stone Theorem}

Let $(\Omega, \mathcal{F}, \mu)$ be a measure space, on which we can
define integration. One essential properties of integration is its
linearity, thus it can be seen as a linear functional on
$L^1(\Omega, \mathcal{F}, \mu)$. This idea leads to another approach
to define integral--Daniell's integral.

\begin{definition}
Let $\Omega$ be an abstract set and $\mathcal{H}$ be a linear space
formed by a family of real valued functions. $\mathcal{H}$ is called
a {\textbf{vector lattice}}\index{Vector lattice} if
$$f\in \mathcal{H}\Rightarrow |f|\in \mathcal{H}, f\wedge 1 \in \mathcal{H}.$$
\end{definition}

\begin{definition}
Suppose that $\mathcal{H}$ is a vector lattice on $\Omega$ and $I$
is a positive linear functional on $\mathcal{H}$, i.e.,
\begin{align*}&f, g\in \mathcal{H}, \alpha, \beta\in \mathbb{R}\Rightarrow
I(\alpha f+\beta g)=\alpha I(f)+\beta I(g);\\ &f\in \mathcal{H},
f\geq 0 \Rightarrow I(f)\geq 0.\end{align*} If $I$ satisfies the
following condition:
$$f_n\in \mathcal{H}, f_n \downarrow 0 \Rightarrow I(f_n)\rightarrow 0,$$
or equivalently,
$$f_n\in \mathcal{H}, f_n \uparrow f\in \mathcal{H}\Rightarrow I(f)=\lim_{n \rightarrow \infty} I(f_n),$$
then $I$ is called a {\textbf{Daniell's integral}}\index{Daniell's
integral} on $\mathcal{H}$.
\end{definition}

\begin{theorem}\label{TheoremDL} {(\textbf{Daniell-Stone theorem})}\index{Daniell-Stone theorem}
Suppose that $\mathcal{H}$ is a vector lattice on $\Omega$ and $I$
is a Daniell integral on $\mathcal{H}$. Then there exists a measure
$\mu\in \mathcal{F}$, where $\mathcal{F}:=\sigma(f: f\in
\mathcal{H})$, such that $\mathcal{H}\subset L^1(\Omega,
\mathcal{F}, \mu)$ and $I(f)=\mu(f)$, $\forall f\in \mathcal{H}$.
Furthermore, if $1\in \mathcal{H}_+^*$, where $\mathcal{H}_+^*:=\{f:
\exists f_n \geq 0, f_n\in \mathcal{H}\ such\ that\ f_n \uparrow
f\}$, then this measure $\mu$ is unique and is $\sigma$-finite.
\end{theorem}

For the proof see Dellacherie and Meyer \cite{DM} (1978, p.59),
Dudley \cite{Dudley} (1995, p.142), or  Yan \cite{Yan} (1998, p.74).

\chapter*{Appendix\  C\newline Solutions of Parabolic Partial Differential Equation}\label{C}
\addcontentsline{toc}{chapter}{Appendix C\ Solutions of Parabolic
Partial Differential Equation} \setcounter{section}{0}

\section{The Definition of Viscosity Solutions}

The notion of viscosity solutions was firstly introduced by Crandall
and Lions (1981) \cite{CL81} and (1983) \cite{CrandallL} (see also
Evans's contribution (1978) \cite{Evans78} and (1980)
\cite{Evans80}) for the first-order Hamilton-Jacobi equation, with
uniqueness proof given in \cite{CrandallL}. The the proof of
second-order case for Hamilton-Jacobi-Bellman equations was firstly
developed by Lions (1982) \cite{Lions82} and (1983) \cite{Lions83}
using stochastic control verification arguments. A breakthrough was
achieved in the second-order PDE theory by  Jensen (1988)
\cite{Jensen}. For all other important contributions in the
developments of this theory we refer to the well-known user's guide
by Crandall, Ishii and Lions (1992) \cite{CIL}. For reader's
convenience, we systematically interpret some parts of \cite{CIL}
required in this book into it's parabolic version. However, up to my
knowledge, the presentation and the related proof for the domination
theorems seems to be a new generalization of the maximum principle
presented in \cite{CIL}. Books on this theory are, among others,
Barles (1994) \cite{barles}, Fleming, and Soner (1992) \cite{FS},
Yong and Zhou (1999) \cite{Yong-Zhou}.

Let $T>0$ be fixed and let $\mathcal{O}\subset \lbrack0,T]\times
\mathbb{R}^{N}$.
We set%
\[
USC(\mathcal{O})=\{ \text{upper semicontinuous functions}\ u:\mathcal{O}%
\rightarrow \mathbb{R}\},
\]%
\[
LSC(\mathcal{O})=\{ \text{lower semicontinuous functions}\ u:\mathcal{O}%
\rightarrow \mathbb{R}\}.
\]
Consider the following parabolic PDE:%
\begin{equation}
\left \{
\begin{array}
[c]{l}%
\text{(E)}\  \partial_{t}u-G(t,x,u,Du,D^{2}u)=0\  \text{on}\
(0,T)\times
\mathbb{R}^{N},\\
\text{(IC)}\ u(0,x)=\varphi(x)\  \text{for}\ x\in \mathbb{R}^{N},
\end{array}
\right.  \label{PE-G}%
\end{equation}
where $G:[0,T]\times \mathbb{R}^{N}\times \mathbb{R}\times \mathbb{R}^{N}%
\times \mathbb{S}(N)\rightarrow \mathbb{R}$, $\varphi \in
C(\mathbb{R}^{N})$. We always suppose that $G$ is continuous and
satisfies the following degenerate elliptic condition:
\begin{equation}
G(t,x,r,p,X)\geq G(t,x,r,p,Y)\  \  \text{whenever}\ X\geq Y. \label{DE}%
\end{equation}
Next we recall the definition of viscosity solutions from Crandall, Ishii and Lions \cite{CIL}. Let $u:(0,T)\times \mathbb{R}%
^{N}\rightarrow \mathbb{R}$ and $(t,x)\in(0,T)\times
\mathbb{R}^{N}$. We denote by $\mathcal{P}^{2,+}u(t,x)$ (the
\textquotedblleft \textbf{parabolic superjet}\index{Parabolic
superjet}\textquotedblright \ of $u$ at $(t,x)$) the set of triples
$(a,p,X)\in \mathbb{R}\times \mathbb{R}^{N}\times \mathbb{S}(N)$
such
that%
\begin{align*}
u(s,y)  &  \leq u(t,x)+a(s-t)+\langle p,y-x\rangle \\
& \ \ \ +\frac{1}{2}\langle X(y-x),y-x\rangle+o(|s-t|+|y-x|^{2}).
\end{align*}
We define
\begin{align*}
\mathcal{\bar{P}}^{2,+}u(t,x):=  &  \{(a,p,X)\in \mathbb{R}\times \mathbb{R}%
^{N}\times \mathbb{S}(N):\exists(t_{n},x_{n},a_{n},p_{n},X_{n})\\
&  \  \  \text{such that}\ (a_{n},p_{n},X_{n})\in \mathcal{P}^{2,+}u(t_{n}%
,x_{n})\  \text{and}\  \\
&  \  \ (t_{n},x_{n},u(t_{n},x_{n}),a_{n},p_{n},X_{n})\rightarrow
(t,x,u(t,x),a,p,X)\}.
\end{align*}
Similarly, we define $\mathcal{P}^{2,-}u(t,x)$ (the
\textquotedblleft
\textbf{parabolic subjet}\index{Parabolic subjet}\textquotedblright \ of $u$ at $(t,x)$) by $\mathcal{P}%
^{2,-}u(t,x):=-\mathcal{P}^{2,+}(-u)(t,x)$ and
$\mathcal{\bar{P}}^{2,-}u(t,x)$ by
$\mathcal{\bar{P}}^{2,-}u(t,x):=-\mathcal{\bar{P}}^{2,+}(-u)(t,x)$.

\begin{definition}
\textup{(i)} A \textbf{viscosity subsolution}\index{Viscosity
subsolution} of (E) on $(0,T)\times \mathbb{R}^{N}$ is a function
$u\in USC((0,T)\times \mathbb{R}^{N})$ such that for each
$(t,x)\in(0,T)\times \mathbb{R}^{N}$,
\[
a-G(t,x,u(t,x),p,X)\leq0 \  \text{for}\ (a,p,X)\in \mathcal{P}
^{2,+}u(t,x);
\]
likewise, a \textbf{viscosity supersolution}\index{Viscosity
supersolution} of (E) on $(0,T)\times \mathbb{R}^{N}$ is a function
$v\in LSC((0,T)\times \mathbb{R}^{N})$ such that for each
$(t,x)\in(0,T)\times \mathbb{R}^{N}$,
\[
a-G(t,x,v(t,x),p,X)\geq0 \ \ \text{for}\ (a,p,X)\in \mathcal{P}
^{2,-}v(t,x);
\]
and a \textbf{viscosity solution}\index{Viscosity solution} of (E)
on $(0,T)\times \mathbb{R}^{N}$ is a function that is simultaneously
a viscosity subsolution and a viscosity supersolution of (E) on
$(0,T)\times \mathbb{R}^{N}$.\newline \indent\textup{(ii)} A
function $u\in USC([0,T)\times \mathbb{R}^{N})$ is called a
\textbf{viscosity subsolution} of \textup{(\ref{PE-G})} on
$[0,T)\times \mathbb{R}^{N}$ if $u$ is a viscosity subsolution
of (E) on $(0,T)\times \mathbb{R}^{N}$ and $u(0,x)\leq \varphi(x)\ $%
for$\ x\in \mathbb{R}^{N}$; the appropriate notions of a viscosity
supersolution and a viscosity solution of \textup{(\ref{PE-G})} on
$[0,T)\times \mathbb{R}^{N}$ are then obvious.
\end{definition}
We now give the following equivalent definition (see Crandall, Ishii
and Lions \cite{CIL}).

\begin{definition}
A viscosity subsolution \index{Viscosity subsolution}of (E), or
$G$-subsolution, on $(0,T)\times \mathbb{R}^{N}$ is a function $u\in
USC((0,T)\times \mathbb{R}^{N})$ such that for all $(t,x)\in
(0,T)\times \mathbb{R}^{N}$, $\phi \in C^{2}((0,T)\times
\mathbb{R}^{N})$
such that $u(t,x)=\phi (t,x)$ and $u<\phi $ on $(0,T)\times \mathbb{R}%
^{N}\backslash (t,x)$, we have%
\begin{equation*}
\partial _{t}\phi (t,x)-G(t,x,\phi (t,x),D\phi (t,x),D^{2}\phi (t,x))\leq 0;
\end{equation*}%
likewise, a viscosity supersolution\index{Viscosity supersolution}
of (E),  or $G$-supersolution,  on $(0,T)\times \mathbb{R}^{N}$
is a function $v\in LSC((0,T)\times \mathbb{R}^{N})$ such that for all $%
(t,x)\in (0,T)\times \mathbb{R}^{N}$, $\phi \in C^{2}((0,T)\times \mathbb{R}%
^{N})$ such that $u(t,x)=\phi (t,x)$ and $u>\phi $ on $(0,T)\times \mathbb{R}%
^{N}\backslash (t,x)$, we have%
\begin{equation*}
\partial _{t}\phi (t,x)-G(t,x,\phi (t,x),D\phi (t,x),D^{2}\phi (t,x))\geq 0;
\end{equation*}%
and a viscosity solution \index{Viscosity solution}of (E) on
$(0,T)\times \mathbb{R}^{N}$ is a function that is simultaneously a
viscosity subsolution and a viscosity supersolution of (E) on
$(0,T)\times \mathbb{R}^{N}$. The definition of a viscosity solution
of \textup{(\ref{PE-G})} on $[0,T)\times \mathbb{R}^{N}$ is the same
as the above definition.
\end{definition}

\section{Comparison Theorem}

We will use the following well-known result in viscosity solution
theory (see Theorem 8.3 of Crandall, Ishii and Lions \cite{CIL}).

\begin{theorem}
\label{Thm-8.3} Let $u_{i}\in$USC$((0,T)\times \mathbb{R}^{N_{i}})$
for
$i=1,\cdots,k$. Let $\varphi$ be a function defined on $(0,T)\times \mathbb{R}%
^{N_{1}+\cdots+N_{k}}$ such that $(t,x_{1},\ldots,x_{k})\rightarrow
\varphi(t,x_{1},\ldots,x_{k})$ is once continuously differentiable
in $t$ and
twice continuously differentiable in $(x_{1},\cdots,x_{k})\in \mathbb{R}%
^{N_{1}+\cdots+N_{k}}$. Suppose that $\hat{t}\in(0,T)$, $\hat{x}_{i}%
\in \mathbb{R}^{N_{i}}$ for $i=1,\cdots,k$ and
\begin{align*}
w(t,x_{1},\cdots,x_{k})  &  :=u_{1}(t,x_{1})+\cdots+u_{k}(t,x_{k}%
)-\varphi(t,x_{1},\cdots,x_{k})\\
&  \leq w(\hat{t},\hat{x}_{1},\cdots,\hat{x}_{k})
\end{align*}
for $t\in(0,T)$ and $x_{i}\in \mathbb{R}^{N_{i}}$. Assume, moreover,
that there exists $r>0$ such that for every $M>0$ there exists
constant $C$ such that for
$i=1,\cdots,k$,%
\begin{equation}%
\begin{array}
[c]{ll}
& b_{i}\leq C\text{ whenever \ }(b_{i},q_{i},X_{i})\in \mathcal{P}^{2,+}%
u_{i}(t,x_{i}),\\
& |x_{i}-\hat{x}_{i}|+|t-\hat{t}|\leq r\text{ and }|u_{i}(t,x_{i}%
)|+|q_{i}|+\left \Vert X_{i}\right \Vert \leq M.
\end{array}
\label{eq8.5}%
\end{equation}
Then for each $\varepsilon>0$, there exist $X_{i}\in
\mathbb{S}(N_{i})$ such
that\newline \textup{(i)} $(b_{i},D_{x_{i}}\varphi(\hat{t},\hat{x}_{1}%
,\cdots,\hat{x}_{k}),X_{i})\in \overline{\mathcal{P}}^{2,+}u_{i}(\hat{t}%
,\hat{x}_{i}),\  \ i=1,\cdots,k,$\newline \textup{(ii)}
\[
-(\frac{1}{\varepsilon}+\left \Vert A\right \Vert )I\leq \left[
\begin{array}
[c]{ccc}%
X_{1} & \cdots & 0\\
\vdots & \ddots & \vdots \\
0 & \cdots & X_{k}%
\end{array}
\right]  \leq A+\varepsilon A^{2},
\]
\textup{(iii)} $b_{1}+\cdots+b_{k}=\partial_{t}\varphi(\hat{t},\hat{x}%
_{1},\cdots,\hat{x}_{k}),$ \newline where
$A=D_{x}^{2}\varphi(\hat{t},\hat {x})\in
\mathbb{S}(N_{1}+\cdots+N_{k})$.
\end{theorem}

Observe that the above condition (\ref{eq8.5}) will be guaranteed by
having each $u_{i}$ be a subsolution of a parabolic equation given
in the following two theorems.

In this section we will give comparison theorem for $G$-solutions
with different functions $G$. \newline \textbf{(G)} We assume that
\[
G:[0,T]\times \mathbb{R}^{N}\times \mathbb{R}\times \mathbb{R}^{N}%
\times \mathbb{S(}N)\rightarrow \mathbb{R},\  \  \ i=1,\cdots,k,
\]
are continuous in the following sense: for each $t\in \lbrack0,T)$,
$v\in \mathbb{R}$, $x$, $y$, $p\in \mathbb{R}^{N}$ and $X\in
\mathbb{S}(N)$,
\begin{align*}
&  |G_{i}(t,x,v,p,X)-G_{i}(t,y,v,p,X)|\\
&  \leq
\bar{\omega}(1+(T-t)^{-1}+|x|+|y|+|v|)\omega(|x-y|+|p|\cdot|x-y|),
\end{align*}
where $\omega$, $\bar{\omega}:\mathbb{R}^{+}\rightarrow
\mathbb{R}^{+}$ are given continuous functions with $\omega(0)=0$.

\begin{theorem}
\label{Thm-dom1}(\textbf{Domination Theorem}\index{Domination
Theorem}) We are given constants $\beta_{i}>0$,
$i=1,\cdots,k$. Let $u_{i}\in$\textrm{USC}%
$([0,T]\times \mathbb{R}^{N})$ be subsolutions of
\begin{equation}
\partial_{t}u-G_{i}(t,x,u,Du,D^{2}u)=0,\  \  \  \ i=1,\cdots,k, \label{visPDE0}%
\end{equation}
on $(0,T)\times \mathbb{R}^{N}$ such that $\left(\sum_{i=1}^{k}  \beta_i u_{i}%
(t,x)\right)^{+}\rightarrow0$, uniformly as $|x|\rightarrow \infty$.
 We assume that the functions $\{G_{i}\}_i=1^k $ satisfies assumption (G) and
 the following domination condition holds for :
\begin{equation}
\sum_{i=1}^{k}\beta_{i}G_{i}(t,x,v_{i},p_{i},X_{i})\leq0,\  \ \text{
}
\label{domi-cond}%
\end{equation}
for each $(t,x)\in(0,T)\times \mathbb{R}^{N}$ and $(v_{i},p_{i},X_{i}%
)\in \mathbb{R}\times \mathbb{R}^{N}\times \mathbb{S}(N)$ such that\
$\sum _{i=1}^{k}\beta_{i}v_{i}\geq0,\
\sum_{i=1}^{k}\beta_{i}p_{i}=0,\  \sum
_{i=1}^{k}\beta_{i}X_{i}\leq0.$

Then a domination also holds for the solutions: if the sum of
initial values $\sum_{i=1}^{k}\beta_{i}u_{i}(0,\cdot)$ is a
non-positive function on $\mathbb{R}^{N}$, then
$\sum_{i=1}^{k}\beta_{i}u_{i}(t,\cdot)\leq0$, for all $t>0$.
\end{theorem}

\begin{proof}
\medskip \medskip We first observe that for $\bar{\delta}>0$ and for each
$1\leq i\leq k$, the functions defined by
$\tilde{u}_{i}:=u_{i}-\bar{\delta }/(T-t)$ is a subsolution of
\[
\partial_{t}\tilde{u}_{i}-\tilde{G}_{i}(t,x,\tilde{u}_{i},D\tilde{u}_{i}%
,D^{2}\tilde{u}_{i})\leq-\frac{\bar{\delta}}{(T-t)^{2}},%
\]
where
$\tilde{G}_{i}(t,x,v,p,X):=G_{i}(t,x,v+\bar{\delta}/(T-t),p,X)$. It
is easy to check that the functions $\tilde{G}_{i}$ satisfy the same
conditions as $G_{i}$. Since $\sum_{i=1}^{k}\beta_{i}u_{i}\leq0$
follows from $\sum _{i=2}^{k}\beta_{i}\tilde{u}_{i}\leq0$ in the
limit $\bar{\delta}\downarrow0$, it suffices to prove the theorem
under the additional assumptions:
\begin{equation}%
\begin{array}
[c]{c}%
\partial_{t}u_{i}-G_{i}(t,x,u_{i},Du_{i},D^{2}u_{i})\leq-c,\  \  \text{where}\ c=\bar
{\delta}/T^{2},\\
\  \  \  \  \  \  \  \  \  \  \  \  \  \text{and }\lim_{t\rightarrow
T}u_{i}(t,x)=-\infty \  \text{uniformly on }[0,T)\times
\mathbb{R}^{N}.
\end{array}
\label{ineq-c}%
\end{equation}
\newline To prove the theorem, we assume to the contrary that
\[
\sup_{(t,x)\in \lbrack0,T)\times \mathbb{R}^{N}}\sum_{i=1}^{k}\beta_{i}%
u_{i}(t,x)=m_{0}>0
\]
We will apply Theorem \ref{Thm-8.3} for $x=(x_{1},\cdots,x_{k})\in
\mathbb{R}^{k\times N}$ and%
\[
w(t,x):=\sum_{i=1}^{k}\beta_{i}u_{i}(t,x_{i}),\  \
\varphi(x)=\varphi_{\alpha
}(x):=\frac{\alpha}{2}\sum_{i=1}^{k-1}|x_{i+1}-x_{i}|^{2}.
\]
For each large $\alpha>0$, the maximum of $w-\varphi_{\alpha}$
achieves at some $(t^{\alpha},x^{\alpha})$ inside a compact subset
of $[0,T)\times \mathbb{R}^{k\times N}$. Indeed, since
\[
M_{\alpha}=\sum_{i=1}^{k}\beta_{i}u_{i}(t^{\alpha},x_{i}^{\alpha}%
)-\varphi_{\alpha}(x^{\alpha})\geq m_{0},
\]
we conclude $t^{\alpha}$ must be inside an interval $[0,T_{0}]$,
$T_{0}<T$ and
$x^{\alpha}$ must be inside a compact set $\{x\in \mathbb{R}^{k\times N}%
:\sup_{t\in \lbrack0,T_{0}]}w(t,x)\geq \frac{m_{0}}{2}\}$. We can
check that (see \cite{CIL} Lemma 3.1)
\begin{equation}
\left \{
\begin{array}
[c]{l}%
\text{(i) }\lim_{\alpha \rightarrow
\infty}\varphi_{\alpha}(x^{\alpha
})=0\text{,}\\
\text{(ii)\ }\lim_{\alpha \rightarrow \infty}M_{\alpha}=\lim_{\alpha
\rightarrow
\infty}\beta_{1}u_{1}(t^{\alpha},x_{1}^{\alpha})+\cdots+\beta
_{k}u_{k}(t^{\alpha},x_{k}^{\alpha}))\\
\; \; \; \; \; \; \; \;\;\quad\quad\qquad\;\;\;=\sup_{(t,x)\in \lbrack0,T)\times \mathbb{R}^{N}}%
[\beta_{1}u_{1}(t,x)+\cdots+\beta_{k}u_{k}(t,x)]\\
\  \  \  \  \  \  \  \  \  \  \  \ \qquad\quad\
=[\beta_{1}u_{1}(\hat{t},\hat{x})+\cdots+\beta
_{k}u_{k}(\hat{t},\hat{x})]=m_{0},
\end{array}
\right.  \label{limit}%
\end{equation}
where $(\hat{t},\hat{x})$ is a limit point of
$(t^{\alpha},x^{\alpha})$. Since $u_{i}\in \mathrm{USC}$, for
sufficiently large $\alpha$, we have
\[
\beta_{1}u_{1}(t^{\alpha},x_{1}^{\alpha})+\cdots+\beta_{k}u_{k}(t^{\alpha
},x_{k}^{\alpha})\geq \frac{m_{0}}{2}.
\]
If $\hat{t}=0$, we have $\limsup_{\alpha \rightarrow \infty}\sum_{i=1}^{k}%
\beta_{i}u_{i}(t^{\alpha},x_{i}^{\alpha})=\sum_{i=1}^{k}\beta_{i}u_{i}%
(0,\hat{x})\leq0$. We know that $\hat{t}>0$ and thus $t^{\alpha}$
must be strictly positive for large $\alpha$. It follows from
Theorem \ref{Thm-8.3} that, for each $\varepsilon>0$ there exist
$b_{i}^{\alpha}\in \mathbb{R}$, $X_{i}\in \mathbb{S}(N)$ such that
\begin{equation}
(b_{i}^{\alpha},\beta_{i}^{-1}D_{x_{i}}\varphi(x^{\alpha}),X_{i}%
)\in \mathcal{\bar{P}}^{2,+}u_{i}(t^{\alpha},x_{i}^{\alpha}),\  \  \sum_{i=1}%
^{k}\beta_{i}b_{i}^{\alpha}=0\ \text{for }i=1,\cdots,k, \label{b-p-eqn}%
\end{equation}
and such that
\begin{equation}
-(\frac{1}{\varepsilon}+\left \Vert A\right \Vert )I\leq \left(
\begin{array}
[c]{cccc}%
\beta_{1}X_{1} & \ldots & 0 & 0\\
\vdots & \ddots & \vdots & \vdots \\
0 & \ldots & \beta_{k-1}X_{k-1} & 0\\
0 & \ldots & 0 & \beta_{k}X_{k}%
\end{array}
\right)  \leq A+\varepsilon A^{2}, \label{ine-matrix}%
\end{equation}
where $A=D^{2}\varphi_{\alpha}(x^{\alpha})\in \mathbb{S}(kN)$ is
explicitly
given by%
\[
A=\alpha J_{kN},\; \text{where}\;J_{kN}=\left(
\begin{array}
[c]{ccccc}%
I_{N} & -I_{N} & \cdots & \cdots & 0\\
-I_{N} & 2I_{N} & \ddots &  & \vdots \\
\vdots & \ddots & \ddots & \ddots & \vdots \\
\vdots &  & \ddots & 2I_{N} & -I_{N}\\
0 & \cdots & \cdots & -I_{N} & I_{N}%
\end{array}
\right)  .
\]
The second inequality of (\ref{ine-matrix}) implies
$\sum_{i=1}^{k}\beta _{i}X_{i}\leq0$. Set
\begin{align*}
p_{1}^{\alpha}  &  =\beta_{1}^{-1}D_{x_{1}}\varphi_{\alpha}(x^{\alpha}%
)=\beta_{1}^{-1}\alpha(x_{1}^{\alpha}-x_{2}^{\alpha}),\\
p_{2}^{\alpha}  &  =\beta_{2}^{-1}D_{x_{2}}\varphi_{\alpha}(x^{\alpha}%
)=\beta_{2}^{-1}\alpha(2x_{2}^{\alpha}-x_{1}^{\alpha}-x_{3}^{\alpha}),\\
&  \vdots \\
p_{k-1}^{\alpha}  &
=\beta_{k-1}^{-1}D_{x_{k-1}}\varphi_{\alpha}(x^{\alpha
})=\beta_{k-1}^{-1}\alpha(2x_{k-1}^{\alpha}-x_{k-2}^{\alpha}-x_{k}^{\alpha
}),\\
p_{k}^{\alpha}  &  =\beta_{k}^{-1}D_{x_{k}}\varphi_{\alpha}(x^{\alpha}%
)=\beta_{k}^{-1}\alpha(x_{k}^{\alpha}-x_{k-1}^{\alpha}).
\end{align*}
Thus $\sum_{i=1}^{k}\beta_{i}p_{i}^{\alpha}=0$. From this together
with (\ref{b-p-eqn}) and (\ref{ineq-c}), it follows that
\[
b_{i}^{\alpha}-G_{i}(t^{\alpha},x_{i}^{\alpha},u_{i}(t^{\alpha},x_{i}^{\alpha
}),p_{i}^{\alpha},X_{i})\leq-c,\  \ i=1,\cdots,k.\
\]
By (\ref{limit}) (i), we also have $\lim_{\alpha \rightarrow \infty}%
|p_{i}^{\alpha}|\cdot|x_{i}^{\alpha}-x_{1}^{\alpha}|\rightarrow0$.
This, together with the domination condition (\ref{domi-cond}) of
$G_{i}$, implies
\begin{align*}
-c\sum_{i=1}^{k}\beta_{i}  &  =-\sum_{i=1}^{k}\beta_{i}b_{i}^{\alpha}%
-c\sum_{i=1}^{k}\beta_{i}\geq-\sum_{i=1}^{k}\beta_{i}G_{i}(t^{\alpha}%
,x_{i}^{\alpha},u_{i}(t^{\alpha},x_{i}^{\alpha}),p_{i}^{\alpha},X_{i})\\
&  \geq-\sum_{i=1}^{k}\beta_{i}G_{i}(t^{\alpha},x_{1}^{\alpha},u_{i}%
(t^{\alpha},x_{i}^{\alpha}),p_{i}^{\alpha},X_{i})\\
&  \  \  \ -\sum_{i=1}^{k}\beta_{i}|G_{i}(t^{\alpha},x_{i}^{\alpha}%
,u_{i}(t^{\alpha},x_{i}^{\alpha}),p_{i}^{\alpha},X_{i})-G_{i}(t^{\alpha}%
,x_{1}^{\alpha},u_{i}(t^{\alpha},x_{i}^{\alpha}),p_{i}^{\alpha},X_{i})|\\
&
\geq-\sum_{i=1}^{k}\beta_{i}\bar{\omega}(1+(T-T_{0})^{-1}+|x_{1}^{\alpha
}|+|x_{i}^{\alpha}|+|u_{i}(t^{\alpha},x_{i}^{\alpha})|)\cdot \omega
(|x_{i}^{\alpha}-x_{1}^{\alpha}|\\
&
\  \  \  \  \  \  \  \  \  \  \  \  \  \  \  \  \  \  \  \  \  \  \  \  \  \  \  \  \  \  \  \  \  \  \  \  \  \  \  \  \  \  \  \  \  \  \  \  \  \  \  \  \  \  \  \  \  \  \  \  \  \  \  \  \  \  \  \ +|p_{i}%
^{\alpha}|\cdot|x_{i}^{\alpha}-x_{1}^{\alpha}|).
\end{align*}
The right side tends to zero as $\alpha \rightarrow \infty$, which
induces a contradiction. The proof is complete.
\end{proof}


\begin{theorem}\label{Thm-dom}
We assume that the functions $G_{i}=G_{i}(t,x,v,p,X)$, $i=1,\cdots ,k$ and $%
G=G(t,x,v,p,X)$ satisfy assumption \textbf{(G)}. We also assume $G$
dominates $\{G_{i}\}_{i=1}^{k}$ in the following sense: for each
$(t,x)\in
\lbrack 0,\infty )\times \mathbb{R}^{N}$ and $(v_{i},p_{i},X_{i})\in \mathbb{%
R\times R}^{N}\times \mathbb{S}(N)$,
\begin{equation}
\sum_{i=1}^{k}G_{i}(t,x,v_{i},p_{i},X_{i})\leq
G(t,x,\sum_{i=1}^{k}v_{i},\sum_{i=1}^{k}p_{i},\sum_{i=1}^{k}X_{i}).\text{
}\ \   \label{App-DominationCon}
\end{equation}%
Moreover there exists a constant $\bar{C}$ such that, for each
$(t,r,x,p)\in
\lbrack 0,\infty )\times \mathbb{R}^{N}\mathbb{\times R}^{N}$ and, $%
Y_{1},Y_{2}\in \mathbb{S}(N)$ such that $Y_{2}\geq Y_{1}$
\begin{equation*}
G(t,x,r,p,Y_{2})-G(t,x,r_{1},p,Y_{1})\geq 0
\end{equation*}%
and there exists a constant $\bar{C}$ such that
\begin{equation*}
|G(t,x,r,p,X_{1})-G(t,x,r,p,X_{2})|\leq
\bar{C}(|r_{1}-r_{2}|+|X_{1}-X_{2}|).
\end{equation*}%
Let $u_{i}\in $\textrm{USC}$([0,T]\times \mathbb{R}^{N})$ be a $G_{i}$%
-subsolution and $u\in $\textrm{LSC}$([0,T]\times \mathbb{R}^{N})$ be a $G$%
-supersolution such that $u_{i}$ and $u$ satisfy polynomial growth
condition. Then $\sum_{i=1}^{k}u_{i}(t,x)\leq u(t,x)$ on
$[0,T)\times \mathbb{R}^{N}$ provided that
$\sum_{i=1}^{k}u_{i}|_{t=0}\leq u|_{t=0}$.
\end{theorem}

\begin{proof}
For a fixed and large constant $\lambda >\bar{C}+C$ we set $\xi
(x):=(1+|x|^{2})^{l/2}$ and
\begin{equation*}
\tilde{u}_{i}(t,x):=u_{i}(t,x)\xi (x)e^{-\lambda t},\ i=1,\cdots
,k,\ \ \tilde{u}_{k+1}(t,x):=-u(t,x)\xi (x)e^{-\lambda t},
\end{equation*}%
where $l$ is chosen to be large enough such that $\sum |\tilde{u}%
_{i}(t,x)|\rightarrow 0$ uniformly as $|x|\rightarrow \infty $. It
is easy to check that, for each $i=1,\cdots ,k+1$, $\tilde{u}_{i}$
is a subsolution of
\begin{equation*}
\partial _{t}\tilde{u}_{i}+\lambda \tilde{u}_{i}-\tilde{G}_{i}(t,x,\tilde{u}%
_{i},D\tilde{u}_{i},D^{2}\tilde{u}_{i})=0,
\end{equation*}%
where, for each $i=1,\cdots ,k$, the function
$\tilde{G}_{i}(t,x,v,p,X)$ is given by
\begin{equation*}
e^{-\lambda t}\xi ^{-1}G_{i}(t,x,e^{\lambda t}v\xi ,e^{\lambda
t}(p+vD\xi ),e^{\lambda t}(X\xi +p\otimes D\xi +D\xi \otimes
p+vD^{2}\xi )),\
\end{equation*}%
and $\tilde{G}_{k+1}(t,x,v,p,X)$ is given by%
\begin{equation*}
-e^{-\lambda t}\xi ^{-1}G(t,x,-e^{\lambda t}v\xi ,-e^{\lambda
t}(p+vD\xi ),-e^{\lambda t}(X\xi +p\otimes D\xi +D\xi \otimes
p+vD^{2}\xi )).
\end{equation*}%
Observe that
\begin{align*}
D\xi (x)& =l\xi (x)(1+|x|^{2})^{-1}x,\ \ \ \  \\
D^{2}\xi (x)& =\xi
(x)[l(1+|x|^{2})^{-1}I+l(l-2)(1+|x|^{2})^{-2}x\otimes x].
\end{align*}%
Thus both $\xi ^{-1}(x)|D\xi (x)|$ and $\xi ^{-1}(x)|D^{2}\xi (x)|$
converges to zero uniformly as $|x|\rightarrow \infty $.

From the domination condition (\ref{App-DominationCon}), for each $%
(v_{i},p_{i},X_{i})\in \mathbb{R\times R}^{N}\times \mathbb{S}(N)$, $%
i=1,\cdots ,k+1$, such that $\sum_{i=1}^{k+1}v_{i}=0$, $%
\sum_{i=1}^{k+1}p_{i}=0$, and $\sum_{i=1}^{k+1}X_{i}=0$, we have
\begin{equation*}
\sum_{i=1}^{k+1}\tilde{G}_{i}(t,x,v_{i},p_{i},X_{i})\leq 0.
\end{equation*}%
For $v$, $r\in \mathbb{R}$, $p\in \mathbb{R}^{N}$ and $X$, $Y\in \mathbb{S}%
(N)$ such that $r\geq 0$, $Y\geq 0$ and $r>0$, since $\tilde{G}$ is
still
monotone in $X$,%
\begin{eqnarray*}
&&\tilde{G}_{k+1}(t,x,v,p,X)-\tilde{G}_{k+1}(t,x,v-r,p,X+Y) \\
&\leq &\tilde{G}_{k+1}(t,x,v,p,X)-\tilde{G}_{k+1}(t,x,v-r,p,X) \\
&\leq &c(\hat{C}+C_{1})r,\text{ }
\end{eqnarray*}%
where the constant $C_{1}$ does not depend on $(t,x,v,p,X)$. We then
apply
the above theorem by choosing $\beta _{i}=1$, $i=1,\cdots ,k+1$. Thus $%
\sum_{i=1}^{k+1}\tilde{u}_{i}|_{t=0}\leq 0$. Moreover for each
$v_{i}\in \mathbb{R}$, $p_{i}\in \mathbb{R}^{N}$ and $X_{i}\in
\mathbb{S}(N)$ such
that $\hat{v}=\sum_{i=1}^{k+1}v_{i}\geq 0$, $\sum_{i=1}^{k+1}p_{i}=0$ and $%
\hat{X}=\sum_{i=1}^{k+1}X_{i}\leq 0$ we have
\begin{eqnarray*}
&&-\lambda \sum_{i=1}^{k+1}v_{i}+\sum_{i=1}^{k+1}\tilde{G}%
_{i}(t,x,v_{i},p_{i},X_{i}) \\
&=&-\lambda \hat{v}+\sum_{i=1}^{k}\tilde{G}_{i}(t,x,v_{i},p_{i},X_{i})+%
\tilde{G}_{k+1}(t,x,v_{k+1}-\hat{v},p_{k+1},X_{k+1}-\hat{X}) \\
&&+\tilde{G}_{k+1}(t,x,v_{k+1},p_{k+1},X_{k+1})-\tilde{G}_{k+1}(t,x,v_{k+1}-%
\hat{v},p_{k+1},X_{k+1}-\hat{X}) \\
&\leq &-\lambda \hat{v}+\tilde{G}_{k+1}(t,x,v_{k+1},p_{k+1},X_{k+1})-\tilde{G%
}_{k+1}(t,x,v_{k+1}-\hat{v},p_{k+1},X_{k+1}) \\
&\leq &-\lambda \hat{v}+c(\bar{C}+C)\hat{v}\leq 0
\end{eqnarray*}%
It follows that all conditions in Theorem \ref{Thm-dom1} are
satisfied. Thus
we have $\sum_{i=1}^{k+1}\tilde{u}_{i}\leq 0$, or equivalently, $%
\sum_{i=1}^{k+1}u_{i}(t,x)\leq u(t,x)$ for $(t,x)\in \lbrack
0,T)\times \mathbb{R}^{N}$.
\end{proof}

The following comparison theorem is a direct consequence of the
above domination theorem.

\begin{theorem}(Comparison Theorem) \label{Comparison}
We are given two functions $G=G(t,x,v,p,X)$ and
$G_{1}=G_{1}(t,x,v,p,X)$ satisfying condition (G). We also assume
that, for each $(t,x,v,p,X)\in
\lbrack 0,\infty )\times \mathbb{R}^{N}\times \mathbb{R\times R}^{N}$ and $%
Y\in \mathbb{S}(N)$ such that $X\geq Y$,
\begin{eqnarray}
G(t,x,v,p,X) &\geq &G_{1}(t,x,v,p,X),  \label{App-Conparison} \\
G(t,x,v,p,Y) &\geq &G(t,x,v,p,X).\ \   \label{GYGX}
\end{eqnarray}%
We also assume that $G$ is a uniform Lipschitz function in $v$ and
$X$,
namely, for each $(t,x)\in \lbrack 0,\infty )\times \mathbb{R}^{N}$ and $%
(v,p,X)$, $(v^{\prime },p,X^{\prime })\in \mathbb{R\times
R}^{N}\times \mathbb{S}(N)$,
\begin{equation*}
|G(t,x,v,p,X)-G(t,x,v^{\prime },p^{\prime },X^{\prime })|\leq \bar{C}%
(|v-v^{\prime }|+|X-X^{\prime }|).
\end{equation*}%
Let $u_{1}\in $\textrm{USC}$([0,T]\times \mathbb{R}^{N})$ be a $G_{1}$%
-subsolution and $u$ $\in $\textrm{LSC}$([0,T]\times
\mathbb{R}^{N})$ be a $G $-supersolution on $(0,T)\times
\mathbb{R}^{N}$ satisfying the polynomial growth condition. Then
$u\geq u_{1}$ on $[0,T)\times \mathbb{R}^{N}$ provided that
$u|_{t=0}\geq u_{1}|_{t=0}$. In particular this comparison
holds for the case where $G\equiv G_{1}$, which is a Lipschitz function in $%
(v,X)$ and satisfies the elliptic condition (\ref{GYGX}).
\end{theorem}

The following special case of the above domination theorem is also
very useful.

\begin{theorem}(Domination Theorem)\label{Domination}
We assume that $G_{1}$ and $G$ satisfy the same conditions given in
the previous theorem except that the condition
(\ref{App-Conparison}) is replaced by the following one: for each
$(t,x)\in \lbrack 0,\infty )\times \mathbb{R}^{N}$ and $(v,p,X)$,
$(v^{\prime },p^{\prime },X^{\prime })\in \mathbb{R\times
R}^{N}\times \mathbb{S}(N)$,
\begin{equation*}
G_{1}(t,x,v,p,X)-G_{1}(t,x,v^{\prime },p^{\prime },X^{\prime })\geq
G(t,x,v-v^{\prime },p-p^{\prime },X-X^{\prime })
\end{equation*}%
Let $u\in $\textrm{USC}$([0,T]\times \mathbb{R}^{N})$ be a $G_{1}$%
-subsolution and $v$ $\in $\textrm{LSC}$([0,T]\times \mathbb{R}^{N})$ be a $%
G_{1}$-supersolution on $(0,T)\times \mathbb{R}^{N}$ and $w$ is a $G$%
-supersolution. They satisfy the polynomial growth condition. If $%
(u-v)|_{t=0}=w|_{t=0}$ then $u-v\leq w$ on $[0,T)\times
\mathbb{R}^{N}$.
\end{theorem}

The following theorem will be frequently used in this book. Let
$G:\mathbb{R}^{N}\times \mathbb{S}(N)\rightarrow \mathbb{R}$ be a
given continuous sublinear function monotonic in $A\in
\mathbb{S}(N)$. Obviously, $G$ satisfies conditions (G) of Theorem
\ref{Thm-dom}. We consider the following $G$-equation:
\begin{equation}
\partial_{t}u-G(Du,D^{2}u)=0,\  \ u(0,x)=\varphi(x).\label{eq-G-heat}%
\end{equation}

\begin{theorem}
\label{Com-G}Let $G:\mathbb{R}^{N}\times \mathbb{S}(N)\rightarrow
\mathbb{R}$ be a given continuous sublinear function monotonic in
$A\in \mathbb{S}(N)$. Then we have

\begin{description}
\item[\textup{(i)}] If $u\in \mathrm{USC}([0,T]\times \mathbb{R}^{N})$ with polynomial
growth is a viscosity subsolution of \textup{(\ref{eq-G-heat})} and
$v\in \mathrm{LSC}([0,T]\times \mathbb{R}^{N})$ with polynomial
growth is a viscosity supersolution of \textup{(\ref{eq-G-heat})},
then $u\leq v$.

\item[\textup{(ii)}] If\ $u^{\varphi}\in C([0,T]\times \mathbb{R}^{N})$ denotes the
polynomial growth solution of \textup{(\ref{eq-G-heat})} with
initial condition $\varphi$, then $u^{\lambda \varphi}=\lambda
u^{\varphi}$ for each $\lambda \geq0$ and $u^{\varphi+\psi}\leq
u^{\varphi}+u^{\psi}$.

\item[\textup{(iii)}] If a given function
$\widetilde{G}:\mathbb{R}^{N}\times
\mathbb{S}(N)\mapsto \mathbb{R}$ is dominated by $G$, i.e.%
\begin{equation*}
\widetilde{G}(p,X)-\widetilde{G}(p^{\prime },X^{\prime })\leq
G(p-p^{\prime
},X-X^{\prime }),\ \ \text{for }p,p^{\prime }\in \mathbb{R}\text{, }%
X,X^{\prime }\in \mathbb{S}(N),
\end{equation*}%
then for each $\varphi \in C(\mathbb{R}^{N})$ satisfying polynomial
growth
condition, there exists a unique $\widetilde{G}$-solution $\tilde{u}%
^{\varphi }(t,x)$ on $[0,\infty )\times \mathbb{R}^{N}$ with
initial condition $\tilde{u}^{\varphi }|_{t=0}=\varphi $ (see the
next section for
the proof of existence), i.e.,%
\begin{equation*}
\partial _{t}\tilde{u}^{\varphi }-\widetilde{G}(D\tilde{u}^{\varphi },D^{2}%
\tilde{u}^{\varphi })=0,\ \ \text{ }\tilde{u}^{\varphi
}|_{t=0}=\varphi .
\end{equation*}%
Moreover%
\begin{equation*}
\tilde{u}^{\varphi }(t,x)-\tilde{u}^{\psi }(t,x)\leq u^{\varphi
-\psi }(t,x),\ \ \text{for }t\geq 0,\ x\in \mathbb{R}^{N}\text{.}
\end{equation*}
Consequently, the following comparison holds: $\psi \geq \varphi$
implies $\tilde{u}^{\psi}\geq \tilde{u}^{\varphi}$.
\end{description}
\end{theorem}

\begin{proof}
By the above corollaries, it is easy to obtain the results.
\end{proof}


\section{Perron's Method and Existence}

The combination of Perron's method and viscosity solutions was
introduced by H. Ishii \cite{Ishii}. For the convenience of readers,
we interpret the proof provided in Crandall, Ishii and Lions
\cite{CIL} into its parabolic situation.

We consider the following parabolic PDE:%
\begin{equation}
\left \{
\begin{array}
[c]{l}%
\partial_{t}u-G(t,x,u,Du,D^{2}u)=0\  \text{on}\ (0,\infty)\times \mathbb{R}%
^{N},\\
u(0,x)=\varphi(x)\  \text{for}\ x\in \mathbb{R}^{N},
\end{array}
\right.  \label{Cauchy-PDE}%
\end{equation}
where $G:[0,\infty)\times \mathbb{R}^{N}\times \mathbb{R}\times \mathbb{R}%
^{N}\times \mathbb{S}(N)\rightarrow \mathbb{R}$, $\varphi \in
C(\mathbb{R}^{N})$.

To discuss Perron's method, we will use the following notations: if
$u:\mathcal{O}\rightarrow \lbrack-\infty,\infty]$ where $\mathcal{O}%
\subset \lbrack0,\infty)\times \mathbb{R}^{N}$, then%
\begin{equation}
\left \{
\begin{array}
[c]{c}%
u^{\ast}(t,x)=\lim_{r\downarrow0}\sup \{u(s,y):(s,y)\in \mathcal{O}%
\  \text{and}\  \sqrt{|s-t|+|y-x|^{2}}\leq r\},\\
u_{\ast}(t,x)=\lim_{r\downarrow0}\inf \{u(s,y):(s,y)\in \mathcal{O}%
\  \text{and}\  \sqrt{|s-t|+|y-x|^{2}}\leq r\}.
\end{array}
\right.  \label{envelope}%
\end{equation}
One calls $u^{\ast}$ the \textbf{upper semicontinuous
envelope}\index{Upper semicontinuous envelope} of $u$; it is the
smallest upper semicontinuous function (with values in
$[-\infty,\infty]$) satisfying $u\leq u^{\ast}$. Similarly,
$u_{\ast}$ is the \textbf{lower semicontinuous envelope}\index{Lower
semicontinuous envelope} of $u$.

\begin{theorem}
\label{Perron}(Perron's Method) Let comparison hold for
\textup{(\ref{Cauchy-PDE})}, i.e., if $w$ is a viscosity subsolution
of \textup{(\ref{Cauchy-PDE})} and $v$ is a viscosity supersolution
of \textup{(\ref{Cauchy-PDE})}, then $w\leq v$. Suppose also that
there is a viscosity subsolution $\underline{u}$ and a viscosity
supersolution $\bar{u}$ of \textup{(\ref{Cauchy-PDE})} that satisfy
the condition
$\underline{u}_{\ast}(0,x)=\bar{u}^{\ast}(0,x)=\varphi(x)$ for $x\in
\mathbb{R}^{N}$. Then
\[
W(t,x)=\sup \{w(t,x):\underline{u}\leq w\leq \bar{u}\  \text{and}\
w\  \text{is a viscosity subsolution of
\textup{(\ref{Cauchy-PDE})}}\}
\]
is a viscosity solution of \textup{(\ref{Cauchy-PDE})}.
\end{theorem}

The proof consists of two lemmas. For the proof of the following two
lemmas, we also see \cite{AT}. The first one is

\begin{lemma}
\label{one-step}Let $\mathcal{F}$ be a family of viscosity
subsolution of \textup{(\ref{Cauchy-PDE})} on $(0,\infty )\times
\mathbb{R}^{N}$. Let $w(t,x)=\sup
\{u(t,x):u\in \mathcal{F}\}$ and assume that $w^{\ast }(t,x)<\infty $ for $%
(t,x)\in (0,\infty )\times \mathbb{R}^{N}$. Then $w^{\ast }$ is a
viscosity subsolution of \textup{(\ref{Cauchy-PDE})} on $(0,\infty
)\times \mathbb{R}^{N}$.
\end{lemma}

\begin{proof}
Let $(t,x)\in (0,\infty )\times \mathbb{R}^{N}$ and consider a sequence $%
s_{n}$, $y_{n}$, $u_{n}\in \mathcal{F}$ such that
$\lim_{n\rightarrow \infty
}(s_{n},y_{n},u_{n}(s_{n},y_{n}))=(t,x,w^{\ast }(t,x))$. There
exists $r>0$
such that $N_{r}=\{(s,y)\in (0,\infty )\times \mathbb{R}^{N}:\sqrt{%
|s-t|+|y-x|^{2}}\leq r\}$ is compact. For $\phi \in C^{2}$ such that
$\phi (t,x)=w^{\ast }(t,x)$ and $w^{\ast }<\phi $ on $(0,\infty
)\times
\mathbb{R}^{N}\backslash (t,x)$, let $(t_{n},x_{n})$ be a maximum point of $%
u_{n}-\phi $ over $N_{r}$, hence $u_{n}(s,y)\leq
u_{n}(t_{n},x_{n})+\phi (s,y)-\phi (t_{n},x_{n})$ for $(s,y)\in
N_{r}$. Suppose that (passing to a
subsequence if necessary) $(t_{n},x_{n})\rightarrow (\bar{t},\bar{x})$ as $%
n\rightarrow \infty $. Putting $(s,y)=(s_{n},y_{n})$ in the above
inequality and taking the limit inferior as $n\rightarrow \infty $,
we obtain
\begin{align*}
w^{\ast }(t,x)& \leq \liminf_{n\rightarrow \infty
}u_{n}(t_{n},x_{n})+\phi
(t,x)-\phi (\bar{t},\bar{x}) \\
& \leq w^{\ast }(\bar{t},\bar{x})+\phi (t,x)-\phi (\bar{t},\bar{x}).
\end{align*}%
From the above inequalities and the assumption on $\phi $, we get $%
\lim_{n\rightarrow \infty
}(t_{n},x_{n},$$u_{n}(t_{n},x_{n}))=(t,x,w^{\ast }(t,x))$ (without
passing to a subsequence). Since $u_{n}$ is a viscosity
subsolution of (\ref{Cauchy-PDE}), by definition we have%
\begin{equation*}
\partial _{t}\phi (t_{n},x_{n})-G(t_{n},x_{n},u_{n}(t_{n},x_{n}),D\phi
(t_{n},x_{n}),D^{2}\phi (t_{n},x_{n}))\leq 0.
\end{equation*}%
Letting $n\rightarrow \infty $, we conclude that%
\begin{equation*}
\partial _{t}\phi (t,x)-G(t,x,w^{\ast }(t,x),D\phi (t,x),D^{2}\phi
(t,x))\leq 0.
\end{equation*}%
Thus $w^{\ast }$ is a viscosity subsolution of (\ref{Cauchy-PDE}) by
definition.
\end{proof}

The second step in the proof of Theorem \ref{Perron} is a simple
\textquotedblleft bump\textquotedblright \ construction that we now
describe.
Suppose that $u$ is a viscosity subsolution of (\ref{Cauchy-PDE}) on $%
(0,\infty )\times \mathbb{R}^{N}$ and $u_{\ast }$ is not a viscosity
supersolution of (\ref{Cauchy-PDE}), so that there is $(t,x)\in
(0,\infty
)\times \mathbb{R}^{N}$ and $\phi \in C^{2}$ with $u_{\ast }(t,x)=\phi (t,x)$%
, $u_{\ast }>\phi $ on $(0,\infty )\times \mathbb{R}^{N}\backslash (t,x)$ and%
\begin{equation*}
\partial _{t}\phi (t,x)-G(t,x,\phi (t,x),D\phi (t,x),D^{2}\phi (t,x))<0.
\end{equation*}%
The continuity of $G$ provides $r,\delta _{1}>0$ such that $N_{r}=\{(s,y):%
\sqrt{|s-t|+|y-x|^{2}}\leq r\}$ is compact and%
\begin{equation*}
\partial _{t}\phi -G(s,y,\phi +\delta ,D\phi ,D^{2}\phi )\leq 0
\end{equation*}%
for all $s,y,\delta \in N_{r}\times \lbrack 0,\delta _{1}]$. Lastly,
we obtain $\delta _{2}>0$ for which $u_{\ast }>\phi +\delta _{2}$ on
$\partial N_{r}$. Setting $\delta _{0}=\min (\delta _{1},\delta
_{2})>0,$ we define
\begin{equation*}
U=\left\{\begin{array} {ll}\max (u,\phi +\delta _{0})\ \ &\text{on}\
N_{r}\\
u\ \ &\text{elsewhere.}
\end{array}\right.
\end{equation*}
By the above inequalities and Lemma \ref{one-step}, it is easy to
check that $U$ is a viscosity subsolution of (\ref{Cauchy-PDE}) on
$(0,\infty )\times \mathbb{R}^{N}.$ Obviously, $U\geq u$. Finally,
observe that $U_{\ast }(t,x)\geq \max (u_{\ast }(t,x),\phi
(t,x)+\delta _{0})>u_{\ast }(t,x);$ hence there exists $(s,y)$ such
that $U(s,y)>u(s,y).$ We summarize the above discussion as the
following lemma.

\begin{lemma}
\label{two-step}Let $u$ be a viscosity subsolution of \textup{(\ref{Cauchy-PDE})} on $%
(0,\infty )\times \mathbb{R}^{N}$. If $u_{\ast }$ fails to be a
viscosity supersolution at some point $(s,z)$, then for any small
$\kappa >0$ there is a viscosity subsolution $U_{\kappa }$ of
\textup{(\ref{Cauchy-PDE})} on $(0,\infty
)\times \mathbb{R}^{N}$ satisfying%
\begin{equation*}
\left \{
\begin{array}{l}
U_{\kappa }(t,x)\geq u(t,x)\  \text{and}\  \sup (U_{\kappa }-u)>0, \\
U_{\kappa }(t,x)=u(t,x)\  \text{for}\  \sqrt{|t-s|+|x-z|^{2}}\geq \kappa .%
\end{array}%
\right.
\end{equation*}
\end{lemma}

\noindent\textbf{Proof of Theorem \ref{Perron}.} With the notation
of the theorem
observe that $\underline{u}_{\ast}\leq W_{\ast}\leq W\leq W^{\ast}\leq \bar {%
u}^{\ast}$ and, in particular,
$W_{\ast}(0,x)=W(0,x)=W^{\ast}(0,x)=\varphi (x)$ for $x\in
\mathbb{R}^{N}$. By lemma \ref{one-step}, $W^{\ast}$ is a
viscosity subsolution of (\ref{Cauchy-PDE}) and hence, by comparison, $%
W^{\ast}\leq \bar{u}$. It then follows from the definition of $W$ that $%
W=W^{\ast}$ (so $W$ is a viscosity subsolution). If $W_{\ast}$ fails
to be a
viscosity supersolution at some point $(s,z)\in(0,\infty)\times \mathbb{R}%
^{N}$, let $W_{\kappa}$ be provided by Lemma \ref{two-step}. Clearly $%
\underline{u}\leq W_{\kappa}$ and $W_{\kappa}(0,x)=\varphi(x)$ for
sufficiently small $\kappa$. By comparison, $W_{\kappa}\leq \bar{u}$
and
since $W$ is the maximal viscosity subsolution between $\underline{u}$ and $%
\bar{u}$, we arrive at the contradiction $W_{\kappa}\leq W$. Hence
$W_{\ast}$ is a viscosity supersolution of (\ref{Cauchy-PDE}) and
then, by comparison for (\ref{Cauchy-PDE}), $W^{\ast}=W\leq
W_{\ast}$, showing that $W$ is continuous and is a viscosity
solution of (\ref{Cauchy-PDE}). The proof is complete. \hfill$\Box$

Let $G:\mathbb{R}^{N}\times \mathbb{S}(N)\rightarrow \mathbb{R}$ be
a given continuous sublinear function monotonic in $A\in
\mathbb{S}(N)$. We now consider the existence of viscosity solution
of the following $G$-equation:
\begin{equation}
\partial_{t}u-G(Du,D^{2}u)=0,\  \ u(0,x)=\varphi(x).  \label{G-eq-ap}
\end{equation}

Case 1: If $\varphi \in C_{b}^{2}(\mathbb{R}^{N})$, then $\underline {u}%
(t,x)=\underline{M}t+\varphi(x)$ and
$\bar{u}(t,x)=\bar{M}t+\varphi(x)$ are respectively the classical
subsolution and supersolution of (\ref{G-eq-ap}),
where $\underline{M}=\inf_{x\in \mathbb{R}^{N}}G(D\varphi(x),D^{2}%
\varphi(x)) $ and $\bar{M}=\sup_{x\in \mathbb{R}^{N}}G(D\varphi(x),D^{2}%
\varphi(x))$. Obviously, $\underline{u}$ and $\bar{u}$ satisfy all
the conditions in Theorem \ref{Perron}. By Theorem \ref{Com-G}, we
know the comparison holds for (\ref{G-eq-ap}). Thus by Theorem
\ref{Perron}, we obtain that $G$-equation (\ref{G-eq-ap}) has a
viscosity solution.

Case 2: If $\varphi \in C_{b}(\mathbb{R}^{N})$ with
$\lim_{|x|\rightarrow \infty }\varphi (x)=0$, then we can choose a
sequence $\varphi _{n}\in
C_{b}^{2}(\mathbb{R}^{N})$ which uniformly converge to $\varphi $. For $%
\varphi _{n}$, by Case 1, there exists a viscosity solution
$u^{\varphi _{n}} $. By comparison theorem, it is easy to show that
$u^{\varphi _{n}}$ is uniformly
convergent, the limit denoted by $u$. Similar to the proof of Lemma \ref%
{one-step}, it is easy to prove that $u$ is a viscosity solution of $G$%
-equation (\ref{G-eq-ap}) with initial condition $\varphi $.

Case 3: If $\varphi \in C(\mathbb{R}^{N})$ with polynomial growth,
then
we can choose a large $l>0$ such that $\tilde{\varphi}(x)=\varphi(x)%
\xi^{-1}(x)$ satisfies the condition in Case 2, where $%
\xi(x)=(1+|x|^{2})^{l/2}$. It is easy to check that $u$ is a
viscosity
solution of $G$-equation (\ref{G-eq-ap}) if and only if $\tilde{u}%
(t,x)=u(t,x)\xi^{-1}(x)$ is a viscosity solution of the following
PDE:
\begin{equation}
\partial_{t}\tilde{u}-\tilde{G}(x,\tilde{u},D\tilde{u},D^{2}\tilde {u})=0,\
\tilde{u}(0,x)=\tilde{\varphi},  \label{G-trans}
\end{equation}
where $\tilde{G}(x,v,p,X)=G(p+v\eta(x),X+p\otimes
\eta(x)+\eta(x)\otimes p+v\kappa(x))$. Here
\begin{align*}
\eta(x) & :=\xi^{-1}(x)D\xi(x)=l(1+|x|^{2})^{-1}x,\  \  \  \  \\
\kappa(x) &
:=\xi^{-1}(x)D^{2}\xi(x)=l(1+|x|^{2})^{-1}I+l(l-2)(1+|x|^{2})^{-2}x\otimes
x.
\end{align*}
Similar to the above discussion, we obtain that there exists a
viscosity solution of (\ref{G-trans}) with initial condition
$\tilde{\varphi}$. Thus there exists a viscosity solution of
$G$-equation (\ref{G-eq-ap}).

We summarize the above discussions as a theorem.

\begin{theorem}
Let $\varphi \in C(\mathbb{R}^{N})$ with polynomial growth. Then
there exists a viscosity solution of $G$-equation
\textup{(\ref{G-eq-ap})} with initial condition $\varphi$.
\end{theorem}

\begin{theorem}
\label{G-Tilde}Let the function $G$ be given as in the previous
theorem and
let $\tilde{G}(t,x,p,X):[0,\infty ]\times \mathbb{R}^{N}\times \mathbb{R}%
^{N}\times \mathbb{S}(N)\mapsto \mathbb{R}$ be a given function
satisfying condition (G) in Theorem \ref{Thm-dom}. We assume that
$\tilde{G}$ is
dominated by $G$ in the sense of%
\begin{eqnarray*}
\tilde{G}(t,x,p,X)-\tilde{G}(t,x,p^{\prime },X^{\prime }) &\leq
&G(p-p^{\prime },X-X^{\prime }), \\
&&\text{for each }t,x\text{, }p,p^{\prime }\text{ and }X\text{,
}X^{\prime }.
\end{eqnarray*}%
Then for each given $\varphi \in C(\mathbb{R}^{N})$ with polynomial
growth
condition the viscosity solution of $\partial _{t}\tilde{u}-\tilde{G}(t,x,D%
\tilde{u},D^{2}\tilde{u})=0$ with $\tilde{u}|_{t=0}=\varphi $ exists
and is unique. Moreover the comparison property also holds.
\end{theorem}

\begin{proof}
It is easy to check that a $G$-solution with $u|_{t=0}=\varphi $ is a $%
\tilde{G}$-supersolution. Similarly denoting $G_{\ast }(p,X):=-G(-p,-X)$, a $%
G_{\ast }$-solution with $u|_{t=0}=\varphi $ is a
$\tilde{G}$-subsolution.

We now prove that comparison holds for $\tilde{G}$-solutions. Let
$u_{1}$ be
a $\tilde{G}$-supersolution and $u_{2}$-be a $\tilde{G}$-subsolution with $%
u_{1}|_{t=0}=\varphi _{1}\in \mathrm{SC}([0,T)\times \mathbb{R}^{N})$, $%
u_{2}|_{t=0}=\varphi _{2}\in \mathrm{PSC}([0,T)\times
\mathbb{R}^{N})$, then by the above domination theorem we can obtain
$u_{1}-u_{2}\leq \hat{u}$ where $\hat{u}$ is a $G$-supersolution
with $\hat{u}|_{t=0}=\varphi _{1}-\varphi _{2}$. On the other hand
it is easy to prove that, in the case when $\varphi _{1}\leq \varphi
_{2}$,  $\hat{u}=(\varphi _{1}(x)-\varphi
_{2}(x))\mathbf{1}_{\{0\}}(t)$ such type of $G$-supersolution.
Consequently
we have $u_{1}\leq u_{2}$. This implies that the comparison holds for $%
\tilde{G}$-equations. We then can apply Theorem \ref{Perron} to prove that $%
\tilde{G}$-solution $u$ with $u|_{t=0}=\varphi $ exists.
\end{proof}

\section{Krylov's Regularity Estimate for Parabolic PDE }

The proof of our new central limit theorem is based on a powerful
$C^{1+\alpha/2,2+\alpha}$-regularity estimates for fully nonlinear
parabolic PDE obtained in Krylov \cite{Krylov1}. A more recent
result of Wang \cite{WangL} (the version for elliptic PDE was
initially introduced in Cabre and Caffarelli \cite{Caff1997}), using
viscosity solution arguments, can also be applied.

For simplicity, we only consider the following type of PDE:%
\begin{equation}
\label{eq-Kr}\partial_{t}u+G(D^{2}u,Du,u)=0,\  \ u(T,x)=\varphi(x),
\end{equation}
where $G:\mathbb{S}(d)\times \mathbb{R}^{d}\times
\mathbb{R}\rightarrow \mathbb{R}$ is a given function and $\varphi
\in C_{b}(\mathbb{R}^{d})$.

Following Krylov {\cite{Krylov1}}, we fix constants $K\geq
\varepsilon>0$, $T>0$ and set $Q=(0,T)\times \mathbb{R}^{d}$. Now we
give the definition of $\mathcal{G}(\varepsilon,K,Q) $ and
$\bar{\mathcal{G}}(\varepsilon,K,Q)$.

The following definition is according to Definition 5.5.1 in Krylov
 {\cite{Krylov1}.}

\begin{definition}
Let $G:\mathbb{S}(d)\times \mathbb{R}^{d}\times
\mathbb{R}\rightarrow \mathbb{R}$ be given, written it as
$G(u_{ij},u_{i},u)$, $i,j=1,\ldots,d$. We denote $G\in
\mathcal{G}(\varepsilon,K,Q)$ if $G$ is twice continuously
differentiable with respect to $(u_{ij},u_{i},u)$ and, for each
real-valued $u_{ij}=u_{ji}$, $\tilde{u}_{ij}=\tilde{u}_{ji}$,
$u_{i}$, $\tilde{u}_{i}$, $u $, $\tilde{u}$ and $\lambda^{i}$, the
following inequalities hold:
\[
\varepsilon|\lambda|^{2}\leq \sum_{i,j}\lambda^{i}\lambda^{j}\partial_{u_{ij}%
}G\leq K|\lambda|^{2},
\]%
\[
|G-\sum_{i,j}u_{ij}\partial_{u_{ij}}G|\leq M_{1}^{G}(u)(1+\sum_{i}|u_{i}%
|^{2}),
\]%
\[
|\partial_{u}G|+(1+\sum_{i}|u_{i}|)\sum_{i}|\partial_{u_{i}}G|\leq M_{1}%
^{G}(u)(1+\sum_{i}|u_{i}|^{2}+\sum_{i,j}|u_{ij}|),
\]%
\begin{align*}
\lbrack M_{2}^{G}(u,u_{k})]^{-1}G_{(\eta)(\eta)}  & \leq
\sum_{i,j}|\tilde
{u}_{ij}|\Big[\sum_{i}|\tilde{u}_{i}|+(1+\sum_{i,j}|u_{ij}|)|\tilde{u}|\Big]\\
& \ \ \ +\sum_{i}|\tilde{u}_{i}|^{2}(1+\sum_{i,j}|u_{ij}|)+(1+\sum_{i,j}|u_{ij}%
|^{3})|\tilde{u}|^{2},
\end{align*}
where the arguments $(u_{ij},u_{i},u)$ of $G$ and its derivatives
are omitted,
$\eta=(\tilde{u}_{ij},\tilde{u}_{i},\tilde{u})$, and%
\begin{align*}
G_{(\eta)(\eta)}:=  & \sum_{i,j,r,s}\tilde{u}_{ij}\tilde{u}_{rs}%
\partial_{u_{ij}u_{rs}}^{2}G+2\sum_{i,j,r}\tilde{u}_{ij}\tilde{u}_{r}%
\partial_{u_{ij}u_{r}}^{2}G+2\sum_{i,j}\tilde{u}_{ij}\tilde{u}\partial
_{u_{ij}u}^{2}G\\
&+\sum_{i,j}\tilde{u}_{i}\tilde{u}_{j}\partial_{u_{i}u_{j}}^{2}G+2\sum
_{i}\tilde{u}_{i}\tilde{u}\partial_{u_{i}u}^{2}G+|\tilde{u}|^{2}\partial
_{uu}^{2}G,
\end{align*}
$M_{1}^{G}(u)$ and $M_{2}^{G}(u,u_{k})$ are some continuous
functions which grow with $|u|$ and $u_{k}u_{k}$ and
$M_{2}^{G}\geq1$.
\end{definition}

\begin{remark}
\label{Kr-re2} Let $\varepsilon I\leq A=(a_{ij})\leq KI$. It is easy
to check
that%
\[
G(u_{ij},u_{i},u)=\sum_{i,j}a_{ij}u_{ij}+\sum_{i}b_{i}u_{i}+cu
\]
belongs to $\mathcal{G}(\varepsilon,K,Q)$.
\end{remark}

The following definition is Definition 6.1.1 in Krylov
{\cite{Krylov1}.}

\begin{definition}
Let a function $G=G(u_{ij},u_{i},u):\mathbb{S}(d)\times \mathbb{R}^{d}%
\times \mathbb{R}\rightarrow \mathbb{R}$ be given. We write $G\in \bar{\mathcal{G}%
}(\varepsilon,K,Q)$ if there exists a sequence $G_{n}\in \mathcal{G}%
(\varepsilon,K,Q)$ converging to $G$ as $n\rightarrow \infty$\ at
every point $(u_{ij},u_{i},u)\in \mathbb{S}(d)\times
\mathbb{R}^{d}\times \mathbb{R}$ such that

\begin{description}
\item[\textup{(i)}] $M_{i}^{G_{1}}=M_{i}^{G_{2}}=\cdots=:M_{i}^{G}$, $i=1,2;$

\item[\textup{(ii)}] for each $n=1,2,\ldots$, the function $G_{n}$ is infinitely
differentiable with respect to $(u_{ij},u_{i},u);$

\item[\textup{(iii)}] there exist constants $\delta_{0}=:\delta_{0}^{G}>0$ and
$M_{0}=:M_{0}^{G}>0$ such that
\[
G_{n}(u_{ij},0,-M_{0})\geq \delta_{0},\ G_{n}(-u_{ij},0,M_{0})\leq-\delta_{0}%
\]
for each $n\geq1$ and symmetric nonnegative matrices $(u_{ij})$.
\end{description}
\end{definition}

The following theorem is Theorem 6.4.3 in Krylov  {\cite{Krylov1} },
which plays important role in our proof of central limit theorem.

\begin{theorem}
\label{Kr-th4} Suppose that $G\in
\bar{\mathcal{G}}(\varepsilon,K,Q)$ and
$\varphi \in C_{b}(\mathbb{R}^{d})$ with $\sup_{x\in \mathbb{R}^{d}}%
|\varphi(x)|\leq M_{0}^{G}$. Then PDE \textup{(\ref{eq-Kr})} has a
solution $u$ possessing the following properties:

\begin{description}
\item[\textup{(i)}] $u\in C([0,T]\times \mathbb{R}^{d})$, $|u|\leq M_{0}^{G}$ on $Q$;

\item[\textup{(ii)}] there exists a constant $\alpha \in(0,1)$ only depending on
$d,K,\varepsilon$ such that for each $\kappa>0$,
\begin{equation}
||u||_{C^{1+\alpha/2,2+\alpha}([0,T-\kappa]\times
\mathbb{R}^{d})}<\infty.
\end{equation}

\end{description}
\end{theorem}

\bigskip

Now we consider the $G$-equation. Let $G:\mathbb{R}^{d}\times \mathbb{S}%
(d)\rightarrow \mathbb{R}$ be a given continuous sublinear function
monotonic in $A\in \mathbb{S}(d)$. Then there exists a bounded,
convex and closed subset $\Sigma \subset \mathbb{R}^{d}\times
\mathbb{S}_{+}(d)$ such that
\begin{equation}
G(p,A)=\sup_{(q,B)\in \Sigma}[\frac{1}{2}\mathrm{tr}[AB]+\left
\langle p,q\right \rangle ]\  \  \  \text{for}\ (p,A)\in
\mathbb{R}^{d}\times
\mathbb{S}(d).\label{re-G}%
\end{equation}
The $G$-equation is%
\begin{equation}
\partial_{t}u+G(Du,D^{2}u)=0,\  \ u(T,x)=\varphi(x).\label{eq-G}%
\end{equation}
We set
\begin{equation}
\tilde{u}(t,x)=e^{t-T}u(t,x).\label{translation}%
\end{equation}
It is easy to check that $\tilde{u}$ satisfies the following PDE:%
\begin{equation}
\partial_{t}\tilde{u}+G(D\tilde{u},D^{2}\tilde{u})-\tilde{u}=0,\  \  \tilde
{u}(T,x)=\varphi(x).
\end{equation}
Suppose that there exists a constant $\varepsilon>0$ such that for each $A,\bar{A}%
\in \mathbb{S}(d)$ with $A\geq \bar{A}$, we have%
\begin{equation}
G(0,A)-G(0,\bar{A})\geq \varepsilon \mathrm{tr}[A-\bar{A}].\label{nonde}%
\end{equation}
Since $G$ is continuous, it is easy to prove that there exists a
constant $K>0$ such
that for each $A,\bar{A}\in \mathbb{S}(d)$ with $A\geq \bar{A}$, we have%
\begin{equation}
G(0,A)-G(0,\bar{A})\leq K\mathrm{tr}[A-\bar{A}].
\end{equation}
Thus for each $(q,B)\in \Sigma$, we have%
\[
2\varepsilon I\leq B\leq2KI.
\]
By Remark \ref{Kr-re2}, it is easy to check that $\tilde{G}(u_{ij}%
,u_{i},u):=G(u_{i},u_{ij})-u\in \bar{\mathcal{G}}(\varepsilon,K,Q)$
and $\delta_{0}^{G}=M_{0}^{G}$ can be any positive constant. By
Theorem \ref{Kr-th4} and (\ref{translation}), we have the following
regularity estimate for $G$-equation (\ref{eq-G}).

\begin{theorem}
Let $G$ satisfy \textup{(\ref{re-G})} and \textup{(\ref{nonde})}, $\varphi \in C_{b}%
(\mathbb{R}^{d})$ and let $u$ be a solution of $G$-equation
\textup{(\ref{eq-G})}. Then there exists a constant $\alpha
\in(0,1)$ only depending on $d,G,\varepsilon$ such that for each
$\kappa>0$,
\begin{equation}
||u||_{C^{1+\alpha/2,2+\alpha}([0,T-\kappa]\times
\mathbb{R}^{d})}<\infty.
\end{equation}
\end{theorem}

%
%

\pagestyle{bibliography}

\chapter*{Index of Symbols}
\addcontentsline{toc}{chapter}{Index of Symbols} \pagestyle{empty}
\normalsize \noindent\begin{tabular}{ll} \vspace{0.3cm}

$\mathcal{A}$ & Coherent acceptable set\quad\pageref{mathcalA}\\


$B_b(\Omega)$ & Space of bounded functions on $\Omega$\quad \pageref{bbomega}\\

$B$ & $G$-Brownian motion\quad\pageref{bbb}\\
\vspace{0.3cm}

$\langle B\rangle$ & Quadratic variation process of $B$\quad\pageref{qvb}\\


$C_b(\Omega)$ & Space of bounded and continuous functions on
$\Omega$\quad \pageref{cbomega}\\

$C_{b,Lip}(\mathbb{R}^n)$ \quad\quad\quad& Space of bounded and Lipschitz
continuous functions on $\mathbb{R}^n$\quad\pageref{cbl}\\
\vspace{0.3cm}
$C_{l,Lip}(\mathbb{R}^n)$ & Space of locally Lipschitz functions\quad\pageref{cllip}\\


$\mathbb{E}$ & Sublinear expectation\quad\pageref{sube}\\
\vspace{0.3cm}
$\hat{\mathbb{E}}$ & $G$-expectation\quad\pageref{ge}\\


\vspace{0.3cm} $\mathcal{H}$ & Space of random variables\quad
\pageref{huah}\\

$L^0(\Omega)$& Space of all $\mathcal{B}(\Omega)$-measurable real
functions\quad \pageref{l0omega}\\

$\mathbb{L}^p_b$ & Completion of $B_b(\Omega)$ under norm
$||\cdot||_p$\quad\pageref{Lpb}\\

\vspace{0.3cm}
$\mathbb{L}^p_c$ & Completion of $C_b(\Omega)$ under norm $||\cdot||_p$\quad\pageref{Lpc}\\


$M^{p,0}_G(0,T)$ & Space of simple processes\quad\pageref{mp0}\\

$M^p_G(0,T)$ & Completion of $M^{p,0}_G(0,T)$ under norm $\{\hat{\mathbb{E}}[\int_0^T|\eta_t|^p dt]\}^{1/p}$\quad\pageref{mgp}\\
\vspace{0.3cm}
$\bar{M}^p_G(0,T)$ & Completion of $M^{p,0}_G(0,T)$ under norm $\{\int_0^T\hat{\mathbb{E}}[|\eta_t|^p dt]\}^{1/p}$\quad\pageref{mbargp}\\


\vspace{0.3cm}
q.s. & Quasi-surely\quad\pageref{qs}\\


$\mathbb{S}(d)$ & Space of $d\times d$ symmetric matrices\quad\pageref{sd}\\

\vspace{0.3cm}
$\mathbb{S}_+(d)$ & Space of non-negative $d\times d$ symmetric
matrices\quad\pageref{splusd}\\


$\rho$ & Coherent risk measure\quad\pageref{rho}\\

$(\Omega,\mathcal{H},{\mathbb{E}})$ & Sublinear expectation
space\quad\pageref{subls}\\

$\overset{d}{=}$ & Identically distributed\quad\pageref{iidd}\\

$\langle x,y \rangle$ & Scalar product of $x,y\in \mathbb{R}^n$\\

$|x|$ & Euclidean norm of $x$\\

$(A,B)$ & Inner product $(A,B):=tr[AB]$\\
\end{tabular}

\newpage

\addcontentsline{toc}{chapter}{Index} \printindex \pagestyle{index}





\end{document}